\documentclass[11pt]{preprint}
\usepackage[left=2.7cm, right=2.7cm, bottom=3cm, top=3.5cm]{geometry}
\usepackage[full]{textcomp}
\usepackage[osf]{newtxtext}
\usepackage{comment}
\usepackage{amssymb,amsmath}
\usepackage{mathtools}
\usepackage{mhenvs,mhsymb,mhequ}
\usepackage{booktabs}
\usepackage{tikz}
\usepackage{mathrsfs}
\usepackage{longtable}
\usepackage{microtype}
\usepackage{wasysym}

\usepackage{float}

\usepackage[normalem]{ulem}

\usepackage{centernot}
\usepackage{enumitem}
\usepackage{breakurl}
\usepackage{hyperref}
\usepackage{color,xcolor}
\usepackage{cases}
\usepackage{xsavebox}
\usepackage{tikz-cd}

\mathtoolsset{showonlyrefs=true}

\def\cX{\mathcal{X}}
\def\cF{\mathcal{F}}
\def\cI{\mathcal{I}}
\def\bZ{\mathbf{Z}}
\def\bR{\mathbf{R}}
\def\cP{\mathcal {P}}

\newif\ifdark
\IfFileExists{dark}{\darktrue}{\darkfalse}

\ifdark

\definecolor{darkred}{rgb}{0.9,0.2,0.2}
\definecolor{darkblue}{rgb}{0.7,0.3,1}
\definecolor{darkgreen}{rgb}{0.1,0.9,0.1}
\definecolor{pagebackground}{rgb}{.15,.21,.18}
\definecolor{pageforeground}{rgb}{.84,.84,.85}
\pagecolor{pagebackground}
\AtBeginDocument{\globalcolor{pageforeground}}

\else

\definecolor{darkred}{rgb}{0.7,0.1,0.1}
\definecolor{darkblue}{rgb}{0.4,0.1,0.8}
\definecolor{darkgreen}{rgb}{0.1,0.7,0.1}
\definecolor{pagebackground}{rgb}{1,1,1}
\definecolor{pageforeground}{rgb}{0,0,0}

\fi

\definecolor{newred}{RGB}{208,16,76}
\definecolor{newgreen}{RGB}{34,125,81}

\definecolor{LB}{rgb}{0.29, 0.63, 0.73}

\makeatletter
\newcommand{\globalcolor}[1]{%
	\color{#1}\global\let\default@color\current@color
}
\makeatother

\DeclareSymbolFont{timesoperators}{T1}{ptm}{m}{n}
\SetSymbolFont{timesoperators}{bold}{T1}{ptm}{b}{n}

\makeatletter
\renewcommand{\operator@font}{\mathgroup\symtimesoperators}
\makeatother

\DeclareMathAlphabet{\mathbbm}{U}{bbm}{m}{n}
\DeclareFontFamily{U}{BOONDOX-calo}{\skewchar\font=45 }
\DeclareFontShape{U}{BOONDOX-calo}{m}{n}{
	<-> s*[1.05] BOONDOX-r-calo}{}
\DeclareFontShape{U}{BOONDOX-calo}{b}{n}{
	<-> s*[1.05] BOONDOX-b-calo}{}
\DeclareMathAlphabet{\mcb}{U}{BOONDOX-calo}{m}{n}
\SetMathAlphabet{\mcb}{bold}{U}{BOONDOX-calo}{b}{n}

\makeatletter 

\newcommand*{\fat}{}
\DeclareRobustCommand*{\fat}{%
	\mathbin{\mathpalette\bigcdot@{}}}
\newcommand*{\bigcdot@scalefactor}{.5}
\newcommand*{\bigcdot@widthfactor}{1.15}
\newcommand*{\bigcdot@}[2]{%
	\sbox0{$#1\vcenter{}$}
	\sbox2{$#1\cdot\m@th$}%
	\hbox to \bigcdot@widthfactor\wd2{%
		\hfil
		\raise\ht0\hbox{%
			\scalebox{\bigcdot@scalefactor}{%
				\lower\ht0\hbox{$#1\bullet\m@th$}%
			}%
		}%
		\hfil
	}%
}

\DeclareRobustCommand{\TitleEquation}[2]{\texorpdfstring{\StrLeft{\f@series}{1}[\@firstchar]$\if b\@firstchar\boldsymbol{#1}\else#1\fi$}{#2}}

\makeatother

\numberwithin{equation}{section}

\theoremnumbering{Alph}
\renewtheorem{theorem*}{Theorem}
\newtheorem{propositionA}[theorem*]{Proposition}

\newtheorem{convention}[lemma]{Convention}

\def\slash{\leavevmode\unskip\kern0.18em/\penalty\exhyphenpenalty\kern0.18em}
\def\dash{\leavevmode\unskip\kern0.18em--\penalty\exhyphenpenalty\kern0.18em}

\def\emptyset{{\centernot\Circle}}

\def\part{\mathrm {part}}

\def\i{\mathrm{i}}

\def\C{\mathbf{C}}
\def\cC{\mathcal C}
\def\D{\mathcal{D}}

\def\cF{{\mathcal F}}
\def\F{\mathcal F}

\def\H{\mathcal{H}}
\def\HH{\mathbb H}
\def\I{\mathcal{I}}
\def\L{\mathcal{L}}
\def\J{\mathcal{J}}
\def\X{\mathcal{X}}

\def\N{\mathbf{N}}
\def\R{\mathbf{R}}
\def\T{\mathbf{T}}
\def\Z{\mathbf{Z}}
\def\E{\mathbb{E}}

\def\f{\frac}
\def\cS{\mathcal {S}}

\def\B{\mathcal{B}}

\def\cJ{\mathcal {J}}
\def\cF{\mathcal {F}}
\def\B{\mathcal{B}}
\def\P{\mathbf{P}}
\def\W{\mathcal {W}}

\def\d{\mathrm{d}}

\def\pt{\partial_{t}}

\def\bn{\boldsymbol{n}}
\def\emptyset{{\centernot\Circle}}

\renewcommand{\leq}{\leqslant}
\renewcommand{\geq}{\geqslant}

\newcommand{\brt}[1]{\left[ #1 \right]}
\newcommand{\pa}[1]{\left( #1 \right)}
\newcommand{\jpb}[1]{\langle #1 \rangle}
\newcommand{\norm}[1]{\left\| #1 \right\|}
\newcommand{\ab}[1]{\left| #1 \right|}

\newcommand{\maxcurly}[1]{\max \{ #1 \}}
\newcommand{\mincurly}[1]{\min \{ #1 \}}
\newcommand{\medcurly}[1]{\mathrm{med}\{ #1 \}}

\def\Sym{\mathrm {Sym}}

\def\1{{\mathbf 1}}

\makeatletter
\DeclareRobustCommand{\cev}[1]{
  {\mathpalette\do@cev{#1}}%
}
\newcommand{\do@cev}[2]{%
  \vbox{\offinterlineskip
    \sbox\z@{$\m@th#1 x$}%
    \ialign{##\cr
      \hidewidth\reflectbox{$\m@th#1\vec{}\mkern7mu$}\hidewidth\cr
      \noalign{\kern-\ht\z@}
      $\m@th#1#2$\cr
    }%
  }%
}
\makeatother

\usepackage{tikz}
\definecolor{gr}{rgb}   {0.,   0.69,   0.23 }
\definecolor{bl}{rgb}   {0.,   0.5,   1. }
\definecolor{mg}{rgb}   {0.85,  0.,    0.85}
\definecolor{yl}{rgb}   {0.8,  0.7,   0.}
\definecolor{or}{rgb}  {0.7,0.2,0.2}

\usepackage{marginnote}
\usepackage{scalerel} 
\usetikzlibrary{shapes.misc}
\usetikzlibrary{shapes.symbols}
\usetikzlibrary{decorations}
\usetikzlibrary{decorations.markings}

\tikzset{
	dot/.style={circle,fill=black,draw=black,inner sep=0pt,minimum size=0.5mm},
	>=stealth,
}
\tikzset{
	dot2/.style={circle,fill=black,draw=black,inner sep=0pt,minimum size=0.2mm},
	>=stealth,
}
\tikzset{
	ddot/.style={circle,fill=black,draw=black,inner sep=0pt,minimum size=0.8mm},
	>=stealth,
}
\tikzset{decision/.style={ 
		draw,
		diamond,
		aspect=1.5
}}
\tikzset{dia2/.style
	={diamond,fill=white,draw=black,inner sep=0pt,minimum size=1mm},
	>=stealth,
}
\tikzset{dia/.style
	={star,fill=black,draw=black,inner sep=0pt,minimum size=1mm},
	>=stealth,
}
\tikzset{dia/.style
	={diamond,fill=black,draw=black,inner sep=0pt,minimum size=1.3mm},
	>=stealth,
}

\makeatletter
\def\DeclareSymbol#1#2#3{\xsavebox{#1}{\tikz[baseline=#2,scale=0.15]{#3}}}
\def\<#1>{\xusebox{#1}}
\makeatother
\newcommand{\pe}{\mathbin{\scaleobj{0.7}{\tikz \draw (0,0) node[shape=circle,draw,inner sep=0pt,minimum size=8.5pt] {\scriptsize  $=$};}}}

\newcommand{\pez}{\mathbin{\scaleobj{0.7}{\tikz \draw (0,0) node[shape=circle,draw,
			fill=white, 
			inner sep=0pt,minimum size=8.5pt]{} ;}}}

\usetikzlibrary{decorations.pathmorphing, patterns,shapes}
\DeclareSymbol{dot}{-2.7}{\draw (0,0) node[dot]{};}

\DeclareSymbol{1}{0}{\draw[white] (-.4,0) -- (.4,0); \draw (0,0)  -- (0,1.2) node[dot] {};}
\DeclareSymbol{2}{0}{\draw (-0.5,1.2) node[dot] {} -- (0,0) -- (0.5,1.2) node[dot] {};}
\DeclareSymbol{3}{0}{\draw (0,0) -- (0,1.2) node[dot] {}; \draw (-.7,1) node[dot] {} -- (0,0) -- (.7,1) node[dot] {};}
\DeclareSymbol{30}{-3}{\draw (0,0) -- (0,-1); \draw (0,0) -- (0,1.2) node[dot] {}; \draw (-.7,1) node[dot] {} -- (0,0) -- (.7,1) node[dot] {};}
\DeclareSymbol{31}{-3}{\draw (0,0) -- (0,-1) -- (1,0) node[dot] {}; \draw (0,0) -- (0,1.2) node[dot] {}; \draw (-.7,1) node[dot] {} -- (0,0) -- (.7,1) node[dot] {};}
\DeclareSymbol{32}{-3}{\draw (0,0) -- (0,-1) -- (1,0) node[dot] {}; \draw (0,0) -- (0,-1) -- (-1,0) node[dot] {}; \draw (0,0) -- (0,1.2) node[dot] {}; \draw (-.7,1) node[dot] {} -- (0,0) -- (.7,1) node[dot] {};}
\DeclareSymbol{20}{-3}{\draw (0,0) -- (0,-1);\draw (-.7,1) node[dot] {} -- (0,0) -- (.7,1) node[dot] {};}
\DeclareSymbol{22}{-3}{\draw (0,0.3) -- (0,-1) -- (1,0) node[dot] {}; \draw (0,0.3) -- (0,-1) -- (-1,0) node[dot] {};\draw (-.7,1) node[dot] {} -- (0,0.3) -- (.7,1) node[dot] {};}
\DeclareSymbol{31p}{-3}{\draw (0,0) -- (0,-1) -- (1,0) node[dot] {}; \draw (0,0) -- (0,1.2) node[dot] {}; \draw (-.7,1) node[dot] {} -- (0,0) -- (.7,1) node[dot] {}; 
	\draw (0,-1) node{\scaleobj{0.5}{\pez}}; 
	\draw (0,-1) node{\scaleobj{0.5}{\pe}}; }
\DeclareSymbol{32p}{-3}{\draw (0,0) -- (0,-1) -- (1,0) node[dot] {}; \draw (0,0) -- (0,-1) -- (-1,0) node[dot] {}; \draw (0,0) -- (0,1.2) node[dot] {}; \draw (-.7,1) node[dot] {} -- (0,0) -- (.7,1) node[dot] {}; 
	\draw (0,-1) node{\scaleobj{0.5}{\pez}};
	\draw (0,-1) node{\scaleobj{0.5}{\pe}};}
\DeclareSymbol{22p}{-3}{\draw (0,0.3) -- (0,-1) -- (1,0) node[dot] {}; \draw (0,0.3) -- (0,-1) -- (-1,0) node[dot] {};\draw (-.7,1) node[dot] {} -- (0,0.3) -- (.7,1) node[dot] {}; 
	\draw (0,-1) node{\scaleobj{0.5}{\pez}};
	\draw (0,-1) node{\scaleobj{0.5}{\pe}};}
\DeclareSymbol{21p}{-3}{\draw (0,0.3) -- (0,-1) -- (1,0) node[dot] {}; 
	\draw (-.7,1) node[dot] {} -- (0,0.3) -- (.7,1) node[dot] {}; 
	\draw (0,-1) node{\scaleobj{0.5}{\pez}};
	\draw (0,-1) node{\scaleobj{0.5}{\pe}};}

\DeclareSymbol{21}{-3}{\draw (0,0.3) -- (0,-1) -- (1,0) node[dot] {}; 
	\draw (-.7,1) node[dot] {} -- (0,0.3) -- (.7,1) node[dot] {}; 
}

\DeclareSymbol{70}{-3}{
	\draw (0,-0.4) -- (0,-1.2);
	\draw (0,0.8) node[dot] {} -- (0,-0.4) -- (0.8,0.6);
	\draw (0.8, .6) -- (1.4, 1.35)  node[dot] {};
	\draw (0.8, .6) -- (1.8, .8)  node[dot] {};
	\draw (0.8, .6) -- (0.85, 1.6)  node[dot] {};
	\draw (0,.6) -- (0,-0.4) -- (-0.8,0.6) ;
	\draw (-0.8, 0.6) -- (-1.4, 1.35)  node[dot] {};
	\draw (-0.8, .6) -- (-1.8, .8)  node[dot] {};
	\draw (-0.8, .6) -- (-0.85, 1.6)  node[dot] {};
}

\DeclareSymbol{7}{-3}{
	\draw (0,0.2) node[dot] {} -- (0,-1) -- (0.8,0);
	\draw (0.8, 0) -- (1.4, 0.75)  node[dot] {};
	\draw (0.8, 0) -- (1.8, 0.2)  node[dot] {};
	\draw (0.8, 0) -- (0.85, 1)  node[dot] {};
	\draw (0,0) -- (0,-1) -- (-0.8,0) ;
	\draw (-0.8, 0) -- (-1.4, 0.75)  node[dot] {};
	\draw (-0.8, 0) -- (-1.8, 0.2)  node[dot] {};
	\draw (-0.8, 0) -- (-0.85, 1)  node[dot] {};
}

\DeclareSymbol{11p}{-3}{\draw (0,0) node[dot2] {} -- (0,-1) -- (1,0) node[dot] {}; 
	\draw (0,0) -- (0,1.2) node[dot] {};  
	\draw (0,-1) node{\scaleobj{0.5}{\pez}}; 
	\draw (0,-1) node{\scaleobj{0.5}{\pe}}; }
\DeclareSymbol{100}{-3}
{\draw[white] (-.4,0) -- (.4,0); 
	\draw (0,0) node[dot2] {} -- (0,-1)  ; 
	\draw (0,0) -- (0,1.2) node[dot] {};}  
\DeclareSymbol{320}{-3}{
	\draw (0,0.4) -- (0,-1.2);
	\draw (0,0.6) -- (0,-0.4) -- (1,0.6) node[dot] {};
	\draw (0,0.6) -- (0,-0.4) -- (-1,0.6) node[dot] {};
	\draw (0,0.6) -- (0,1.8) node[dot] {};
	\draw (-.7, 1.6) node[dot] {} -- (0,0.6) -- (.7,1.6) node[dot] {};}

\DeclareSymbol{90}{-3}{%
  \draw (0,-1.2) -- (0, 0.7);
  \draw (0, 0.7) -- (-0.6, 1.7)  node[dot] {};
  \draw (0, 0.7) -- (0.6, 1.7)  node[dot] {};
  \draw (0, 0.7) -- (0, 1.9)  node[dot] {};
  \draw (0,-0.4) -- (0.8, 0.2);
  \draw (0.8, 0.2) -- (1.1, 1.2)  node[dot] {};
  \draw (0.8, 0.2) -- (1.7, 0.9)  node[dot] {};
  \draw (0.8, 0.2) -- (1.8, 0.3)   node[dot] {}; 
  \draw (0,-0.4) -- (-0.8, 0.2);
  \draw (-0.8, 0.2) -- (-1.1, 1.2)  node[dot] {};
  \draw (-0.8, 0.2) -- (-1.7, 0.9)  node[dot] {};
  \draw (-0.8, 0.2) -- (-1.8, 0.3)   node[dot] {}; 
  }
 
\DeclareSymbol{9}{-3}{%
  \draw (0,-0.44) -- (0, 0.7);
  \draw (0, 0.7) -- (-0.6, 1.7)  node[dot] {};
  \draw (0, 0.7) -- (0.6, 1.7)  node[dot] {};
  \draw (0, 0.7) -- (0, 1.9)  node[dot] {};
  \draw (0,-0.4) -- (0.8, 0.2);
  \draw (0.8, 0.2) -- (1.1, 1.2)  node[dot] {};
  \draw (0.8, 0.2) -- (1.7, 0.9)  node[dot] {};
  \draw (0.8, 0.2) -- (1.8, 0.3)   node[dot] {}; 
  \draw (0,-0.4) -- (-0.8, 0.2);
  \draw (-0.8, 0.2) -- (-1.1, 1.2)  node[dot] {};
  \draw (-0.8, 0.2) -- (-1.7, 0.9)  node[dot] {};
  \draw (-0.8, 0.2) -- (-1.8, 0.3)   node[dot] {}; 
  }

\usetikzlibrary{mindmap,backgrounds,trees,arrows,decorations.pathmorphing,decorations.pathreplacing,positioning,shapes.geometric,shapes.misc,decorations.markings,decorations.fractals,calc,patterns}

\tikzset{>=stealth',
	cvertex/.style={circle,draw=black,inner sep=1pt,outer sep=3pt},
	vertex/.style={circle,fill=black,inner sep=1pt,outer sep=3pt},
	star/.style={circle,fill=yellow,inner sep=0.75pt,outer sep=0.75pt},
	tvertex/.style={inner sep=1pt,font=\scriptsize},
	gap/.style={inner sep=0.5pt,fill=white}}
\tikzstyle{mybox} = [draw=black, fill=blue!10, very thick,
rectangle, rounded corners, inner sep=10pt, inner ysep=20pt]
\tikzstyle{boxtitle} =[fill=blue!50, text=white,rectangle,rounded corners]

\tikzstyle{decision} = [diamond, draw, fill=blue!20,
text width=4.5em, text badly centered, node distance=2.5cm, inner sep=0pt]
\tikzstyle{block} = [rectangle, draw, fill=blue!20,
text width=2.7cm, text centered, rounded corners, minimum height=3em]

\tikzstyle{line} = [draw, very thick, color=black!50, -latex']

\tikzstyle{cloud} = [draw, ellipse,fill=red!40, 
node distance=2.5cm,
minimum height=1em]

\tikzstyle{cloud2} = [draw, ellipse,fill=red!30, text=white,text width=10em, node distance=2.5cm, text centered, minimum height=4em]

\tikzstyle{cloud3} = [draw, ellipse, fill=cyan!30, 
node distance=2.5cm,
minimum height=3em, 
text width=1.7cm]

\tikzstyle{cloud4} = [draw, ellipse,fill=orange!70, node distance=2.5cm,
text width=2.1cm,  
minimum height=3em]

\tikzstyle{cloud5} = [draw, ellipse,fill=red!20, node distance=2.5cm,
text centered, 
text width=3.0cm, 
minimum height=3em]

\tikzstyle{cloud6} = [draw, ellipse,fill=red!20, node distance=2.5cm,
text width=1.5cm, 
text centered, 
minimum height=3em]

\tikzset{
	position/.style args={#1:#2 from #3}{
		at=(#3.#1), anchor=#1+180, shift=(#1:#2)
	}
}

\begin{document}
	\title{Stochastic nonlinear wave equation with rougher than white noise}
	\author{Xue-Mei Li$^{1,2}$ and Xianfeng Ren$^2$}
	\institute{EPFL, Switzerland. \email{xue-mei.li@epfl.ch} \and Imperial, U.K. \email{xue-mei.li@imperial.ac.uk, x.ren@imperial.ac.uk}} 
	\maketitle

	\begin{abstract}
We study the singular stochastic wave equation on $\T^2$, with a cubic nonlinearity and Gaussian rough `Matérn' forcing (a Fourier multiplier of order $\alpha>0$ applied to space-time white noise) and establish local well-posedness for $\alpha < \tfrac{3}{8}$. This extends \cite{GKO18} beyond white noise and strengthens the quadratic-case result \cite{OO21} ($\alpha<\f 12$). Our argument develops new trilinear estimates in Bourgain spaces together with  case-specific cubic counting estimates. 
\end{abstract}



	\tableofcontents

\section{Introduction}\label{sec:introduction}

 We study the stochastic wave equation with cubic non-linearity (SNLW) in the two dimensional torus.
\begin{equ}\label{eq:SNLW}
\begin{cases}
(\partial^2_t+1-\Delta)u=-u^3+  \jpb{\nabla}^{\alpha} \xi,\\
(u,  \pt u)|_{t=0}= (u_0,u_1),
\end{cases}
\end{equ}
where $\xi$ denotes mean-zero space-time white noise and $\alpha>0$.

The cubic nonlinearity $u^3$ represents nonlinear self-interaction of the wave field, while the stochastic forcing $\jpb{\nabla}^{\alpha} \xi$ models a rough external input, capturing highly oscillatory or turbulent behaviour. Our focus is on the local well-posedness theory— it is well known that wave equations with nonlinearities of degree three or higher may develop shocks or blow up in finite time. Classical counterexamples for deterministic equations were given in \cite{Reed, John};  for stochastic equations with sufficiently regular noise, see \cite{Chow02}.  The central question is:  how rough can the noise be while still permitting a local well-posed solution theory?

Stochastic perturbations of low regularity of mathematical naturally arise in mathematical models of turbulence and complex microscopic interactions. Of particular interest are Fourier multipliers of the form  $F(\Delta)\xi$, which preserve the spatial homogeneity of the white noise~$\xi$.
 In the time-independent negative $\alpha$ setting, the corresponding covariance functions are known as Mat\'ern covariance functions, which admit explicit expressions in terms of Bessel functions \cite{GH90, Ste70, LRL11, Whi53, BW16, LSEO17}.  

Let $\Delta$ denote the Laplacian and $\jpb{\nabla} := \sqrt{1 - \Delta}$  the Bessel potential operator.  For $\alpha>0$, the Gaussian noise
\[
\jpb{\nabla}^{\alpha} \xi(x,t), \quad t \geq 0, \, x \in \mathbb{T}^2,
\]
takes values in the space of distributions. Indeed, for any $\epsilon > 0$,  almost surely,
\[
\jpb{\nabla}^{\alpha} \xi \in W^{-1-\alpha-\epsilon,\infty}_x(\T^2).
\]
This places the equation within the scope of singular stochastic partial differential equations (SPDEs).

A central challenge in singular SPDEs is the definition of products of distributions, often guided by informal regularity heuristics. Typically, a product of functions with regularities $s_1$ and $s_2$ is well-defined if $s_1 + s_2 > 0$. When $s_1 > 0$ and $s_1 \geq s_2$, the product typically inherits the lower regularity $s_2$. In the stochastic setting, such products, even when requiring renormalisation, can often be given a consistent meaning, with the resulting regularity given by $\min(s_1, s_2, s_1+s_2)$.
This perspective is closely related to techniques developed for parabolic singular SPDEs, which have recently seen tremendous advances—most notably through the theory of regularity structures~\cite{Hai14} and the paracontrolled calculus~\cite{GIP15}. Related ideas also appear in the classical renormalisation group approach; see, e.g.,~\cite{Kup16}.

In the stochastic wave equation setting,  the linear propagator $\mathcal{J} :=(\partial_t^2 + 1 - \Delta)^{-1}$  improves the regularity of the forcing term by at most one order.  

For the corresponding deterministic linear equation, the solution is given explicitly \begin{equ}\label{St}
\boldsymbol{S}(t)(u_0, u_1)=\cos (t\sqrt{1 - \Delta})u_0+\f{\sin (t\sqrt{1 - \Delta})}{\sqrt{1 - \Delta}}u_1.
\end{equ}
 If the initial displacement satisfies $u_0\in H^{s}_x(\T^2)$ and $u_1 \in H^{s-1}_x(\T^2)$, then the solution can at best be expected to lie in $H^{s}_x(\T^2)$- the same regularity as the initial data.

For the stochastic linear equation, the mild solution takes the form
$$\boldsymbol{S}(t)(u_0, u_1) + \sum_{n \in \Z^2} e^{in\cdot x}\int_0^t \f{\sin ((t-s)\sqrt{1 + |n|^2})}{(1 + |n|^2)^{\f {1-\alpha}{2} }} \d \beta_n (s).$$
where $\{\beta_n\}_{n\in \Z^2}$ are standard complex-valued Brownian motions satisfying the symmetric condition described in \S\ref{sec:Wiener}. With this formulation, the real-valued Wiener integral improves the regularity of the driving Brownian motions by at best $1-\alpha$, reflecting the action of the fractional operator $\jpb{\nabla}^\alpha$ encoded in the stochastic forcing term. Consequently, any solution $u$ to~\eqref{eq:SNLW} can at best lie in $W^{-\alpha-\epsilon,\infty}_x(\T^2)$. For $\alpha > 0$, $u$ is expected to be a distribution rather than a function, and the cubic nonlinearity $u^3$ is ill defined in the classical sense. This makes~\eqref{eq:SNLW} a genuinely singular SPDE, and any meaningful solution theory necessarily requires both regularisation and renormalisation.

 A natural first step is to regularise the noise. Let $\xi_{N} := \P_{N} \xi$ denote a spectrally truncated version of $\xi$, where
$
\widehat{\P_{N} \xi}(n) := \mathbf{1}_{\{|n| \le N\}}\, \widehat{\xi}(n)$,
with $\widehat{\xi}(n)$ the spatial Fourier coefficients of $\xi$. We then consider the regularised equation:
\begin{equation}\label{eq:truncated-SNLW}
\begin{cases}
(\partial^2_t + 1 - \Delta) \tilde u = -\tilde u^3 + \jpb{\nabla}^{\alpha} \xi_{N},\\
(\tilde u, \pt \tilde u)|_{t=0} = (u_0, u_1),
\end{cases}
\end{equation}
which admits smooth solutions for each fixed $N$.

However, as shown in~\cite{AHR96},  for certain nonlinearities $F$, solutions to
$(\partial^2_t+1-\Delta)u=F(u)+\xi_{\le N}$
may converge to solutions of the linear equation, with the nonlinear effects disappearing entirely in the limit—a phenomenon known as \emph{triviality}. This highlights the inadequacy of regular approximations unless the nonlinearity is properly renormalised to account for its interaction with the singular noise. See also \cite{OOR20}.

Before introducing renormalization, we recall a common strategy in the solution theory of singular SPDEs: to seek a solution of the form
\[
u = \Psi + v,
\]
where $\Psi$ solves the linear stochastic equation with additive rough noise, and $v$ is a regular remainder. This decomposition isolates the singular effects in $\Psi$, while $v$ satisfies a nonlinear but better-posed equation. The underlying principle—common to both nonlinear ODEs and PDEs—is to separate the linear and nonlinear components. This strategy was successfully applied by Da Prato and Debussche~\cite{DD03} for parabolic SPDEs and by Bourgain~\cite{Bou96} in dispersive settings; see also~\cite{Mc95}.

For the non-linear stochastic wave equation, this strategy was implemented for additive white noise and Wick-renormalized nonlinearities on $\T^2$ in the work of Gubinelli, Koch, and Oh ~\cite{GKO18}. In the white noise case,  the linear solution belongs to the space $C([0,T]; W^{-\epsilon,\infty}_x(\T^2))$, ensuring that any power of the solution remains within the same space after suitable renormalization.  On $\T^3$, however, the linear solution lies in $C([0,T]; W^{-\frac{1}{2}-\epsilon,\infty}_x(\T^3))$, making the nonlinear analysis more delicate. The quadratic SNLW  $$ (\pt^{2}+1-\Delta)u= -u^2+  \xi $$ was resolved in the seminal work ~\cite{GKO24} using paracontrolled calculus. The paracontrolled structure, together with the tensor structure in the ansatz, was also used in \cite{Bri20b} to prove global existence for the cubic Hartree-type equation
$ (\pt^{2}+1-\Delta)u = -(V * u)u^2$
with random initial data, where $V$ is the kernel of  $\jpb{\nabla}^\alpha$ with $\alpha <0$,  extending results of Okamoto-Oh-Tolomeo \cite{OOT20}. In particular, the random tensor component was inspired by the random tensor operator developed by Deng, Nahmod, and Yue in \cite{DNY22, DNY23}, which was settled for the random initial problem for Schr\"{o}dinger equations. Also see for the wave equation, c.f. \cite{Bri20b}.

Stochastic wave equations have also been studied under random initial data~\cite{BT08a, BT08b, OPT22, ORT23, Bri20a, Bri20b, BDNY24},  fractional Gaussian noise~\cite{De19, De20}. Focusing on quadratic nonlinearities, Oh and Okamoto, \cite{OO21} established well-posedness for the equation $(\pt^{2}+1-\Delta)u = -u^2 +   \jpb{\nabla}^\alpha \xi $ in the range $\alpha \in (0, \f 12)$ on $\mathbb{T}^2$. The problem of cubic non-linearities remains open.

Without using paracontrolled calculus, Oh, Wang, and Zine \cite{OWZ22} proved local well-posedness for the cubic stochastic wave equation on $3d$ torus driven by noise smoother than white space-time
 noise: $$(\pt^{2}+1-\Delta)u = -u^3  +  \jpb{\nabla}^{-\alpha} \xi,$$
 where $\alpha>0$ and $-$ sign gives smoothness of the noise.
 
For cubic nonlinearities in two dimensions, the critical regularity of the initial data is $s = \tfrac{1}{4}$, corresponding to conformal invariance. In this work, we focus on determining how rough the driving noise $\jpb{\nabla}^{\alpha} \xi$ can be taken to be while still allowing for meaningful well-posedness theory.

\subsection{Overview of the analysis}
We begin by analyzing the regularized equation~\eqref{eq:truncated-SNLW}, writing $\tilde u_N = \<1>_N + \tilde v_N$, where $\<1>_N$ solves the linear stochastic wave equation with smooth noise $\jpb{\nabla}^\alpha \xi_{N}$ and zero initial condition. Explicitly,
\[
\<1>_N(t,x) := \sum_{|n|\le N} e_n(x) \int_0^t \frac{\sin((t-s)\jpb{n})}{\jpb{n}^{1-\alpha}}\, \d \beta_n(s),
\]
which converges in $C([0,T]; W^{-\alpha-,\infty}_x)$ almost surely 
\footnote{By $a+$, we refer to $a+\epsilon$ for arbitrarily small $\epsilon>0$. The same rule applies to $a-$.}. Expanding the nonlinearity yields terms involving powers of $\<1>_N$, which diverge as $N \to \infty$. These divergences are renormalised via Wick ordering:
\[
\<2>_N := (\<1>_N)^2 - \sigma_N(t), \qquad \<3>_N := (\<1>_N)^3 - 3 \sigma_N(t) \<1>_N,
\]
where $\sigma_N(t) = \mathbb{E}[|\<1>_N(t)|^2] \sim tN^{2\alpha}$ diverges as $N \to \infty$.

This leads to a renormalized equation for the regular component $v_N$:
\begin{equ}\label{eq:1st-re-vN}\begin{cases}
(\partial_t^2 + 1 - \Delta) v_N &= -(v_N)^3 - 3\<1>_N (v_N)^2 - 3\<2>_N v_N - \<3>_N\\
(v_N, \partial_t v_N)|_{t=0}& =\; (u_0,u_1).
\end{cases}
\end{equ}
Under suitable regularity assumptions on the stochastic objects and initial data, a local solution theory is developed for the limit equation via deterministic fixed-point arguments in Bourgain-type spaces.

The term $\<3>_N$ gives rise to the stochastic object $\<30> = \mathcal{J}(\<3>)$, which captures additional nonlinear interactions and satisfies
 $$(\partial_t^2+1 - \Delta) \<30>_N = \<3>_N, \quad \<30>_N(0)=0.$$ 
 The regularity of $\<30>$ plays a crucial role in determining the appropriate rough/regular decomposition. In particular, it must belong to the solution space in order for the fixed-point argument to apply to~(\ref{eq:1st-re-vN}).

The decomposition~\eqref{eq:1st-re-vN} is sufficient to treat the case $\alpha < \tfrac{1}{8}$, but it breaks down for rougher noise. When $\alpha \in [\tfrac{1}{8}, \tfrac{3}{10})$, multilinear smoothing effects—stemming from dispersive cancellations—still permit us to close the fixed-point argument. In this regime, nonlinear products such as $\cJ(\<2> v)$ are best interpreted as random linear operators acting on $v$, denoted by $\mathcal{J}^{\<2>} (v)$.
More precisely, for $\alpha \in [\f 1 8,\f 1 4)$, we have that  $\fJ^{\<2>}\in \L^{\frac{1}{4}+,\frac{1}{2}+}_T$ almost surely, while for $  \alpha \in [\f 1 4, \f 3 8)$,  the operator lies almost surely in the class $ \L^{2\alpha - \f  14+, \frac{1}{2}+}_T$ almost surely.

\begin{table}[H]
\centering
\small
\renewcommand{\arraystretch}{1.3}
\setlength{\tabcolsep}{5pt}
\begin{tabular}{|c|c|c|c|c|}
\hline
\textbf{Range of $\alpha$} & \textbf{$\Psi$} & \textbf{Stochastic Symbols} & \textbf{regularity of $v$} & \textbf{Remarks} \\
\hline
$(0, \tfrac{1}{8})$ & $\<1>$ & $\<2>, \<3>$ & $\mathcal{X}^{\frac{1}{4}+, \frac{1}{2}+}_T$ & 1sr order expansion, use (\ref{eq:1st-re-vN}) \\
\hline
$[\tfrac{1}{8}, \tfrac{1}{4})$ & $\<1>$ & $\<2>, \<3>$ & $\mathcal{X}^{\frac{1}{4}+, \frac{1}{2}+}_T$ & use (\ref{eq:1st-re-vN}) and random operator $\mathcal{J}^{\<2>}$ \\
\hline
$[\tfrac{1}{4}, \tfrac{3}{10})$ & $\<1>$ & $\<2>,\<30>$ & $\mathcal{X}^{2\alpha - \frac{1}{4}+, \frac{1}{2}+}_T$ & use (\ref{eq:1st-re-vN}) and $\mathcal{J}^{\<2>}$\\
\hline
$[\tfrac{3}{10}, \tfrac{3}{8})$ & $\<1> + \<30>$ & $\<2>, \<3>, \<30>, \<31>, \<320>, \<70>, \<90>$ & $\mathcal{X}^{2\alpha - \frac{1}{4}+, \frac{1}{2}+}_T$ & 2nd order expansion, use (\ref{eq:2nd-re-vN})\ \\
\hline
\end{tabular}
\caption{Decomposition strategy and stochastic symbols by range of $\alpha$}
\end{table}

By Proposition~\ref{prop:reg-of-cube}, for $\alpha < \frac{1}{4}$, the term $\<30>$, which lies in $\X_T^{1-2\alpha-, \f 12 +}$,
is sufficiently smooth to be absorbed into the smooth component of the solution. For $\alpha>\f 14$, Corollary~\ref{cor-object-3} shows that $\<30>$ gains only one order $\alpha$ from the multilinear smooth effect and  that $$\<30>\in \X_T^{\f 54-3\alpha, \f 12+},$$
see also  (\ref{eq:s-alpha}. Up to the threshold $\alpha=\f 3{10}$, the crucial term $\<30>$ continues to belong to the appropriate solution space.
Moreover, for $\alpha<\f 14$, it already suffices to control $\<3>$, while for
for $\alpha\in [\frac{1}{4}, \frac{3}{10})$, we would need to use the regularity of  $\<30>$, which belongs to $ C\bigl([0,T]; W^{2\alpha-\f 14+, \infty}(\mathbf{T}^2)\bigr)$. In contrast, for $\alpha > \frac{3}{10}$, $\<30>$ is too rough and must instead be treated as part of the rough component.

Below, we provide a more detailed analysis of the regularity properties and derive the renormalized equation.

\subsection{The renormalized SNLW }
In all regimes $\alpha \in (0, \tfrac{3}{8})$, we work with the same renormalized equation
\begin{equ} 
(\partial_t^2 + 1 - \Delta) u_N = -u_N^3 - 3\sigma_N u_N + \jpb{\nabla}^\alpha \xi_{\le N}
\end{equ}
and we write $u_N=\Psi_N+v_N$. Heuristically, as $N\to\infty$,  the renormalized equation takes the following unified form:
\begin{equ}\label{eqn:renormalized-SNLW}
(\partial_t^2 + 1 - \Delta) u = -u^3 + 3 \cdot \infty \cdot u + \jpb{\nabla}^\alpha \xi,
\end{equ}
where the divergent term $3 \cdot \infty \cdot u$ is interpreted through Wick renormalization. It arises as a counterterm needed to cancel divergences in the product $u^3$ due to the interaction between the rough stochastic component and the nonlinearity.

 The goal of this article is to provide a pathwise solution theory for~\eqref{eqn:renormalized-SNLW}, constructing all required stochastic objects and establishing local well-posedness in an appropriate function space for each regime of $\alpha$.

\subsection{Analysis of the rough/regular decomposition} 
To implement the fixed-point argument to~\eqref{eq:1st-re-vN}, we work in Bourgain-type spaces $\mathcal{X}^{s,b}$, with initial data in $\mathcal{H}^s$. Observe that if the solution space for~\eqref{eq:1st-re-vN} is taken to be $\mathcal{X}^{s,b}$, then the corresponding initial data must lie in $\mathcal{H}^s$. 

As shown later in Corollary~\ref{cor-object-3}, for any $\epsilon > 0$, we have
 \[
\<30> \in C([0,T]; W^{s_\alpha - \epsilon, \infty}_x( \mathbf{T}^2))
\] almost surely, where the exponent $s_\alpha$ is defined in~\eqref{eq:s-alpha}. In particular, when $\alpha < \tfrac{1}{4}$, we have $s_\alpha > \tfrac{1}{4}$. For $\alpha > \tfrac{1}{4}$, we have $s_\alpha \ge 2\alpha - \tfrac{1}{4}$ if and only if $\alpha < \tfrac{3}{10}$. 

For larger values of $\alpha$, the object $\<30> = \mathcal{J}(\<3>)$ no longer belongs to the solution space and must instead be incorporated into $\<1> $--- the rougher component of the solution decomposition--prompting a second-order decomposition:
\[
u_N =( \<1>_N + \<30>_N )+ v_N,
\]
and a different remainder equation given below in ~\eqref{eq:2nd-re-vN}.

\subsection{Acknowledgments}
Li acknowledges support by Swiss National Science Foundation project mint 10000849 and  UKRI  grant
EP/V026100/1, Ren acknowledges a Doris Chen Mobility Award from Imperial.  The authors benefited from the support from Mathematics of Physics NCCR.

\section{Main results}\label{sec:main-results}
Our main result is a pathwise local well-posedness theorem for (\ref{eqn:renormalized-SNLW}). To explain this further, we recall that the strategy is to isolate the singular stochastic component and decompose the solution of the renormalized smooth equation as:
$u_N=\Psi_N+v_N$,
the term $\Psi_N$ built from $\<1>_N$ and possibly higher-order terms, depending on the roughness index $\alpha$. We consider two regimes:

\begin{itemize}
\item[ (1)] High regularity regime: $\alpha \in (0, \tfrac{3}{10})$.
We define $\Psi_N := \<1>_N$ and $v$ solves (\ref{eq:1st-re-vN})

\item [(2)] Low regularity regime: $\alpha \in [\tfrac{3}{10}, \tfrac{3}{8})$.
The decomposition is refined to
$\Psi_N = \<1>_N + \<30>_N$ and 
\begin{equ}\label{eq:2nd-re-vN}
\begin{cases}
(\partial_t^2+1 -\Delta) v_N
 &= \;- 3\<1>_N( \<30>_N)^2 - 3\<1>_N (v_N)^2+ 3\<2>_N\<30>_N-3\<2>_N v_N + 6\<1>_N\<30>_N v_N\\
&\;\quad +(\<30>_N)^3+3\<30>_N v^2_N -3(\<30>_N)^2  v_N-(v_N)^3 \\
(v_N, \partial_t v_N)|_{t=0}& =\; (u_0,u_1)
\end{cases}
\end{equ}

\end{itemize}
We show that all stochastic objects converge in suitable distributional spaces as $N \to \infty$, and that the limit remainder equations are locally well-posed. This yields a rigorous solution theory for the original renormalized equation.

\subsection{Main theorem}
\begin{theorem*}[Pathwise Local Well-Posedness]
Let $\alpha \in (0, \tfrac{3}{8})$. Then for every initial data $(u_0, u_1)$ in a suitable Sobolev space $\mathcal{H}$, the renormalized equation~\eqref{eqn:renormalized-SNLW} admits a unique local-in-time solution of the form:
$u=\Psi+v$
where $\Psi$ is a rough stochastic distribution constructed from the linear stochastic wave equation, and $v$ is a more regular correction in a space $\X$ satisfying a nonlinear remainder equation. The precise form of $\Psi$, the function spaces involved, and the well-posedness thresholds depend on the value of $\alpha$:

\begin{itemize}
\item For $\alpha\in (0, \f  3 {10})$,  we take $\Psi=\<1>$, and $v$  solves(\ref{eq:1st-v-});
\item For $\alpha \in ( \f  3 {10}, \f 3 8)$, we take $\Psi=\<1>+\<30>$,   and $v$  solves (\ref{eq:2nd-v-}).
\end{itemize}

For $\alpha \in (0, \tfrac{1}{4})$, we obtain $v$ in  $ \X^{\f 1 4 + \epsilon, \f 1 2 + \delta}_T$ for initial data in $\mathcal{H}^{\f 1 4 + \epsilon}$ for some $0<\epsilon,\delta \ll 1$.
For $\alpha \in [\tfrac{1}{4}, \tfrac{3}{8})$, the critical regularity depends on $\alpha$ and we have $v\in\mathcal{X} = \X^{2\alpha - \f 1 4 + \epsilon, \f 1 2 + \delta}_T$ for some other $0<\epsilon,\delta \ll 1$.

In each case, the solution map is continuous with respect to the enhanced data (initial conditions and renormalized stochastic symbols).
\end{theorem*}
\begin{remark}
By the uniqueness we mean that if $\Psi_i+v_i$ are two solutions then $\Psi_1-\Psi_2 \in \X$, and that the remainder equations are well posed.
\end{remark}

We also establish the convergence of the Wick symbols and their derivatives, the regularity of associated random operators, and the continuity of the solution map in the corresponding enhanced topologies.

\subsection{Limit theorems}
For any real numbers $s$ and $b$,  we denote by $\X^{s,b}_T$ the restricted Bourgain space, over a finite time interval 
$[0,T]$;  see Section~\ref{sec:preliminaries} for definitions and conventions.
We have already some understanding of the two key stochastic objects $\<1>,\<2>,$ that appear in the renormalized dynamics. In fact, according to  \cite{OO21}, the following holds if $\alpha \in (0, \f 38)$ and $T>0$: for any $\epsilon>0$, 
  $$(\<1>,\<2>):=\lim_{N\to \infty}(\<1>_N,\<2>_N)$$
exists in the product space
 $$C([0,T];W^{-\alpha-\epsilon,\infty}_x) \times C([0,T];W^{-2\alpha-\epsilon,\infty}_x).$$

In particular, for $\alpha \in (0, \frac{1}{8})$, we expect $v$ to lie in $\X_T^{\frac{1}{4}+, \frac{1}{2}+}$, so the product $\<2>\,v$ is well defined as a distribution. In this regime, the mollified products converge, and the limit is independent of the mollification. As a guiding principle, this is possible whenever the sum of regularities is positive.

For rougher noise, where the product $\<2>\, v$ cannot be defined directly due to insufficient regularity, we instead apply the Duhamel operator $\mathcal{J}$ to the regularized product before taking the limit of the mollification. Rather than forming the product first—which would be ill posed—we take into account the smoothing effect of the Duhamel operator by interpreting it as a random operator on the solution space. Namely, we define $\fJ^{\<2>}v := \mathcal{J}(\<2>\,v)$ and show that it is well defined as a bounded linear operator acting on the space of solutions.

The space of bounded operators on $\X_T^{s,b}$—equipped with the norm introduced in (\ref{ro-norm}) in Section~\ref{sec:local-wellposedness}—is denoted by $\mathcal{L}_T^{s, b}$. This operator norm allows us to gain a time factor, which will be crucial in proving that the solution map is a contraction. In the notation of (\ref{ro-norm}), the symbol $\mathcal{L}_T^{s,s, b}$ is used to handle the case when the domain and target spaces differ.

\begin{propositionA}\label{propb}
\begin{enumerate}
\item [(i)] Let  $\alpha \in [\f 1 8,\f 1 4)$. There exists $\epsilon>0$ and $ \delta>0$ such that for any real number $T>0$,  the  operator 
$$\fJ^{\<2>}\in \L^{\frac{1}{4}+\epsilon, \frac{1}{4}+\epsilon,\frac{1}{2}+\delta}_T $$ almost surely.
\item [(ii)] 
Let  $  \alpha \in [\f 1 4, \f 3 8)$. There exists some $\epsilon>0, \delta>0$ such that for any real number $T>0$,
$$\fJ^{\<2>}\in \L^{2\alpha - \f  14+\epsilon, 2\alpha - \f  14+\epsilon,\frac{1}{2}+\delta}_T$$ almost surely.
\end{enumerate}
\end{propositionA}

\subsection{Limit remainder equations}

We now formulate the renormalized remainder equations, arising from the limits of the regularized system. The form of these equations depends on the noise regularity parameter 
$\alpha$

Let $\alpha\in (0, \f 3 {10})$.  We already know that the limit
 $(\<1>,\<2>):=\lim_{N\to \infty}(\<1>_N,\<2>_N)$
exists in $$C([0,T];W^{-\alpha-\epsilon,\infty}_x) \times C([0,T];W^{-2\alpha-\epsilon,\infty}_x).$$
 With these two terms, we obtain the first-order remainder equation:
\begin{equ}\label{eq:1st-v-}
(\partial_t^2 + 1 - \Delta) v = \;-v^3-3\<1>v^2-3\<2> \,v-\<3>.
\end{equ}
For $\alpha\in [\f 3 {10}, \f 38)$, Proposition~\ref{propb} ensures the existence of the stochastic objects, and we formulate the second-order equation:
\begin{equs}\label{eq:2nd-v-}
(\partial_t^2+1 -\Delta) v
 &= -3\<1>( \<30>)^2 - 3\<1> v^2+ 3\<2>\<30>-3\<2> v+6\<1>\<30>v\\
&\quad+(\<30>)^3+3\<30> v^2 -3(\<30>)^2  v-v^3.
\end{equs}

\subsubsection{Local well-posedness}
The proof of the following well-posedness theorem is given by Proposition \ref{prop:LWP-18} to Proposition \ref{prop:LWP-38}.

\begin{theorem*}[Local Well-Posedness]\label{thm:well-posed}
Let $\alpha \in (0, \f 3 {8})$. Then:
\begin{enumerate}
\item [(1)] For $\alpha \in (0, \f 3 {10})$, for any initial condition in $\H$, the first-order remainder equation~\eqref{eq:1st-v-} admits a unique local solution $v\in \X$, depending on the subrange of $\alpha$:
\begin{itemize}
\item [(1a)]If $ \alpha \in (0,  \f 1 4)$, then for all $\epsilon>\epsilon_0$ and $\delta>\delta_0$,  we can take
\begin{equ}
 \qquad \X= \X^{\f 1 4+\epsilon,\f 12 +\delta}_T \subset C([0,T];H^{\f 1 4+\epsilon}_x), \qquad  \H= \mathcal{H}^{\f 1 4+\epsilon};
\end{equ}
where $\epsilon_0$ and $\delta_0$ are some positive numbers.
\item [(1b)] If $ \alpha \in (\f 14,  \f 3 {10})$, then for all $\epsilon>\epsilon_0$ and $\delta>\delta_0$,  
 \begin{equ}
\X= \X^{2\alpha - \f 14 +\epsilon,\f 12 +\delta}_T \subset C([0,T];H^{2\alpha - \f  14 +\epsilon}_x), \qquad \mathcal{H} = \mathcal{H}^{2\alpha - \f 14 +\epsilon}, 
\end{equ} 
where $\epsilon_0$ and $\delta_0$ are some positive numbers.
\end{itemize}
\item [(2)] For $\alpha\in (\f 3 {10},  \f 3 8)$ and for any initial data in $\H$,  the remainder equation  \eqref{eq:2nd-v-} admits a unique local solution in $\X$ for all $\epsilon>\epsilon_0$ and $\delta>\delta_0$,  where
\begin{equ}
 \X= \X^{2\alpha - \f 14 +\epsilon,\f 12 +\delta}_T \subset C([0,T];H^{2\alpha - \f  14 +\epsilon}_x), \qquad \mathcal{H} = \mathcal{H}^{2\alpha - \f 14 +\epsilon}.
\end{equ} 
and $\epsilon_0$ and $\delta_0$ are some positive numbers.
\end{enumerate}
In each case, the solution exists up to a (random) stopping time $T>0$.
\end{theorem*}

\subsubsection{Continuous dependence on enhanced data}
We now describe the continuous dependence of the solution map on enhanced data, which includes both the initial conditions and the renormalized stochastic terms.

Let 
\begin{equ}\label{Theta1}
\begin{aligned}
\Theta &= (u_0,u_1,\<1>,\<2>,\<3>), \qquad \qquad  \hbox{if }\alpha \in (0, \f 1 8);\\
\Theta &= (u_0,u_1,\<1>,\fJ^{\<2>},\<3>), \qquad \qquad \hbox{if } \alpha \in [\f 1 8,  \f 3 {10}),
\end{aligned} \end{equ}
for (\ref{eq:1st-v-}), and
 \begin{equ}\label{Theta2}
\Theta = (u_0,u_1,\<1>,\fJ^{\<2>},\<30>,\<31>,\<320>,\<70>), 
\end{equ}
 for (\ref{eq:2nd-v-}) and $\alpha \in [\f 3 {10},  \f 38)$.

\begin{theorem*}[Continuity of the Solution Map]\label{thm:cty}
Let $ \alpha \in (0, \f 3 {8})$ and define the enhanced parameter set $\Theta$ as above. Then
there exists some $\epsilon_0,\delta_0>$ such that for all $\epsilon>\epsilon_0$ and $\delta>\delta_0$, the map
 \begin{equ}
\Theta \in  \mathcal{Y}^{\alpha,\epsilon,\delta}_T  \to v\in  \X, 
\end{equ}
is continuous almost surely, where the  function space $\mathcal{Y}^{\alpha,\epsilon,\delta}_T $ depends on the range of 
$\alpha$. The full break down is as follows:

\begin{enumerate}
\item [(1)] For  $ \alpha \in (0, \f 1 8)$,  $\X= \X^{\f 1 4+\epsilon,\f 12 +\delta}_T$ and
 $$ \mathcal{Y}^{\alpha,\epsilon,\delta}_T = \mathcal{H}^{\f 1 4+\epsilon} \times C([0,T];W^{-\alpha-\epsilon,\infty}_x) \times C([0,T];W^{-2\alpha-\epsilon,\infty}_x) \times \X^{-3\alpha-\epsilon,\delta}_T$$

\item [(2)] For  $  \alpha \in [\f 1 8,  \f 1 4)$,  $\X= \X^{\f 1 4+\epsilon,\f 12 +\delta}_T$ and $$ \mathcal{Y}^{\alpha,\epsilon,\delta}_T  = \mathcal{H}^{\f 1 4+\epsilon} \times C([0,T];W^{-\alpha-\epsilon,\infty}_x) \times \L^{\frac{1}{4}+\epsilon, \frac{1}{4}+\epsilon,\frac{1}{2}+\delta}_T \times \X^{-3\alpha-\epsilon,\f 12+\delta}_T.
$$

\item [(3)]For $\alpha \in [\f 1 4, \f 3 {10})$, $\X= \X^{2\alpha - \f 14 +\epsilon,\f 12 +\delta}_T $  and
 \begin{equ}
\mathcal{Y}^{s,\alpha,\epsilon,\delta}_T = \mathcal{H}^s \times C([0,T];W^{-\alpha-\f \epsilon 2,\infty}_x) \times \L^{2\alpha - \f  14+\epsilon, 2\alpha - \f  14+\epsilon,\frac{1}{2}+\delta}_T \times \X^{\f 54-3\alpha-\epsilon,\f 12+\delta}_T,
\end{equ}
\item [(4)] For $ \alpha \in [\f 3{10}, \f 38)$,  $\X=\X^{2\alpha - \f 14 +\epsilon,\f 12 +\delta}_T$
and 
\begin{equs}
\mathcal{Y}^{s,\alpha,\epsilon,\delta}_T = &\mathcal{H}^s \times C([0,T];W^{-\alpha-\f \epsilon 2,\infty}_x) \times \L^{2\alpha - \f  14+\epsilon, 2\alpha - \f  14+\epsilon,\frac{1}{2}+\delta}_T \times \X^{\f 54-3\alpha-\epsilon,\f 12+\delta}_T \\
&\quad \times C([0,T];W^{-\alpha- \epsilon,\infty}_x) \times  \X^{\f 56-\alpha-\epsilon,\f 12+\delta}_T.
\end{equs}
\end{enumerate}
\end{theorem*}

\section{The renormalized model and regularity analysis} \label{sec:deduction}
In this section we derive the remainder term $v_N$ for various ranges of the roughness parameter $\alpha > 0$, along with the corresponding equations for the regular part of the solution. By classical theory for deterministic wave equations, for each $N$ the regular residual equations admits a unique solution in the space $C_tH_x^s([0,T]\times \T^2)$ for $s > \frac{1}{4}$. 
Under the regularity conditions
\[
-2\alpha - \epsilon + \frac{1}{4} > 0, \quad -\alpha-\epsilon +\f 14 >0
\]  
the nonlinear terms  
\[
3\<1>_N (v_N)^2 \quad \text{and} \quad 3\<2>_N v_N
\]  
are expected to be deterministically well-defined. Consequently,  the existence and uniqueness of solutions follow from a fixed-point arguments in the corresponding function space. These conditions imply $\alpha \in (0, \frac{1}{8})$.

In regimes $\alpha\ge  \f 18$, however, this decomposition becomes insufficient.
 At first glance, one may expect that a more refined expansion of the stochastic nonlinear wave equation is needed to accommodate the loss of regularity. 
Yet, for $\alpha\in [\f 18, \f 3{10})$, a key mechanism of dispersive PDEs comes into play—namely, multilinear smoothing effects, which arise from cancellation in frequency interactions within multilinear wave propagations.  In more concrete terms, these estimates will be obtained by trilinear estimates to be given later.
  
To make this precise, consider the integral formulation of the remainder equation~\eqref{eq:1st-re-vN}:  
\begin{equ}\label{vn-1}
v_N(t) = \cos(t\jpb{\nabla})u_0  +  \frac{\sin(t\jpb{\nabla})}{\jpb{\nabla}}u_1 - 3\cJ(v_N^3) - 3\cJ(\<1>_N v_N^2) - 3\cJ(\<2>_N v_N) -\cJ(\<3>_N),
\end{equ}
where $\mathcal{J} = (\partial_t^2 + 1 - \Delta)^{-1}$ denotes the Duhamel operator associated with the wave equation. To ensure that the right-hand side lies in a space of spatial regularity at least $\tfrac{1}{4} + \epsilon$ for sufficiently small $\epsilon > 0$, each of the following terms must possess regularity exceeding this threshold:
\[
\cJ(\<3>_N), \quad \cJ(\<2>_N v_N), \quad \cJ(\<1>_N v_N^2).
\]  
  We now analyze these terms below. 
  
 \subsection{The regularity of the trilinear term for  $\alpha\in [\f 14, \f 5{12})$}
For $\alpha$ in this range, we extend the multilinear smoothing strategy to the term $\cJ(\<3>_N)$ and refine the regularity estimate for the stochastic object
\begin{equ}
\<30>_N := \mathcal{J}(\<3>_N).
\end{equ}  originally obtained in Case 1. Specifically, we define
\begin{equ}\label{eq:s-alpha}
s_\alpha =
\begin{cases}
&1-2\alpha, \quad 0<\alpha <\f 14,\\
&\f 54-3\alpha,\quad \f14 \le \alpha<\f 5{12},
\end{cases}
\end{equ}
and show in Corollary~\ref{cor-object-3} that for any $\epsilon > 0$,
 \[
\<30> \in { C\bigl([0,T]; W^{s_\alpha - \epsilon, \infty}( \mathbf{T}^2)\bigr)}
\] 
almost surely.

\begin{remark}
The refinement for $\<30>_N $ yields two key insights. First, it improves the previously known regularity $1 - 3\alpha - \epsilon$ for $\<3>_N$, revealing a phase transition at $\alpha = \tfrac{1}{4}$. That is, for $\alpha < \tfrac{1}{4}$, the regularity gain depends on $\alpha$, whereas for $\alpha \ge \tfrac{1}{4}$, the gain plateaus at a fixed $\tfrac{1}{4}$. This behaviour mirrors the quartic case discussed in~\cite{OO21}.

Secondly, this phase transition has an analogue in the action of the random operator $\mathcal{J}^{\<2>_N}$, which we describe following its introduction in the next section.
\end{remark}

\subsection{Range 1: $\alpha \in (0,\f 14)$}
\begin{itemize}
\item [(i)] The term $\cJ(\<3>_N)$ gains one derivative in spatial regularity via the operator $\cJ$. Since $\<3>_N$ has spatial regularity $-3\alpha - \epsilon$ almost surely, the term, for any $\epsilon>0$, the term $\cJ(\<3>_N)$ lies in 
$C([0,T];W^{1 - 3\alpha - \epsilon, \infty}_x)$, which exceeds $\frac{1}{4}$ for $\alpha<\f 14$.  
\item [(ii)] The term $\cJ(\<2>_N v_N)$ initially appears more problematic. The product $\<2>_N v_N$ inherits the lower regularity of $\<2>_N$, which is $-2\alpha - \epsilon$.
Thus, $\cJ(\<2>_N v_N)$ naively has regularity $1 - 2\alpha - \epsilon$, suggesting admissibility for $\alpha < \tfrac{3}{8}$. However, high-high to low frequency interactions in bilinear estimates (see Lemma 8.1) restrict this threshold to $\alpha < \tfrac{1}{8}$. To overcome this, we reinterpret the term as a multilinear stochastic operator:
 \[
   \cJ^{\<2>_N} \colon v_N \mapsto \cJ(\<2>_N v_N).
 \]  
 Due to cancellation effects between $\<2>_N$ and $v_N$, this operator becomes a bounded random linear map from $\mathcal{X}^{\tfrac{1}{4}+\epsilon, \tfrac{1}{2}+\delta}_T$ to itself, for $\alpha \in (0, \tfrac{1}{4})$. Here, $\mathcal{X}^{s,b}_T$ denotes a Bourgain space adapted to the wave equation (see Section 3).

\item[(iii)] The term $\mathcal{J}(\<1>_N v_N^2)$ is less restricted. Standard trilinear estimates suffice to guarantee regularity above $\tfrac{1}{4}$ when $\alpha < \tfrac{1}{4}$.
 \end{itemize}

In summary, by employing the first-order expansion $u_N = \<1>_N + v_N$ and leveraging multilinear smoothing, we establish local well-posedness in the extended regime $\alpha \in (0, \tfrac{1}{4})$. This framework balances noise roughness against the smoothing effects of dispersive interactions, demonstrating how cancellations in frequency interactions can enlarge the admissible range of $\alpha$.

\subsection{Range 2: $\alpha \in (\f 14, \f 3{10})$}
 For $\alpha \ge \tfrac{1}{4}$, we show (see Lemma~9.1) that
\[
   \cJ^{\<2>_N} \colon \mathcal{X}^{2\alpha - \frac{1}{4} + \epsilon, \frac{1}{2} + \delta}_T \to \mathcal{X}^{2\alpha - \frac{1}{4} + \epsilon, \frac{1}{2} + \delta}_T,
\]  
is bounded as a linear operator almost surely, where $\mathcal{X}^{s, b}_T$ denotes a Bourgain-type space adapted to the wave equation. This matches the expected regularity reduction: while the product $\<2>_N v_N$ formally requires $v_N$ to lie in $\mathcal{X}^{2\alpha - \frac{1}{4} + \epsilon, \frac{1}{2} + \delta}_T$, multilinear smoothing relaxes this requirement to $2\alpha - \tfrac{1}{4} + \epsilon$, reflecting a fixed $\tfrac{1}{4}$ gain instead of $\alpha$-dependent scaling.

To ensure the random operator operates in the solution space, we now require $v_N \in \mathcal{X}^{2\alpha - \tfrac{1}{4} + \epsilon, \f 12+\delta}_T$. 
For the stochastic term  $\<30>_N$ to be compatible with this regularity
the condition  
\[
\frac{5}{4} - 3\alpha > 2\alpha - \frac{1}{4}
\]  
must hold, leading to the constraint $\alpha < \frac{3}{10}$.  A rigorous formulation and proof of local well-posedness in this regime is provided in Proposition~\ref{prop:LWP-310}.

\subsection{Range 3.  $\alpha \in [ \frac{3}{10}, \f 3 8)$}
The combined regularity requirements for $ \cJ^{\<2>_N}$ and for $\<30>$ in the rougher noise regimes $\alpha \ge \frac{3}{10}$
necessitate a refinement of the solution framework.  In particular, we introduce a second-order expansion that isolates the lower-regularity component
 $\<30>$ from the residual term $v_N$. 

Importantly, no further renormalisation beyond Wick ordering is required. We return to the renormalised equation:
\begin{equ} 
\begin{cases}
(\partial_t^2 + 1 - \Delta) u_N = -u_N^3 - 3\sigma_N u_N + \jpb{\nabla}^\alpha \xi_{\le N}, \\
(u_N, \partial_t u_N)\big|_{t=0} = (u_0, u_1).
\end{cases}
\end{equ}
which retains its counter-term $3\sigma_N u_N$.  

We now decompose the solution as:
\[
u_N = \<1>_N - \<30>_N + w_N,
\]  
where the remainder $w_N$ absorbs the residual nonlinear interactions.

The remainder term $w_N$ then satisfies:  
\[
(\partial_t^2 + 1 - \Delta) w_N = -(\<1>_N - \<30>_N + w_N)^3 + \<3>_N.
\]  
Expanding the cubic term yields 
\begin{equs}
-(\<1>_N - \<30>_N + w_N)^3 + \<3>_N &= - \<3>_N - 3\<1>_N \<30>^2_N - 3\<1>_Nw^2_N + 3\<2>_N\<30>_N-3\<2>_N w_N + 6\<1>_N\<30>_Nw_N\\
&\quad +\<30>^3_N+3\<30>_Nw^2_N -3\<30>^2_Nw_N-w^3_N +\<3>_N\\
&\quad -3\sigma_N(\<1>_N-\<30>_N+w_N)
\end{equs}
which simplifies the equation for $w_N$ to:
\begin{equ}
\begin{cases}
(\partial_t^2+1 -\Delta) w_N
 &= \;- 3\<1>_N \<30>^2_N - 3\<1>_Nw^2_N + 3\<2>_N\<30>_N-3\<2>_N w_N + 6\<1>_N\<30>_Nw_N\\
&\quad +\<30>^3_N+3\<30>_Nw^2_N -3\<30>^2_Nw_N-w^3_N \\
(w_N, \partial_t w_N)|_{t=0}& =\; (u_0,u_1)
\end{cases}
\end{equ}

The corresponding mild solution formulation reads:
\begin{equs}
w_N(t) &= \cos(t\jpb{\nabla})u_0  +  \frac{\sin(t\jpb{\nabla})}{\jpb{\nabla}}u_1 - 3\cJ(\<1>_Nw^2_N)+ 3\CJ(\<30>_Nw^2_N)-3\CJ(\<30>^2_Nw_N)-\cJ(w^3_N)\\
&\quad +3\cJ(\<2>_N\<30>_N)-3\cJ(\<1>_N \<30>^2_N)+\cJ(\<30>^3_N)+6\cJ(\<1>_N\<30>_Nw_N)-3\cJ(\<2>_N w_N).
\end{equs}

To ensure that the term $\mathcal{J}(\<2>_N w_N)$—ill-defined in the classical sense—is meaningful, we once again reinterpret it as a random operator acting on $\mathcal{X}^{2\alpha - \tfrac{1}{4} + \epsilon, \tfrac{1}{2} + \delta}_T$. This necessitates that $w_N$, along with every term on the right-hand side, belongs to the same space.

The terms
$$\mathbf S(t)(u_0,u_1) - 3\cJ(\<1>_Nw^2_N)+ 3\CJ(\<30>_Nw^2_N)-3\CJ(\<30>^2_Nw_N)-\cJ(w^3_N)$$
are shown to lie in this space via trilinear estimates developed in Section~\ref{sec:deterministic-estimates}. The remaining terms require more delicate treatment.

\begin{itemize}
\item [(1)] For $\mathcal{J}(\<30>_N^3)$, since $\<30>_N$ is less regular than $w_N$, direct trilinear control in $\mathcal{X}^{2\alpha - \tfrac{1}{4} + \epsilon, \tfrac{1}{2} + \delta}_T$ is unavailable. Instead, we define this as a stochastic symbol:
\begin{equ}
\<90>_N := \cJ(\<30>^3_N).
\end{equ}
{\bf Conjecture:}  We conjecture that, with full exploitation of multilinear smoothing effects, one can show
 \[
\<90>_N \in C_tW^{\frac{5}{2} - 3\alpha - \epsilon, \infty}_x([0, T] \times \mathbf{T}^2)
\]  
almost surely.

Even without such effects, the Duhamel operator provides a one-derivative gain. We prove in Proposition~\ref{prop:reg-cubic-cubic-cubic} that
$$\<90>_N\in\X^{\frac{9}{4} - 3\alpha - \epsilon}_T$$ almost surely, for any $\epsilon > 0$, which suffices for fixed-point theory when
\[
\frac{9}{4} - 3\alpha - \epsilon > 2\alpha - \frac{1}{4},
\]  
i.e., $\alpha < \tfrac{1}{2}$ for small $\epsilon > 0$.
\item [(2)] For the term $\cJ(\<1>_N\<30>_N w_N)$, which imposes stronger regularity demands on $w_N$, we again reinterpret it as a random operator:
\[
\cJ^{\<31>_N} \colon v_N \mapsto \cJ(\<31>_N v_N),
\]  
with  $\<31>_N := \<1>_N\<30>_N$ defined as a stochastic symbol, facilitated with stochastic analytic techniques.

We show that $\<31>_N$ inherits the worst regularity from $\<1>_N$, namely
\[
\<31>_N \in C_tW^{-\alpha -, \infty}_x([0, T] \times \mathbf{T}^2)
\]  
almost surely.

\item [(3)]Similarly, the term $\mathcal{J}(\<2>_N \<30>_N)$ is interpreted as a quintic stochastic symbol:

\[\<320>_N := \cJ(\<2>_N \<30>_N).
\]  
For $\alpha < \frac{3}{8}$, we establish via multilinear smoothing that
 $$\<320>_N\in \X^{\frac{5}{6} - \alpha -, \f 12+}$$ 
 almost surely  for any $\epsilon>0$, which is higher than $2\alpha - \frac{1}{4}$ and thus sufficient for our fixed-point argument (see Proposition~10.5).

\item [(4)] For  $\cJ(\<1>_N\<30>^2_N)$, which also defies direct trilinear control, we define the septic stochastic symbol
\[
\<70>_N := \cJ(\<1>_N\<30>^2_N).
\]  

Unlike $\<90>_N$  and $\<320>_N$, establishing the regularity of $\<70>_N$ requires partial multilinear smoothing.
We prove that $$\<70>_N\in \X^{1-\alpha-\epsilon, \f 12+\epsilon} \subset \X^{2\alpha-\f 14+\epsilon, \f 12+\epsilon}$$
for $\alpha <\f 5{12}$ and sufficiently small $\epsilon>0$. 
\item [(5)] The terms  $\cJ(\<2>_N w_N)$ can be interpreted as a random operator as in case $\alpha \in [\f 14, \f 3{10})$.
\end{itemize}

With all terms in the mild formulation possessing regularity strictly greater than $2\alpha - \tfrac{1}{4} + \epsilon$, we conclude by applying a deterministic fixed-point argument in the Bourgain space $\mathcal{X}^{2\alpha - \tfrac{1}{4} + \epsilon, \tfrac{1}{2} + \delta}_T$, establishing local existence and uniqueness of solutions.

\section{Preliminaries}\label{sec:preliminaries}

In the rest of the article, we denote the probability space by $(\Omega,\cF,\P)$, and $\xi$ refers to the space-time white noise on the two torus $\T^2$.
If $\eta$ is a tempered distribution on $\T^2$ 
we denote by $\{\hat \eta(n): n\in \Z^2\}$ its set of Fourier coefficients.   
Given any function $m: \Z^2\to \C$, we define the multiplier operator $m(\Delta)$ by its action on $\eta$ via
 $\widehat{m(\Delta)\eta}(n)=m(n)\hat \eta(n)$.
In particular,  setting $m(n)=\jpb{n}= \sqrt{1+|n|^2}$ (the Japanese bracket),we define the noise in the stochastic nonlinear wave equation (SNLW) as
 $$ \widehat{\jpb{\nabla}^\alpha  \xi}(n)= \jpb{n}^\alpha \hat \xi(n).$$
We write $\eta\in \cC^{\beta-}$ to mean that $\eta \in \cC^{\beta-\epsilon}$ for any $\epsilon>0$.
 For example, the sample paths of $\jpb{\nabla}^\alpha  \xi$ almost surely belong to $\cC^{-1-\alpha-}$.

\subsection{Conventions and notations}

Since this paper involves a large number of symbols and notations, we collect the most frequently used ones in this subsection for the reader’s convenience.
The notation $f \lesssim g$  indicates that there exists a constant $c$, independent of the data in $f$ and $g$
such that $f\lesssim c g$.
Similarly, $f \lesssim_{\gamma} g$ means that the proportionality constant may depend on $\gamma$.

\begin{itemize}
\item $\| \cdot \|_{L^p}$ denotes the norm in $L^p(\Omega)$.
\item  $\jpb{n} = \sqrt{1+|n|^2}$, where $|\cdot|$ is the Euclidean norm and $n \in \Z^2$.
\item For any $n_1, \cdots, n_k \in \Z^2$, we introduce the symbol 
\begin{equation}
n_{12\cdots k} := n_1 +n_2 + \cdots +n_k.
\end{equation}
For example, $n_{123} = n_1+n_2+n_3$ and $n_{23}=n_2+n_3$.

\item $\jpb{\nabla}$ denotes $(1-\Delta)^{\f 12}$ where $\Delta$ is the Laplacian.
\item For $s \in \R$ and $ 1 \le p \le \infty$, $W^{s,p}_x$ denotes the $L^p$-based Sobolev space in the spatial variable. If $f \colon T^2 \to \R$, its norm is given by
\begin{equation}
\|f\|_{W^{s,p}_x} := \| \jpb{\nabla}^s f \|_{L^p_x},
\end{equation}
with the standard modification when $p=\infty$.  For $p=2$, we dentoe $H^s := W^{s,2}$.
\item $\mathcal{H}^s=H^{s}_x \times H^{s-1}_x$ denotes the initial function space for the renormalised SNLW wave, where  $s \in \R$ is typically taken to be greater than $\f 14$.
\item $\X^{s,b},\X^{s,b}_I$ are defined in  Subsection \ref{Bourgain-space}.
 They denotes respectively the Bourgain-weave space and restricted Bourgain-wave space. 
 \item  The notation $A_j = \{n \in \mathbb{Z}^2 : 2^{j-1} \le |n| \le 2^{j+1}\}$ denotes the dyadic annulus in frequency space. If $\psi: \mathbb{Z}^2 \to \mathbb{R}$ is a bump function supported in $\{n \in \mathbb{Z}^2 : \tfrac{1}{2} \le |n| \le 2\}$, then $\P_j u$ denotes the Littlewood–Paley projection of a tempered distribution $u$, with the decomposition $u = \sum_{j \in \mathbb{N}_0} \P_j u$. See \S\ref{sec:tri-linear}.
 
\item We keep the notation $\d x$ for the normalized Lebesgue measure on $\T^2$ in the sense that
\begin{equation}
\int_{\T^2} 1 \, \d x =1.
\end{equation}
\item We use the notation $\widehat{f}(\cdot)$ to denote the spatial Fourier transform while use $\widetilde{f}(\cdot,\cdot)$ to denote the space-time Fourier transform.
\end{itemize}

\subsection{Wave Sobolev spaces (Bourgain space) and estimates}\label{Bourgain-space}
The $\X^{s,b}$ space provides a suitable solution space for our equation. In the context of nonlinear Schrödinger equations, it is often referred to as the Bourgain space or hyperbolic Sobolev space \cite{Bou93a}; while for nonlinear wave equations, it is also known as the $H^{s,\delta}$ space or wave-Sobolev space \cite{Klainerman-Machedon-94}. We list below some basic properties of Bourgain spaces that will be frequently used in this paper.

We use the following notation for the Fourier transform on $\T^d$ and $\R \times \T^d$, respectively:
\begin{equs}
\widehat{g}(\cdot,n)&= (\mathcal{F}_{x} g)(\cdot,n) := \int_{\T^d} g(\cdot, x) e^{\i n \cdot x}\,\d x, \\
\tilde{g}(\lambda,n) &= (\mathcal{F}_{t,x} g)(\lambda,n) := \int_{\R \times \T^d} g(t,x) e^{-\i n\cdot x - \i \lambda t} \, \d x \d t.
\end{equs}
We denote by $L^qL^r$ or $L^q_tL_x^r$ the space of functions with finite norm:
$$\|g\|_{q,r}\equiv \|g\|_{L^q_tL_x^r}=\Bigl(\int_{\T^2}\bigl( \int_\R |g(t,x)|^q dt\bigr)^{\f  r q}dx\Bigr)^{\f 1r}.$$

\begin{definition}
Let $s,b \in \R$. The space $\X^{s,b}(\R \times \T^d)$, associated with the wave equation on $\T^d$, is defined as the closure of the Schwartz space
 $\mathcal{S}(\R \times \T^d)$ under the norm
\begin{equation}\label{eq:Xsb-norm}
\|u\|_{\X^{s,b}} :=\Bigl \| \jpb{n}^s \jpb{|\lambda| - \jpb{n}}^b (\cF_{t,x} u)(\lambda,n)\bigr\|_{L^2_\lambda \ell^2_n(\R \times \Z^d)}.
\end{equation}
We abbreviate this space as $\X^{s,b}$. 

For any interval $I \subset \R$, we define the restricted norm by
\begin{equation}\label{Xsb-norm-T}
\|u\|_{\X^{s,b}(I)} := \inf_v \{ \|v\|_{\X^{s,b}} \colon v|_I = u \}.
\end{equation}
We denote by $\X^{s,b}_I$  the corresponding closure.  In the special case $I= [0,T]$, we write $\X^{s,b}_T$.
\end{definition}

The spaces $\X^{s,b}$ form Banach spaces and satisfy the basic inclusion property
$$\X^{s,b} \subset \X^{\tilde s,\tilde b}, \qquad \hbox{whenever $s> \tilde s$ and $b> \tilde b$}$$
together with the natural duality inherited from the $L^2$ setting. In particular, for any 
$s,b\in \R$,
$$\X^{s,b}=(\X^{-s, -b})^*$$

The $\X^{s,b}$ space are naturally related to the spatial Sobolev spaces $H^s_x$ on the torus $\T^2$. In fact, for any time interval $I\subset\R$,
 $$\X^{s,0}(I)=L^2_t(I; H^s_x(\T^2)),$$
where the norm is given by 
$$\|u\|_{L^2_t(I; H^s_x(\T^2))}
=\bigl(\sum_n \int_I  \bigl \| \jpb{n}^s  (\cF_{x} u(t, \cdot))(n)\bigr\|\bigr)^{\f 12}.$$
 As an immediate consequence, for every $b\ge 0$,   
$$\|u\|_{L^2_t(I; H^s_x(\T^2))} \le \|u\|_{\X^{s,b}}.$$

We next define the Duhamel operator $\J$:
\begin{equ}\label{Duhamel}
\J(F)(t,x)=\int_0^t \frac{\sin\pa{(t-s)\jpb{\nabla}}}{\jpb{\nabla}} F(s,x)\, \d s.
\end{equ}

\begin{lemma}\cite[Lemma 2.5]{OWZ22} \label{lem:nlXsb}
Let $s \in \R$, $b'\in (-\f 12, 0]$ and $b \in [0, b'+1]$. Let $T \in (0, 1]$.  For $F\in \X^{s-1,b'}_T$, one has
	\begin{equs}\label{nl-Xsb}
	 \norm{\J(F)}_{\X^{s,b}_T} \lesssim T^{1-b+b'} \norm{F}_{\X^{s-1,b'}_T}.
	\end{equs}
\end{lemma}

\begin{definition}
Let $s\in [0,1)$. A pair of numbers $(q,r)$, where $q \in [2,\infty]$ and $r\in [2, \infty)$,  is said to be (wave) $s$-admissible if they satisfy the following two conditions:
\begin{itemize}
\item (scaling)
\begin{equation}\label{eqn:scaling-condition}
\frac{1}{q} + \frac{2}{r} = 1 - s ,
\end{equation}
\item  (wave admissibility)
\begin{equation}\label{eqn:wave-admissiblity-condition}
\frac{2}{q} + \frac{1}{r}\leq \frac{1}{2}.
\end{equation}
\end{itemize}
\end{definition}

\begin{remark}\label{remark-admissible-pair}
 Note that $(q,r)$ is wave $s-$admissible pair if and only if $\f 1r \ge \f 12 -\frac{2}{3}s$. For any $r \ge 2$,  $(q,r)$ is $s$-admissible if $\f 1q\le \f s 3$.
\end{remark}

The following useful estimates is taken from \cite[Theorem 2.6, pp.77]{Tao06}, c.f. \cite{Ginibre-Velo-95, Shatah-Struwe98, LS95},
\begin{lemma}[Strichartz estimates for wave]\label{Stri}
Let $u$ denote the solution to the deterministic nonlinear wave equation on $\T^2$:
\begin{equation}\label{eq:nlw}
\begin{cases}
(\partial_{tt}+(1-\Delta))u = F(t,x), \\
u(0)=u_0, \, \partial_t u(0) =u_1
\end{cases}.
\end{equation}
Let $ s \geq 0$, $(q,r)$ be wave $s$-admissible, and $(\tilde q, \tilde r)$ be wave $1-s$-admissible.
Then for any $T>0$, we have the following estimates on $[0,T]\times \T^2$:
\begin{align}
\|u\|_{L^q_tL^r_x} + \|u\|_{C^0_tH^s_x} + \| \partial_t u\|_{C^0_tH^{s-1}_x} 
\lesssim \pa{\|(u_0,u_1)\|_{\H^s}+\|F\|_{L^{\tilde{q}}_tL^{\tilde{r}}_x}}
\end{align}
where the proportional constants depend only on  ${q,r,s}$.
\end{lemma}
The Sobolev norms refers to  the homogeneous norm. These estimates are sharp for the homogeneous linear equation.

One special case of the Strichartz estimates for wave equation is the following energy estimate:
\begin{lemma}\cite[Lemma 4.8]{Bri20b} \label{energy-est}
Let $s \ge 0, 1/2 < b < 1$, and $I \subset \R$  be a compact interval, let $t_0 \in I$, and let $u$ the solution to \eqref{eq:nlw}. Then,
\begin{equation}
\|u\|_{\X^{s,b}(I)} \lesssim (1+|I|)^2 (\|(u_0,u_1)\|_{\H^s} + \|F\|_{\X^{s-1,b-1}(I)}).
\end{equation}
\end{lemma}
\begin{remark}
The square in the pre-factor can likely be improved but it is good enough for our application.
\end{remark}

 The following lemmas are taken from  Lemmas 4.2- 4.4 in \cite{Bri20b}, \cite[Lemma 2.6]{Tao06}, \cite[section3]{Erdogan-Tzirakis}.

 \begin{proposition} \label{character-Xsb}
 \begin{enumerate}
 \item { \bf Characterization of $\X^{s,b}$ .} Let  for any $s, b \in \R, \sigma \in \{ \pm 1\}$.  Then the following holds for any complex-valued function $u$ on $\R \times \T^2$:
\begin{equation}\label{eq:Xsb-upper}
\|u\|_{\X^{s,b}} \lesssim \min_{\sigma \in \{ \pm \}} \|\jpb{n}^s \jpb{\lambda}^{b} \cF_{t,x}(u)(\lambda - \sigma \jpb{n}, n)\|_{\ell^2_nL^2_\lambda}.
\end{equation}
The proportional is universal, independent of $u$, $s,b, \sigma$. Moreover, we have the equivalence:
\begin{equation}\label{eq:Xsb-equiv}
\|u\|_{\X^{s,b}} \sim \min_{\substack{u=u_+ + u_- \\ u_+, u_- \in \X^{s,b}}} \max_{\sigma \in \{ \pm \}} \|\jpb{n}^s \jpb{\lambda}^{b} \cF_{t,x}(u_{\pm})(\lambda - \sigma \jpb{n}, n)\|_{\ell^2_nL^2_\lambda}.
\end{equation}
\item {\bf Time-localization lemma. }
Let  $b_1, b_2 \in (-\frac 12, \frac 12)$ with $b_1 \le b_2$.
For every $F \in \X^{s,b}$,  $\psi \in\cS(\R)$,  and  $\tau \in (0, 1]$,  the following inequalities hold: 
\begin{align}
\|\psi\bigl( t/\tau\bigr) F\|_{\X^{s,b_1}} &\lesssim \tau^{b_2 - b_1} \|F\|_{\X^{s,b_2}}, \qquad
\|F\|_{\X^{s,b_1}([0,\tau])} &\lesssim \tau^{b_2 - b_1} \|F\|_{\X^{s,b_2}([0,\tau])}.
\end{align}
The proportionality constants depend only on $s, b_1$, and $ b_2$.
\item Let $s \in \R$, $b\in (-1/2,1/2)$, and $I \subset \R$ an interval. Let $F \in \X^{s,b'}$ and $G \in \X^{s,b}(I)$, we have
\begin{equation}
\|\1_{I} F\|_{\X^{s,b}} \lesssim \| F\|_{\X^{s,b}}, \qquad \|G\|_{\X^{s,b}_I} \sim \|\1_I G\|_{\X^{s,b}}.
\end{equation}
 \end{enumerate}
\end{proposition}

Another useful lemma vital  in application is the linear transference principle of $\X^{s,b}$ space, taken from \cite[Lemma 2.9]{Tao06}.
\begin{lemma}[Linear transference] 
Let $b >1/2$ and $s \in \R$. Assume $Y$ be a Banach space of functions on $\R \times \T^d$ with the property that
\begin{equation}
\|e^{\i t \tau_0} e^{ \pm \i t \jpb{\nabla}} f \|_{Y} \lesssim \|f\|_{H^{s}}, \qquad \forall \; f \in H^s(\T^d),\; \forall \tau_0 \in \R.
\end{equation}
Then we have the embedding
\begin{equation}
\|u \|_{Y} \lesssim \|u\|_{\X^{s,b}}, \qquad \forall \;u \in \X^{s,b}.
\end{equation}
\end{lemma}

To go from frequency space $\X^{s,b}$ to physical space $L^q_tH^{s}_x$, we need the following inhomogeneous Strichartz estimate taken from \cite[Lemma 4.9]{Bri20b}:
\begin{lemma}[Inhomogeneous Strichartz estimate in $\X^{s,b}$]\label{lem:inhomo-Stri}
Let $1/2 < b< 1$ and $s \in \R$. Assume $I \subset \R$ be a compact interval and $F \colon I \times \T^2 \to \R$. Then, we have 
\begin{equ}\label{eq:inhomo-Stri}
\|F\|_{\X^{s-1,b-1}(I)} \lesssim \|F\|_{L^{2b}_tH^{s-1}_x(I \times \T^2)}.
\end{equ}
\end{lemma}

By this lemma and the Strichartz estimates (see Lemma \ref{Stri}), one obtains the embedding of an $L^qL^r$ space in a $\X^{s,b}$.

\begin{lemma}[Transference Principle]\label{lem:transference-principle}
Let $(q,r)$ be a wave $s-$admissible pair where $s\in (0,1)$.  Then, for any $b > \frac{1}{2}$ and any $ u \in \X^{s,b}$, we have
\begin{equation}
\|u\|_{L^q_t L^r_x(\R \times \T^2)} \lesssim \|u\|_{\X^{s,b}(\R \times \T^2)}.
\end{equation}
\end{lemma}

\section{Interpolation and non-linear estimates}\label{sec:deterministic-estimates}

Throughout the section let $T>0$ and $s\in (0,1)$. For brevity, we sometimes omit the subscript $T$ in the norm on $\X^{1-s, 1-b}_T$. 

\subsection{Interpolation  between $L^p$ and $\X^{s,b}$ spaces}	
The following lemma is key to our estimates. Since we do not find a proof in the literature, a proof is included for reader's convenience.
\begin{lemma}[Interpolation Inequality] \label{interpolation-lemma}
Let  $s\in (0,1)$, and $T\subset \R$ an interval. Let  $q>2$ and $r\ge 2$ be a wave admissible pair. 
Then, any $b>\f 12$, and $\theta \in (0,1)$, we have
$$\norm{v}_{L^{q_\theta}_TL^{r_\theta}_x} \lesssim \norm{v}_{\X^{s \theta, b \theta}_T},$$
where $q_\theta$ and $r_\theta$ satisfy the admissible conditions:
\begin{equ}\label{theta}
\frac{1}{q_\theta}=\frac{1- \theta} {2}+\frac{\theta}{q}, \qquad \frac{1}{r_\theta}=\frac{1- \theta} {2}+\frac{\theta}{r}.
\end{equ}
 \end{lemma}
\begin{proof}
 Set $w_0=1$ and $w_1=\jpb{\xi}^s\jpb{|\tau|-\jpb{\xi}}^b$, we interpret the $\X^{s,b}$ spaces as weighted $L^2$ spaces on the Fourier transforms of functions on $[0,T]\times \T^2$. Specifically, a function $u$ is in $ \X^{s,b}$ if and only if its Fourier transform $\hat u\in L^2([0,T]\times \T^2, w_0 dtdx)$. Observe that $\X^{0,0}=L^{2,2}$. 
 For  $b>\f 12$, by Lemma \ref{lem:transference-principle},  we also have $\norm{v}_{L_T^qL_x^r} \lesssim  \norm{v}_{\X_T^{s,b}}$.

According to Theorem 5.5.3 \cite[pp.120]{BL76}, the interpolation spaces for the weighted $L^2$ space, obtained via the complex method with parameter $\theta \in (0,1)$, are given by
$$(\X^{0,0}, \X^{s,b})_{[\theta]}=(\X^{s\theta, b\theta}).$$
Furthermore, also using the complex interpolation method, Theorem 5.1.2 \cite[pp.107]{BL76} determines the interpolation space of $L^2_T(L^2_x)$ and $L^q_T(L^r_x)$,  for vector space valued functions, with parameter $\theta$, as
\begin{equation}
(L^2_T(L^2_x),L^q_T(L^r_x))_{[\theta]}=L^{q_\theta}_x(L^{r_\theta}_x),
\end{equation}
where $q_\theta, r_\theta$ are given by (\ref{theta}).

By the interpolation theorem in \cite[Theorem 4.1.2]{BL76}, the inclusion map is continuous from $\X^{s\theta, b\theta}$ to $L^{q_\theta}_x(L^{r_\theta}_x)$, consequently we have the bound $\norm{v}_{L^{q_\theta}_TL^{r_\theta}_x} \lesssim \norm{v}_{\X^{s \theta, b \theta}_T}$, completing the proof for the statement.
\end{proof}

An immediate application allow us to obtain the following lemma, which outlines the philosophy for obtaining embedding theorem of the form 
$\X^{s-1, b-1}\subset L^qL^r$. The difficulty in applying this lemma is that there is no obvious choice of a number $s'>\f 12$ together with a wave $s'$-admissible pair $(q', r')$, and an interpolation constant $\theta$  that satisfies the restrictions in the lemma.
As a rule of thumb, we choose $\theta$ close $0$ and $\tilde b$ as close to $\f 12$ as possible.

\begin{lemma}\label{lem:Xsb-to-LqLr}
Let $s\in (0,1)$, $b > \frac{1}{2}$. Let  $\theta\in (0,1)$ be a number satisfying that
$$\theta \tilde s\le 1-s, \qquad \theta \tilde b \le 1-b,$$
where $\tilde b>\f 12$  and $\tilde b\in (0,1)$ are any numbers.
Let $(\tilde q,\tilde r)$ be a wave $\tilde s$-admissible pair. Then  we have
$$\|u\|_{\X^{s-1, b-1}}\le \|u\|_{q_\theta ,r_\theta } \le T^{\f 1 {q_\theta}-\f 1q }\|u\|_{q,r}$$
for any $(q,r)$ satisfying the following inequality. 
$$\f 1 q\le \f {1+\theta} 2 -\f \theta {\tilde q}, \qquad \f 1r \le  \f {1+\theta} 2 -\f \theta {\tilde r}.$$
\end{lemma}
	\begin{proof} 	
	Let $\f 1 { q}+\f 1 { q'}=1$ and $\f 1{ r}+\f 1 { r'}=1$. Then,
	\begin{equs}
\|u\|_{\X^{s-1,b-1}}
&= \sup_{\|w\|_{\X^{1-s, 1-b}}=1}  \Bigl | \int_{ [0,T]\times \T^2} u w \, \d x\d t  \Bigr|
 \le  \sup_{\|w\|_{\X^{1-s, 1-b}}=1}
\|u\|_{q,r}  \|w\|_{ q',r'}.
\end{equs}
Let $b'>\f 12$. Since $(\tilde q,\tilde r)$ a wave $s'$-admissible pair, we apply  Lemma \ref{interpolation-lemma} to obtain
$ \|w\|_{\tilde q_\theta, \tilde r_\theta}  \le \|w\|_{\X^{s'\theta ,b'\theta}}$. 
Consequently, if furthermore $q'\le q_\theta$ and  $r'\le r_\theta$,   $\tilde s\theta \le 1-s$ and $\tilde b\theta \le 1-b$ then
$$\|w\|_{q',r'}\le \|w\|_{\tilde q_\theta, \tilde r_\theta}  \le \|w\|_{\X^{\tilde s\theta ,\tilde b\theta}}\le \|w\|_{\X^{1-s, 1-b}}.$$

It follows that
\begin{equs}
\|u\|_{\X^{s-1,b-1}}\le \sup_{\|w\|_{\X^{1-s, 1-b}}=1}
\|u\|_{q,r}  \|w\|_{\X^{1-s, 1-b}}
 \le \|u\|_{ q,r},
 \end{equs}
as required.
	\end{proof}

\subsection{Tri-linear estimates}\label{sec:tri-linear}
We establish some trilinear estimate which is indispensable for solving the remainder equation in the appropriate function spaces by Banach fixed-point argument.

Let $\psi: \mathbb{Z}^2 \to \mathbb{R}$ be a smooth bump function supported in the annulus $\{n \in \mathbb{Z}^2 : \tfrac{1}{2} \le |n| \le 2\}$ of the frequency plane. In fact which is identically equal to $1$ in 
$\{n \in \mathbb{Z}^2 : \tfrac{1}{2} \le |n| \le 1\}$. For each $j \in \mathbb{N}$, define the dyadic frequency cutoff
$$
\psi_j(n) := \psi(n\, 2^{-j}) \qquad \text{and} \qquad \psi_0(n) := 1 - \sum_{j \in \mathbb{N}} \psi_j(n).
$$
Then, for $j \in \mathbb{N}$, the function $\psi_j$ is supported in the discrete annulus
\begin{equation}\label{annulus}
A_j := \{n \in \mathbb{Z}^2 : 2^{j-1} \le |n| \le 2^{j+1} \},
\end{equation}
and serves as a smooth approximation of the indicator function $\1_{A_j}$.

Let $P_j u$ denote the Littlewood–Paley projection of a tempered distribution $u$: $$
\widehat{P_j u}(n) := \psi_j(n)\, \hat{u}(n), \qquad \text{for } n \in \mathbb{Z}^2.
$$
initially defined on $L^2$ extending to distributions. Each operator $P_j$ acts as a frequency localization to the annulus $A_j$, effectively isolating the Fourier modes of $u$ at scale $\sim 2^j$.

The family $(P_j)_{j \in \mathbb{N}_0}$ yields the Littlewood–Paley decomposition
$$
u = \sum_{j \in \mathbb{N}_0} P_j u,
$$
where the sum converges in the sense of tempered distributions.

For the sum of the two vectors $n_1+n_2$ equals a third vector $n_3$, with values in, respectively,  the annulus $A_j, A_k,A_l$,  there are two ways for the vector sum $n_1+n_2$ equal $n$. (1) the size of $n$ is about the same as one of the two. (2) The size of $n_1$ and $n_2$ are about the same and much larger than that of $n$, which happens when the angle between $n_1$ and $n_2$ are close to each other, with small angle $\theta$. This follows from the Cosine Rule
$|n|^2=|n_1|^2+|n_2|^2-2|n_1||n_2| \cos\theta$  that $|n_1+n_2|$ is about the size  $ 2^{k \vee l}$.
These same relations translate into the Littlewood-Paley blocks. From this we have the following lemma:
\begin{lemma}\label{vanishing}
For any $j,k,l\in \N_0$ and  for any functions $f,g$,  $\P_j(\P_k f\P_lg)=0$ if $j,k,\ell$ satisfy one of the following relations:
\begin{itemize}
\item $j> (k\vee \ell)+2$;
\item  $|k-\ell |> 2$ and $j\le (k\vee \ell)-2$.
\end{itemize}
\end{lemma}

\begin{proof}
Consider the frequencies of the  Fourier transform of the projection $\P_j(\P_k f \P_l g)$, supported in~$A_j$.
\begin{equs}
\P_j(\P_k f\P_\ell g)&=\sum_{n\in A_j} \widehat{\P_j(\P_k f\P_ \ell g)}(n)e^{in\cdot x}
=\sum_{n\in \N_0}\psi_j(n) \widehat{\P_k f\P_\ell g}(n)e^{in\cdot x}\\
&=\sum_{n\in A_j}\psi_j(n) \sum_{n=n_1+n_2} \psi_k(n_1) \hat f(n_1)\psi_\ell(n_2)\hat g(n_2)e^{in\cdot x}.
\end{equs}
We have the further restrictions $n_1\in A_k$ and $n_2\in A_\ell$. 
If $j> (k\vee \ell)+2$, then
$$2^{j-1} > 2 \cdot 2^{k\vee l}\ge 2^{k+1}+2^{\ell+1}.$$
And $A_j\cap (A_k+A_l)=0$. Similarly  this is so if  $|k-\ell| >2$ and $j\le (k\vee \ell)-2$.
\end{proof}

\begin{proposition}\label{prop:tri1}
Let $0 <\alpha<\f 14$ and $ 0 < T \le 1$. If $$ \epsilon <\mincurly{ \f 1{12}, \f 14 -\alpha} \quad \hbox{and} \quad 
 0<\delta<\f 23 \epsilon,$$
 then there exists some $\eta=\eta(\epsilon)>0$ such that  the following holds
\begin{equ}
\Bigl\|uvw \Bigr\|_{\X^{-\f 3 4 +\epsilon, -\f 12 +\delta}_T}
 \lesssim T^\eta  \bigl\| u \bigr\|_{L^\infty_TW^{-\alpha-\epsilon,\infty}_x} \bigl\|v\bigr\|_{\X^{\f 1 4 +\epsilon,\f 12+\delta}_T}
\bigl \|  w \bigr\|_{\X^{\f 14 +\epsilon,\f 12+\delta}_T}\label{eq:tri1}
\end{equ}
for all functions $u,v,w$. 
\end{proposition}

\begin{proof} 
For $n_2,n_3\in \Z^2$,  let $n_{23}=n_2+n_3$, $N_{23}=|n_{23}|$.  This notation is inspired by: 
$$\widehat{vw}(n_{23})=\sum_{\stackrel {n_2, n_3\in \Z^2} 
{n_2+n_3=n_{23}} }\hat v(n_2) \hat w(n_3).$$
Similarly one defines  $N_1=|n_1|$, $N_{123}=|n_{123}|$ and $n_{123}=n_1+n_2+n_3$ etc.  With these notations we have the following Littlewood-Paley decomposition with respect to frequencies $n_{123},n_1$ and $n_{23}$,
\begin{equ}
\|uvw\|_{\X^{-\f 34+\epsilon, -\f 12 +\delta}_T} 
\le \sum_{N_{123},N_1,N_{23} \ge 1} 
\Bigl \| P_{N_{123}} \Big( P_{N_1}u \cdot P_{N_{23}}(vw)\Big) \Bigr\|
_{\X^{-\f 34+\epsilon, -\f 12 +\delta}_T},
\end{equ}
where $N_k=2^k$, $N_{23}=2^{k_{23}}$, and $N_{123}=2^{k_{123}}$ are dyadic numbers.

Since $N_{123}  \le |n_1| + |n_{23}| \le 2 \max (N_1, N_{23})$, we have two cases: either 
$N_{123}$ is much smaller than $\max (N_1, N_{23})$ or it is about the same size as $ \max (N_1, N_{23})$. More specifically:
$$\begin{aligned}
 &{}N_{123}  < \f 12 \max\{N_1,N_{23} \}\\
  \f 12 \max\{N_1,N_{23} \}& \le N_{123} \le 2\max\{N_1,N_{23} \}.\end{aligned}
  $$
 We denote respectively the two cases as $N_{123} \ll \maxcurly{N_1,N_{23}}$,
 $N_{123} \sim \maxcurly{N_1,N_{23}}$.

we have the following estimate
\begin{equ}\label{cubic-estimate}
\|uvw\|_{\X^{-\f 34+\epsilon, -\f 12 +\delta}_T}
 \lesssim  T^{\eta'(\epsilon)} \| u \|_{L^\infty_TW^{-\alpha-\epsilon,\infty}_x} \|v\|_{\X^{\f 1 4 +\epsilon,b}_T} \|  w \|_{\X^{\f 14 +\epsilon,b}_T} C(N_1,N_2, N_3, N_{123}).
\end{equ}
where $$C(N_1,N_2, N_3, N_{123})=\min \Bigl\{ \sum_{\substack{N_{123},N_{23} \\ N_1,N_2,N_3 \ge 1}} (N_{123})^{-\f 3 4+\epsilon} (N_1)^{\alpha+\epsilon} (N_3)^{\f 1 4 -\epsilon }, \sum_{\substack{N_{123},N_{23} \\ N_1,N_2,N_3 \ge 1}}  N^{\alpha+\epsilon}_1 N^{-\f 1 4 -\epsilon }_3
\Bigr\}.$$
Before proving this estimate we first show that the summation in $C(N_1,N_2, N_3, N_{123})$ is finite.

By symmetry, we assume $ N_2\le N_3$. By Lemma \ref{vanishing}, we have the following cases. 

\begin{itemize}
\item[(a)]   Case $|k_1-k_{23}|\le 2$ and $k_{123} \le \max (k_1,k_{23})+2 $. The first constraint 
translats into $\f 14 N_{23} \le N_1 \le  4N_{23}$ and in particular  $N_1\le 4N_3$.
Consequently, 
$$ \sum_{\substack{N_{123},N_{23} \\ N_1,N_2,N_3 \ge 1}}  N^{\alpha+\epsilon}_1 N^{-\f 1 4 -\epsilon }_3 \|   \lesssim  \sum_{\substack{N_{123},N_{23} \\ N_1,N_2,N_3 \ge 1}} N_3^{\alpha-\f 14}.$$
Since we summing up only over dyadic numbers of the form $N_1 =2^{k_1}$ and $N_3 =2^{k_3}$,
\begin{equ}
\label{dya-sum}
 \sum_{\substack{N_{123},N_{23} \\ N_1,N_2,N_3 \ge 1}} N_3^{\alpha-\f 14}\le \sum_{k_3}(k_3)^3 2^{k_3(\alpha-\f 14)},
\end{equ}
which is finite since $\alpha<\f 14$.

\item [(b)] Case $\max (k_1,k_{23} )-2 \le k_{123} \le \max (k_1,k_{23})+2 $ which translates into $\f 14 \max (N_1,N_{23}) \le  N_{123} \le 4\max (N_1,N_{23})$.  If furthermore $N_1\le N_{23}\le 2N_3$, then  $N_{123}$ and $N_{23}$ are controlled by $N_3$ and part (a) applies. 

If on the other hand, $N_{23}\le N_1$, $N_{123}$ is about the same order as $ N_1$, and $N_3$ is controlled by $N_1$.
we use the bound
\begin{equ}
 \sum_{\substack{N_{123},N_{23} \\ N_1,N_2,N_3 \ge 1}} N^{-\f 3 4+\epsilon}_{123} N^{\alpha+\epsilon}_1 N^{\f 1 4 -\epsilon }_3 \lesssim \sum_{\substack{N_{123},N_{23} \\ N_1,N_2,N_3 \ge 1}} N_1^{-\f 1 2+\alpha+\epsilon} \end{equ}
 is finite under the condition 
$$0<\alpha < \f 14 \quad \hbox{ and } 0<\epsilon <\f 14 - \alpha,$$
\end{itemize}
Assertion \eqref{eq:tri1} now follows. We now go ahead proving \eqref{eq:tri1}.

{\bf Bound 1.} 
Let $0< \epsilon < \f 1{12}$ and $0<\delta<\f 2 3\epsilon$. Then by the embedding lemma, \eqref{eq:aux-to-tri1}, we have
\begin{equ}\label{basic-5}
\begin{aligned}
\|P_{N_{123}}\Big( P_{N_1}u \cdot P_{N_{23}}(vw)\Big)\|_{\X^{-\f 34+\epsilon, -\f 12 +\delta}_T} \lesssim &T^{\eta'(\epsilon)}\|P_{N_{123}}\Big( P_{N_1}u \cdot P_{N_{23}}(vw)\Big)\|_{L^{\f 65}_TL^{\f 6 5}_x}\\
\lesssim & T^{\eta'(\epsilon)}\|P_{N_1}u \|_{L^\infty_TL^\infty_x} \|vw\|_{L^{\f 65}_TL^{\f 6 5}_x}.
\end{aligned}
\end{equ}

Next, we make use of the dyadic decomposition for frequencies $n_2$ and $n_3$, which leads to 
\begin{align}
&\|P_{N_{123}}\Big( P_{N_1}u \cdot P_{N_{23}}(vw)\Big)\|_{\X^{-\f 34+\epsilon, -\f 12 +\delta}_T} 
\lesssim T^{\eta'(\epsilon)} \|P_{N_1}u \|_{L^\infty_TL^\infty_x}\|vw\|_{L^{\f 6 5}_TL^{\f 6 5}_x} \\
\lesssim &T^{\eta'(\epsilon)} \sum_{N_2,N_3 \ge 1} \|P_{N_1}u \|_{L^\infty_TL^\infty_x}\|P_{N_2}v \cdot P_{N_3}w\|_{L^{\f 6 5}_TL^{\f 6 5}_x} \\
\le &  T^{\eta'(\epsilon)} N^{\alpha+\epsilon}_1  \sum_{N_2,N_3 \ge 1} \|u \|_{L^\infty_TW^{\alpha+\epsilon,\infty}_x}  \|P_{N_2}v\|_{L^3_TL^3_x} \|P_{N_3}w\|_{L^2_TL^2_x} \\
\le & T^{\eta'(\epsilon)} N^{\alpha+\epsilon}_1  \sum_{N_2,N_3 \ge 1} \|u \|_{L^\infty_TW^{\alpha+\epsilon,\infty}_x} \|P_{N_2}v \|_{L^{\f {12}{1+4\epsilon}}_TL^{\f 3{1-2\epsilon}}_x} \|P_{N_3}w\|_{L^2_TL^2_x}
\end{align}
where we utilize H\"{o}lder's inequality repeatly and $0 < T \le 1$. In the third line we have used the fact that 
\begin{equ}
\label{Payley-projection}
 \|P_{N_1}u \|_{L^\infty_TL^\infty_x} \sim N^{\alpha+\epsilon}_1 \|\jpb{\nabla}^{-\alpha-\epsilon}P_{N_1}u \|_{L^\infty_TL^\infty_x} \sim  N^{\alpha+\epsilon}_1 \|u \|_{L^\infty_TW^{\alpha+\epsilon,\infty}_x}.
\end{equ}

Since $\Big( \f {12}{1+4\epsilon}, \f 3{1-2\epsilon} \Bigr)$ is $(\f 1 4 +\epsilon)$-admissible, we apply Lemma \ref{lem:transference-principle} to obtain, for any $b>\f 12$:
\begin{equ}\label{cubic-estimate}
\|uvw\|_{\X^{-\f 34+\epsilon, -\f 12 +\delta}_T} \lesssim T^{\eta'(\epsilon)} \| u \|_{L^\infty_TW^{-\alpha-\epsilon,\infty}_x} \|v\|_{\X^{\f 1 4 +\epsilon,b}_T}  \|  w \|_{\X^{\f 14 +\epsilon,b}_T}  \sum_{\substack{N_{123},N_{23} \\ N_1,N_2,N_3 \ge 1}}  N^{\alpha+\epsilon}_1 N^{-\f 1 4 -\epsilon }_3
\end{equ}

\medskip

{\bf Bound 2.}  If we directly start from the definition of $\X^{s,b}_T$-spaces, we would have
\begin{align}
\|uvw\|_{\X^{-\f 34+\epsilon, -\f 12 +\delta}_T}  \le & \sum_{N_{123},N_1,N_{23} \ge 1}\| \jpb{\nabla} ^{-\f 3 4+\epsilon}P_{N_{123}} \Big( P_{N_1}u \cdot P_{N_{23}}(vw) \Big)\|_{L^2_TL^2_x} \\
\lesssim &\sum_{N_{123},N_1,N_{23} \ge 1} (N_{123})^{-\f 3 4+\epsilon} \|(P_{N_1}u \cdot P_{N_{23}}(vw) )\|_{L^2_TL^2_x} \\
\lesssim &  \sum_{N_{123},N_1,N_{23} \ge 1} (N_{123})^{-\f 3 4+\epsilon}( N_1)^{\alpha+\epsilon} \| u \|_{L^\infty_TW^{-\alpha-\epsilon,\infty}_x} \|P_{N_{23}}(vw)\|_{L^2_TL^2_x}.
\end{align}
We further apply the Littlewood-Paley decomposition to  frequencies $n_2$ and $n_3$.  Let $\epsilon<\f 12$.  We first apply H\"{o}lder's inequality, then go down to a lower H\"older space to gain a factor $\eta'(\epsilon) >0$:
\begin{align}
\|P_{N_{23}}(vw)\|_{L^2_TL^2_x} \le &\sum_{N_2,N_3 \ge 1} \|P_{N_2} v P_{N_3} w\|_{L^2_TL^2_x} 
\le  \sum_{N_2,N_3 \ge 1} \|P_{N_2} v \|_{L^4_TL^{\f3 {1-2\epsilon}}_x} \|P_{N_3} w\|_{L^4_TL^{\f 6 {1+2\epsilon}}_x} \\
\le & T^{\eta'(\epsilon)}  \sum_{N_2,N_3 \ge 1} \|P_{N_2} v \|_{L^{\f {12}{1+4\epsilon}}_TL^{\f3 {1-2\epsilon}}_x} \|P_{N_3} w\|_{L^6_TL^{6}_x} .
\end{align}

Next we  apply Lemma \ref{lem:transference-principle} to bound the $L^pL^q$ norm by the $\X^{s,b}$ norm. For any $b >1/2$,
\begin{equs}
\|P_{N_2} v \|_{L^{\f {12}{1+4\epsilon}}_TL^{\f3 {1-2\epsilon}}_x} \|P_{N_3} w\|_{L^6_TL^{6}_x} 
&\le \|v \|_{L^{\f {12}{1+4\epsilon}}_TL^{\f3 {1-2\epsilon}}_x} \|P_{N_3}w\|_{L^6_TL^{6}_x} \\
&\le N^{\f 1 4 -\epsilon }_3 \|v\|_{\X^{\f 1 4 +\epsilon,b}_T} \|  w \|_{\X^{\f 14 +\epsilon,b}_T}
\end{equs}
We have used the transfer principle and the fact that $\big( \f {12}{1+4\epsilon}, \f 3{1-2\epsilon} \bigr)$ is a wave-admissible pair for $\f 1 4 +\epsilon$, and that $\big( 6,6 \big)$ is a admissible pair for the regularity $s=\f 12$.
Putting everything together,  we obtain 
\begin{equ}\label{cubic-estimate2}
\|uvw\|_{\X^{-\f 34+\epsilon, -\f 12 +\delta}_T}
 \lesssim  T^{\eta'(\epsilon)} \| u \|_{L^\infty_TW^{-\alpha-\epsilon,\infty}_x} \|v\|_{\X^{\f 1 4 +\epsilon,b}_T} \|  w \|_{\X^{\f 14 +\epsilon,b}_T} \sum_{\substack{N_{123},N_{23} \\ N_1,N_2,N_3 \ge 1}} N^{-\f 3 4+\epsilon}_{123} N^{\alpha+\epsilon}_1 N^{\f 1 4 -\epsilon }_3.
\end{equ}
This completes the proof of the bound \eqref{cubic-estimate} and the proposition.
\end{proof}

\begin{lemma}\label{lem:aux-to-tri1}
Let $0< \epsilon < \f 1{12}$ and $0<\delta<\f 2 3\epsilon$. Then there exists a constant $\eta=\eta(\epsilon)>0$ such that, for every function~$z$:
\begin{equ}
\| z \|_{\X^{-\f 34 +\epsilon, - \f 12+\delta}_T} \lesssim T^{\eta} \;\|z\|_{L^{6/5}_TL^{6/ 5}_x}.
 \label{eq:aux-to-tri1}
\end{equ}
\end{lemma}
\begin{proof}
Recall that $\X^{s,b}=(\X^{-s, -b})^*$. For brevity, write 
 $$\X=\X^{-\f 34 +\epsilon, - \f 12+\delta}_T, \qquad \X^*=\X^{\f 34 -\epsilon, \f 12-\delta}_T.$$
Using this duality relation and Hölder's inequality, we obtain
\begin{align}
\| z \|_{\X}  = &\sup_{ \|\tilde{z}\|_{\X^*} = 1} \Big| \int_{[0,T] \times \T^2} z\tilde{z} \, \d x \d t \Big| 
\le  \|z\|_{L^{\f 65}_TL^{\f 6 5}_x} \|\tilde{z}\|_{L^{6}_TL^6_x}.
\end{align}

Next, we apply Lemma \ref{lem:Xsb-to-LqLr} to bound $\|\tilde{z}\|_{L^{\tilde q_\theta}_TL^{\tilde r_\theta}_x} $ in terms of its norm in $\X^*$. To this end,  we set 
\begin{equs}
(q,r)&=(\f 56,\f 56), \qquad (s,b)=(\f 34 -\epsilon,\f 12-\delta)\\
\tilde s&=\f 34, \quad (\tilde q, \tilde r)=\big(8,16\big)
\end{equs}
One verifies that $(8,16)$ is a $\tilde{s}$-admissible pair.  

It remains to identify the interpolation constant $\theta$, which we take to be of the form $1-\f {4}{3}\epsilon$, where $\epsilon$ will be chosen to satisfy
 $$ \f {\theta}{8} +\f {1-\theta}{2} < \f 16, \qquad \f {\theta}{16} +\f {1-\theta}{2} = \f 16.
$$
For simplicity, the second condition is taken as an equality, while the first must hold strictly since it concerns the time integrability in $L^q$.  

Finally, we choose $\tilde b$ sufficiently close to $\frac{1}{2}$ and $\delta$ so that
$$ \f1 2-\delta > \f 12\Big(\f {3-4\epsilon}{3}\Big).$$
This yields the admissible range
$$0< \epsilon < \f 1{12}, \qquad 0<\delta<\f 2 3\epsilon. $$
Under these conditions, Lemma~\ref{lem:Xsb-to-LqLr} applies to give, $\|\tilde{z}\|_{L^{6}_TL^6_x} \le T^{\eta(\epsilon)} \|\tilde{z}\|_{L^{q_\theta}_TL^{r_\theta}_x}$, which completes the proof.
\end{proof}

\begin{proposition}\label{prop:tri2}
Let $ \alpha\in (0,\f 14)$ and $T\in (0, 1]$.  For any $0< \epsilon < \f 1{12}$, $0<\delta<\f 2 3\epsilon$, there exists $\eta(\epsilon)>0$ such that
\begin{equ}\label{eq:tri2}
\|uvw\|_{\X^{-\f 3 4 +\epsilon, -\f 12 +\delta}_T} 
\lesssim T^{\eta(\epsilon)}  
\| u \|_{\X^{\f 1 4 +\epsilon,\f 12+\delta}_T} \|v\|_{\X^{\f 1 4 +\epsilon,\f 12+\delta}_T} \|  w \|_{\X^{\f 14 +\epsilon,\f 12+\delta}_T}.\end{equ}
\end{proposition}

\begin{proof}
By duality and H\"{o}lder's inequality,  we have 
\begin{align}
&\|u v w \|_{\X^{-\f 34 +\epsilon, -\f 12 +\delta}_T} 
{=  \sup_{\|z\|_{\X^{\f 34-\epsilon,\f 12 -\delta}_T} =1 } \Big| \int_{[0,T]\times \T^2} u v w z \, \d x \d t\Big|}
 \\
& \le  \sup_{ \|z\|_{\X^{\f 34-\epsilon,\f 12 -\delta}_T} =1 } 
\Big( \|u\|_{L^{\f {12}{1+4\epsilon}}_TL^{\f 3{1-2\epsilon}}_x} \|v\|_{L^{\f {12}{1+4\epsilon}}_TL^{\f 3{1-2\epsilon}}_x} \|w\|_{L^{\f {12}{1+4\epsilon}}_TL^{\f 3{1-2\epsilon}}_x}  \|z\|_{L^{\f 4{3-4\epsilon}}_TL^{\f 1{2\epsilon}}_x}  \Big).
\end{align}
Since the pair $\big( \f {12}{1+4\epsilon}, \f 3{1-2\epsilon} \bigr)$ is wave-admissible for $s=\f 1 4 +\epsilon$,  by the transfer principle Lemma \ref{lem:transference-principle},
$\X^{\f 14+\epsilon,b} \subset L^{\f {12}{1+4\epsilon}}_TL^{\f 3{1-2\epsilon}}_x $. 
Finally since $\|z\|_{\X^{\f 34-\epsilon,\f 12 -\delta}_T}=1$, we only need to show 
\begin{equ}
 \|z\|_{L^{\f 4{3-4\epsilon}}_TL^{\f 1{2\epsilon}}_x}  \lesssim T^{\eta(\epsilon)}  \|z\|_{\X^{\f 34-\epsilon,\f 12 -\delta}_T}.
\end{equ}
This can be obtained using the $\f 34$-admissible
 pair $\big( \f 4{1-8\epsilon}, \f 1{\epsilon}\big)$. 
 Letting  $\theta = \f {3-4\epsilon}3$, one can solve
 $$
q_\theta= \f {(1-8\epsilon)\theta}4+ \f {1-\theta}2 <\f{3-4\epsilon} 4,  \qquad \, r_\theta=\epsilon\theta +\f {(1-\theta)}2  \le 2\epsilon.$$
Finally one choose $\delta$ such that $ \f1 2-\delta > \f 12 \theta$. It follows from the interpolation inequality that $ \|z\|_{L^{\f 4{3-4\epsilon}}_TL^{\f 1{2\epsilon}}_x}  \le T^{\eta(\epsilon)} \|z\|_{L^{q_\theta}_TL^{r_\theta}_x}$, which concludes the proof.
\end{proof}

\subsection{Trilinear-estimates  on  $\X^{2\alpha-\f 14, b}$ spaces}
\begin{proposition}\label{prop:tri3}
Let $\f 14 \le \alpha<\f 3 8$. Let $\epsilon, \delta$ be positive numbers such that 
$$\epsilon < \mincurly{\f {3-8\alpha} 4, \f {(5-8\alpha)(24\alpha-5)}{12(1+8\alpha)} }, \quad \delta< \f {2\epsilon}{5-8\alpha}.$$Then, there exists some $\eta=\eta(\epsilon,\alpha)>0$ such that
\begin{equ}\label{eq:tri3}
\|uvw\|_{\X^{2\alpha - \f 14 +\epsilon-1, -\f 12 +\delta}_T} \lesssim T^{\eta} \| u \|_{L^\infty_TW^{-\alpha-\f \epsilon 2,\infty}_x} \|v\|_{\X^{2\alpha - \f 14+\epsilon,\f 12+\delta}_T} \|  w \|_{\X^{2\alpha - \f 14 +\epsilon,\f 12+\delta}_T},\label{eq:tri3}
\end{equ}
for all $u,v,w$.
\end{proposition}
\begin{proof}
The proof is analogous to that for Proposition \ref{prop:tri1}.  The main difference is to choose, and work with, a $2\alpha -\f 14 +\epsilon$-admissible pair,  versus a $\f 14$-admissible pair.

In the estimate below, by symmetry, we assume that $N_3 \ge N_2$. Using the Littlewood-Paley decomposition, and since we have:
\begin{align}
\|uvw\|_{\X^{2\alpha - \f 14 +\epsilon-1, -\f 12 +\delta}_T} \le &\sum_{\substack{N_{123},N_{23} \\ N_1,N_2,N_3 \ge 1}} \|P_{N_{123}}\Big( P_{N_1}u \cdot P_{N_{23}}(P_{N_2}v \cdot P_{N_3}w)\Big)\|_{\X^{2\alpha - \f 14 +\epsilon-1, -\f 12 +\delta}_T} \\
\lesssim & T^{\eta'(\epsilon,\alpha)} \sum_{\substack{N_{123},N_{23} \\ N_1,N_2,N_3 \ge 1}}  \|P_{N_1}u \cdot P_{N_{23}}(P_{N_2}v \cdot P_{N_3}w)\|_{L^{\f {18}{13}}_TL^{\f 6 5}_x}
\end{align}
where we have applied \eqref{eq:aux-to-tri3} from Lemma \ref{lem:aux-to-tri3}, which is applicable for the range of $\epsilon$ and $\delta$ specified in the statement.

It suffices to consider two cases: $N_{123} \ll \maxcurly{N_1,N_{23}}$ or $N_{123} \sim \maxcurly{N_1,N_3}$.

If $N_{123} \ll \maxcurly{N_1,N_{23}}$, we apply H\"{o}lder's inequality to obtain 
\begin{align}
&  \|P_{N_1}u \cdot P_{N_{23}}(P_{N_2}v \cdot P_{N_3}w)\|_{L^{\f {18}{13}}_TL^{\f 65}_x} \\
\le &\|P_{N_1}u\|_{L^\infty_TL^\infty_x} \|P_{N_2}v\|_{L^{\f 9 2}_TL^3_x} \|P_{N_3}w\|_{L^2_TL^2_x} \\
 \lesssim & N^{\alpha+\f \epsilon 2}_1 \|u\|_{L^\infty_TW^{-\alpha-\f \epsilon 2}_x} \|P_{N_2}v\|_{L^{\f 9 2}_TL^3_x} \|P_{N_3}w\|_{L^2_TL^2_x} \\
 \le  &T^{\eta''(\epsilon,\alpha)} N^{\alpha+\f \epsilon 2}_1 \|u\|_{L^\infty_TW^{-\alpha-\f \epsilon 2}_x} \|v\|_{L^{\f {12}{8\alpha-1+4\epsilon}}_TL^{\f 3{2-4\alpha-2\epsilon}}_x} N^{-2\alpha+ \f 14-\epsilon}_3\|w\|_{\X^{2\alpha-\f 14+\epsilon,\f 12+\delta}_T},
\end{align}
The above H\"{o}lder's inequality is applicable when $\f 14 \le \alpha <\f 3 8$ with $0< \epsilon < \mincurly{\f {11}3-8\alpha, 1-2\alpha}$.

Note $\Big(\f {12}{8\alpha-1+4\epsilon}, \f 3{2-4\alpha-2\epsilon}\Big)$ is $2\alpha -\f 14 +\epsilon$-admissible and we must have $N_1 \sim N_{23}$, it yields,
\begin{align}
\|uvw\|_{\X^{\alpha+\epsilon-1, -\f 12 +\delta}_T} \lesssim &\sum_{\substack{N_{123},N_{23} \\ N_1,N_2,N_3 \ge 1}}  N^{\alpha+\f \epsilon 2}_1 \|u\|_{L^\infty_TW^{-\alpha-\f \epsilon 2}_x} \|v\|_{\X^{2\alpha-\f 14+\epsilon,\f 12+\delta}_T} N^{-2\alpha+ \f 14-\epsilon}_3\|w\|_{\X^{2\alpha-\f 14+\epsilon,\f 12+\delta}_T}
\end{align} 
which is summable if $\alpha+\f \epsilon 2 -2\alpha+ \f 14-\epsilon <0$ for any and $\epsilon >0$, which is equivalent to $\alpha \ge \f 14$.

The remaining case $N_{123} \sim \maxcurly{N_1,N_3}$ can be treated exactly as in the proof for Proposition \ref{prop:tri1}.
\end{proof}

\begin{lemma}\label{lem:aux-to-tri3}
Let $\f 14\le \alpha<\f 3 8$ and $0<T\le 1$. Assume that
$$0<\epsilon <\f {3-8\alpha} 4 \qquad \hbox{and} \quad 0<\delta< \f {2\epsilon}{5-8\alpha}.$$
Then there exists a constant $\eta=\eta(\epsilon,\alpha)>0$ such that for all functions $z$,
\begin{equ}
\| z \|_{\X^{2\alpha - \f 14 +\epsilon -1, - \f 12+\delta}_T} \lesssim T^{\eta}  \|z\|_{L^{\f {18}{13}}_TL^{\f 65}_x}. \label{eq:aux-to-tri3}
\end{equ}
\end{lemma}
\begin{proof}
We begin as in Lemma~\ref{lem:aux-to-tri1}, writing $\X^*=\X^{1-2\alpha+\f 14 -\epsilon,\f 12-\delta}_T$ for the dual space.
\begin{align}
\| z \|_{\X^{2\alpha -\f 14 +\epsilon -1, - \f 12+\delta}_T}  = &\sup_{ \|\tilde{z}\|_{\X^*} = 1 }\Big| \int_{[0,T] \times \T^2} z\tilde{z} \, \d x \d t \Big| 
\le \|z\|_{L^{\f {18}{13}}_TL^{\f 65}_x} \|\tilde{z}\|_{L^{\f {18}5}_TL^6_x}.
\end{align}

Choose the regularity parameter $s=\f 54 - 2\alpha$.  Then the pair $\bigl(\f {12}{5-8\alpha}, \f 3{4\alpha-1}\bigr)$ is $\f 5 4 - 2\alpha$-admissible because the condition
\begin{equ}
\f 1 r \ge \f {4\alpha-1}3, \qquad \f 1 q \le \f {5-8\alpha}{12}
\end{equ}
holds.
(Also note that $(2,2)$ is $s$-admissible for $s=0$.) Choose the interpolation constant $$ \theta = \f{5-8\alpha-4\epsilon}{5-8\alpha}.$$
Assume:
 \begin{equ}\label{eq:interpolation-condition-54-2alpha}
 \f {(5-8\alpha)\theta}{12} + \f {1-\theta}2 < \f 5 {18},  \qquad
\f {(4\alpha-1)\theta}3 + \f {1-\theta}2 \le \f 1 6
 \end{equ}
 and that
 \begin{equ}
 \f 12-\delta  > \f 12 \, \f{5-8\alpha-4\epsilon}{5-8\alpha}.\end{equ}
Under these conditions, we may apply the interpolation Lemma, Lemma~\ref{lem:Xsb-to-LqLr} to obtain
$$\|\tilde{z}\|_{L^{q_\theta}_rL^{r_{\theta}}_x} \lesssim \|\tilde{z}\|_{\X^{\f 54-2\alpha-\epsilon,\f 12-\delta}_T}$$ 
From the above conditions, we compute that
$$0<\epsilon < \mincurly{\f {3-8\alpha} 4, \f {(5-8\alpha)(24\alpha-5)}{12(1+8\alpha)} },\qquad 0<\delta< \f {2\epsilon}{5-8\alpha}.$$
Finally, by Lemma \ref{lem:transference-principle}, $ \|\tilde{z}\|_{L^{\f {18} 5}_TL^6_x} \le T^{\eta(\epsilon,\alpha)} \|\tilde{z}\|_{L^{q_\theta}_TL^{r_\theta}_x}$. Substituting back, this completes the proof. \end{proof}

Similarly, we have the following trilinear estimate:
\begin{proposition}\label{prop:tri4}
Let $\f 14  \le \alpha< \f 38$ and $ 0 < T \le 1$.  Assume that
 $$0< \epsilon <  \f {(5-8\alpha)^2}{32(1-\alpha)} \quad 
  \hbox{and} \quad  0<\delta< \f {2\epsilon}{5-8\alpha}.$$ 
  Then there exists a constant $\eta=\eta(\epsilon,\alpha)>0$ such that for all $u,v,w$,
\begin{equ}\label{eq:tri4}
\|uvw\|_{\X^{2\alpha -\f 14 +\epsilon -1, -\f 12 +\delta}_T} 
\lesssim T^{\eta} \;\| u \|_{\X'} \|v\|_{\X'} \|  w \|_{\X'},
\end{equ}
 where $\X'=\X^{2\alpha -\f 14 +\epsilon,\f 12+\delta}_T$. 
\end{proposition}
\begin{proof}
We argue by duality. Setting,  $\X^*=\X^{1-2\alpha + \f 14-\epsilon,\f 12 -\delta}_T$ and apply Hölder’s inequality in mixed time–space Lebesgue norms, introduce the exponents
$$\tilde q=\f {12}{8\alpha-1+4\epsilon}, \quad \tilde r= \f 3{2-4\alpha-2\epsilon}.$$ 
With this choice, we obtain
$$
\norm{u v w }_{\X^{2\alpha- \f 14 +\epsilon -1, -\f 12 +\delta}_T} 
\le  \sup_{\|z\|_{\X^*} =1} 
\Big( \bigl\|u \bigr\|_{\tilde q, \tilde r} \bigl\|v\bigr\|_{\tilde q,\tilde r} \bigl \|w\bigr\|_{\tilde q,\tilde r} \bigl  \|z\bigr\|_{{\f 3 {4-8\alpha-4\epsilon}}, {\f 1{4\alpha-1+2\epsilon}}} \Big)
$$
The pair $\big(\tilde q,\tilde r \big)$ is wave‑admissible for $s'=2\alpha - \f 14+\epsilon$, so that for any $\delta>0$, 
$$\bigl\|u \bigr\|_{\tilde q, \tilde r} \le \bigl\|u \bigr\|_{\X^{s', \f 12+\delta} }, $$
and similarly for $v$ and $w$.

The remaining factor is the term involving $z$. We claim that
$$\|z\|_{L^{\f 3 {4-8\alpha-4\epsilon}}_TL^{\f 1{4\alpha-1+2\epsilon}}_x}  \lesssim T^{\eta(\epsilon,\alpha)}\|z\|_{\X^*}.$$
To verify this, introduce $\tilde s:=\f 5 4 - 2\alpha$.  Then the pair 
$\big(\f {12}{5-8\alpha}, \f 3{4\alpha-1}\big)$ is wave -admissible for regularity $\tilde s$.
For $\f 14  \le \alpha< \f 38$ we observe that  $$\f {12}{5-8\alpha} > \f 4{5-8\alpha-4\epsilon}, \qquad \f 3{4\alpha-1} > \f 1{4\alpha-1+2\epsilon}.$$
Thus, we interpolate between admissible pairs by setting $\theta = \f{5-8\alpha-4\epsilon}{5-8\alpha}$.
Interpolation yields conditions:
\begin{align}
  q_\theta\f {(5-8\alpha)\theta}{12} + \f {1-\theta}2 < \f {5-8\alpha-4\epsilon}4, \qquad r_\theta =\f {(4\alpha-1)\theta}3 + \f {1-\theta}2 \le 4\alpha-1+2\epsilon, \end{align}
and we also require
 $$   \f 12-\delta  > \f 12\theta.$$ 
 These inequalities are consistent provided
 $$0< \epsilon<\f {(5-8\alpha)^2}{32(1-\alpha)}, \qquad 0<\delta<\f {2\epsilon}{5-8\alpha}.$$
 With these parameter restrictions in place, Lemma~\ref{lem:Xsb-to-LqLr} applies to the dual factor 
$z$ and gives the desired bound with a time-gain $T^{\eta(\epsilon, \alpha)}$.
 Combining all the above  establishes the desired inequality. 
    \end{proof}

\section{Local well-posedness and stability of the remainder equation}\label{sec:local-wellposedness}

In this section we show that the limit equation is well posed. Throughout the section we take $T\in (0, 1]$. We first recall the classical result.
Denote by  $\mathbf S(t)$ the wave propagator, i.e. the solution map to the linear wave equation. Explicitly, $$\mathbf S(t)(u_0,u_1)=\cos(t\jpb{\nabla})u_0  +  \frac{\sin(t\jpb{\nabla})}{\jpb{\nabla}}u_1$$
with $t\in [0,1]$. 
Let $\mathcal S$ denote the set of symmetric Schwartz functions $\eta :\R\to \R$ with $\eta=1$ on $[-1,1]$ and with support in $(-2,2)$. We will use the following facts:
\begin{lemma}
\begin{itemize}
\item [(1)] For any $s>0$, $\X^{s,b} \subset C(\R, H^s_x)$;
\item[(2)]  For any $\eta \in \mathcal S$, 
$\|\eta u\|_{\X^{s,b}}\lesssim_{\eta, b} \|u\|_{\X^{s,b}}$ and 
consequently $\X_T^{s,b} \subset C([0,T], H^s_x)$
\item  [(3)]
For any $\epsilon,\delta>0$, the solution map $\mathbf S(t)$ is a contraction on $\X^{\f 14+\epsilon,\f 12 +\delta}_T$:
\begin{align}
\| \mathbf S(t)(u_0,u_1)\|_{\X^{\f 14+\epsilon,\f 12 +\delta}_T} \le  \|(u_0,u_1) \|_{\mathcal{H}^{\f 14+\epsilon}}.
\end{align}
\end{itemize}
\end{lemma}

We first prove the auxiliary bilinear estimate which will be useful in proving the local well-posedness.

\begin{lemma}\label{lem:bi1}
\begin{itemize}
\item [(1)] Let $0<\alpha<\f 18$. If $0<\epsilon < \f 3 8 - \alpha$ and $\delta>0$, then for all functions $u,v$,\begin{equ}\label{eq:bi1}
\|uv\|_{\X^{-\f 34+\epsilon, -\f 12+2\delta}_T} \lesssim \|u\|_{L^\infty_TW^{-2\alpha-\epsilon,\infty}_x} \|v\|_{\X^{\f 1 4+\epsilon, \f 12+\delta}_T}.
\end{equ}
\item [(2)]Let $\f 3{10} \le \alpha<\f 38$. If $0<\epsilon < \f 3 8 - \alpha$ and $\delta>0$, we have
\begin{equ}\label{eq:bi1}
\|uv\|_{\X^{2\alpha -\f 54+\epsilon, -\f 12+2\delta}_T} \lesssim \|u\|_{L^\infty_TW^{-\alpha-\epsilon,\infty}_x} \|v\|_{\X^{2\alpha-\f 1 4+\epsilon, \f 12+\delta}_T} \label{eq:bi-for-38}
\end{equ}
for any $u,v$.

\end{itemize}
\end{lemma}

\begin{proof}
Inserting the Littlewood-Paley decomposition  we have:
\begin{align}
\|uv\|_{\X^{-\f 34+\epsilon, -\f 12+2\delta}_T} \le & \sum_{N_{12},N_1,N_2 \ge 1} \|P_{N_{12}}( P_{N_1}u \cdot P_{N_2} v)\|_{\X^{-\f 34+\epsilon, -\f 12+2\delta}_T}\\
\le & \sum_{N_{12},N_1,N_2 \ge 1} N^{-\f 34+\epsilon}_{12} \|P_{N_1}u \cdot P_{N_2} v \|_{L^2_TL^2_x}.
\end{align}
Without loss of generality, we assume $N_1 \ge N_2$. There are two subcases.
\begin{itemize}
\item $N_1\sim N_2$, i.e.  $\f 12 N_1 \le N_2\le N_1$
\item $N_1 \gg N_2$,  i.e. $N_2< \f 12 N_1$ in which case $N_{12} \sim N_1$.
\end{itemize}

If  $N_1 \gg N_2$, then for any $\epsilon>0$ and $b > \f 12$, by Lemma \ref{lem:transference-principle}, one has
\begin{align}
\|P_{N_1}u \cdot P_{N_2} v \|_{L^2_TL^2_x} 
&\le \|P_{N_1} u \|_{L^\infty_TL^\infty_x} \|P_{N_2} v \|_{L^2_TL^2_x}
\le N^{2\alpha+\epsilon}_1 \|u \|_{L^\infty_TW^{-2\alpha-\epsilon,\infty}_x} \;\|v \|_{L^2_TL^2_x} \\
&\lesssim N^{2\alpha+\epsilon}_1 \|u\|_{L^\infty_TW^{-2\alpha-\epsilon,\infty}_x} \; \|v \|_{\X^{\f 1 4+\epsilon,b}_T} 
\end{align}
Consequently, 
\begin{align}
\|uv\|_{\X^{-\f 34+\epsilon, -\f 12+2\delta}_T} \le & \sum_{N_{12},N_1,N_2 \ge 1} N^{-\f 34+\epsilon}_{12}N^{2\alpha+\epsilon}_1\|u \|_{L^\infty_TW^{-2\alpha-\epsilon,\infty}_x}\|v \|_{\X^{\f 1 4+\epsilon,b}_T} 
\end{align}

Note that in this subcase, we have $N_{12} \sim N_1$. Then, for  $\alpha < \f 3 8$ and $\epsilon < \f 3 8 -\alpha$, the sum
$$\sum_{N_{12},N_1,N_2 \ge 1} N^{-\f 34+\epsilon}_{12}N^{2\alpha+\epsilon}_1
\le \sum_{N_{12},N_1,N_2 \ge 1} N_1^{2\alpha-\f 34+2\epsilon}<\infty.$$

If $N_1 \sim N_2$, we have $N_{12}\sim N_1$ or $N_{12} \lesssim N_1$. The former case is the same as the case $N_1 \gg N_2$ above. In the latter case, 
we cannot make use of the decay in $N_{12}$. To this end, we have
\begin{align}
\|uv\|_{\X^{-\f 34+\epsilon, -\f 12+2\delta}_T} \le & \sum_{N_{12},N_1,N_2 \ge 1} N^{-\f 34+\epsilon}_{12} N^{2\alpha+\epsilon}_1 \|u \|_{L^\infty_TW^{-2\alpha-\epsilon,\infty}_x} \;\|P_{N_2}v \|_{L^2_TL^2_x}\\
\le & \sum_{N_{12},N_1,N_2 \ge 1} N^{2\alpha+\epsilon}_1 N^{-\f 1 4 -\epsilon}_2 \|v \|_{L^\infty_TW^{-2\alpha-\epsilon,\infty}_x}\|v \|_{\X^{\f 1 4+\epsilon,b}_T} 
\end{align}
Hence, it is summable only when $0 < \alpha < \f 18$.
\end{proof}

\subsection{Case $0<\alpha <\f 14$}
Let  $0<\alpha <\f 14$. As discussed in \S\ref{sec:deduction}, the limiting remainder equation
\begin{equ}\label{eq:1st-v}
\begin{cases}
(\partial^2_t+1-\Delta) v= -3\<3> +3 \<2>v-3\<1>v^2-3v^2\\
(v,\partial_t v)|_{t=0} =(u_0,v_0)
\end{cases}
\end{equ}
is expected to be locally well-posed in $\X^{\f 1 4+\epsilon,\f 12 +\delta}_T$.
\begin{proposition}\label{prop:LWP-18}
Let $0 < \alpha < \f 1 8$. There exist thresholds $\epsilon_0>0$ and $\delta_0>0$ such that, for all $\epsilon>\epsilon_0$ and $\delta >\delta_0$, and for every initial data $(u_0,u_1) \in \mathcal{H}^{\f 1 4+\epsilon}$, there exist a (random time) $0<T<1$ for which the cubic nonlinear SNLW \eqref{eq:1st-v} admits, almost surely,  a unique solution on $[0, T]$ in
\begin{equ}
\X^{\f 1 4+\epsilon,\f 12 +\delta}_T \subset C([0,T];H^{\f 1 4+\epsilon}_x).
\end{equ}
Furthermore,  for almost surely all $\omega$,
the solution $v$ depends continuously on the enhanced data set
\begin{equ}
\Theta = (u_0,u_1,\<1>,\<2>,\<3>)
\end{equ}
where $\Theta $ ranges in the space
\begin{equ}
\mathcal{Y}^{\alpha,\epsilon,\delta}_T = \mathcal{H}^{\f 1 4+\epsilon} \times C([0,T];W^{-\alpha-\epsilon,\infty}_x) \times C([0,T];W^{-2\alpha-\epsilon,\infty}_x) \times \X^{-3\alpha-\epsilon,\delta}_T,
\end{equ}
with the norm
$$\|\Theta\|_{\mathcal{Y}^{\alpha,\epsilon}_T} = \|(u_0,u_1)\|_{\mathcal{H}^{\f 1 4+\epsilon}}+ \| \<1> \|_{L^\infty_TW^{-\alpha-\epsilon,\infty}_x}+ \|\<2>\|_{L^\infty_TW^{-2\alpha-\epsilon,\infty}_x}+\|\<3>\|_{\X^{-3\alpha-\epsilon,\delta}_T}.$$
\end{proposition}

\begin{proof}
To establish well-posedness of ~\eqref{eq:1st-sol}, we apply the contraction mapping principle. 
 There $T\in (0, 1)$ and constants $K_i$ such that 
\begin{equ}\label{chicken-fixed-1}
\| \<30>_N\|_{\X^{\f 14+\epsilon,b}_T}<K_4, \qquad  a.s.
\end{equ}
(Proposition \ref{prop:reg-of-cube}),
$\| \<1> \|_{L^\infty_TW^{-\alpha-\epsilon,\infty}_x}\le K_2$ and 
$ \|\<2>\|_{L^\infty_TW^{-2\alpha-\epsilon,\infty}_x}\le K_3$, almost surely, see Lemma \ref{reg-cherry}. We further assume that 
$ \|(u_0,u_1)\|_{\mathcal{H}^{\f 1 4+\epsilon}}\le  K_1$.
Set $K=\sum_{j=1}^4 K_j$.

Define the solution map
\begin{equ}\label{eq:1st-sol}
\Gamma_{\Theta} (v)(t) =  \mathbf S(t)(u_0,u_1) -3 \cJ(\<3>) +3 \cJ(\<2>v)-3\cJ(\<1>v^2)-3\cJ(v^3)
\end{equ}

We first show that $\Gamma_\Theta$ maps a suitable closed ball in  $\bar{B}_R \subset \X^{\f 1 4+\epsilon,\f 12 +\delta}_T$ to itself.

For $\alpha\in (0, \f 38)$, the linear estimate holds $$\|\mathbf S_t(u_0,u_1)\|_{\X^{\f 14+\epsilon,\f 12 +\delta}_T} \le  \|(u_0,u_1) \|_{\mathcal{H}^{\f 14+\epsilon}}$$
holds.
When $\alpha\in (0, \f 14)$, we fix parameters  $$0< \epsilon < \mincurly{ \f 1{12}, \f 14 -\alpha}, \qquad 0< \delta<\f 2 3\epsilon$$
and apply Lemma~\ref{lem:nlXsb} together with the trilinear estimates of  Lemma \ref{prop:tri1} and Lemma \ref{prop:tri2} to obtain 
\begin{equ}\label{fixed-point-2}\begin{aligned}
\|\cJ(\<1>v^2)\|_{\X^{\f 14+\epsilon,\f 12 +\delta}_T} \lesssim &\|\<1>v^2\|_{\X^{-\f 34+\epsilon,-\f 12 +\delta}_T} \le C_1T^\eta  \| \<1> \|_{L^\infty_TW^{-\alpha-\epsilon,\infty}_x} \|v\|^2_{\X^{\f 14+\epsilon,\f 12 +\delta}_T}\\
\|\cJ(v^3)\|_{\X^{\f 14+\epsilon,\f 12 +\delta}_T} \lesssim &\|v^3\|_{\X^{-\f 34+\epsilon,-\f 12 +\delta}_T} \le C_2 T^\eta \|v\|^3_{\X^{\f 14+\epsilon,\f 12 +\delta}_T}.
\end{aligned}
\end{equ}

When $\alpha \in (0, \f 18)$,  we take $\delta >0$ and $0<\epsilon < \f 38$, and invoke the bilinear estimate \eqref{eq:bi1} to handle the mixed term:
\begin{align}
\|\cJ(\<2>v)\|_{\X^{\f 14+\epsilon,\f 12 +\delta}_T}  \lesssim T^\delta \|\<2>v\|_{\X^{-\f 34+\epsilon,-\f 12 +2\delta}_T} \le  C_3T^\delta \|\<2>\|_{L^\infty_TW^{-2\alpha-\epsilon,\infty}_x} \|v\|_{\X^{\f 1 4+\epsilon, \f 12+\delta}}. \label{eq:cherry-v}
\end{align}
 
Combining these estimates, we set  $\theta =\mincurly{\eta,\delta}, C=3\maxcurly{C_1,C_2,C_3,1}$. This yields
$$\begin{aligned}
\|\Gamma_{\Theta}(v)\|_{\X^{\f 14+\epsilon,\f 12 +\delta}_T}
& \le CK +CT^\theta( \|v\|^3_{\X^{\f 14+\epsilon,\f 12 +\delta}_T}\\
&\quad +K\|v\|^2_{\X^{\f 14+\epsilon,\f 12 +\delta}_T}+K\|v\|_{\X^{\f 1 4+\epsilon, \f 12+\delta}}) 
\end{aligned}$$
Let $R = 2 C K$. For any $v\in \bar B_R$, we can choose $T$ sufficiently small so that $\Gamma_{\Theta}(v) \in \bar{B}_R$, showing that $\Gamma_{\Theta}$ maps $\bar B_R$ to itself.

Next we prove that  $\Gamma_\Theta$ is a contraction on $\bar{B}_R$.
The initial condition $(u_0,u_1)$ enters only through the linear term $S_t$,  and $\cJ(\<2>v)$ is linear. Thus the nonlinear part is  $$\tilde \Gamma_\theta (v)=-\cJ(\<1>v^2)+\cJ(v^3).$$
The standard multilinear estimate, as used above,  gives
\begin{equs}
\|\tilde \Gamma_{\Theta}(v) - \tilde \Gamma_{\Theta}(w)\|_{\X^{\f 14+\epsilon,\f 12 +\delta}_T} 
&\le C_1 T^\eta \| \<1> \|_{L^\infty_TW^{-\alpha-\epsilon,\infty}_x} \|v-w\|^2_{\X^{\f 14+\epsilon,\f 12 +\delta}_T}
+C_2 T^\eta \|v-w\|^3_{\X^{\f 14+\epsilon,\f 12 +\delta}_T}\\
&\lesssim (1+R)^2T^\eta  \|v-w\|_{\X^{\f 14+\epsilon,\f 12 +\delta}_T}.
\end{equs}
Since $v,w\in \Bar B_R$, the terms in parentheses are uniformly bounded by $ (1+R)^2$.

By the Banach fixed‑point theorem, we conclude that there exists $0<T' \le T \le 1$ and a unique solution $v$ on $[0, T']$ solving \eqref{eq:1st-v} in the class $$\X^{\f 14+\epsilon,\f 12 +\delta}_T \subset C([0, T']; H^{\f 14 +\epsilon}_x ).$$

Finally, the continuity of the solution with respect to the enhanced data ( $u_0,u_1,\<1>,\<2>,\<3>) $
follows from the bounds obtained above and the fact that the solution map $\Gamma_\theta$  depends linearly on each component of~$\Theta$.
\begin{equ}
\|u-v\|_{\X^{\f 14+\epsilon,\f 12 +\delta}_T} = \| \Gamma_\Theta(u) - \Gamma_{\overline{\Theta}}(v) \|_{\X^{\f 14+\epsilon,\f 12 +\delta}_T} \leq C \| \Theta - \overline{\Theta} \|_{\mathcal{Y}^{\alpha,\epsilon}_T},
\end{equ}
which completes the proof for the proposition.
\end{proof}

\subsection{Case $\f 1 8 \le \alpha <\f 14$ and random tensor}
To overcome the upper bound restriction for $\alpha>\f 18$, we will exploit the smoothing effect of the Duhamel operator and multi-linear effects, introducing the  random tensor operator $\fJ^{\<2>}$ with an appropriate operator norm adapted to the this setting.

Let $\L(B_1;B_2)$ denote the space of bounded linear operators between Banach spaces $B_1$ and $B_2$, and let $I\subset \R$ be a closed interval. For any $s_1,s_2 \in \R$, we define
\begin{equation}\label{eq:operator-norm}
\L^{s_1,s_2,b}_I = \bigcap_{I' \subset I} \L(\X^{s_1,b}_{I'};\X^{s_2,b}_{I'}),
\end{equation}
where the intersection runs over all nonempty closed sub-intervals $I'\subset I$. 

Let $\ell>0$ be a small parameter. Then $\L^{s_1,s_2,b}_I $ becomes a Banach space when equipped with the norm \begin{equs}\label{ro-norm}
	\norm{T}_{\L^{s_1,s_2;b}_I} = \sup_{I' \subset I} \f{\norm{T}_{\L(\X^{s_1,b}_{I'};\X^{s_2,b}_{I'})}}{|I'|^{\ell} }.
\end{equs}
where $|I'|$ denotes the length of the interval $I'$.

Th introduction of  he weight $|I'|^{\ell}$ in (\ref{ro-norm}) is intentional: this factor will later allow us to gain an extra smallness in key multilinear estimates, which is essential for applying a fixed‑point argument.

Define the enhanced initial data set space, in which $
\Theta = ((u_0,u_1),\<1>,\fJ^{\<2>},\<3>) $ takes values.

\begin{equ}
\mathcal{Y}^{\alpha,\epsilon,\delta}_T = \mathcal{H}^{\f 1 4+\epsilon} \times C([0,T];W^{-\alpha-\epsilon,\infty}_x) \times \L^{\frac{1}{4}+\epsilon, \frac{1}{4}+\epsilon,\frac{1}{2}+\delta}_T \times \X^{-3\alpha-\epsilon,\f 12+\delta}_T,
\end{equ}
with the associated norm 
$$ \|\Theta\|_{\mathcal{Y}^{\alpha,\epsilon}_T} = \|(u_0,u_1)\|_{\mathcal{H}^{\f 1 4+\epsilon}}+ \| \<1> \|_{L^\infty_TW^{-\alpha-\epsilon,\infty}_x}+ \|\fJ^{\<2>}\|_{\L^{\f 1 4+\epsilon, \f 1 4+\epsilon; \f 1 2+\delta}_T}+\|\<3>\|_{\X^{-3\alpha-\epsilon,\delta}_T}.$$

\begin{proposition}\label{prop:LWP-14}
Let $\f 1 8 \le  \alpha < \f 1 4$. Then there exist threshold $\epsilon_0>0$ and $\delta_0>0$ such that, for every $\epsilon > \epsilon_0$ and $\delta > \delta_0$, the remainder equation \eqref{eq:1st-v}  is locally well-posed in the space $$\X^{\f 1 4+\epsilon,\f 12 +\delta}_T.$$

More precisely,  there exist $\epsilon_0>0$ and $\delta_0>0$ with the following property.
for any initial data $(u_0,u_1) \in \mathcal{H}^{\f 1 4+\epsilon}$ where $\epsilon > \epsilon_0$, 
there exists a random time $0<T<1$ and a unique solution $v$ to \eqref{eq:1st-v} on $[0, T]$ in the class
\begin{equ}
\X^{\f 1 4+\epsilon,\f 12 +\delta}_T \subset C([0,T];H^{\f 1 4+\epsilon}_x)
\end{equ}
Furthermore,  for almost surely all $\omega$, the solution $v$ depends continuously on the enhanced data set
\begin{equ}
\Theta = (u_0,u_1,\<1>,\fJ^{\<2>},\<3>)
\end{equ}
where $\Theta$ range through the space $\mathcal{Y}^{\alpha,\epsilon,\delta}_T$.
\end{proposition}

\begin{proof}
We know that there exists $T\in (0,1)$ such that $K := \|\Theta\|_{\mathcal{Y}^{\alpha,\epsilon}_T}$ is finite almost surely.  Define the solution map \begin{equ}
\Gamma_{\Theta} (v)(t) = \mathbf S_t(u_0,u_1) -3 \cJ(\<3>) +3 \fJ^{\<2>}(v)-3\cJ(\<1>v^2)-3\cJ(v^3).
\end{equ}
By Proposition \ref{tensors}, the object $\fJ^{\<2>}$ defines a random operator, 
 and we have the estimate
\begin{equ}\label{eq:tensor-inequ}
\| \fJ^{\<2>} (v) \|_{\X^{\f 1 4+\epsilon, \f 12+\delta}_T} \le C_4 T^\ell \|\fJ^{\<2>}\|_{\L^{\f 1 4+\epsilon, \f 1 4+\epsilon; \f 1 2+\delta}_T} \| v\|_{\X^{\f 1 4+\epsilon, \f 12+\delta}_T} = C_4 T^\ell K \| v\|_{\X^{\f 1 4+\epsilon, \f 12+\delta}_T}.
\end{equ}

Combining the bound \eqref{eq:tensor-inequ} with the regularity estimate (\ref{chicken-fixed-1}) and the fixed‑point inequality (\ref{fixed-point-2}), we obtain the necessary contraction estimate for the solution map $\Gamma_{\Theta} (v)(t) $. The remaining steps — verifying the contraction mapping argument and the continuity with respect to the enhanced data — follow exactly as in the case $\alpha\in (0, \f 18)$.
\end{proof}

\subsection{Case $\f 1 4 \le \alpha <\f 3{10}$}
In the naiive counting of regularities, one postulates that the solution can be established only when $0<\alpha<\f 14$. However, the multilinear smoothing effect helps to improve the regularity of $\<30>$ from $1-3\alpha-\epsilon$ to $\f 54 -3\alpha-\epsilon$ for any $\epsilon >0$ when $\alpha \ge \f 14$. Consequently, even without introducing a second‑order expansion of the ansatz, we can still establish local well‑posedness of \eqref{eq:1st-v} up to $\alpha < \f 3{10}$, since 
\begin{equ}
\f 54 -3\alpha-\epsilon > 2\alpha -\f 14 +\epsilon,
\end{equ}
where the latter inequality matches the regularity requirement for the solution space $\X^{2\alpha-\f 14+\epsilon,b}$
and the threshold $2\alpha -\f 14 +\epsilon$ arises from applying the tensor estimate for $\fJ^{\<2>}$.

To state the theorem clearly, we define the solution map
\begin{equ}
\Gamma_{\Theta} (v)(t) =  \mathbf S_t(u_0,u_1) -3\<30> +3 \fJ^{\<2>}(v)-3\cJ(\<1>v^2)-3\cJ(v^3).
\end{equ}
and introduce the enhanced initial data set $$\Theta := ((u_0,u_1),\<1>,\fJ^{\<2>},\<30>)$$ which takes values in the product space
\begin{equ}
\mathcal{Y}^{s,\alpha,\epsilon,\delta}_T = \mathcal{H}^s \times C([0,T];W^{-\alpha-\f \epsilon 2,\infty}_x) \times \L^{2\alpha - \f  14+\epsilon, 2\alpha - \f  14+\epsilon,\frac{1}{2}+\delta}_T \times \X^{\f 54-3\alpha-\epsilon,\f 12+\delta}_T,
\end{equ}
 equipped with the natural product norm.

\begin{proposition}\label{prop:LWP-310}
Let $\f 1 4 \le  \alpha < \f 3 {10}$. Then there exist $\epsilon_0$ and $\delta_0>0$, such that
for any initial data $(u_0,u_1) \in \mathcal{H}^{2\alpha -\f 14+\epsilon}$ where $\epsilon>\epsilon_0$, and for any $\delta<\delta_0$, for almost surely all $\omega$, 
there exists a time $0<T(\omega)<1$ and a unique solution $v$ to~\eqref{eq:1st-v} on $[0, T]$ in the class
\begin{equ}
\X:=\X^{2\alpha - \f 14 +\epsilon,\f 12 +\delta}_T \subset C([0,T];H^{2\alpha - \f  14 +\epsilon}_x).
\end{equ}
Furthermore, the data-solution map 
$$\Theta \to \Gamma_\Theta $$  is continuous  from $ \mathcal{Y}^{s,\alpha,\epsilon,\delta}_T$ to $\X$.
\end{proposition}

\begin{proof}
By Proposition \ref{prop:reg-of-cube}, there exists a random time $T$ such that $\<30>_N\in\X^{\f 54-3\alpha-\epsilon,\f 12 +\delta}_T$. So for almost surely all $\omega$, there exists $T>0$ such that
$$ \| \<30>_N\|_{\X^{2\alpha -\f 14+\epsilon,b}_T}<\infty.$$
Let $\Omega_0$ be the set of full measure such that the random objects in $\Theta$ and $(u_0,u_1)$ are bounded in their respective distribution space. Let $T>0$ be such a random time and let $K$, which denotes the $\omega$,  denotes the bound.

Again we can ignore the linear term $\mathbf S(t)$. For $\fJ^{\<2>} $, we apply the second part of Proposition \ref{tensors}, which is applicable to the range of $\alpha\in [\f 14, \f 38)$. It shows that $\fJ^{\<2>} $ is a bounded random operator on $\X^{\f 1 4+\epsilon, \f 12+\delta}_T$
with finite $L_p$ bounded linear operator norm, c.f. Proposition \ref{prop:LWP-14}.

Progressing to the next term, let $\epsilon < \f 12 \big(\f 3 2 - 5\alpha \big)$ (i.e.
$2\alpha-\f 14 \ge \f 54-3\alpha$) and further assume that
$$0< \epsilon <  \f {(5-8\alpha)^2}{32(1-\alpha)}, \quad  \quad 0<\delta< \f {2\epsilon}{5-8\alpha}.$$

We apply Lemma \ref{lem:nlXsb}, followed by Proposition \ref{prop:tri3}, to obtain
$$\| \J(\<1>v^2)\|_{ \X^{2\alpha - \f 14 +\epsilon,\f 12 +\delta}}
\le \| \<1>v^2\|_{\X^{2\alpha - \f 14 +\epsilon-1, -\f 12 +\delta}_T} 
\lesssim T^{\eta} \| \<1> \|_{L^\infty_TW^{-\alpha-\f \epsilon 2,\infty}_x} (\|v\|_{\X})^2 .$$
For the triple $v$ term, we apply Lemma \ref{lem:nlXsb} and Proposition \ref{prop:tri4}:
\begin{equ}
\|J(v^3)\|_{\X^{2\alpha - \f 14 +\epsilon, \f 12 +\delta}_T} 
\le \|v^3\|_{\X^{2\alpha - \f 14 +\epsilon-1, -\f 12 +\delta}_T} \le T^\eta \bigl(\|v\|_{\X^{2\alpha - \f 14 +\epsilon, \f 12 +\delta}_T} \bigr)^3.
\end{equ}
 
For  the random tensor term, we apply Proposition  \ref{tensors}, 
 random tensor inequality \begin{equ}\label{eq:tensor-inequ}
\| \fJ^{\<2>} (v) \|_{\X^{2\alpha-\f 1 4+\epsilon, \f 12+\delta}_T} \lesssim  T^\ell \|\fJ^{\<2>}\|_{\L^{2\alpha-\f 1 4+\epsilon, 2\alpha-\f 1 4+\epsilon; \f 1 2+\delta}_T} \| v\|_{\X^{2\alpha-\f 1 4+\epsilon, \f 12+\delta}_T},
\end{equ}
which gives the restriction $s>2\alpha-\f 14$ when $\alpha \ge \f 14$. This completes the proof for the proposition.
\end{proof}

\subsection{Case $\f 3 {10} \le \alpha <\f 3 8$}
In this range of $\alpha$, as noted in \S\ref{sec:deduction}, we need the second order expansion. Hence, the limiting remainder equation becomes   
\begin{equ}\label{eq:2nd-v}
\begin{cases}
(\partial^2_t + 1 -\Delta) v = - v^3 + 3(\<30>-\<1>)v^2 - 3(\<30>^2-2\<30>\<1>)v-3\<2>v\\
\qquad \quad +\<30>^3-3\<30>^2\<1> + 3\<30>\<2>\\
(v,\partial_t v)|_{t=0} =(u_0,v_0)
\end{cases},
\end{equ}
for which the solution map is:
\begin{equ}\label{solution-map-3}
\begin{aligned}\Gamma_{\Theta} (v)(t) = & \mathbf S_t(u_0,u_1) -3 \cJ(v^3) +3\cJ(\<30>v^2)-3\cJ(\<1>v^2)-3\cJ(\<30>^2v)\\
&\quad +6\cJ(\<30>\<1>v)-3 \fJ^{\<2>}(v) + \cJ(\<30>^3)-3\<70>+3\<320>
\end{aligned}
\end{equ}
for which the enhanced initial data -parameter becomes
\begin{equ}
\Theta = ((u_0,u_1),\<1>,\fJ^{\<2>},\<30>,\<31>,\<320>,\<70>,\<90>).
\end{equ}
We recall that
 $$\<31>=\<30>\<1>, \quad \<320>=\J(\<30>\<2>), \quad \<70>=\J((\<30>)^2\<1>), \quad \<90>=\J((\<30>)^3).$$
 Note that the last term involves products of functions, the product in the first three terms are the form of 
product of a function with a distribution. After applying the Duhamel operator, the second and the third term take values in a function space.

  The initial data-parameter takes values in the following product space with the product norm
\begin{equs}
\mathcal{Y}^{s,\alpha,\epsilon,\delta}_T = &\mathcal{H}^s \times C([0,T];W^{-\alpha-\f \epsilon 2,\infty}_x) \times \L^{2\alpha - \f  14+\epsilon, 2\alpha - \f  14+\epsilon,\frac{1}{2}+\delta}_T \times \X^{\f 54-3\alpha-\epsilon,\f 12+\delta}_T \\
&\quad \times C([0,T];W^{-\alpha- \epsilon,\infty}_x) \times  \X^{\f 56-\alpha-\epsilon,\f 12+\delta}_T
\times  \X_T^{1-\alpha-\epsilon, \f 12+\delta}.
\end{equs}

\begin{proposition}\label{prop:LWP-38}
Let $\f 3 {10} \le  \alpha < \f 3 8$. Then there exist $\epsilon_0,\delta_0>0$, such that for any initial data $(u_0,u_1) \in \mathcal{H}^{2\alpha -\f 14+\epsilon}$ where $\epsilon>\epsilon_0$, and for any  $\delta > \delta_0$, there exist a random time $0<T<1$ and a unique solution $v$ to \eqref{eq:2nd-v} on $[0, T]$ in the class
\begin{equ}
\X^{2\alpha - \f 14 +\epsilon,\f 12 +\delta}_T \subset C([0,T];H^{2\alpha - \f  14 +\epsilon}_x).
\end{equ}
Furthermore, almost surely the map $\Theta\to \Gamma_\Theta$
is continuous from $\mathcal{Y}^{s,\alpha,\epsilon,\delta}_T$ to $\X^{2\alpha - \f 14 +\epsilon,\f 12 +\delta}_T$.
\end{proposition}

\begin{proof}
We emphasize that we treat the cubic product of the cubic wick product of $\<1>$ (cubic-cubic-cubic) as a new stochastic symbol $\<90>$. By Proposition \ref{prop:reg-cubic-cubic-cubic}, the symbol $\<90>$ has the regularity is $\f 94 -3\alpha -\epsilon$ for any $\epsilon>0$. This regularity is larger than $2\alpha-\f 1 4$.

Similarly, with the help of Proposition \ref{prop:reg-of-quartic}, the term $\cJ(\<30>\<1>v)$ is understood as $\cJ(\<31>v)$. Using the symbols,
\begin{equ}\label{solution-map-4}
\begin{aligned}\Gamma_{\Theta} (v)(t) = & \mathbf S_t(u_0,u_1) -3 \cJ(v^3) +3\cJ(\<30>v^2)-3\cJ(\<1>v^2)-3\cJ(\<30>^2v)\\
&\quad +6\cJ(\<31>v)-3 \fJ^{\<2>}(v) + \<90>-3\<70>+3\<320>
\end{aligned}
\end{equ}

Since the solution space here is the same as that in Proposition \ref{prop:LWP-310}, the terms 
$S_t$, $\cJ(v^3)$, and $\cJ(\<1>v^2)$ were already presented in the first order expansion, and shown to be in the solution space, with $L_p$ bounds.
We have already seen in Proposition~\ref{prop:LWP-310} that $\fJ^{\<2>} $ is a bounded random operator on $\X^{\f 1 4+\epsilon, \f 12+\delta}_T$ (the proof there is valid for the range $\alpha\in [\f 14, \f 38)$).

To the rest of the terms we apply the trilinear estimates from \S\ref{sec:deterministic-estimates} and tensor estimate from \S\ref{sec:deterministic-tensor}. 

By Proposition \ref{prop:reg-of-quartic}, $\<31>$ belongs to  $L^\infty_t W^{-\alpha-\epsilon,\infty}_x ([0,T] \times \T^2)$ almost surely, so we would need to take into account of the Duhamel operator before taking the product with $v$, and show that it is a bounded random operator on the solution space.  Indeed, by the bilinear estimate below, see \eqref{eq:bi-for-38}, we have
\begin{equ}
\|\cJ(\<31>v)\|_{\X^{2\alpha-\f 14 +\epsilon,\f 12+\delta}_T} \lesssim T^\delta \|\<31>v\|_{\X^{2\alpha -\f 54+\epsilon, -\f 12+2\delta}_T} \lesssim T^\delta \|\<31>\|_{L^\infty_TW^{-\alpha-\epsilon,\infty}_x} \|v\|_{\X^{2\alpha-\f 1 4+\epsilon, \f 12+\delta}_T}
\end{equ}
for any $\delta>0$.
Finally, 
the terms $\<70>$ and $\<320>$ are treated in Proposition \ref{prop:reg-septic} and Proposition \ref{prop:reg-of-quintic}.
This completes the proof that $\Gamma_\theta$ maps a ball $B_K$ in the initial data-parameter space into itself. It is now trivial to see that it is a contraction map, since by using the bounds on $v$, $\Gamma_\theta$ is controlled by a linear map.
We have completed the proof for the stated result.  \end{proof}

\section{Malliavin calculus on complex-valued symmetric sequences }
\label{sec:Wiener}

There are potential sources of confusion when defining multiple Wiener–Itô integrals associated with a sequence of complex-valued Brownian motions.
The purpose of this section is  precise and self-contained construction of real-valued multiple Wiener integrals on a Hilbert space of $L^2$ sequences satisfying a prescribed symmetry condition, with respect to a sequence of complex valued Brownian motions,  and to state the key properties that will be used throughout the article.

Typically, we construct the Gaussian space $\H$ from a sequence of independent standard  real-valued Gaussian random variables $\{W_n\}_{n=1}^\infty$. The space $\H$ is defined as the $L^2$-closure of the set of all finite linear combinations $\sum_{i=1}^m c_i W_{n_i}$ where $c_i\in \R$ and $n_i\in \N$. For each $n$, we define $\H_n$ (the space of polynomials of degree at most $n$) as the $L^2$-closure of the span of random variables of the form 
$$P(W_{n_1}, \dots, W_{n_k}),$$ where $P$ is a real polynomial of degree at most $n$ in $k$ variables, and $k, n_1, \dots, n_k$ are finite integers.

The confusion surrounding the Malliavin calculus associated with complex-valued the sequence 
 $\{\beta_n\}$ arises from two related issues:
\begin{itemize}
\item [(1)] there are two possible choices of inner product, and
\item[(2)] there are two ways to regard $V= \ell^2(L^2(\R;\C))$ as the complexification of real spaces.
\end{itemize}
The standard approach is to view $V$ as the complexification of $ \ell^2(L^2(\R;\R))$. However, there is also a non-standard complexification obtained by decomposing each element $\{v_n\}\in V$ as
$$v_n=\f {v_n+\bar v_{-n}} 2 +\i \f {v_n-\bar v_{-n}} {2i}.$$
In addition, one encounters a second, “natural” inner product on $V$ given by
$$\langle(a_n), (b_n)\rangle^{\tilde{}}=\sum_{n\in \Lambda} a_n b_{n}.
$$

{\bf Convention:}
 $\Lambda\subset \Z^2$ denotes a symmetric subset of $\Z^2$. The integration will be on  $[0,T] \subset\R$, an interval.

\subsection{ Real Hilbert space $\HH$ of complex-valued symmetric sequences}
Given a symmetric set $\Lambda\subset \Z^2$, and $T\subset \R$ an interval, we now define a real  Hilbert of functions with a built-in symmetry constraint. Consider
$$\{(h_n)_{n\in \Lambda} \;|\; h_n\in L^2([0,T];\C): h_{-n}=\bar h_n, n\in \Lambda\}.$$
with following real inner product:
\begin{equ}\label{inner-product}
\langle h,k\rangle_{\HH}=\sum_{n\in \Lambda}\int_0^T h_n(s)  k_{-n} (s)ds.
\end{equ}
To see this inner product always takes real values, we check that
$$\overline{\langle h,k\rangle}_{\HH}=\sum_{n\in \Lambda}\int_0^T  \bar h_n(s) \bar  k_{-n} (s)ds=\langle h,k\rangle_{\HH}.$$
$$\HH:= \{(h_n)_{n\in \Lambda} \;|\; h_n\in L^2([0,T];\C): h_{-n}
=\bar h_n, n\in \Lambda\, \quad \|h\|_{\HH}<\infty \}.$$
We further observe that  $\HH=L^2([0,T], H)$, where the “fiber space” $H$ is defined by
\begin{equ}\label{H0}
 H=\{(a_n)_{n\in \Lambda} \in \ell^2 \;|\; a_n\in \C,  a_{-n}=\bar a_n, \forall n \in \Lambda, \qquad \|h\|_H<\infty. \},
\end{equ}
The inner product on $H$ is :
$$\langle(a_n), (b_n)\rangle_{H}=\sum_{n\in \Lambda} a_n b_{-n}.$$
Consequently, for $h,k\in \HH$,
\begin{equ}
\label{inner-product-H}
\langle h,k\rangle_{\HH}=\int_T \langle h(s), k(s) \rangle_{H} ds=\sum_{n\in \Lambda} \int_T h_n(s) k_{-n}(s) ds,
\end{equ}
whenever the integral is well-defined.

We take the standard inner product on the tensor space $\otimes^m H_0$. Let, for $i=1, \dots, k$,  $a_i=\{a_i(n_i), n_i\in \Lambda\}\in \HH_0$ and $b_i=\{b_i(n_i), n_i\in \Lambda\}\in \HH_0$. Then, on elementary tensors, the inner product is given by
$$\langle \otimes^m a_i, \otimes^m b_i\rangle=\Pi_{i=1}^k  \langle a_i, b_i\rangle_{H}=\Pi_{i=1}^k \sum_{n_i\in \Lambda} a_i(n_i)b_i(-n_i).$$

To conclude, we define the Hilbert space on which the Wiener integration is defined to be:
\begin{equ}\label{Hilbert}
 \HH:= \left\{  (h_n(s))_{n \in \Lambda} : \hbox{ for a.e. $s$, }  \{h_n(s): n\in \Lambda\} \in H_0 , \text{ and } \int_T \|h(s)\|_{H_0}^2 \mathrm{d}s < \infty \right\}. 
\end{equ}

\subsection{Wiener process on $\HH$}
Let  $(\Omega, \mathcal{F}, \mathbb{P})$ be a probability space rich enough to support the following objects. Consider two independent families of two-sided real-valued standard Wiener processes
$\{W^{(n,1)}(t), t \ge 0\}_{n\in \mathbb{Z}^2}$ and $\{W^{(n,2)}(t), t \ge 0\}_{n\in \mathbb{Z}^2}$, each indexed by $n \in \mathbf{Z}^2$. We define a sequence of complex-valued processes $\{\beta_n(t), t \ge 0\}_{n \in \mathbf{Z}^2}$ by:
$$ \beta_n(t) = \frac{1}{\sqrt{2}} W^{(n,1)}(t) + \frac{\i}{\sqrt{2}} W^{(n,2)}(t). $$
Furthermore, we impose the following crucial symmetry condition relating the processes for $n$ and $-n$:
$$ W^{(-n,1)}(t) = W^{(n,1)}(t) \quad \text{and} \quad W^{(-n,2)}(t) = -W^{(n,2)}(t), \qquad  n \in \mathbf{Z}^2, t \ge 0. $$
This symmetry implies: $$\beta_{-n}(t) = \overline{\beta_n(t)} \qquad \hbox{ for all $n \in \mathbf{Z}^2$}.$$

Let $f\in L^2([0, t]; \mathbf{C})$ be a complex-valued deterministic function. For each $n \in \Z^2$ the Wiener integral of $f$ with respect to $\beta_n$ is defined via complex linearity:
$$ \int_0^t f(s) \mathrm{d} \beta_n(s) := \frac{1}{\sqrt{2}} \int_0^t f(s) \mathrm{d} W_s^{(n,1)} + \frac{\i}{\sqrt{2}} \int_0^t f(s) \mathrm{d} W_s^{(n,2)}. $$
Under the assumption $\beta_{-n} = \overline{\beta_n}$, the processes $\{\beta_n\}$ satisfy the following covariance relations (Itô isometry): for any $f, g \in L^2([0, t]; \mathbf{C})$ and any $m, n \in \mathbf{Z}^2$,
\begin{equs} \label{basic-isometry}
\mathbb{E} \left[ \int_0^t f(s) \mathrm{d} \beta_n(s) \overline{\int_0^t g(s) \mathrm{d} \beta_m(s)} \right] &= \delta_{n,m} \int_0^t f(s) \overline{g(s)} \mathrm{d}s, \\
\mathbb{E} \left[ \int_0^t f(s) \mathrm{d} \beta_n(s) \int_0^t g(s) \mathrm{d} \beta_m(s) \right] &= \delta_{n,-m} \int_0^t f(s) g(s) \mathrm{d}s.
\label{basic-isometry2}
\end{equs}
Note: $\delta_{n,m}$ is the Kronecker delta, equal to $1$ if $n=m$ and $0$ otherwise. The second identity follows directly from the symmetry $\beta_{-m} = \overline{\beta_m}$.

\subsection{Real-valued Wiener Integrals on $\HH$}

Taking the interval $[0,T]=\R$ for simplicity, it can be established verbatim on any interval.  We may now introduce the real-valued Wiener integral operator $\mathcal{W}: \HH \to L^2(\Omega)$:
$$ \mathcal{W}(h) := \sum_{n \in \Lambda} \int_\R h_n(s) \mathrm{d}\beta_n(s), \qquad h \in \HH. $$
By the definition of $\HH$, this sum converges in $L^2(\Omega)$. For any $h \in \HH$, $\mathcal{W}(h)$ is a real-valued Gaussian random variable. Indeed,
$$ \overline{\mathcal{W}(h)} = \sum_{n \in \Lambda} \int_\R \overline{h_n(s)} \mathrm{d}\overline{\beta_n(s)} = \sum_{n \in \Lambda} \int_\R h_{-n}(s) \mathrm{d}\beta_{-n}(s). $$
Since $\Lambda$ is symmetric, changing variable $m = -n$ yields that
  $$\mathcal{W}(h)=\sum_{m \in \Lambda} \int_\R h_{m}(s) \mathrm{d}\beta_{m}(s)$$
showing that $\mathcal{W}(h)$ is real.
\medskip

The operator $\mathcal{W}$ is an isometry between the real Hilbert space $\HH$ and $L^2(\Omega, \mathbf{R})$, Using (\ref{basic-isometry2}),
\begin{align*}
\mathbb{E}[\mathcal{W}(h) \mathcal{W}(k)] 
&= \sum_{n, m \in \Lambda} \mathbb{E} \left[ \int_\R h_n(s) \mathrm{d}\beta_n(s) \int_\R k_m(s) \mathrm{d}\beta_m(s) \right] \\
&= \sum_{n, m \in \Lambda} \delta_{n, -m} \int_\R h_n(s) k_m(s) \mathrm{d}s  \\
&= \sum_{n \in \Lambda} \int_\R h_n(s) k_{-n}(s) \mathrm{d}s = \langle h, k \rangle_{\HH}.
\end{align*}
Moreover, since $\mathcal{W}(k)$ is real, we also have 
\begin{equation}
\mathbb{E}[\mathcal{W}(h) \overline{\mathcal{W}(k)}] = \mathbb{E}[\mathcal{W}(h) \mathcal{W}(k)].
\end{equation}
 
 Finally, observe that the set of functions of the form
$$h: t\mapsto \sum_{i=1}^k  a^{(i)}\1_{A_i}(t), \qquad \qquad k\in \N, A_i \in \B(\R), a^{(i)}\in H_0 $$
 is dense in $\HH$. For such simple functions, the Wiener integral takes the explicit form::
$$I_1(h) := \mathcal{W}\Bigl( \sum_{i=1}^k  a^{(i)}\1_{A_i}(t) \Bigr)=\sum_{i=1}^k\sum_{n\in \Lambda} a^{(i)}_n \int_\R \1_{A_i}(s) \d \beta_n(s).$$

We now consider the $m$-th tensor product of $\HH$, denoted $\HH^{\otimes m}$. This is the real Hilbert space spanned by elementary tensors of the form
$$h^{(1)} \otimes \dots \otimes h^{(m)}, \qquad h^{(i)} \in \HH.$$
An element  $h\in\HH^{\otimes m}$ can be viewed as a collection of $L^2$- functions on $\R^m$,
indexed by multi-indices $(n_1, \dots, n_m)\in \Lambda^m$. We may write this either as 
$$\{ h_{n_1, \dots, n_k}(t_1, \dots, t_k): n_1, \dots, n_k \in \Lambda\}$$  
or,  more conveniently,  as
 $$\{ (h((n_1,t_1), \dots, (n_k,t_k)): n_1, \dots, n_k \in \Lambda\},$$
 with the symmetric restriction that the values are consistent when the multi-index is negated:
 $$(h((n_1,t_1), \dots, (n_k,t_k))=(\bar h((-n_1,t_1), \dots, \bar h(-n_k,t_k)).$$
 The space $\HH^{\otimes m}$ 
 can be identified with $$L^2(\mathbf{R}^m; \HH_0^{\otimes m}),$$
 where $\HH_0^{\otimes m}$ is the m-th tensor power of $H_0$.
  The inner product on the tensor space on $\otimes^k \HH$ is given by
\begin{equ}
\langle \otimes^m h_i, \otimes^m k_i\rangle_{\otimes \HH^k}=\sum_{n_1, \dots, n_m\in \Lambda} \int_\R \dots \int_\R h_{n_1,\dots, n_m}(s_1, \dots s_k)k_{-n_1,\dots, -n_m}(s_1, \dots s_k) \Pi_{i=1}^k ds_i.
\end{equ}

The space $\HH^{\otimes m}$ is generated by the linear combinations of simple functions of the form:
$$h(t_1,\dots, t_m)= \otimes_{k=1}^ m  a^{(k)}  \1_{A_k}(t_k), \qquad a^{(k)}\in \otimes^m H_0, A_n\in \B(\R) $$
where $\B(\R)$ denotes Borel subset of $\R$ and the set $\{A_n\}$ are assumed to be disjoint and identified by symmetry $A_{n}=A_{-n}$.  

For such simple functions, we define the multiple Wiener integral by
\begin{equ}\label{multiple-Wiener-Ito}
I_m(h)=  \Pi_{k=1}^m \W(a^{(k)}\1_{A_k}(t_i)),
\end{equ}
which extends by linearity and density to all of $\HH^{\otimes m}$. In particular for elementary tensors
$h=\otimes_{i=1}^mh^{(i)}$, this yields the natural identification
$$I_m(\otimes_{i=1}^mh^{(i)})=\mathcal W(h^{(1)} \cdots \mathcal W(h^{(m)})).$$

Let $\HH^{\odot m}$ denote the symmetric tensor product of $\HH$, which is a subspace of $\HH^{\otimes m}$. This space can be identified with $L^2_\text{sym}(\mathbf{R}^m; H_0^{\otimes m})$, the set of symmetric functions in $L^2(\mathbf{R}^m; H_0^{\otimes m})$ that satisfies the required symmetry condition. 

The symmetrization operator $\text{Sym}: \mathcal{H}^{\otimes m} \to \mathcal{H}^{\odot m}$ is defined by:
\begin{equ}
 \text{Sym}(f)((s_1, n_1), \dots, (s_m, n_m)) = \frac{1}{m!} \sum_{\sigma \in S_m} f((s_{\sigma(1)}, n_{\sigma(1)}), \dots, (s_{\sigma(m)}, n_{\sigma(m)})), 
\end{equ}
where $S_m$ is the permutation group of $\{1, \dots, m\}$. 

The multiple Wiener–Itô integrals are then defined for $f \in \HH^{\odot m}$
as mappings: $$I_m: \HH^{\odot m} \to L^2(\Omega, \mathbf{R})$$
and they satisfy several key structural properties, which we summarize below. 

\begin{lemma} \label{le:basic-wiener}
Let $m, k \in \mathbb{N}$.
\begin{itemize}
    \item[(1)] Linearity: $I_m$ is a linear map from $\HH^{\odot m}$ to $L^2(\Omega)$.
    \item[(2)] Symmetry: For any $f \in \HH^{\otimes m}$, we have $I_m(f) = I_m(\text{Sym}(f))$.
    \item[(3)] Isometry: For $f \in \HH^{\odot m}$ and $g \in \HH^{\odot k}$,
    \begin{equation} \label{isometry}
    \mathbb{E} [ I_{m}(f) I_{k}(g) ] = \delta_{m,k} m! \langle f, g \rangle_{\HH^{\otimes m}}.
    \end{equation}
    \item[(4)] Iterated Integral Representation: For $f \in \HH^{\odot m}$,
    \begin{equation} \label{iterated-integral-representation}
    I_m(f) = m! \sum_{n_1,\dots, n_m \in \Lambda} \int_0^T \int_0^{s_m} \dots \int_0^{s_2} f_{n_1, \dots, n_m}(s_1,\dots, s_m) \mathrm{d}\beta_{n_1}(s_1)\dots \mathrm{d}\beta_{n_m}(s_m)
    \end{equation}
    with integration over the simplex $0 < s_1 < s_2 < \dots < s_m < T$.
\end{itemize}
\end{lemma}

We also need the product formula, and for this we first define the notion of contraction.
\begin{definition}\label{def:tensor}
If $f \in \HH^{\odot m}$ and $g \in \HH^{\odot k}$, their contraction $f \otimes_r g \in \HH^{\otimes (m+k-2r)}$ for $r \in \{0, 1, \dots, \min(m, k)\}$ is defined by its components:
\begin{align*}
(f \otimes_r g)_{p_1, \dots, p_{m+k-2r}} &(t_1, \dots, t_{m+k-2r}) \\
= \sum_{q_1, \dots, q_r \in \Lambda} \int_{T^r} & f_{p_1, \dots, p_{m-r}, q_1, \dots, q_r}(t_1, \dots, t_{m-r}, u_1, \dots, u_r) \\
&\times g_{p_{m-r+1}, \dots, p_{m+k-2r}, -q_1, \dots, -q_r}(t_{m-r+1}, \dots, t_{m+k-2r}, u_1, \dots, u_r) \mathrm{d}u_1 \dots \mathrm{d}u_r.
\end{align*}
Note that $f \otimes_r g$ is not necessarily symmetric; in the product formula we always work with its symmetrization $\Sym(f \otimes_r g)$.
\end{definition}

\begin{lemma}[Product Formula]\label{product}
Let $f \in \HH^{\odot m}$ and $g \in \HH^{\odot k}$. Then
\begin{equation} \label{product-formula}
I_m(f) I_k(g) = \sum_{r=0}^{\min(m, k)} r! \binom{m}{r} \binom{k}{r} I_{m+k-2r}(\text{Sym}(f \otimes_r g)).
\end{equation}
In particular, for $k=1$, $g \in \HH$,
$$ I_m(f) I_1(g) = I_{m+1}(\text{Sym}(f \otimes g)) + m I_{m-1}(\text{Sym}(f \otimes_1 g)). $$
\end{lemma}

We next define multi-index Hermite polynomials. Let $\ell = (\ell_1, \dots, \ell_k) \in \mathbb{N}_0^k$ be a multi-index. For $x = (x_1, \dots, x_k) \in \mathbb{R}^k$, set
$$ H_\ell(x) = \prod_{j=1}^k H_{\ell_j}(x_j).$$

\begin{lemma}[Wick's Formula] \label{Wick}
Let $f_1, \dots, f_k \in \mathcal{H}$ be an orthonormal set, i.e., $\langle f_i, f_j \rangle_{\mathcal{H}} = \delta_{i,j}$. Let $\ell = (\ell_1, \dots, \ell_k) \in \mathbb{N}_0^k$ be a multi-index with total order $|\ell| = \sum \ell_j = m$. 
Define the tensor $$f^{\otimes \ell} = \text{Sym}(f_1^{\otimes \ell_1} \otimes \dots \otimes f_k^{\otimes \ell_k}) \in \HH^{\odot m}.$$ Then
$$ I_m(f^{\otimes \ell}) = H_\ell(\mathcal{W}(f_1), \dots, \mathcal{W}(f_k)). $$
\end{lemma}

Finally, we record a stochastic Fubini-type result, which is standard in the literature.
\begin{lemma}[Stochastic Fubini theorem] \label{commute}
Let $g: [0, T] \to \R$ be deterministic and bounded. Let $f: [0, T] \times [0,T]^k$ be such that for each $r\in [0,T]$, the map $r  \mapsto f(r, \cdot)$ takes values in $\HH^{\odot k}$.

Then
\begin{align*}
&\int_0^T g(r) I_k(f(r, \cdot)) \mathrm{d}r \\
&\quad = \int_0^T g(r) \left( k! \sum_{n_1,\dots, n_k} \int_0^T \dots \int_0^{s_2} f_{n_1,\dots,n_k}(r, s_1, \dots, s_k) \mathrm{d}\beta_{n_1}(s_1) \dots \mathrm{d}\beta_{n_k}(s_k) \right) \mathrm{d}r \\
&\quad = k! \sum_{n_1,\dots, n_k} \int_0^T \dots \int_0^{s_2} \left( \int_{s_k}^T g(r) f_{n_1,\dots,n_k}(r, s_1, \dots, s_k) \mathrm{d}r \right) \mathrm{d}\beta_{n_1}(s_1) \dots \mathrm{d}\beta_{n_k}(s_k).
\end{align*}
\end{lemma}
To justify the interchange of integration, it suffices to show that for any Brownian motion $B_t$
  and any adapted stochastic process $g(s,\cdot)$, the following identity holds:
  $$ \int_0^t \int_0^t \1_{[0,r]}(s) g(s,r) dB_s dr=\int_0^t \int_0^t \1_{[s,t]}(r)  g(s,r) dr dB_s. $$
This equality can first be verified for simple functions by exchanging summations, and then extended to the general case via limits.

\subsection{Chaos expansion omn $\HH$}
The image $\H=  \{\mathcal{W}(h) : h \in \HH \} $ is called the first Wiener chaos.
In general, let $(\Omega, \F, \P)$ be a given probability space and let $\H$ be a Hilbert sub-space of $L^2(\Omega, \F, \P)$ consisting of mean-zero Gaussian random variables. Let $\H_n$ denote the subspace of $L^2(\Omega, \F, \P)$ generated by polynomials of degree $n$ in variables from $\H$. Specifically, 
$$\H_n=\overline{ \Bigl\{P(\xi_1, \dots, \xi_k): k\in \mathbb N, \xi_1, \dots, \xi_k \in \H\Bigr\}}$$
where the closure is taken in $L^2$ $P$ is a polynomial in $k$ variables with real coefficients and of degree $n$. 

By an orthogonalization procedure, we obtain mutually orthogonal subspaces $\H^n$ as follows: Let 
$$\H^0=\H_0, \quad \H^1 = \H,$$ 
and and inductively define
$$\H^n=\H_n\cap (\H_0\oplus\dots \oplus \H^{n-1})^\perp$$
so that 
$$\H_n=\H^0\oplus\dots \oplus \H^{n-1}.$$

For example, $\H_0$ is the space of constant random variables. The space $\H_1 = \H_0 \oplus \H$ consists of random variables of the form $c + \xi$ where $c$ is a constant and $\xi \in \H$. For $n \ge 2$, $\H_n$ contains non-Gaussian random variables.

Let $$\mathcal{K} = \oplus_{n=0}^\infty \H^n$$ be the orthogonal direct sum.
A fundamental result — known as the Wiener Chaos Expansion or the Cameron–Martin Theorem — asserts that
$$L^2(\Omega, \sigma(\H), \P)= \oplus_{n=0}^\infty \H^n,$$ 
the closed subspace of $L^2$ consisting of all random variables measurable with respect to the 
$\sigma$-algebra generated by $\H$. Each $\H_n$ is referred to as the 
n-th Wiener chaos. In fact, $\H_n$ is the closed linear span of (appropriately normalized) Hermite polynomials of degree  $n$ applied to elements of $\H$.
\subsection{Hypercontractivity inequality}
Given the Gaussian Hilbert space $\H (=\H^1)$ and its chaos expansion, the following holds for any $p \ge 2$ and any $Q \in \H^k$:
$$ \norm{Q}_{L^p(\Omega)} \leq (p-1)^{\frac{k}{2}}\norm{Q}_{L^2(\Omega)} $$
This inequality is known as Nelson's hypercontractivity inequality.

For $Q \in \H_k$ (the inhomogeneous chaos space up to order $k$), write the orthogonal decomposition $$Q=\sum_{j=0}^k Q_j, \quad Q_j \in \H^j.$$
Assuming $Q$ has mean zero (so that $Q_0 = E[Q] = 0$), we obtain, for $p \ge 2$:
$$ \|Q\|_{L^p} \le \sum_{j=1}^k \|Q_j\|_{L^p} \le \sum_{j=1}^k (p-1)^{\frac{j}{2}} \|Q_j\|_{L^2} \le \sqrt{\sum_{j=1}^k (p-1)^j} \|Q\|_{L^2}, $$
where the last inequality follows from the Cauchy-Schwarz inequality and the orthogonality of the chaos decomposition, which implies $\|Q\|_{L^2}^2 = \sum_{j=1}^k \|Q_j\|_{L^2}^2$.

This framework applies naturally to multiple Wiener integrals, in the setting where the Gaussian Hilbert space $\H=\HH$, consisting of first-order Wiener integrals of symmetric kernels (see Section~\ref{sec:Wiener}). It also extends to finite complex linear combinations of products of complex Gaussian random variables, provided the resulting sum is real-valued. In this context, the Gaussian Hilbert space $\H$ can be realized as the space of first-order Wiener integrals $I_1(f)$ for $f$ in a suitable Hilbert space (e.g., $L^2(\mathbf{R}_+)$), and the $k$-th chaos $\H^k$ corresponds to $k$-th order multiple Wiener-Itô integrals $I_k(f_k)$ with symmetric kernels $f_k$. The result also applies to finite complex linear combinations of products of complex Gaussian random variables, when the sum is real valued.

We can apply this result to multiples Wiener integrals in the setting where the Gaussian Hilbert space are Wiener integrals of symmetric sequences, \S\ref{sec:Wiener}. It also applies to finite complex linear combinations of products of complex Gaussian random variables, when the sum is real valued. In this context, the Gaussian Hilbert space $\H$ can be realized as the space of first-order Wiener integrals $I_1(f)$ for $f$ in a suitable Hilbert space (e.g., $L^2(\mathbf{R}_+)$), and the $k$-th chaos $\H^k$ corresponds to $k$-th order multiple Wiener-Itô integrals $I_k(f_k)$ with symmetric kernels $f_k$.

\section{Convergence of stochastic symbols}\label{sec:stochastic-symbols}

To study the regularity of the stochastic integrals, we first represent them in terms of Wiener integrals.
Let $\Lambda:=\Z_N^2=\{n\in \Z^2: |n|\le N\}$ be the symmetric lattice set.

\subsection{Stochastic symbols as iterated Wiener integrals on $\HH$}

  For every $t, x, N$, define, for $n\in \Z^2$,
\begin{equ}
h^{(1),t,x,N}_n(s) := e^{\i n \cdot x}\cdot \1_{[0,t]}(s) \1_{\Z^2_N}(n) \f {\sin((t-s) \jpb{n})}{\jpb{n}^{1-\alpha}}.
\end{equ}
where $\jpb{n}=\sqrt{1+|n|^2}$. 
Then we can write \begin{equ}
\<1>_N(t)  = I_1\Bigl(h^{(1),t,x,N}\Bigr)= \sum_{n \in \Z^2} e^{\i n \cdot x} \1_{\Z^2_N}(n) \int_0^\infty \f {\sin((t-s) \jpb{n})}{\jpb{n}^{1-\alpha}} \, \d \beta^n_s, \label{eq:lolipop-w}
\end{equ}
where $h^{(1),t,x,N}=(h^{(1),t,x,N}_n, n \in \Z^2_N)\in \HH$.

The wick products $\<2>_N(t)$ and $\<3>_N(t)$ are then given by Hermite polynomial via Wick's formula:
$$\<2>_N(t) = I_2[\otimes^2 h^{(1),t,x,N}], \qquad \<3>_N(t) = I_3[\otimes^3 h^{(1),t,x,N}].$$
We now denote:
$$\left\{
\begin{aligned}
h^{(2),t,x,N} &= \otimes^2 h^{(1),t,x,N}, \qquad  h^{(3),t,x,N} = \otimes^3 h^{(1),t,x,N},\\
h^{(2),t,x,N}_{n_1,n_2}(s_1,s_2) &= e^{\i (n_1+n_2)\cdot x} \cdot \prod_{j=1}^2 \1_{[0,t]}(s_j) \1_{\Z^2_N}(n_j) \f {\sin((t-s_j) \jpb{n_j})}{\jpb{n_j}^{1-\alpha}},\\
 h^{(3),t,x,N}_{n_1,n_2,n_3}(s_1,s_2,s_3) &= e^{\i (n_1+n_2+n_3)\cdot x} \cdot \prod_{j=1}^3 \1_{[0,t]}(s_j) \1_{\Z^2_N}(n_j) \f {\sin((t-s_j) \jpb{n_j})}{\jpb{n_j}^{1-\alpha}}.
\end{aligned}\right.
$$
 By Lemma~\ref{le:basic-wiener}, these quantities can be expressed in terms of iterated Wiener–Itô integrals. Specifically, 
\begin{equ}
\label{eq:cherry-w}
\left\{\begin{aligned}
\<2>_N(t) &:= \sum_{n_1,n_2 \in \Z^2} e^{\i (n_1+n_2)\cdot x} \int_{[0,\infty)^2} \prod_{j=1}^2 \1_{[0,t]}(s_j) \1_{\Z^2_N}(n_j) \f {\sin((t-s_j) \jpb{n_j})}{\jpb{n_j}^{1-\alpha}} \, \d \beta^{n_1}_{s_1} \d \beta^{n_2}_{s_2},  
\\
 \<3>_N(t) &:= \sum_{n_1,n_2,n_3 \in \Z^2} \!\!e^{\i (n_1+n_2+n_3)\cdot x} 
 \int_{[0,\infty)^3} \! \prod_{j=1}^3 \1_{[0,t]}(s_j) \1_{\Z^2_N}(n_j)  \f {\sin((t-s_j) \jpb{n_j})}{\jpb{n_j}^{1-\alpha}}\, \d \beta^{n_1}_{s_1} \d \beta^{n_2}_{s_2}   \d \beta^{n_3}_{s_3}.
\end{aligned}\right.
\end{equ}

To simplify the notation, we  define
\begin{equ}
\left\{\begin{aligned}
h^{(1),t}_n &=  \1_{[0,t]}(s) \f {\sin((t-s) \jpb{n})}{\jpb{n}^{1-\alpha}}, \\ 
h^{(2),t}_{n_1,n_2} &= \otimes^2 h^{(1),t}_{\cdot} = \prod_{j=1}^2 \1_{[0,t]}(s_j) \f {\sin((t-s_j) \jpb{n_j})}{\jpb{n_j}^{1-\alpha}}, \\
h^{(3),t}_{n_1,n_2,n_3} &= \otimes^3 h^{(1),t}_{\cdot} = \prod_{j=1}^3 \1_{[0,t]}(s_j) \f {\sin((t-s_j) \jpb{n_j})}{\jpb{n_j}^{1-\alpha}},\\
h^{(3),t,\cJ}_{n_1,n_2,n_3}&:= \CJ_t(h^{(3),\cdot}_{n_1,n_2,n_3})=\int_0^t \f {\sin((t-s)\jpb{n_{123}})}{\jpb{n_{123}}} h^{(3),s}_{n_1,n_2,n_3}(s_1,s_2,s_3) \, \d s,
\end{aligned}\right.
\label{define-h}
\end{equ}
The Duhamel operator acting on symmetric sequences from $\HH$ by its action on each component, thus as a sequence, $$h^{(3),t,\cJ}=\CJ_t(h^{(3),\cdot}).$$

Recall that we use $I$ for the Wiener integral of sequences. We shall use the standard notation, $\I$,  for iterated Wiener integral for a real or complex valued function $f$:
\begin{equ}\label{wiener-k}
\I_k[ f]: = \int_{[0,\infty)^k} f(s_1,\cdots,s_k) \, \d \beta^{n_1}_{s_1} \cdots \d \beta^{n_k}_{s_k},
\end{equ}
Then we can write,  where we omit the sup-script on the first order tensor $h^{t,x,N}=h^{(1),t, x,N}$ whenever it does not cause confusion,
\begin{equs}
\left\{\begin{aligned}
\<1>_N(t) (x)&=I_1(h^{t, x,N})= \sum_{n \in \Z^2_N} e^{\i n \cdot x} \I_1[h^{(1),t}_n],\\
 \<2>_N(t)(x) &= I_2(h^{(2),t, x})=\sum_{n_1,n_2 \in \Z^2_N} e^{\i (n_1+n_2)\cdot x}  \I_2[h^{(2),t}_{n_1,n_2}] \\
\<3>_N(t)(x) &= I_3[h^{(3),t,x,N}]= \sum_{n_1,n_2,n_3 \in \Z^2_N} e^{\i (n_1+n_2+n_3)\cdot x} \I_3[h^{(3),t}_{n_1,n_2,n_3}].
\end{aligned}\right.
\label{eq:lolipop-cherry-cube}
\end{equs}

We can now use the stochastic Fubini theorem, Proposition \ref{commute}, to represent $\<30>_N:=\cJ(\<3>)$ as follows:
\begin{align}
\<30>_N(t) &= \sum_{n_1,n_2,n_3 \in \Z^2_N} e^{\i (n_1+n_2+n_3)\cdot x} \cJ\Big( \I_3[h^{(3),s}_{n_1,n_2,n_3}] \Big)(t) \\
&=\sum_{n_1,n_2,n_3 \in \Z^2_N} e^{\i (n_1+n_2+n_3)\cdot x} \I_3[h^{(3),t,\cJ}_{n_1,n_2,n_3}].\label{eq:chicken-foot}
\end{align}
$$=\sum_n e^{\i n\cdot x} \sum_{n_1,n_2,n_3 : n_1+n_2+n_3=n}  \I_3[h^{(3),t,\cJ}_{n_1,n_2,n_3}].$$

The regularities of $\<1>_N$ and $\<2>_N$ have been proved in \cite{OO21} which is stated as following:
\begin{lemma}\label{reg-cherry}
We have
\begin{itemize}
\item Given $\alpha >0$, for any $s<-\alpha$, $\{ \<1>_N\}_{N \in \N}$ is a Cauchy sequence in $C([0,T];W^{s,\infty}_x)$ almost surely. In particular, denoting the limit by $\<1>$, we have
\begin{equ}
\<1> \in C([0,T]; W^{-\alpha-\epsilon,\infty}_x)
\end{equ}
for any $\epsilon>0$ almost surely. 
\item Given $0<\alpha <\f12$, for any $s<-2\alpha$, $\{ \<2>_N\}_{N \in \N}$ is a Cauchy sequence in $C([0,T];W^{s,\infty}_x)$ almost surely. In particular, denoting the limit by $\<2>$, we have
\begin{equ}
\<2> \in C([0,T]; W^{-2\alpha-\epsilon,\infty}_x)
\end{equ}
for any $\epsilon>0$ almost surely. 
\end{itemize}
\end{lemma}

\subsection{Elementary estimates on kernels}
\begin{lemma}\label{lem:symbols-in-Xsb}
Let $k\in N$, $J\subset \R$ an interval, and $\Lambda\subset \Z^2$ a symmetric set. 
Define the symmetric subsets 
$$\Lambda_n=\{ (n_1, \dots, n_k):  n_1+\dots+n_k=n, n_j\in \Lambda  \}.$$
For each $t\in J$ fixed , consider a family of tensors in $ \otimes^k\HH$: 
$G(t,x)\in L^2([0,T]^k; \otimes^kH)$,  where $G=\{G_{n_1, \dots, n_k}(t,x)\}$ with
$$G_{n_1, \dots, n_k}(t,x)=e^{i(n_1+\dots +n_k)\cdot x} g(t)_{n_1, \dots, n_k}.$$
Then we have: \begin{equs} \label{Xsb-stochastic}
{}&\Bigl\| \|I_k(G)\|_{\X_J^{s-1,b-1}}\Bigr \|_{L^p(\Omega)} \\
 &\lesssim p^{\f k2} \|  \jpb{\lambda}^{b-1} \jpb{n}^{s-1} 
\Bigl\| \Bigl( \int_{[0, \infty)^k} \sum_{(n_1, \dots, n_k)\in\Lambda_n}
\Bigl | \F_t (\1_Jg) (\cdot)_{n_1, \dots, n_k}(\lambda-\sigma \jpb{n}, n)\Bigr|^2 \d s\Bigr)^{\f 12}\Bigr\|_{\ell^2_nL^2_\lambda}
\end{equs}
where $\d s=\Pi_{j=1}^k ds_j$. Note that $$I_k[G(t,x)]=\sum_n e^{i n\cdot x}\I_k[\sum_{ (n_1, \dots, n_k)\in \Lambda_n}  g(t)_{n_1, \dots, n_k} ],$$
where the first integral $I_k$ indicates the $k$ iterated Wiener integral over the symmetric sequences, and the second integral $\I_k$ denotes Wiener integral over complex number valued functions. 

Similarly, for fixed $t \in [0,T]$, we have
\begin{equ}\label{sobolev-bound}
\Bigl\| \|I_k[G]\|_{H_x^{s-1}}\Bigr \|_{L^p(\Omega)} \lesssim p^{\f k2} \|  \jpb{m}^{s-1} 
 \Bigl\|   \sum_{(m_1, \dots, m_k)\in \Lambda_m}\Bigl( \int_{[0, \infty)^k} |g_{m_1,\dots, m_k}(\cdot)|^2 \d s\Bigr)^{\f 12}\Bigr\|_{\ell^2_m}
\end{equ}
\end{lemma}

\begin{proof}
Consider the component of the space-time Fourier transform of $G$ with index $m$: 
\begin{equs}
\widetilde{G}(\lambda,m) 
&=  \cF_t (\widehat{G}(\cdot,m))
=\F_t\Bigl( I_k[ g_m(t) \1_{J}(t)]\Bigr)
=I_k\Bigl[  \F_t(g_m(t) \1_{J}(t))\Bigr].
\end{equs}
where the second equality follows from the stochastic Fubini theorem, c.f. Lemma \ref{commute}.

For $\sigma \in \{ \pm 1\}$, by a change of variable, we have \begin{equs}
\Bigl\| \|G\|_{\X_J^{s-1,b-1}} \Bigr\|_{L^p(\Omega)} 
&\lesssim \Bigl\| \|\jpb{\lambda}^{b-1} \jpb{m}^{s-1}\widetilde {G\1_J} (\lambda - \sigma \jpb{m}, m) \|_{\ell^2_m L^2_\lambda }\Bigr \|_{L^p(\Omega)} \\
&=  \Bigl\| \|\jpb{\lambda}^{b-1} \jpb{m}^{s-1}  I_k\Bigl[\F_t \bigl(g_m(\cdot)\1_{J}(\cdot)\bigr)(\lambda - \sigma \jpb{m}, m) \Bigr] \|_{\ell^2_mL^2_\lambda }\Bigr \|_{L^p(\Omega)} \\
\end{equs}
We applied Minkowski's integral inequality to switch the $L^p(\Omega)$ and $\ell^2_nL^2_\lambda$-norms
\begin{equs}
\Bigl\| \|G\|_{\X^{s-1,b-1}} \Bigr\|_{L^p(\Omega)} 
& \lesssim \Bigl \|  \jpb{\lambda}^{b-1}  \jpb{n}^{s-1} 
I_k\Bigl[\F_t \bigl(g_m(\cdot)\1_{J}(\cdot)\bigr)(\lambda - \sigma \jpb{m}, m) \Bigr]
  \|_{L^p(\Omega) } \Bigr\|_{\ell^2_mL^2_\lambda} \\
 &\lesssim  p^{\f k2} \|  \jpb{\lambda}^{b-1} \jpb{m}^{s-1}  I_k[\F_t \bigl(g_m(\cdot)\1_{J}(\cdot)\bigr)(\lambda - \sigma \jpb{m}, m) ]\|_{L^2(\Omega)} \|_{\ell^2_mL^2_\lambda} \\
 &=p^{\f k2} \|  \jpb{\lambda}^{b-1} \jpb{m}^{s-1} 
 \Bigl( \int_{[0, \infty)^k} |\F_t \bigl(g_m(\cdot)\1_{J}(\cdot)\bigr)(\lambda - \sigma \jpb{m}, m) ]|_{\otimes H^k}^2 \d s\Bigr)^{\f 12}\|_{\ell^2_mL^2_\lambda} ,
\end{equs}
where  in the penultimate line we used the Gaussian hyper-contractivity, followed by It\^{o} isometry (Lemma \ref{le:basic-wiener}).  This concludes the first inequality. Finally for $t \in [0,T]$,
\begin{align}
\| \|G\|_{H^s_x} \|_{L^p(\Omega)} 
&\lesssim p^{k/2}  \Bigl\|  \|\jpb{m}^{s} I_k[g_m(t)] \|_{L^2(\Omega)} \Bigr\|_{\ell^2_m} = \|\jpb{m}^{s}  \Bigl( \int_{[0, \infty)^k} |g_m(\cdot)|_{\otimes H^k}^2 \d s\Bigr)^{\f 12}\|_{\ell^2_m}.
\end{align}
This complete the proof.
\end{proof}

\textbf{Convention:}\, In the following subsections, we set $0 < T \le 1$ and $p \ge 2$.

\subsection{The cubic symbols $\<3>$ and $\<30>$}\label{sec-symbols}

Let $k\in \mathbf{N}$, we define  \begin{equs}
h^{(k),t}_{n_1,\cdots, n_k} & =  \prod_{j=1}^k \1_{[0,t]}(s_j)   \f {\sin((t-s_j)\jpb{n_j})}{\jpb{n_j}^{1-\alpha}},\\
h^{(k),t,x,N}_{n_1,\cdots,n_k} &= e^{\i  \sum_{j=1}^k n_j \cdot x} \cdot \prod_{j=1}^k  \1_{\Z^2_N}(n_j)
h^{(k),t}_{n_1,\cdots, n_k}, \\
h^{(k),t,x,N} &= (h^{(k),t,x,N}_{n_1,\cdots,n_k}).
\end{equs}

By Euler identity $\sin \theta =\f {e^{\i \theta} -e^{-\i \theta}}{2\i}$,  we have:
\begin{lemma}\label{Euler}
For any $k \ge 1$, we have
\begin{equ}\label{eq:sin-to-e}
\prod_{j=1}^k \sin((t-s_j)\jpb{n_j}) = \f 1{(2\i)^k} \sum_{\substack{ \sigma_j \in \{\pm 1\}, j=1, \dots, k}}  \Pi_{j=1}^k \sigma_j  \; e^{\i t  \sum_{\ell=1}^k \sigma_\ell \jpb{n_\ell}} e^{\i  -\sum_{\ell =1}^k s_\ell \sigma_\ell \jpb{n_\ell}}.
\end{equ}
\end{lemma}

\begin{proposition}\label{prop:reg-of-cube}
Let $s_\alpha$ be defined by~(\ref{eq:s-alpha}). Then, for any  $s <s_\alpha$, there exists $b_0>\f 12$ ( depending on $(\alpha, s)$) such that for all $b>b_0$, \begin{equ}
\Bigl\| \| \<3>_N\|_{\X^{s-1,b-1}_T} \Bigr\|_{L^p(\Omega)}  \lesssim p^{\f 32} T^{\f 32}.
\end{equ}
Moreover, the sequence $\{ \<3>_N \}$ is Cauchy in $L^p(\Omega; \X^{s-1,b-1}_T)$.

In addition, by Duhamel's formula, \begin{equ}
 \Bigl \| \| \<30>_N\|_{\X^{s,b}_T}\Bigr\|_{L^p(\Omega)} \lesssim p^{\f32} T^{\f 32}.
\end{equ}
for the same values of $s$ and $b$ as above. Furthermore, the sequence $\{ \<30>_N \}$ is Cauchy in $L^p(\Omega; \X^{s,b}_T)$.
\end{proposition}

\begin{proof} 
As an application of Lemma \ref{lem:nlXsb}, the bound for $\<30>_N$ follows immediately from the corresponding result for $\<3>_N$.  It therefore remains to establish the upper bound for
$\| \|\<3>_N\|_{\X^{s-1,b-1}_T} \|_{L^p(\Omega)}$. 

We begin by applying the Littlewood–Paley decomposition to $\<3>_N$ with respect to the frequencies $n_1,n_2,n_3\in \Z^2_N$:
\begin{equ}
\<3>_N(t,x) = \sum_{N_1,N_2,N_3 \ge 1} \<3>^{\mathrm{dyadic}}_{N,N_1,N_2,N_3}(t,x).
\end{equ}
where
\begin{equ}
\widehat{\<3>^{\mathrm{dyadic}}}_{N,N_1,N_2,N_3}(t,n)
 =  \sum_{\substack{n=n_1+n_2+n_3\\n_1,n_2,n_3 \in\Z^2_N}} \prod_{j=1}^3 \1_{|n_j| \sim N_j} \widehat{\<3>}_N(t,n).
 \end{equ}
We analyze these contributions in four steps.

{\bf Step 1.} 
We begin by expanding the norm  
$\| \| \<3>^{\mathrm{dyadic}}_{N,N_1,N_2,N_3}\|_{\X^{s-1,b-1}_T}  \|_{L^p(\Omega)}$-norm. 

By the definition (see \eqref{eq:lolipop-cherry-cube}), $
\widehat{\<3>}_N(t,n) = \sum_{\substack{n_{123}=n \\n_1,n_2,n_3 \in\Z^2_N}} \I_3[h^{(3),t}_{n_1,n_2,n_3}]$, and therefore:
\begin{equs}
\widehat{\<3>^{\mathrm{dyadic}}}_{N,N_1,N_2,N_3}(t,n)
 &=\sum_{\substack{n=n_{123}\\n_1,n_2,n_3 \in\Z^2_N}} \prod_{j=1}^3 \1_{|n_j| \sim N_j} 
 \I_3[h^{(3),t}_{n_1,n_2,n_3}](t).
  \label{dyadic-cube}
\end{equs}

Moreover, for $\sigma \in \{ \pm \}$, 
 \begin{equs}
\widetilde{\1_{[0,T]}(\cdot)\<3>^{\mathrm{dyadic}}_N}(\lambda - \sigma \jpb{n},n) 
= &  \sum_{\substack{n=n_{123}\\n_1,n_2,n_3 \in\Z^2_N}}  \prod_{j=1}^3 \1_{|n_j| \sim N_j}   \I_3 \bigl( \cF_t(\1_{[0,T]}h^{(3),t}_{n_1,n_2,n_3})\bigr)(\lambda - \sigma \jpb{n},n),
\end{equs}
By Lemma \ref{lem:symbols-in-Xsb},

\begin{equs} 
{}&\| \| \<3>^{\mathrm{dyadic}}_{N,N_1,N_2,N_3}\|_{\X^{s-1,b-1}_T}  \|_{L^p(\Omega)}\\
 &\lesssim p^{\f k2} \|  \jpb{\lambda}^{b-1} \jpb{n}^{s-1} 
\Bigl\| \bigl( \int_{[0, \infty)^k}  \sum_{\substack{n=n_{123}\\n_1,n_2,n_3 \in\Z^2_N}}  \prod_{j=1}^3 \1_{|n_j| \sim N_j} 
\Bigl | \F_t (\1_{[0,T]}h^{(3),t}_{n_1,n_2,n_3})(\lambda-\sigma \jpb{n}, n)\Bigr|^2 \d s\bigr)^{\f 12}\Bigr\|_{\ell^2_nL^2_\lambda}.
\end{equs}

By \eqref{eq:sin-to-e} and the explicit definition of $h^{(3),t}_{n_1,n_2,n_3}$ from \eqref{define-h}, we see:
\begin{equs}
 \cF_t(\1_{[0,T]}(t)h^{(3),t}_{n_1,n_2,n_3}) 
=  \f 1{(2\i)^3} \sum_{ \sigma_1,\cdots,\sigma_3 \in \{\pm 1\} } (\sigma_1 \sigma_2 \sigma_3) e^{\i \psi(\boldsymbol{s},\bn)}  \int_{\maxcurly{s_1,s_2,s_3}}^T   e^{-\i t (\lambda + \phi(\bn))} \, \d t.
\end{equs}
where  $\phi$ and $\psi$ are the phase functions:
\begin{equs}
\psi(\boldsymbol{s},\bn) =-\sum_{\ell =1}^3 s_\ell \sigma_\ell \jpb{n_\ell}, \qquad \phi(\bn)  = - \sigma \jpb{n} - \sum_{\ell =1}^3 \sigma_\ell \jpb{n_\ell},
\end{equs}
with $\bn=(n_1,n_2,n_3)$, $\boldsymbol{s}=(s_1,s_2, s_3)$.

Integrating in $t$ first and then integrating in $s_1,s_2,s_3$, we arrive at the simplified estimate: 
\begin{equs} 
&{}\int_{[0,T]^3} \Big| e^{\i \psi(\boldsymbol{s},\bn)}  \int_{\maxcurly{s_1,s_2,s_3}}^T   e^{-\i t (\lambda + \phi(\boldsymbol{n}) )} \, \d t \Big|^2 \, \d s_1\d s_2 \d s_3 \d \lambda\\
&\le \int_{[0,T]^3} \f 2{\jpb{\lambda +\phi(\boldsymbol{n})}} \, \d s_1 \d s_2 \d s_3 \\
& \lesssim T^3 \f 2{\jpb{\lambda +\phi(\boldsymbol{n})}}.
\end{equs}

Hence, we have
\begin{equ}\label{eq:stochasti-symbol-in-Xsb}
\begin{aligned}
 \| \| \<3>^{\mathrm{dyadic}}_{N,N_1,N_2,N_3}  \|_{\X^{s-1,b-1}_T} \|^2_{L^p(\Omega)} 
 &\lesssim p^3T^3 
  \sum_{n_1,n_2,n_3 \in \Z^2_N} \jpb{n_{123}}^{2s-2}  \\
 &\quad \prod_{j=1}^3 \1_{|n_j| \sim N_j} \jpb{n_j}^{-2+2\alpha}  \;  \int_{\R} \f {\jpb{\lambda}^{2b-2}}{\jpb{\lambda +\phi(\bn)}^2} \, \d \lambda. 
  \end{aligned}
 \end{equ}
Let $\beta<-1$. Then uniformly over $c\in \R$,
\begin{equ}\label{basic-integration}
\int_\R   \f{\jpb{t}^{\beta}}{ \jpb{t+c}^{2}}dt
\lesssim \jpb{c}^{\beta} 
\end{equ}
 This estimate is easily seen by splitting the integration region. By symmetry, assume $c>0$:
\begin{equs}
\int_{{|t|}\ge 3c}  \f{\jpb{t}^{\beta}}{ \jpb{t+c}^{2}}dt \lesssim \jpb{c}^{\beta-1}, \qquad
\int_{{|t|}\le \f c 2}  \f{\jpb{t}^{\beta}}{ \jpb{t+c}^{2}}dt \lesssim \jpb{c}^{\beta}, \quad 
\int_{\f c 2 \le {|t|}\le 3c}  \f{\jpb{t}^{\beta}}{ \jpb{t+c}^{2}}dt \lesssim \jpb{c}^{\beta-1}.
\end{equs}
 We deduce that
$$ \int_{\R} \f {\jpb{\lambda}^{2b-2}}{\jpb{\lambda +\phi(\bn)}^2} \, \d \lambda 
\lesssim \jpb{\phi(\bn)}^{2b-2}.$$
and hence
\begin{equ}\label{eq:stochasti-symbol-in-Xsb-2}
\begin{aligned}
&{} \| \| \<3>^{\mathrm{dyadic}}_{N,N_1,N_2,N_3}  \|_{\X^{s-1,b-1}_T} \|^2_{L^p(\Omega)} \\
 &\lesssim p^3T^3 
  \sum_{n_1,n_2,n_3 \in \Z^2_N} \jpb{n_{123}}^{2s-2}  
 \prod_{j=1}^3 \1_{|n_j| \sim N_j} \jpb{n_j}^{-2+2\alpha}  \;  \jpb{\phi(\bn)}^{2b-2}. \end{aligned}
 \end{equ}
It remains to show that
$$ \sum_{N_1,N_2,N_3 \ge 1}   \sum_{n_1,n_2,n_3 \in \Z^2_N} \jpb{n_{123}}^{2s-2}  
 \prod_{j=1}^3 \1_{|n_j| \sim N_j} \jpb{n_j}^{-2+2\alpha}  \;  \jpb{\phi(\bn)}^{2b-2}<\infty.$$

{\bf Step 2.}
We deduce the case where $b<\f 12$ and the case where $b=1$,  which by interpolation allows us to conclude for any $b>\f 12$. 

\begin{itemize} \item [\bf (i) ]  
We fix any $b_1<\f12 $ and let $s_1$ be any real number and denoting $B=\{m: |m-\phi(\bn)|\le 1\}$,
\begin{equ}\label{estimate-phase}
|\phi(\bn)|^{2b_1-2}
\lesssim   \sum_{m \in \Z} \f{\1_B(m)} {\jpb{m}^{2-2b_1}}
\lesssim \sup_{m\in \Z}\1_B(m).
\end{equ}
By \eqref{eq:stochasti-symbol-in-Xsb-2}
\begin{equs}
 \| \| \<3>^{\mathrm{dyadic}}_{N,N_1,N_2,N_3} \|_{\X^{s_1 -1,b_1 -1}_T} \|^2_{L^p(\Omega)} 
\lesssim  p^3T^3 \sum_{n_1,n_2,n_3 \in \Z^2_N}  \jpb{n_{123}}^{2\underline  s-2}  \prod_{j=1}^3 \1_{|n_j| \sim N_j} \jpb{n_j}^{-2+2\alpha}  \sup_{m\in \Z}\1_B(m).
\end{equs}
Apply cubic sum estimate \eqref{eq:cubic-sum}, which holds for $0<\alpha <\f 5{12}$, we obtain:
\begin{equ}\label{eq:stochasti-symbol-in-Xsb-3}
 \| \| \<3>^{\mathrm{dyadic}}_{N,N_1,N_2,N_3} \|_{\X^{s_1 -1,b_1 -1}_T} \|^2_{L^p(\Omega)} 
 \lesssim  p^3T^3 \maxcurly{N_1,N_2,N_3}^{2(s_1 -s_\alpha)}.
\end{equ}

\item [\bf (ii) ]
Now take $b=1$ and $s=1$, we make estimates in the $\X^{0,0}_T=L^2([0,T]\times \T^2)$-norm.
By~\eqref{eq:stochasti-symbol-in-Xsb},
\begin{equs}
 \| \| \<3>^{\mathrm{dyadic}}_{N,N_1,N_2,N_3} \|_{\X^{0,0}_T}\|^2_{L^p(\Omega)} \lesssim & p^3T^3 \sum_{n_1,n_2,n_3 \in \Z^2_N}    \prod_{j=1}^3 \1_{|n_j| \sim N_j} \jpb{n_j}^{-2+2\alpha} \\
\lesssim  & p^3 T^3 \Pi_{j=1}^3 N_j^{-2+2\alpha} \sum_{n_1,n_2,n_3 \in \Z^2_N}    \prod_{j=1}^3 \1_{|n_j| \sim N_j}.\end{equs}
Bound the number of lattice point  in $ \1_{|n_j| \sim N_j}$by the trivial bound $N_j^2$,
we obtain
\begin{equ}\label{eq:X00-of-cube}
 \| \| \<3>^{\mathrm{dyadic}}_{N,N_1,N_2,N_3} \|_{\X^{0,0}_T}\|^2_{L^p(\Omega)}
 \lesssim  p^3 T^3 \maxcurly{N_1,N_2,N_3}^{6\alpha}
 \end{equ}

\item [\bf (iii)]  We now consider the case $b>\f 12$. Setting $$\theta = \f {1-b}{1-  b_1},$$
and choose $s_1 $ satisfying that
\begin{equ}
  s-1 = (s_1-1)\theta.
\end{equ}
Interpolating the bounds \eqref{eq:stochasti-symbol-in-Xsb-3} and \eqref{eq:X00-of-cube}, we obtain
\begin{equs}
 \| \| \<3>^{\mathrm{dyadic}}_{N,N_1,N_2,N_3}\|_{\X^{s-1,b-1}_T}\|^2_{L^p(\Omega)}
 & \lesssim T^3 \maxcurly{N_1,N_2,N_3}^{2(s_1-s_\alpha)\theta}\maxcurly{N_1,N_2,N_3}^{6\alpha(1-\theta)}\\
 &=p^3 T^3 \maxcurly{N_1,N_2,N_3}^{2\big(s-s_\alpha + \f {b-\underline  b}{1-  b_1}(s_\alpha-1+3\alpha) \big)}.
 \label{eq:cubic-interpolated}
\end{equs}
The proportional constants are independent of $N_1, N_2,N_3$. For this bound to be summable in $N_i$ we show that it is possible to choose $b_1$ such that 
\begin{equ}\label{constraints-cubic}
 s < s_\alpha - \f {b-b_1 }{1-b_1 }(s_\alpha-1+3\alpha).
 \end{equ}
\end{itemize}

{\bf Step 3.}
Given any $s<s_\alpha$, and set
$$b = \f 12+\delta, \quad b_1 = \f 12 -\delta.$$
The  interpolation constraint \eqref{constraints-cubic} can be rewrite as 
\begin{equ}
\f {4\delta}{1+2\delta}(s_\alpha -1 +3\alpha) \le s_\alpha-s,
\end{equ}
which holds for any $\delta$ sufficiently small. Let  $\delta_0$ be a number, which depends on $\alpha$ and $s$, smaller than this threshold. 

By summing over all dyadic numbers $N_1,N_2,N_3$, it yields, 
 \begin{equ}
\| \| \<3>_N\|_{\X^{s-1,b-1}_T}\|^2_{L^p(\Omega)} \le p^3 T^3 \sum_{1 \le N_1,N_2,N_3 \le N} \maxcurly{N_1,N_2,N_3}^{2\big(s-s_\alpha + \f {b-\underline  b}{1-\underline  b}(s_\alpha-1+3\alpha) \big)},\end{equ}
with which we conclude that for any $s<s_\alpha$, there exists $\delta_0>0$ such that for $\delta<\delta_0$,
 $\{\<3>_N, N\in \N\} $ is bounded in $L^p(\Omega; \X^{s_\alpha -\epsilon, \f 12 +\delta}_T)$,
 
Set $$b_0:= \f 12 +\delta_0,$$
we conclude the first part of the proposition.

{\bf Step 4.}
To  show that $\{ \<3>_N \}_{N \in \N}$ is a Cauchy sequence in $L^p(\Omega; \X^{s_\alpha -\epsilon, \f 12 +\delta}_T)$, we estimate
\begin{equ}
 \| \| \<3>_N - \<3>_M \|_{\X^{s-1,b-1}_T}\|^2_{L^p(\Omega)}.
\end{equ}
The only change in the proof from bounding  $\<3>_N$ itself to bounding $\<3>_N - \<3>_M$ lies in the summation domain: the lattice sum over $Z^2_N$ is replaced by 
\begin{equ}
 Z^2_{M,N} := \{ n \in \Z^2 \colon M < |n| \le N \}.
\end{equ}
Applying the Littlewood–Paley decomposition at this stage simply shifts the dyadic indices from
 $ 1\le N_j  \le N$ to $ M  \le N_j \le  N$ for all $j=1,2,3$, without introducing any further modifications to the argument.
 Using the interpolated bound \eqref{eq:cubic-interpolated}, we find
\begin{align}
 \| \| \<3>_N - \<3>_M \|_{\X^{s-1,b-1}_T}\|_{L^p(\Omega)}  \le & \sum_{\substack{ M  \le N_j \le  N  \\  j=1,2,3 }}  \| \| \<3>^{\mathrm{dyadic}}_{N,N_1,N_2,N_3} - \<3>^{\mathrm{dyadic}}_{M,N_1,N_2,N_3} \|_{\X^{s-1,b-1}_T}\|_{L^p(\Omega)}\\
 \lesssim & p^{\f32} T^{\f 32} \sum_{\substack{ M  \le N_j \le  N  \\  j=1,2,3 }}  \maxcurly{N_1,N_2,N_3}^{s-s_\alpha + \f {b-b_1 }{1-b_1 }(s_\alpha-1+3\alpha)}  \\
 \lesssim & p^{\f32} T^{\f 32}  M^{s-s_\alpha + \f {b-b_1 }{1-b_1 }(s_\alpha-1+3\alpha)} ,
\end{align}
provided that $ s < s_\alpha -  \f {b-b_1 }{1-b_1 }(s_\alpha-1+3\alpha) $. This decay in $M$ shows that
$\{ \<3>_N \}$ is indeed a Cauchy sequence in $L^p$, completing the proof of the required regularity result.
\end{proof}

\begin{remark}\label{rmk:b-vs-b_}
The interpolation method used in the proof of the preceding proposition will also feature in our analysis of other stochastic symbols. Because the arguments are essentially identical and the corresponding
 $L^2$ norm are straightforward, we restrict attention to the case $b<\f 12$ and omit further details for the remaining cases.
 \end{remark}

The next sections apply this streamlined interpolation strategy directly to quartic and higher‑order stochastic terms.

In Proposition \ref{prop:reg-of-cube}, we bounded the $\X^{s,b}$ norm of $\<3>$ and $\<30>$.
However, if we want to determine the regularity of $\<30>$ in the spaces $C([0,T]; W^{s,\infty}_x)$ or $C([0,T]; \mathcal{C}^s_x)$, the translation invariance of the space-time white noise becomes crucial. More precisely, denote by $\mathcal{T}_y$ the translation operator by $y\in \T^2$: 
$$\mathcal{T}_y f(x)=f(x+y).$$ A random function $f$ is spatially invariant if and only if
$$(\hat f(n): n\in \Z^2) \stackrel {\textrm{law}}=(\widehat{ \mathcal{T}_y f}(n), n\in \Z^2).$$ 
Note that for any function $m:\T^2 \to \C$, the following actions commute: 
$$\tau_x m(\jpb{\nabla})=m(\jpb{\nabla})\tau_x,$$
 since both are Fourier multipliers. In fact, for each $n$, we have 
$\widehat{\jpb{\nabla}\tau_x f}(n)=\widehat{\tau_x\jpb{\nabla}f}(n)$ because $e^{-\i n\cdot y} \hat f(n)= \widehat {\mathcal{T}_y f}(n) \stackrel{\textrm{law}} {=} \hat f(n)$. 

Consequently if $f$ is spatially stationary, then so is $m(\jpb{\nabla})f$. Moreover, the same property holds for the projection operator $\Pi_\Lambda f$
$$\Pi_\Lambda f(x)=\sum_{n \in \Lambda} e^{\i n \cdot x} \hat f(n),$$
 provided $\Lambda$ is a symmetric subset of $\Z^2$. It follows that $ \jpb{\nabla}^{s} \1_{[0,T]}\<30>_N(t,x)$ remains stationary.

\begin{lemma}\label{lem:Xsb-to-Wsinfty}
Let $s \in \R$ and let $f, g: \Omega\times [0,T]\times \T^2 \to \C$ be spatially stationary Gaussian random variables. 
Then, for any $\epsilon,\delta>0$, the following uniform estimate holds:
\begin{equ}
\| \|f\|_{L^\infty_TW^{s,\infty}_x }\|_{L^p(\Omega)}  \lesssim p^{3/2} \| \|f\|_{\X^{s+\epsilon+\delta,\frac{1}{2}+\delta}_T}\|_{L^2(\Omega)}. \label{eq:symbols-in-Sobolev-1}
\end{equ}
Consequently, 	\begin{equ}
\| \|\J(g)\|_{L^\infty_TW^{s,\infty}_x}\|_{L^p(\Omega)} \lesssim p^{3/2} \| \|g\|_{\X^{s+\epsilon+\delta-1, \delta -\frac{1}{2}}_T} \|_{L^2(\Omega)}, \label{eq:symbols-in-Sobolev-2}
\end{equ}
where $\cJ$ is the Duhamel operator for the wave equation.
\end{lemma}

\begin{proof}
Fix $0<\epsilon, \delta \ll 1$, and take $q \ge 2$ such that $\epsilon q >2$ (for the spatial variable) and $\delta q >1$ (for the temporal variable). 
Applying Sobolev inequalities both in space and in time, we have
\begin{equ}
\| \|f\|_{L^\infty_TW^{s,\infty}_x} \| _{L^p(\Omega)} \lesssim \| \|f\|_{W^{\delta,q}_tW^{s+\epsilon,q}_x} \| _{L^p(\Omega)}.
\end{equ}

We would like to exchange the order of the $L^p(\Omega)$ norm and the $L^q(\T^2)$ norm in order to eliminate the stochastic integrals. For this purpose we may either take $p$ large enough or simply set $p=q$. 

Observe that $$\jpb{\nabla_t}^{\delta} \jpb{\nabla_x}^{s+\epsilon} f$$ is also a spatial stationary family of Gaussian random variables, and in particular its $L^p(\Omega)$ norm is constant in $x$.
Hence, by Gaussian hypercontractivity we obtain:
\begin{equs}
 \| \|f\|_{W^{\delta,q}_tW^{s+\epsilon,q}_x} \|_{L^p(\Omega)}
&=  \| \| \jpb{\nabla_t}^{\delta} \jpb{\nabla_x}^{s+\epsilon} f \|_{L^q_{t,x}} \|_{L^p(\Omega)}\\
\le & \| \| \jpb{\nabla_t}^{\delta} \jpb{\nabla_x}^{s+\epsilon} f \|_{L^p(\Omega)} \|_{L^q_{t,x}}  
\lesssim  p^{3/2} \| \| \jpb{\nabla_t}^{\delta} \jpb{\nabla_x}^{s+\epsilon} f \|_{L^2(\Omega)} \|_{L^q_{t}}.
\end{equs}
where we have applied Minkowski's integral inequality, and in the last step removed the integration in $x$ by using the fact that
 $$\|x\mapsto  \jpb{\nabla_t}^{\delta} \jpb{\nabla_x}^{s+\epsilon} f \|_{L^2(\Omega)} $$ is constant,
 together with the convention $\int_{\mathbb{T}^2} 1 \, dx = 1$.
 
  Repeating the same argument, we insert an $L^2(\mathbb{T}^2)$ norm, exchange it with the $L^2(\Omega)$ norm, and apply Plancherel’s identity to obtain:
\begin{equs}
 \| \jpb{\nabla_t}^{\delta}\jpb{\nabla_x}^{s+\epsilon}f \|_{L^2(\Omega)} 
&=  \| \| \jpb{\nabla_t}^{\delta}\jpb{\nabla_x}^{s+\epsilon}f \|_{L^2(\Omega)} \| _{L^2_x}\\
 &=  \| \| \jpb{\nabla_t}^{\delta} \jpb{n}^{s+\epsilon} \hat f(t,\cdot)(n) \|_{\ell^2_n} \|_{L^2(\Omega)}.
\end{equs}
Combining these estimates and applying the triangle inequality once more, we find
\begin{equs}
\| \|f\|_{L^\infty_TW^{s,\infty}_x} \| _{L^p(\Omega)} & \le{} \| \| f\|_{W^{\delta,q}_tW^{s+\epsilon,q}_x} \|_{L^p(\Omega)}\\
&\le  p^{3/2}  \Big\|  \| \| \jpb{\nabla_t}^{\delta} \jpb{n}^{s+\epsilon}\hat f(t,\cdot)(n)\|_{\ell^2_n} \|_{L^2(\Omega)} \Big\|_{L^q_t}\\
& \leq p^{3/2}  \| \| \jpb{\nabla_t}^{\delta} \jpb{n}^{s+\epsilon} \hat f(t,\cdot)(n) \|_{L^q([0,T])} \|_{L^2(\Omega) \times \ell^2_n}\\
&\lesssim p^{3/2} \| \| \jpb{\lambda}^{\delta}   \jpb{n}^{s+\epsilon} \widetilde{\1_{[0,T]}f} (\lambda, n)  \|_{L^{q'}_\lambda} \|_{L^2(\Omega)\times\ell^2_n}
\end{equs} 
where $1 < q' \le 2$  and $\f 1 q+\f 1 {q'}=1$. We used the Hausdorff–Young inequality in time to control the $L^q$ norm of a function by the $L^{q'}$ norm of its inverse Fourier transform:
$$\|\tilde{f}\|_{L^{q'}}\le \|f\|_{L^q}.$$

To pass to the $\X^{s,b}$ norm, we work with the $L^2_\lambda$ norm and include the $\jpb{\lambda - \jpb{n}}$ factor.  Using the inequality
   $\jpb{\lambda}^\delta \lesssim \jpb{\lambda - \jpb{n}}^\delta \jpb{n}^\delta$, we deduce 
 \begin{equs}
\| \|f\|_{L^\infty_TW^{s,\infty}_x} \|_{L^p(\Omega)}  \lesssim 
&\lesssim p^{3/2} \| \| \jpb{\lambda}^{\delta}   \jpb{n}^{s+\epsilon} \hat f  (\lambda, n)  \|_{L^{q'}_\lambda} \|_{L^2(\Omega)\times\ell^2_n}\\
&\lesssim p^{3/2}  \|  \| \jpb{\lambda - \jpb{n}}^\delta  \jpb{n}^{s+\epsilon+\delta} \widetilde{\1_{[0,T]}f} (\lambda, n) \|_{L^{q'}_\lambda(\mathbf{R})} \|_{L^2(\Omega)\times\ell^2_n}.
 \end{equs}

We then apply Hölder’s inequality with $$\frac{1}{q'} = \frac{1}{\frac{2q'}{2-q'}} + \frac{1}{2}$$ to obtain
\begin{equs}
 \| \jpb{\lambda - \jpb{n}}^\delta  \jpb{n}^{s+2\epsilon} \widetilde{f} (\lambda, n) \|_{L^{q'}_\lambda}  
 \le \| \jpb{\lambda - \jpb{n}}^{-\frac{1}{2}} \|_{L^{\f {2q'}{2-q'}}_\lambda} \;  \| \jpb{\lambda - \jpb{n}}^{\frac{1}{2}+\delta}   \jpb{n}^{s+\epsilon+\delta} \widetilde{f} (\lambda, n)\|_{L^2_\lambda}.
\end{equs}

Since $\f {q'}{2-q'}>1$, the  $L^{\frac{2q'}{2-q'}}_\lambda$ norm of the function $\lambda \mapsto \jpb{\lambda - \jpb{n}}^{-\frac{1}{2}}$  is finite. Hence, we conclude 
\begin{equs}
\| \|f\|_{L^\infty_TW^{s,\infty}_x} \|_{L^p(\Omega)}  \lesssim & p^{3/2} \| \| \jpb{\lambda - \jpb{n}}^{\frac{1}{2}+\delta}  \jpb{n}^{s+\epsilon+\delta} \widetilde{\1_{[0,T]}f}(\lambda,n)\|_{L^2_\lambda \ell^2_n} \|_{L^2(\Omega)}, 
\end{equs}
where the inner term is controlled by the $\X^{s+\epsilon+\delta,\frac{1}{2}+\delta}_T$ norm of $f$.

For smaller $p$, Hölder’s inequality is first applied in $L^p(\Omega)$
 to ensure that the Sobolev inequality applies in the initial step, thereby completing the proof of \eqref{eq:symbols-in-Sobolev-1}.
Finally, in the case $f = \mathcal{J}(g)$, we apply Lemma~\ref{nl-Xsb} to obtain
\begin{equs}
\|  \| \J(g) \|_{\X^{s+\epsilon+\delta,\frac{1}{2}+\delta}_T} \|_{L^p(\Omega)}
\lesssim  p^{3/2} \| \|\J (g)\|_{\X^{s+\epsilon+\delta,\frac{1}{2}+\delta}_T}\|_{L^2(\Omega)} \lesssim p^{\f 32} \| \|g\|_{\X^{s+\epsilon+\delta-1,\frac{1}{2}+\delta-1}_T} \|_{L^2(\Omega)},
\end{equs}
which completes the proof.
\end{proof}

As a consequence of the proposition and applying Lemma \ref{lem:Xsb-to-Wsinfty}, we have the following. 
\begin{corollary}\label{cor-object-3}
Let $0 < \alpha < \f 5{12}$ and let $s_\alpha$ be defined as \eqref{eq:s-alpha}. For any $\epsilon>0$, and $s < s_\alpha - \epsilon$, we have 
\begin{equ}
\| \sup_{t \in [0,T]} \| \<3>_N(t) \|_{ W^{s-1,\infty}_x} \|_{L^p(\Omega)} \lesssim p^{\f 32} T^{\f 32}.
\end{equ}
\begin{equ}
\| \sup_{t \in [0,T]} \| \<30>_N \|_{W^{s,\infty}_x} \|_{L^p(\Omega)} \lesssim p^{\f 32} T^{\f 32} .
\end{equ}
Furthermore, the sequence $\<30>_N$ is Cauchy in $L^p(C([0,T]);W^{s,\infty}_x)$.
\end{corollary}
\begin{proof}
The claim follows directly from Proposition~\ref{prop:reg-of-cube} by setting $\delta < \f 12\mincurly{\epsilon, \f {4\delta}{1+2\delta}(s_\alpha -1 +3\alpha)}$.
\end{proof}

\subsection{Quintic symbols $\<32>$ and $\<320>$}
Let $\<2>_N$ and $\<30>_N$ be defined by \eqref{eq:cherry-w} and \eqref{eq:chicken-foot}, respectively. In this section we study the quintic term, $$\<32>_N := \<2>_N \<30>_N ,$$ a product of two Wiener chaos elements at different order. 
Applying the product formula \eqref{product} for multiple Wiener integrals yields the decomposition
\begin{equs}
\<32>_N  &= I_2[h^{(2),t,x,N}] I_3[h^{(3),t,x,N}]\\
&= I_5[h^{(2),t,x,N} \otimes_0 h^{(3),t,x,N,\cJ}] + 6I_3[h^{(2),t,x,N} \otimes_1 h^{(3),t,x,N,\cJ}] \\
& \quad + 6I_1[h^{(2),t,x,N} \otimes_2 h^{(3),t,x,N,\cJ}].
\end{equs}

These three terms on the right hand side are referred to, respectably, as the non-resonance term, the single resonance term, and the double resonance term. These labels reflect the types of index pairings that arise in the contractions. For convenience, we denote them by
 $$\<32>^{(5)}, \quad \<32>^{(3)}, \quad \<32>^{(1)}.$$
More explicitly,
 \begin{equ}\label{Qunintic-tems}
 \left\{\begin{aligned}
 \<32>^{(5)}_N(t) =&\sum_{n_1,\cdots, n_5 \in \Z^2_N} e^{\i (n_1+\cdots n_5)\cdot x} \I_5[g^{(5),t}_{n_1,\cdots,n_5}], \\
 \<32>^{(3)}_N(t) =&\sum_{n_3,n_4, n_5 \in \Z^2_N} e^{\i (n_3+ n_4+ n_5)\cdot x} \I_3[g^{(3),t}_{n_3,n_4,n_5}],  \\
 \<32>^{(1)}_N(t) =&\sum_{n_3  \in \Z^2_N} e^{\i n_3 \cdot x} \I_1[g^{(1),t}_{n_3}].
 \end{aligned} \right.
 \end{equ}
 where, by the product formula (\ref{product}), the kernel will be specified below:

 \begin{equ}\label{h1}
  \left\{\begin{aligned}
&{}g^{(5),t}_{n_1,\cdots,n_5}(s_1, \cdots,s_5) =  \prod_{j=1,5} \1_{[0,t]}(s_j) \f {\sin((t-s_j) \jpb{n_j})}{\jpb{n_j}^{1-\alpha}} \\
&\quad \times \int_0^t \f {\sin((t-s)\jpb{n_{234}})}{\jpb{n_{234}} }\prod_{j=2}^4 \1_{[0,s]}(s_j) \f {\sin((s-s_j) \jpb{n_j})}{\jpb{n_j}^{1-\alpha}} \, \d s, \\
&{}g^{(3),t}_{n_3,n_4,n_5}(s_3, s_4, s_5) =  \1_{[0,t]}(s_5) \f {\sin((t-s_5) \jpb{n_j})}{\jpb{n_j}^{1-\alpha}}  \sum_{n_2 \in \Z^2_N} \int_0^t \f {\sin((t-s)\jpb{n_{234}})}{\jpb{n_{234}} } \prod_{j=3}^4 \1_{[0,s]}(s_j)\\
&\quad \times  \f {\sin((s-s_j) \jpb{n_j})}{\jpb{n_j}^{1-\alpha}} \Bigl( \int_0^s  \frac{\sin ((t-s_2)\jpb{n_2}) \sin ((s-s_2)\jpb{n_2})}{ \jpb{n_2}^{2-2\alpha}}  \, \d s_2 \Bigr) \, \d s, 
\\
&{}g^{(1),t}_{n_3} =  \sum_{n_2,n_4 \in \Z^2_N} \int_0^t  \f {\sin((t-s)\jpb{n_{234}})}{\jpb{n_{234}} }  \1_{[0,s]}(s_3) \f {\sin((s-s_3)\jpb{n_3})}{\jpb{n_3}^{1-\alpha}} \\
&\quad \times \Bigl( \int_0^s \int_0^s \prod_{j=2,4}  \f {\sin((t-s_j)\jpb{n_j})\sin((s-s_j)\jpb{n_j})}{ \jpb{n_j}^{2-2\alpha}}  \, \d s_2 \d s_4 \Bigr) \, \d s. 
 \end{aligned} \right.
\end{equ}
Note that the contractions in (\ref{product}) are in both time and in frequencies.

We obtain the following regularity result for the quintic term. We use the fact that Wiener chaoses of different orders, with respect to the same underlying Wiener process, are independent --- here the underlying Wiener process is the symmetric sequence $\{\beta_n: n\in \Z_N^2\}$:
\begin{proposition}\label{prop:reg-of-quintic}
Let $\alpha \in [\f 14 ,\f 38)$. For any  $s <\f 56 -\alpha$,  there exists $b_0>\f 12$ (possibly depending on $(\alpha, s)$) such that for all $b > b_0$ we have
\begin{equ}
\| \| \<32>_N\|_{\X^{s-1,b-1}_T} \|_{L^p(\Omega)}  \lesssim p^{\f 52} T^{\f 52}.
\end{equ}
Moreover, the sequence $\{ \<32>_N \}$ is  Cauchy in $L^p(\Omega; \X^{s-1,b-1}_T)$.
\end{proposition}

\begin{proof}
We estimate each contribution in (\ref{Qunintic-tems}) separately.
\begin{itemize}
\item Non-resonance term $\<32>^{(5)}_N(t)$. 

We analyze it by applying the Littlewood–Paley decomposition to $\<32>^{(5)}_N$ with respect to the frequencies $n_1,n_2,n_3,n_4,n_5\in \Z^2_N$. 
More precisely we make a dyadic decomposition in its $n$-th frequency 
$$\widehat{\<32>^{(5)}_{N}}(t,n)=\sum_{n_i\in \Z_N^2, \& n_1+\dots, n_5=n} \I_5[g^{(5),t}_{n_1,\cdots,n_5}]$$
as follows.
Let $N_j=2^{k_j}$ denote dyadic number and define:
\begin{equ}
\widehat{\<32>^{(5),\mathrm{dyadic}}}_{N,N_1,N_2,N_3,N_4,N_5}(t,n)
 =  \sum_{\substack{n=n_1+\cdots+n_5\\n_1,\cdots,n_5 \in\Z^2_N}} \prod_{j=1}^5 \1_{|n_j| \sim N_j} \widehat{\<32>^{(5)}_N}(t,n).
 \end{equ}
Then,
\begin{equ}
\widehat{\<32>^{(5)}_{N}}(t,n) = \sum_{N_1,N_2,N_3, N_4, N_5 \ge 1}\widehat{\<32>^{(5),\mathrm{dyadic}}}_{N,N_1,N_2,N_3,N_4,N_5}(t,n).
\end{equ}
Here the summation are only over dyadic numbers, bearing in mind all indices are less or equal to $N$. More precisely,
\begin{equ}\label{sum-log}
\sum_{N_1,N_2,N_3, N_4, N_5 \ge 1}:=\sum_{j=1}^5 \sum_{k_j=1}^{\log_2(N)} .
\end{equ}

Hence, the estimate reduces to controlling each dyadic component.

By Lemma \ref{lem:symbols-in-Xsb}, we have
\begin{equs} \label{eq:g5-in-Xsb}
{}& \| \|\<32>^{(5),\mathrm{dyadic}}_N\|_{\X^{s-1,b-1}_T} \|_{L^p(\Omega)}  \\
\lesssim &p^{\f 52}\| \jpb{n}^{s-1} \jpb{\lambda}^{b-1} \Big(\sum_{\substack{n=n_{1\cdots 5}\\ n_1,\cdots, n_5 \in \Z^2_N}} \int_{[0,\infty)^5} |\cF_t(\1_{[0,T]}(t)h^{(5),t}_{n_1,\cdots, n_5})|^2 \, \d s_1 \cdots \d s_5 \Big)^{\f 12}\|_{\ell^2_n L^2_\lambda}.
\end{equs}

To organize the oscillatory structure, for $\sigma \in \{\pm\}$, we introduce the phase functions
\begin{equs}
\phi(\bn) &:= \phi(n_2,n_3,n_4) = \sigma_{234} \jpb{n_{234}} - \sum_{\ell =2}^4 \sigma_\ell \jpb{n_\ell},\\
\psi(\bn) &:= \psi(n_1,\cdots,n_5) = - \sigma \jpb{n_{12345}} - \sigma_{234} \jpb{n_{234}} - \sigma_1 \jpb{n_1} - \sigma_5 \jpb{n_5}, \\
\widetilde{\psi}(\bn) &:= \widetilde{\psi}(n_1,\cdots,n_5) = - \sigma \jpb{n_{12345}} + \sigma_{234} \jpb{n_{234}} - \sum_{\ell =1}^5 \sigma_\ell \jpb{n_\ell}.
\end{equs}

By the definition of $g^{(5),t}_{n_1,\cdots,n_5}$ and applying Lemma \ref{eq:sin-to-e} twice, we have, 
\begin{equs}
& |\cF_t(\1_{[0,T]}(t)g^{(5),t}_{n_1,\cdots,n_5})(\lambda-\sigma \jpb{n},n)| \\
 \lesssim  & \jpb{n_{234}}^{-1}   \prod_{j=1}^5 \jpb{n_j}^{-1+\alpha}\sum_{\sigma_{234}, \sigma_1,\cdots,\sigma_5 \in \{\pm\} } \f1 {\jpb{\phi(n_2,n_3,n_4)}\mincurly{\jpb{\lambda+\psi(n_1,\cdots,n_5)},\jpb{\lambda+\widetilde{\psi}(n_1,\cdots,n_5)}} }.
\end{equs}
Substituting this into \eqref{eq:g5-in-Xsb} and using the restrictio $n=n_{12345}$, we obtain
\begin{equs}
\| \|\<32>^{(5),\mathrm{dyadic}}_N\|_{\X^{s-1,b-1}_T} \|^2_{L^p(\Omega)}  \lesssim &p^5 T^5 \sum_{\substack{ \sigma_{12345}, \sigma_{234}, \\\sigma_1,\cdots,\sigma_5 \in \{\pm\}} }  \sum_{n_1,\cdots, n_5 \in \Z^2_N} \prod_{j=1}^5 \jpb{n_j}^{-2+2\alpha} \jpb{n_{12345}}^{2s-1}  \\
&\quad \times  \f{1}{\jpb{n_{234}}^2\jpb{\phi(\bn)}} \int_\R  \f {\jpb{\lambda}^{2b-2}} {\mincurly{\jpb{\lambda+\psi(\bn)},\jpb{\lambda+\widetilde{\psi}(\bn)}}^2}.
\end{equs}

As noted in Remark~\ref{rmk:b-vs-b_}, it suffices to estimate the $\X^{s,b}$-norm of the dyadic components of $\<32>^{(5)}_N$ for $b<\f 1 2$. We denote such a dyadic block by 
$\<32>^{(5),\mathrm{dyadic}}_N$ and let $b$ denote any number slightly smaller than $\f 12$
For the dyadic piece, we have,\begin{equs}
\| \|\<32>^{(5),\mathrm{dyadic}}_N\|_{\X^{s-1,b-1}_T} \|^2_{L^p(\Omega)} 
 \lesssim & p^5T^5    \sum_{n_1,\cdots, n_5 \in \Z^2_N} \prod_{j=1}^5 \1_{|n_j| \sim N_j} \jpb{n_j}^{-2+2\alpha} \jpb{n_{12345}}^{2s-1}  \jpb{n_{234}}^{-2}  \\
&\quad \times \f 1 {\jpb{\phi(\bn)}^2} \int_\R  \f {\jpb{\lambda}^{2b-2}} {\mincurly{\jpb{\lambda+\psi(\bn)},\jpb{\lambda+\widetilde{\psi}(\bn)}}^2}.
\end{equs}

By applying \eqref{basic-integration} and decomposing the phase contributions, we obtain
\begin{equs}
\f 1 {\jpb{\phi(\bn)}^2} \int_\R  \f {\jpb{\lambda}^{2b-2}} {\mincurly{\jpb{\lambda+\psi(\bn)},\jpb{\lambda+\widetilde{\psi}(\bn)}}^2} & \lesssim \f 1 {\jpb{\phi(\bn)}^2 \mincurly{\jpb{\lambda+\psi(\bn)},\jpb{\lambda+\widetilde{\psi}(\bn)}^{2-2b}} }\\
&\lesssim
\sum_{m,m' \in \Z} \f {\1_{|\phi-m|\le 1}}{\jpb{m}^2} \f {\1_{|\psi-m'|\le 1}+\1_{|\widetilde{\psi}-m'| \le 1}}{\jpb{m'}^{2-2b}}\\
&\lesssim \sup_{m,m' \in \Z} \1_{|\phi-m|\le 1} (\1_{|\psi-m'|\le 1}+\1_{|\widetilde{\psi}-m'| \le 1}).
\end{equs}
where, in the last step, we summed over $m, m'$ and used that $b<\f 12$.

Applying the quintic non‑resonance counting estimate from Lemma~\ref{lem:quintic-non-resonance}, we deduce that, for $s< \f 56 - \alpha$, there exists a small $\epsilon>0$ such that  
\begin{equs}
\| \|\<32>^{(5),\mathrm{dyadic}}_N\|_{\X^{s-1,b_{-}-1}_T} \|^2_{L^p(\Omega)}  
& \lesssim p^5T^5    \sum_{n_1,\cdots, n_5 \in \Z^2_N} \prod_{j=1}^5 \1_{|n_j| \sim N_j} \jpb{n_j}^{-2+2\alpha} \jpb{n_{12345}}^{2s-1}  \jpb{n_{234}}^{-2} \\
& \lesssim  p^5T^5  \maxcurly{N_1,\cdots, N_5}^{-2\epsilon}.
\end{equs}

Summing over all dyadic piece, then, it yields
\begin{equs}
 \| \|\<32>^{(5),\mathrm{dyadic}}_N\|_{\X^{s-1,b-1}_T} \|^2_{L^p(\Omega)} & \le  \sum_{N_1,\cdots,N_5 \ge 1}  \| \|\<32>^{(5),\mathrm{dyadic}}_N\|_{\X^{s-1,b-1}_T} \|^2_{L^p(\Omega)}  \\
& \lesssim  p^5T^5 \sum_{N_1,\cdots,N_5 \ge 1} \maxcurly{N_1,\cdots, N_5}^{-2\epsilon} \\
& \lesssim p^5T^5N^{-2\epsilon+2\kappa},
\end{equs}
where $\kappa>0$ is chosen sufficiently small so that $\kappa \ll \epsilon$.

\item The single resonance term $\<32>^{(3)}_N(t)$ from (\ref{Qunintic-tems}).

 Recall that  
 $$\widehat{ \<32>^{(3)}_N}(t,n) =\sum_{n_3,n_4, n_5 \in \Z^2_N \& n_3+n_4+n_5=n} \I_3[g^{(3),t}_{n_3,n_4,n_5}](n), $$
where $g^{(3),t}_{n_3,n_4,n_5}$  is as in (\ref{h1}) which involves another summation, in $n_2$, indicating the spatial contraction coming from the product formula for iterated Wiener integrals. For $N_j=2^{k_j}$, define:
\begin{equ}
\widehat{\<32>^{(3),\mathrm{dyadic}}}_{N, N_2, N_3, N_4, N_5}(t,n)
 =  \sum_{\substack{n=n_3+n_4+n_5\\ n_3, n_4,n_5 \in\Z^2_N}} \prod_{j=2}^5 \1_{\{|n_j| \sim N_j\}} \widehat{\<32>^{(3)}_N}(t,n).
\end{equ}
where we have incorporated the extra $\1_{ \{|n_2| \sim N_2\}}$ from the formula of $g^{(3),t}_{n_3,n_4,n_5}$.
Denoting $n_{345}=n_3+n_4+n_5$ and  applying \eqref{eq:sin-to-e} together with Lemma \ref{lem:symbols-in-Xsb}, we obtain the bound
\begin{equs}
{} & \| \|\<32>^{(3),\mathrm{dyadic}}_N\|_{\X^{s-1,b-1}_T} \|^2_{L^p(\Omega)} \\
\lesssim &p^{\f 32} \| \jpb{n}^{s-1} \jpb{\lambda}^{b-1} \Big(\sum_{\substack{n=n_{345}\\ n_3,n_4, n_5 \in \Z^2_N}} \prod_{j=2}^5 \1_{|n_j| \sim N_j} \int_{[0,\infty)^3} |\cF_t(\1_{[0,T]}(t)g^{(3),t}_{n_3,n_4,n_5})|^2 \, \d s_3 \d s_4 \d s_5 \Big)^{\f 12}\|_{\ell^2_n L^2_\lambda}
\end{equs}
Here, the resonance structure differs from the non‑resonant case, so we introduce a new set of phase functions:
 \begin{equs}
\phi(\bn) &:= \phi(n_2,n_3,n_4) = \sigma_{234}\jpb{n_{234}} -  \sum_{j=2}^4 \sigma_j\jpb{n_j} , \label{phi-n}\\
\psi(\bn) &:= \psi(n_2,n_3,n_4,n_5) = - \sigma_{345}\jpb{n_{345}} - \sigma_{234}\jpb{n_{234}} - \sigma_5 \jpb{n_5}, \\
\widetilde{\psi}(\bn) &:= \widetilde{\psi}(n_2,n_3,n_4,n_5) = -\sigma_{345}\jpb{n_{345}} + \sigma_{234}\jpb{n_{234}} - \sum_{\ell=2}^5 \sigma_\ell \jpb{n_\ell}.
\end{equs}
With this notation, the Fourier transform satisfies
\begin{equs}
{} &|\cF_t(\1_{[0,T]}(t)g^{(3),t}_{n_3,n_4,n_5})(\lambda - \sigma_{345} \jpb{n_{345}},n_{345})| \lesssim   \prod_{j=3}^5 \jpb{n_j}^{-1+\alpha}  \prod_{j=3}^5 \1_{|n_j| \sim N_j}  \\
& \quad \times \Big| \sum_{n_2 \in \Z^2_N}  \jpb{n_{234}}^{-1} \1_{|n_2| \sim N_2} \jpb{n_2}^{-2+2\alpha} \f1{\jpb{\phi(\bn)}\mincurly{\jpb{\lambda+\psi(\bn)},\jpb{\lambda+\widetilde{\psi}(\bn)}    }} \Big|.
\end{equs}

At first glance, one might try to proceed as in the non‑resonance case, counting under restrictions of the form
$$\1_{|\phi-m| \le 1} \1_{|\psi-m'| \le 1}, \qquad \1_{|\phi-m| \le 1} \1_{|\widetilde{\psi}-m'| \le 1}.$$ 
 However, this direct approach is problematic here because the indices $n_{234}$ and $n_{345}$ are no longer coupled---they appear in separated terms.
 
 In fact we can work with a strong norm $L^p(H^s)$ which allows to overcome this difficulty. Specifically, we switch from $\X^{s,b}$ framework to an $H^s$ -based argument and apply
 the inhomogeneous Strichartz estimate \eqref{eq:inhomo-Stri} and Gaussian hypercontractivity, which yields
\begin{equs}
 \| \|\<32>^{(3),\mathrm{dyadic}}_N(N_2,\cdots,N_5)\|_{\X^{s-1,b-1}_T} \|^2_{L^p(\Omega)} \lesssim & \| \|\<32>^{(3),\mathrm{dyadic}}_{N,N_2,\cdots,N_5}\|_{L^{1+}_TH^{s-1}_x} \|^2_{L^p(\Omega)} \\
 \lesssim &T^{2(\f 1{1+} - \f 12)}\| \|\<32>^{(3),\mathrm{dyadic}}_{N,N_2,\cdots,N_5}\|_{L^2_TH^{s-1}_x} \|^2_{L^p(\Omega)} \\
\lesssim  & p^3T^{2(\f 1{1+} - \f 12)} \| \|\<32>^{(3),\mathrm{dyadic}}_{N,N_2,\cdots,N_5}\|_{H^{s-1}_xL^2_T} \|^2_{L^2(\Omega)} \\
\le & p^3T^\theta \sup_{t \in [0,T]}  \| \|\<32>^{(3),\mathrm{dyadic}}_{N,N_2,\cdots,N_5}\|_{H^{s-1}_x} \|^2_{L^2(\Omega)}.
\end{equs}
for some $\theta >0$, where $1+$ denotes an exponent slightly larger than $1$.

By the $H^s_x$-norm version in Lemma \ref{lem:symbols-in-Xsb}, we have 
\begin{align}\label{eq:h3-in-Sobolev}
&\| \|\<32>^{(3),\mathrm{dyadic}}_N \|_{H^{s-1}_x}\|^2_{L^2(\Omega)} \\
= &\sum_{n_3,n_4,n_5\in\Z^2_N  } \jpb{n_{345}}^{2s-2} \prod_{j=2}^5 \1_{|n_j| \sim N_j}  \int_{[0,\infty)^3} |g^{(3),t}_{n_3,n_4,n_5}|^2 \, \d s_3 \d s_4 \d s_5.
\end{align}
By definition of $g^{(3),t}_{n_3,n_4,n_5}$ and applying the \eqref{eq:sin-to-e}, we have 
\begin{align}
\Big| \int_{[0,\infty)^3} |g^{(3),t}_{n_3,n_4,n_5}|^2 \, \d s_3 \d s_4 \d s_5 \Big| \le &T^3 \prod_{j=3}^5 \jpb{n_j}^{-2+2\alpha} \\
& \quad \times \Big| \sum_{n_2 \in \Z^2_N}  \jpb{n_{234}}^{-1} \1_{|n_2| \sim N_2} \jpb{n_2}^{-2+2\alpha} \f1{\jpb{\phi(\bn)}} \Big|^2
\end{align}
where $\phi(\bn)$ is given by \eqref{phi-n}.
Hence, we needs to control
\begin{align}
\| \|\<32>^{(3),\mathrm{dyadic}}_N \|_{H^{s-1}_x}\|^2_{L^2(\Omega)}  \lesssim &T^3 \sum_{n_3,n_4,n_5\in\Z^2_N  } \jpb{n_{345}}^{2s-2}  \prod_{j=3}^5 \1_{|n_j| \sim N_j|} \jpb{n_j}^{-2+2\alpha} \\
&\quad  \times \Big| \sum_{n_2 \in \Z^2_N}  \jpb{n_{234}}^{-1} \1_{|n_2| \sim N_2} \jpb{n_2}^{-2+2\alpha} \f1{\jpb{\phi(\bn)}} \Big|^2
\end{align}
\textbf{Step 1} We get two ways summing over $n_2$. Denote
\begin{equ}
I:= \Big| \sum_{n_2 \in \Z^2_N}  \jpb{n_{234}}^{-1} \1_{|n_2| \sim N_2} \jpb{n_2}^{-2+2\alpha} \f1{\jpb{\phi(\bn)}} \Big|^2
\end{equ}
By decomposing the phase function $\phi(\bn)$, we have 
\begin{align}
I \le & \Big| \sum_{m \in \Z} \sum_{n_2 \in \Z^2_N}  \1_{|n_2| \sim N_2} \jpb{n_{234}}^{-1} \1_{|n_2| \sim N_2} \jpb{n_2}^{-2+2\alpha} \f{\1_{|\phi-m| \le 1}}{\jpb{m}} \Big|^2 
\end{align}
To estimate the square,   we take supremum in $m$ in one of the factor.  The summation over $m$  of the supremum leads to $\log(m)$, see (\ref{sum-log}),  we note also that $|m| \lesssim \maxcurly{N_2,N_3,N_4}$. For the other factor in the product, we apply basic resonance counting estimate from Lemma \ref{lem:basic-resonant}:
\begin{equs}
\sum_{m \in \Z} \sum_{n_2 \in \Z^2_N}  \1_{|n_2| \sim N_2} \jpb{n_{234}}^{-1} \1_{|n_2| \sim N_2} \jpb{n_2}^{-2+2\alpha} \f{\1_{|\phi-m| \le 1}}{\jpb{m}} \lesssim  \log(2+N_2) \jpb{n_{34}}^{-\f 32+2\alpha}.
\end{equs} We deduce
\begin{align}
I \lesssim  & \log(2+\maxcurly{N_2,N_3,N_4}) \Bigl( \sup_{m \in \Z} \sum_{n_2 \in \Z^2_N}   \jpb{n_{234}}^{-1} \1_{|n_2| \sim N_2} \jpb{n_2}^{-2+2\alpha} \1_{|\phi-m| \le 1} \Bigr)\\
& \quad \times  \Bigl( \sum_{m \in \Z} \sum_{n_2 \in \Z^2_N}  \1_{|n_2| \sim N_2} \jpb{n_{234}}^{-1} \1_{|n_2| \sim N_2} \jpb{n_2}^{-2+2\alpha} \f{\1_{|\phi-m| \le 1}}{\jpb{m}} \Bigr)\\
\lesssim & \log(2+\maxcurly{N_2,N_3,N_4}) \log(2+N_2) \jpb{n_{34}}^{-\f 32+2\alpha} \\
&\quad \times \Bigl( \sup_{m \in \Z} \sum_{n_2 \in \Z^2_N}   \jpb{n_{234}}^{-1} \1_{|n_2| \sim N_2} \jpb{n_2}^{-2+2\alpha} \1_{|\phi-m| \le 1} \Bigr).
\end{align}

\textbf{Step 2}. Substituting back to the dyadic piece and first summing over $n_5$, we have 
\begin{align}
{}& \| \|\<32>^{(3),\mathrm{dyadic}}_{N,N_2,\cdots,N_5}\|_{\X^{s-1,b-1}_T} \|^2_{L^p(\Omega)} \\
\lesssim & p^{\f 32} T^3 \log(2+\maxcurly{N_2,N_3,N_4}) \log(2+N_2) \sum_{n_3,n_4,n_5\in\Z^2_N} \jpb{n_{345}}^{2s-2}  \prod_{j=3}^5 \1_{|n_j|\sim N_j} \jpb{n_j}^{-2+2\alpha}   \\
 &\quad \times \jpb{n_{34}}^{-\f 32+2\alpha} \Bigl( \sup_{m \in \Z} \sum_{n_2 \in \Z^2_N}   \jpb{n_{234}}^{-1} \1_{|n_2| \sim N_2} \jpb{n_2}^{-2+2\alpha} \1_{|\phi-m| \le 1} \Bigr) \\
\lesssim & p^{\f 32} T^3 \log(2+\maxcurly{N_2,N_3,N_4}) \log(2+N_2) N^{-2\kappa}_5 \\
&\quad \times  \sup_{m \in \Z} \sum_{n_2, n_3,n_4 \in\Z^2_N}\jpb{n_{234}}^{-1}  \jpb{n_{34}}^{2s-\f 72+4\alpha+2\kappa}  \prod_{j=2}^4 \1_{|n_j| \sim N_j} \jpb{n_j}^{-2+2\alpha} \1_{|\phi-m| \le 1} 
\end{align}
where $\kappa>0$ is chosen sufficiently small (to be determined) and we used the convolution inequality in the last line and leads to condition:
\begin{equ}
0<s < 1-\alpha-\kappa.
\end{equ}

For any $0<\epsilon<\f 5 6-\alpha$ and $s=\f 5 6-\alpha-\epsilon$, in order to ensure that
\begin{equ}
2s-\f 72 +4\alpha +2\kappa  <-1,
\end{equ}
we require $$\kappa < \epsilon +  \mincurly{\f 5 6-2\alpha, \f 16}.$$

Applying the special cubic sum estimate \eqref{eq:special-cubic-sum}, we further deduce 
\begin{equs}
{} &\| \|\<32>^{(3),\mathrm{dyadic}}_N\|_{\X^{s-1,b-1}_T} \|^2_{L^p(\Omega)}\\
\lesssim & p^3 T^4 \sum_{N_2 \ge 1} \log(2+N_2) \log^2(2+\maxcurly{N_2,N_3,N_4}) N^{-\kappa}_5 \maxcurly{N_2,N_3,N_4}^{-\epsilon}
\end{equs}
which is summable across all dyadic scales $ N_2,N_3,N_4,N_5 \ge 1$ . 

For $s \le 0$, we can get the desired bound by using the embedding $H^{s'}_x \subset H^{s}_x$ for some $s'>0$.

\item The double resonance term $\<32>^{(1)}_N(t)$ from (\ref{Qunintic-tems}). 

The dyadic piece for this term is: 
\begin{equ}
\widehat{\<32>^{(1),\mathrm{dyadic}}}_{N, N_2, N_3, N_4, N_5}(t,n)
 =  \sum_{\substack{n=n_3, n_3 \in\Z^2_N}} \prod_{j=2}^4 \1_{|n_j| \sim N_j} \widehat{\<32>^{(1)}_N}(t,n),
\end{equ}
where $n_2$ and $n_4$ are indices in the contraction where we have contracted  $n_1$ with $n_2$ and also $n_4$ with $n_5$.
Again we switch from $\X^{s,b}$-norm  to $H^s$-norm for $\<32>^{(1)}_N(t)$ by applying the inhomogeneous Strichartz estimate \eqref{eq:inhomo-Stri} and Gaussian hypercontractivity which leads to
\begin{equs}
 \| \|\<32>^{(1),\mathrm{dyadic}}_{N, N_2, N_3, N_4, }\|_{\X^{s-1,b_-1}_T} \|^2_{L^p(\Omega)} \lesssim pT \sup_{t \in [0,T]}  \| \|\widehat{\<32>^{(1),\mathrm{dyadic}}}_{N, N_2, N_3, N_4}\|_{H^{s-1}_x} \|^2_{L^2(\Omega)}.
\end{equs}

Utilizing Lemma \ref{lem:symbols-in-Xsb}, for fixed $t \in [0,T]$, we obtain
\begin{align}
 \| \| \<32>^{(1),\mathrm{dyadic}}_{N, N_2, N_3, N_4, N_5} \|_{H^{s-1}_x} \|^2_{L^2(\Omega)} \le \sum_{n_3 \in \Z^2_N  } \jpb{n_3}^{2s-2} \prod_{j=2}^4 \1_{|n_j| \sim N_j} \int_{[0,\infty)} |g^{(1),t}_{n_3}|^2 \, \d s_3.
\end{align}

Recalling the expression for $g^{(1),t}_{n_3}$, we deduce
\begin{equs}
|g^{(1),t}_{n_3}|^2 & \lesssim \jpb{n_3}^{-2+2\alpha} \Big| \sum_{n_2,n_4 \in \Z^2_N} \jpb{n_2}^{-2+2\alpha}  \jpb{n_4}^{-2+2\alpha} \jpb{n_{234}}^{-1} \\
&\qquad \times\int_{s_3}^t  \sin((t-s)\jpb{n_{234}})  \sin((s-s_3)\jpb{n_3}) \cos((t-s)\jpb{n_2})\cos((t-s)\jpb{n_4}) \, \d s\Big|^2.
\end{equs}
Applying  \eqref{eq:sin-to-e} and integrate with $s$, we have 
we obtain
\begin{equ}
\|g^{(1),t}_{n_3}|^2 \lesssim T \sum_{n_3 \in \Z^2_N  } \jpb{n_3}^{2s-4+2\alpha} \Big| \sum_{n_2,n_4 \in \Z^2_N} \jpb{n_2}^{-2+2\alpha}  \jpb{n_4}^{-2+2\alpha} \jpb{n_{234}}^{-1} \f 1{\jpb{\phi(\bn)}} \Big|^2
\end{equ}
where $\phi(\bn) = \sigma\jpb{n_{234}} + \sigma_2 \jpb{n_2} + \sigma_3 \jpb{n_3} + \sigma_4 \jpb{n_4}$.

Introducing an auxiliary parameter $m \in \Z$ as before and note that $|m| \lesssim \maxcurly{N_2,N_3,N_4}$, we rewrite
\begin{equs}
\| \| \<32>^{(1),\mathrm{dyadic}}_{N, N_2, N_3, N_4, N_5} \|_{\X^{s-1,b_-1}_T} \|^2_{L^2(\Omega)} \lesssim &  pT\log(2+\maxcurly{N_2,N_3,N_4}) \sum_{n_3 \in \Z^2_N  } \1_{|n_3| \sim N_3} \jpb{n_3}^{2s-4+2\alpha} \\
&\;\; \times \Big| \sup_{m \in \Z} \sum_{n_2,n_4 \in \Z^2_N} \prod_{j=2,4} \1_{|n_j| \sim N_j} \jpb{n_j}^{-2+2\alpha} \jpb{n_{234}}^{-1}  \1_{|\phi-m| \le 1} \Big|^2
\end{equs}

By first summing over $n_4$ and using the basic resonant estimate (Lemma \ref{lem:basic-resonant}), for $\alpha \ge  \f 14$, we have
\begin{equs}
{} & \Big| \sup_{m \in \Z} \sum_{n_2,n_4 \in \Z^2_N} \prod_{j=2,4} \1_{|n_j| \sim N_j} \jpb{n_j}^{-2+2\alpha} \jpb{n_{234}}^{-1}  \1_{|\phi-m| \le 1} \Big|\\
\lesssim & \log(2+N_4)\Big|  \sum_{n_2  \in \Z^2_N} \1_{|n_2| \sim N_2} \jpb{n_2}^{-2+2\alpha} \jpb{n_{23}}^{-\f 3 2+2\alpha}  \Big| \lesssim \log(2+N_4)\jpb{n_3}^{-\f 32+4\alpha},
\end{equs}
where the last step follows from the discrete convolution sum inequality under condition $0<\alpha < \f 3 8$.

Therefore, we deduce, for $\f 14 \le \alpha < \f 3 8$,
\begin{equs}
{} & \| \| \<32>^{(1),\mathrm{dyadic}}_{N, N_2, N_3, N_4, N_5} \|_{\X^{s-1,b-1}_T} \|^2_{L^2(\Omega)} \\
\lesssim & pT\log(2+\maxcurly{N_2,N_3,N_4})  \log(2+N_4) \sum_{n_3 \in \Z^2_N  } \1_{|n_3| \sim N_3} \jpb{n_3}^{2s-7+10\alpha} \\
\lesssim & pT\log(2+\maxcurly{N_2,N_3,N_4}) \log(2+N_4)  N^{2s-5+10\alpha}_3, \label{eq:32-1-dispersion}
\end{equs}
where in the last line we sum over $n_3$.

\textbf{Case (a)}. If $\maxcurly{N_2,N_3,N_4} \sim N_3$ and we take $s< \f 5 6 -\alpha$, then summing over all the dyadic pieces yields $$ \| \|\<32>^{(1)}_N\|_{\X^{s-1,b-1}_T} \|^2_{L^p(\Omega)} \lesssim pT .$$

\textbf{Case (b)}. If instead $\maxcurly{N_2,N_3,N_4} = \maxcurly{N_2,N_4}$, we assume without loss of generality that $N_2 \ge N_4$. We have the logarithmic divergence in $N_2$ which causes issues in summing over $n_2$ and $n_3$. To resolve it, we need to explore the sine-cancellation lemma. More precisely, we give another way to bound 
\begin{equs}
{} & \sum_{n_3 \in \Z^2_N  } \1_{|n_3| \sim N_3} \jpb{n_3}^{2s-4+2\alpha}  \Big| \sum_{n_2,n_4 \in \Z^2_N} \prod_{j=2,4} \1_{|n_j| \sim N_j} \jpb{n_j}^{-2+2\alpha}  \jpb{n_{234}}^{-1} \\
&\quad \times\int_{s_3}^t  \sin((t-s)\jpb{n_{234}})  \sin((s-s_3)\jpb{n_3}) \cos((t-s)\jpb{n_2})\cos((t-s)\jpb{n_4}) \, \d s\Big|^2.
\end{equs}

Define
\begin{equ}
f_{N_2,N_4}(t,s,n_2; s_3,n_3,n_4) = \jpb{n_2+n_{34}}^{-1} \jpb{n_2}^{-2} \sin((s-s_3)\jpb{n_3}) \cos((t-s)\jpb{n_4}) 
\end{equ}
and set $A= \maxcurly{N_3,N_4}$. We then observe the uniform bounds
\begin{equs}
|f(t,s,n_2)| \lesssim A N^{-3}_2, \quad  |\partial_s f(t,s,n_2)| \lesssim A N^{-3}_2, \quad  |f(t,s,n_2) -f(t,s,-n_2)|  \lesssim  AN^{-4}_2,
\end{equs}
These estimates verify the assumptions of Lemma ~\ref{lem:sine-cancellation}. Applying this lemma yields
\begin{equs} 
{} & \Big| \sum_{|n_4| \sim N_4} \jpb{n_4}^{-2+2\alpha}  N^{2\alpha}_2 \sum_{|n_2| \sim N_2} \int_{s_3}^t  \sin((t-s)\jpb{n_2+n_{34}})  \cos((t-s)\jpb{n_2}) f(t,s,n_2) \, \d s \Big| \\
 \lesssim & N^{2\alpha}_4 N^{2\alpha}_2 \maxcurly{N_3,N_4}^2 N^{-1}_2
\end{equs}
Substituting this estimate back into the original expression \eqref{eq:h3-in-Sobolev}, for the contribution of the term under consideration yields
\begin{equ}\label{eq:32-1-sine}
\| \|\<32>^{(1),\mathrm{dyadic}}_{N,N_2,N_3,N_4}\|_{\X^{s-1,b-1}_T} \|^2_{L^2(\Omega)} \lesssim T N^{2s-2+2\alpha}_3N^{2\alpha}_4 N^{-1+2\alpha}_2 \maxcurly{N_3,N_4}^2  
\end{equ}

Interpolating \eqref{eq:32-1-dispersion} and \eqref{eq:32-1-sine}, we obtain, for $s<\f 56 -\alpha$, 
\begin{equ}
\| \|\<32>^{(1),\mathrm{dyadic}}_N\|_{\X^{s-1,b-1}_T} \|^2_{L^2(\Omega)} \lesssim T \maxcurly{N_2,N_3,N_4}^{-2\epsilon},
\end{equ}
for some sufficiently small $\epsilon>0$,. This bound is summable over all dyadic scales $N_2,N_3, N_4$, and we have completed the estimates for the third term.
\end{itemize}
This concludes the proof of the proposition.
\end{proof}

\begin{remark}
We emphasize that the restriction $\alpha < \f 3 8$ actually originates from the high-high to low interaction. More precisely, at the dyadic level, applying Lemma~\ref{lem:cubic-sup-count}, we have
\begin{align}
&\Big| \sum_{m \in \Z} \sum_{n_2,n_4 \in \Z^2_N} \jpb{n_2}^{-2+2\alpha}\jpb{n_4}^{-2+2\alpha} \jpb{n_{234}}^{-1}  \f 1 {\jpb{m}} \cdot \1_{|\phi-m| \le 1} \Big|\\
\lesssim & N^{-2+2\alpha}_2N^{-2+2\alpha}_4 N^{-1}_{234} \log(\maxcurly{N_2,N_3,N_4})\Big| \sup_{m \in \Z} \sum_{n_2,n_4 \in \Z^2_N}\1_{|n_2| \sim N_2} \1_{|n_4|\sim N_4|} \1_{|\phi-m| \le 1}\Big| \\
\lesssim & N^{-2+2\alpha}_2N^{-2+2\alpha}_4 N^{-1}_{234} \log(\maxcurly{N_2,N_3,N_4}) \\
&\qquad \times \Big(\mincurly{N_2,N_{34},N_{234}}\Big)^{-\f 12}  \Big(\mincurly{N_2,N_{234}}\Big)^2 \Big(\mincurly{N_4,N_{34}}\Big)^2 \\
\le &\mincurly{N_2,N_{234}} \Big(\mincurly{N_4,N_{34}}\Big)^{\f 32} N^{-\f 12}_{234} N^{-\f 32+2\alpha}_2N^{-2+2\alpha}_4  \log(\maxcurly{N_2,N_3,N_4}).
\end{align}

If we consider the output frequency $n=n_3$ on the dyadic scale, we eventually compute
\begin{align}
N^{s-1+\alpha}_3 \mincurly{N_2,N_{234}} \Big(\mincurly{N_4,N_{34}}\Big)^{\f 32} N^{-\f 12}_{234} N^{-\f 32+2\alpha}_2N^{-2+2\alpha}_4  \log^2(\maxcurly{N_2,N_3,N_4}).
\end{align}

The worst case occurs when $N_3 \ll \mincurly{N_2,N_4}$. In this situation, the best we can obtained is 
\begin{equ}
\maxcurly{N_2,N_4}^{-\f 32 +4\alpha} \log^2(\maxcurly{N_2,N_4}),
\end{equ}
which forces the condition $\alpha < \f 3 8$. Thus, if the output frequency $n$ is low and does not contribute to the summation, we are restricted to the range $0<\alpha <\f 38$.
\end{remark}
\begin{remark}\label{non-wick}
This observation made in the previous remark has an important implication: to investigate  the range $\alpha \ge \f 38$,  it becomes essential to renormalize the term  $\<32>^{(1)}_N$. Indeed, the preceding estimate shows that
\begin{equ}
\lim_{N \to \infty } \mathbb{E}[|\<32>^{(1)}_N|^2] = \infty,
\end{equ}
signalling the divergence of the second moment. To address this, we naturally introduce the non‑Wick ordered renormalization of the single‑resonance term $\<32>^{(1)}_N$:
\begin{equ}
\<32>^{(1)}_N - \mathbb{E}[|\<32>^{(1)}_N|^2] .
\end{equ}
This subtraction removes the diverging component, but its implementation and analysis  require additional work. A detailed treatment of this renormalization will be developed in forthcoming work.
\end{remark}

\subsection{Septic symbols $\<7>$ and $\<70>$}
To identify the regularity of stochastic symbol $$\<7>:=\<30>_N \<1>_N\<30>_N,$$ 
we consider all possible resonance subcases. However, because the required regularity of $ \<30>_N \<1>_N\<30>_N$ is only $-\alpha-$, 
we need only use the partial multilinear smooth effect within $\<30>_N$,
 without fully exploiting the multilinear smooth between $\<30>_N, \<1>_N$ and $\<30>_N$. 
 
As mentioned earlier, we claim that the spatial regularity of $\<7>_N:=\<30>_N \<1>_N\<30>_N$ is $-\alpha-\epsilon$ for any~$\epsilon <0$: 
\begin{proposition}\label{prop:reg-septic}
 Let $T \le 1$ and $p \ge 2$. Assume $s<1-\alpha$ and $\f 14 \le \alpha < \f 3 8$. We have, for any $b > \f 12$, 
\begin{equ}
\| \| \<7>_N\|_{\X^{s-1,b-1}_T} \|_{L^p(\Omega)} \lesssim p^{\f 7 2} T.
\end{equ}
\end{proposition}

\begin{proof}
As with the same of Proposition \ref{prop:reg-of-quintic}, it is easier to explore the Sobolev norm.
By inhomogeneous Strichartz estimate from Lemma \ref{lem:inhomo-Stri} and Gaussian hypercontractivity, we have 
\begin{align}
\| \| \<7>_N\|_{\X^{s-1,b-1}_T} \|_{L^p(\Omega)} \lesssim p^{\f72} \sup_{t \in [0,T]} \| \| \<7>_N\|_{H^{s-1}_x} \|_{L^2(\Omega)}.
\end{align}
It suffices to prove the right hand side is bounded uniformly in $[0,1]$  and uniformly in $N \ge 1$.

{\bf  Bound the Wiener Integrals.}
We first simplify and explore symmetries.
By the product formula to $(\<30>_N )^2\<1>_N$:
\begin{equs}
\<7>_N&=\<1>_N\sum_{r=0}^3 r! \Bigl(\binom{3}{r}\Bigr)^2 I_{6-2r}\bigl[(\text{Sym}(\cJ(h^{(3),t,x,N}) \otimes_r \cJ(h^{(3),t,x,N} )))\bigr]\\
&=\sum_{r=0}^3 r! \Bigl(\binom{3}{r}\Bigr)^2 I_{7-2r} \Bigl[\text{Sym} (h^{(1),t,x,N}\otimes_0 (\cJ(h^{(3),t,x,N}) \otimes_r \cJ(h^{(3),t,x,N}) ))\Bigr]\\
&\qquad +\sum_{r=0}^2 r! (\binom{3}{r})^2 I_{5-2r} \Bigl[\text{Sym} (h^{(1),t,x,N}\otimes_1 (\cJ(h^{(3),t,x,N}) \otimes_r \cJ(h^{(3),t,x,N}) ))\Bigr].
\end{equs}
Since the estimates for the individual terms obtained from symmetrization are of the the same magnitude, we only need to bound the $H^s$ norm  of the form
\begin{equ}\label{terms-7}
\begin{aligned}
A_{7-2r} :={} &I_{7-2r} \Bigl[\text{Sym}(h^{(1),t,x,N} \otimes _0 (\cJ(h^{(3),t,x,N}) \otimes_r \cJ(h^{(3),t,x,N}) ))\Bigr], \qquad r=0,1,2,3\\
B_{5-2r} :={} &I_{5-2r} \Bigl[\text{Sym}(h^{(1),t,x,N}\otimes_1 (\cJ(h^{(3),t,x,N}) \otimes_r\cJ(h^{(3),t,x,N})))\Bigr], \qquad r=0,1,2.
\end{aligned}\end{equ}

The proof is reduced to proving that
$$\sup_{t \in [0,T]} \| \| A_{7-2r}\|_{H^{s-1}_x} \|_{L^2(\Omega)}<\infty,
\qquad \sup_{t \in [0,T]} \| \| B_{5-2r}\|_{H^{s-1}_x} \|_{L^2(\Omega)}<\infty.$$

To this end, we write these It\^o integrals in the following form:
$$I_k[G(t,x)]=\sum_n e^{i n\cdot x}\I_k[\sum_{ (n_1, \dots, n_k)\in \Lambda_n}  g^{(k)}(t)_{n_1, \dots, n_k} ],$$
where $\Lambda_n=\{(n_1, \dots, n_k): n=n_1+\dots+n_k\}$, 
then apply \eqref{sobolev-bound} from Lemma \ref{lem:symbols-in-Xsb}, to remove the iterated Wiener integrals. Namely,
\begin{equs}
{}&\Bigl\| \|\sum_{m \in \Z^2} e^{\i m \cdot x}\I_k[\sum_{ (m_1, \dots, m_k)\in \Lambda_n}  g^{(k)}(t)_{m_1, \dots, m_k} ]\|_{H_x^{s-1}}\Bigr \|_{L^p(\Omega)} \\
&\lesssim
 p^{\f k2} \|  \jpb{m}^{s-1} 
 \Bigl\|  \Bigl(  \sum_{(m_1, \dots, m_k)\in \Lambda_m}\int_{[0, \infty)^k} |g^{(k)}_{m_1,\dots, m_k}(\cdot)|^2 \d s\Bigr)^{\f 12}\Bigr\|_{\ell^2_m}.
\end{equs}
Then it remains to prove the right hand side is finite, after identifying the integrands $g^{(k)}$, i.e. showing  being finite.
These terms in term will be proved to be finite in Lemma \ref{lem:septic-sum} further below.

For the terms $A_{7-2r}$ and $B_{5-2r}$, we denote by $g^{7-2r}(t)$ and $k^{5-2r}(t)$, respectively, their integrands, which are sequences. We are not concerned with the symmetrization. The estimates would be independent of the the choice of pairings associated with the contraction,  The individual terms in the sequences are indexed by $m_k$ as follows:
\begin{equs}\label{g7}
g^{(7-2r)}_{m_1,\cdots m_{7-2r}}(t) &:=  \Bigl(h^{(1),t,N}\otimes_0 (\cJ(h^{(3),t,N}) \otimes_r \cJ(h^{(3),t,N}) \Bigr)_{m_1,\cdots m_{7-2r}},\\
k^{(5-2r)}_{m_1,\cdots m_{5-2r}}(t) &:=  \Bigl(h^{(1),t,N}\otimes_1 (\cJ(h^{(3),t,N}) \otimes_r \cJ(h^{(3),t,N}) \Bigr)_{m_1,\dots m_{5-2r}},
\label{k5}
\end{equs}

Therefore, for example, for $A_1$ when $r=3$, we have 
\begin{equs}
{}&g^{(1)}_{m} (t, s_4)\\
&=\1_{m=n_4} \int_0^t \int_0^t \int_0^t  \cJ(h_{n_1,n_2, n_3}^{(3),t,N}) (s_1,s_2,s_3) h_{n_4}^{t,N}(s_4) \cJ(h_{-n_1,-n_2, -n_3}^{(3),t,N}) )(s_1,s_2,s_3) \d s_1 \d s_2 \d s_3,
\end{equs}
where $n_1$ and $n_5$, $n_2$ and $n_6$, $n_3$ and $n_7$ are contracted.

{\bf Basic estimates.}  Recall~\eqref{eq:lolipop-cherry-cube}
\begin{equs}\<1>_N(t) &=I_1[h^{(1),t,x,N}] = \sum_{n \in \Z^2_N} e^{\i n \cdot x} \I_1[h^{(1),t}_n],\\
 \<3>_N(t) & =I_3[h^{(3),t,x,N}]= \sum_{n_1,n_2,n_3 \in \Z^2} e^{\i (n_1+n_2+n_3)\cdot x} \I_3[h^{(3),t,N}_{n_1,n_2,n_3}].
\end{equs}
where \begin{equs}
h^{(1),t,x,N}_{m}(s_1,s_2,s_3) &= e^{\i m\cdot x} \cdot  \1_{[0,t]}(s) \1_{\Z^2_N}(m) \f {\sin((t-s) \jpb{m})}{\jpb{m}^{1-\alpha}}\\
 h^{(3),t,x,N}_{n_1,n_2,n_3}(s_1,s_2,s_3) &= e^{\i (n_1+n_2+n_3)\cdot x} \cdot \prod_{j=1}^3 \1_{[0,t]}(s_j) \1_{\Z^2_N}(n_j) \f {\sin((t-s_j) \jpb{n_j})}{\jpb{n_j}^{1-\alpha}}.
 \end{equs}
Thus, $$\cJ(h^{(3),t,N})_{n_i,n_j,n_k}  =\cJ(h^{(3),t})_{n_i,n_j,n_k} \Pi_{\ell=i,j,k} \1_{\{n_\ell \in \Z^2_N\}},$$ where,
\begin{equ}
\cJ(h^{(3),t})_{n_i,n_j,n_k} (s_i,s_j,s_k) = \int_0^t \frac{\sin((t-s)\jpb{n_{ijk}})}{\jpb{n_{ijk}}}  \Pi_{\ell =i,j,k} \1_{[0,s]}(s_\ell)  \f {\sin((s-s_\ell)\jpb{n_\ell})}{\jpb{n_\ell}^{1-\alpha}} \, \d s.
\end{equ}
Define $$\phi(\bn) : = \phi(n_i,n_j,n_k) = \sigma \jpb{n_{ijk}} +\sigma_i\jpb{n_i} + \sigma_j \jpb{n_j} + \sigma_k\jpb{n_k}.$$
 Applying Lemma \ref{Euler} and simplify the notation we have:
\begin{equ}
\cJ_t(h^{(3)}_{n_i,n_j,n_k} )=\Psi(t,n_i,n_j,n_k) \Phi(t,n_i,n_j,n_k)
\end{equ}
where 
\begin{equs}\Psi(t,n_i,n_j,n_k) (s_i,s_j, s_k) 
&:=  \sum_{\substack{ \sigma,\sigma_i,\sigma_j \\ \sigma_l \in \{\pm 1\}}} e^{\i \sum_{\ell = i,j,k } s_\ell \sigma_\ell \jpb{n_\ell}} \prod_{\ell =i,j,k} \1_{[0,t']}(s_\ell)\\
\Phi(t,n_i,n_j,n_k) (s_i,s_j, s_k)& := \f 1 {16} \jpb{n_{ijk}}^{-1} \Bigl( \prod_{\ell =j,k,l} \jpb{n_\ell}^{-1+\alpha} \Bigr) e^{-\i t \sigma \jpb{n_{ijk}}} \int_{\max\{s_i,s_j,s_k\}}^t e^{-\i s \phi(\bn)} \, \d s.
\end{equs}
We use the trivial bound $|\Psi|\lesssim 1$,  and 
\begin{equ}
\Big| e^{-\i t \sigma \jpb{n_{ijk}}} \int_{\maxcurly{s_i,s_j,s_k}}^t e^{-\i s \phi(\bn)} \, \d s \Big| 
\lesssim \f 1{\jpb{\phi(\bn)}}, \quad \Big|\f{e^{\i t\sigma_4\jpb{n_4}}}{2\i \jpb{n_4}^{1-\alpha}} \Big| \lesssim \jpb{n_4}^{-1+\alpha}.\end{equ}
In particular,
\begin{equ}\label{h1-brutal-bound}
 |h^{(1),t,N}_{n_4}(s_4) | \le \1_{[0,t]}(s_4)\f 1 {\jpb{m}^{1-\alpha}}\1_{Z_N^2}(m),
 \end{equ}
and
\begin{equ}\label{J-h3}
|\cJ_t(h^{(3)}_{n_i,n_j,n_k} )|\lesssim \Bigl |\Phi(t,n_i,n_j,n_k) (s_i,s_j, s_k)\Bigr|
\lesssim \jpb{n_{ijk}}^{-1}\Bigl( \prod_{\ell =j,k,l} \jpb{n_\ell}^{-1+\alpha} \Bigr) \f 1{\jpb{\phi(\bn)}}.
\end{equ}
By decomposing the phase function, we have 
\begin{equ}\label{J-h3}
\begin{aligned}
{}&|\cJ_t(h^{(3)}_{n_i,n_j,n_k} )| \lesssim \Psi(n_i,n_j,n_k), \\
{}&|\cJ_t(h^{(3),t, N}_{n_1,n_2,n_3} )(s_s, s_2, s_3)| 
\lesssim  \Pi_{j=1}^3\1_{0,t}(s_j) \Psi(n_1,n_2,n_3) \Pi_{\ell=1}^3\1_{\{|n_\ell|\le N\}},
\end{aligned}
\end{equ}
where $$\Psi(n_i,n_j,n_k) :=  \sum_{m \in \Z} \f1{\jpb{m}} \jpb{n_{ijk}}^{-1}\Bigl( \prod_{\ell =j,k,l} \jpb{n_\ell}^{-1+\alpha} \Bigr)\1_{|\phi-m| \le 1}.$$

For clarity, we label the frequencies in the first factor $\<30>$ by $n_1,n_2,n_3$, and those in the second $\<30>$ by $n_5,n_6,n_7$ while the frequencies coming from $\<1>$ is labeled $n_4$. Below, we first use a case-by-case approach to reduce the seven terms appearing in $A_{7-2r}$ and $B_{5-2r}$ to essentially two cases, namely, whether $n_4$ is contracted. We then apply the basic estimate above and carry out a counting analysis to complete the proof.
\medskip

We return to show estimate $\| \|A_{7-2r}\|_{H^{s-1}_x} \|^2_{L^p(\Omega)}$ and the correspondign $B$ terms.

{\bf Case 1.}  We consider $$A_7 = I_7(h^{(1),t,x,N} \otimes_0  (\cJ(h^{(3),t,x,N}) \otimes_0 \cJ(h^{(3),t,x,N}) ).$$
 In this case, there is no contraction.  We have 
\begin{equs}
\| \|A_7\|_{H^{s-1}_x} \|^2_{L^p(\Omega)} \lesssim &p^7 \sum_{n \in \Z^2} \jpb{n}^{2s-2} \Bigl( \sum_{n=n_1+\cdots+n_7} \int_{[0, \infty)^7}  |g^{(7),t,N}_{n_1,n_2,\cdots,n_7}|^2 \, \d \boldsymbol{s} \Bigr) \\
\lesssim & p^7 T^7 \sum_{n_1,\cdots, n_7 \in \Z^2_N} \jpb{n_{1234567}}^{2s-2} \Psi^2(n_1,n_2,n_3) \jpb{n_4}^{-2+2\alpha}  \Psi^2(n_5,n_6,n_7) 
\end{equs}
where in the last step we used the estimate from \eqref{J-h3}. By definition, c.f. \eqref{g7},
\begin{equ}
g^{(7),t,N}_{n_1,n_2,\cdots,n_7} =  \cJ_t(h^{(3),t,N}_{n_1,n_2,n_3}) h^{(1),t,N}_{n_4} \cJ_t(h^{(3),t,N}_{n_5,n_6,n_7}),
\end{equ}
using the above estimate we obtain:
\begin{equ}
|g^{(7),t,N}_{n_1,n_2,\cdots,n_7}| \le \prod_{j=1}^7 \1_{\Z^2_N}(n_j) \jpb{n_4}^{-1+\alpha}\Psi(n_1,n_2,n_3) \Psi(n_5,n_6,n_7).
\end{equ} 
If we insert dyadic decomposition with respect to frequencies $n_{1234567},n_{123},n_4$ and $n_{567}$ and denote the resulting blocks by $A^{\mathrm{dyadic}}_7$, we obtain 
\begin{equs}
\| \|A^{\mathrm{dyadic}}_7\|_{H^{s-1}_x} \|^2_{L^p(\Omega)}
& \lesssim p^7 T^7 \sum_{n_1,\cdots, n_7 \in \Z^2_N} \jpb{n_{1234567}}^{2s-2} \mathbf{1}_{|n_{1234567}| \sim N_{1234567}}\mathbf{1}_{|n_{123} \sim N_{123}|} \1_{|n_{567} \sim N_{567}|} \\
&\quad \1_{|n_4| \sim N_4} \Psi^2(n_1,n_2,n_3) \jpb{n_4}^{-2+2\alpha}  \Psi^2(n_5,n_6,n_7). \label{A7}
\end{equs}
The final estimates on the right hand side will be carried out after completing the case analysis.

\medskip

{\bf Case 2.} We next consider $$A_5 = I_5(h^{(1),t,x,N} \otimes_0  (\cJ(h^{(3),t,x,N}) \otimes_1 \cJ(h^{(3),t,x,N}) ).$$
In this case, there is one contraction between $\<30>_N$ and $\<30>_N$; for definiteness, let  $n_1$ be contracted with $n_5$. By the definition of $g^{(5),t,N}_{n_2,n_3,n_4,n_6,n_7} $ and the definition of tensors, c.f. Definition \ref{def:tensor}, we sum over the contracted index $n_1$, set $s_1=s_5$, and integrate over $s_1$ to arrive at:
 $$g^{(5),t,N}_{n_2,n_3,n_4,n_6,n_7}=\sum_{n_1\in \Z^2_N} \int_0^t \cJ_t(h^{(3),t,N}_{n_1,n_2,n_3}) (s_1, \cdot, \cdot) h^{(1),t,N}_{n_4}(\cdot)  \cJ_t(h^{(3),t,N}_{-n_1,n_6,n_7})(s_1, \cdot, \cdot) \, \d s_1. $$
 We now use the bound (\ref{h1-brutal-bound}) fro $ |h^{(1),t,N}_{n_4}(s_4) | $
  and the bound in \eqref{J-h3} for the $h^{(3)}$ term, to obtain that for $t \in [0,T]$,
 \begin{equs}
{}& | |g^{(5),t,N}_{n_2,n_3,n_4,n_6,n_7}| \\
 &\le \Pi_{j=2,3,4,6,7} \1_{[0,t]}(s_j)\f 1 {\jpb{n_4}^{1-\alpha}}\1_{Z_N^2}(n_4) \sum_{n_1}\Pi_{\ell=2,3,6,7}\1_{Z_N^2}(n_\ell)\Psi(n_1,n_2,n_3)\Psi(n_1,n_6,n_7).
  \end{equs}
 
Taking supremum bound in time and writing $n_{23467}=n_2+n_3+n_4+n_6+n_7$, we obtain:
\begin{equs}\label{A5}
{}& \| \|A_5\|_{H^{s-1}_x} \|^2_{L^p(\Omega)}
 \lesssim p^5 \sum_{n \in \Z^2} \jpb{n}^{2s-2} \Bigl( \sum_{n_2,n_3, n_4, n_6, n_7 | n=n_{23467}} \int_{[0, \infty)^5}  |g^{(5),t,N}_{n_2,n_3,n_4,n_6,n_7}|^2 \, \d \boldsymbol{s} \Bigr) \\
\lesssim & p^5 T^5 \sum_{n_2,n_3,n_4,n_6,n_7 \in \Z^2_N} 
\jpb{n_{23467}}^{2s-2} \jpb{n_4}^{-2+2\alpha}}\Big( {\sum_{n_1}\Psi(n_1,n_2,n_3)\Psi(n_1,n_6,n_7) \Big)^2.
\end{equs}
In the above, $\d \boldsymbol{s} $ denotes $ \d s_2 \d s_3 \d s_4 \d s_6 \d s_7$.

\medskip

{\bf Case 3.}  Consider $A_3 = I_3(h^{(1),t,x,N} \otimes_0  (\cJ(h^{(3),t,x,N}) \otimes_2 \cJ(h^{(3),t,x,N}) )$. We may assume that $n_1$ is contracted with $n_5$ and $n_2$ with $n_6$. Then, for $t \in [0,T]$, we obtain
\begin{equs}
{}&\| \|A_3\|_{H^{s-1}_x} \|^2_{L^p(\Omega)}\\
& \lesssim  p^3T^7 \sum_{n_3,n_4,n_7 \in \Z^2_N} \jpb{n_{347}}^{2s-2} \Big( \sum_{n_1, n_2\in \Z^2_N} \Psi(n_1,n_2,n_3) \jpb{n_4}^{-1+\alpha} \Psi(n_1,n_2,n_7) \Bigr)^2\\
&= p^3T^7 \sum_{n_3,n_4,n_7 \in \Z^2_N} \jpb{n_{347}}^{2s-2}  \jpb{n_4}^{-2+2\alpha} \Big( \sum_{n_1, n_2\in \Z^2_N} \Psi(n_1,n_2,n_3) \Psi(n_1,n_2,n_7) \Bigr)^2.
\label{A3}
\end{equs}
\bigskip

{\bf Case 4.} The final $A$-term is $A_1 = I_1(h^{(1),t,x,N} \otimes_0  (\cJ(h^{(3),t,x,N}) \otimes_3 \cJ(h^{(3),t,x,N}) )$. In this case, $n_4$ is the only index that is not contracted. Hence, for $t \in [0,T]$, we have 
\begin{equs}\label{A1}
\| \|A_3\|_{H^{s-1}_x} \|^2_{L^p(\Omega)} \lesssim & p^1T^7  \sum_{n_4 \in \Z^2_N} \jpb{n_4}^{2s-4+2\alpha}   \Big( \sum_{n_1, n_2,n_3\in \Z^2_N} \Psi^2(n_1,n_2,n_3) \Bigr)^2.
\end{equs}

{\bf Case 5.} The first $B$-term to be considered is $$B_5 = I_5 ( (h^{(1),t,x,N}\otimes_1 (\cJ(h^{(3),t,x,N}) \otimes_0 \cJ(h^{(3),t,x,N}))).$$  In this case, there is one contraction and we may assume that $n_4$ is contracted with $n_1$, i.e. we set  $n_4=-n_1$, the integral kernel $k$ defined by \eqref{k5}
has the form:
\begin{equ}
k^{(5),t,N}_{n_2,n_3,n_5,n_6,n_7} 
= \sum_{n_1 \in \Z^2} \int_0^t \cJ_t(h^{(3),t,N}_{n_1,n_2,n_3}) (s_1, \cdot, \cdot)
h^{(1),t,N}_{-n_1}(s_1) \cJ_t(h^{(3),t,N}_{n_5,n_6,n_7}) (\cdot, \cdot, \cdot)\, \d s_1.
\end{equ}
Apply again (\ref{h1-brutal-bound}) and (\ref{J-h3}) to obtain:
\begin{equs}\label{B5}
\| \|B_5\|_{H^{s-1}_x} \|^2_{L^p(\Omega)} 
\lesssim & p^5T^5 \sum_{n \in \Z^2} \jpb{n}^{2s-2} \Big( \sum_{n=n_2+n_3+n_5+n_6+n_7} \int_{[0,\infty)^5} |k^{(5),t,N}_{n_2,n_3,n_5,n_6,n_7}|^2 \d \boldsymbol{s} \Bigr)\\
\lesssim & p^5T^5 \sum_{n_2,n_3,n_5,n_6,n_7 \in \Z^2_N} \jpb{n_{23567}}^{2s-2}  \Psi^2(n_5,n_6,n_7) \Big( \sum_{n_1\in \Z^2_N} \Psi(n_1,n_2,n_3) \jpb{n_1}^{-1+\alpha} \Bigr)^2
\end{equs}

{\bf Case 6.} Similarly, for $$B_3 = I_3 ( (h^{(1),t,x,N}\otimes_1 (\cJ(h^{(3),t,x,N}) \otimes_1 \cJ(h^{(3),t,x,N}))), $$ we have the bound
\begin{equs}\label{B3}
{}&\| \|B_3\|_{H^{s-1}_x} \|^2_{L^p(\Omega)} \\
&\lesssim p^3T^5 \sum_{n_3, n_6, n_7 \in \Z^2_N} \jpb{n_{367}}^{2s-2}   \Big( \sum_{n_1,n_2\in \Z^2_N} \Psi(n_1,n_2,n_3) \jpb{n_1}^{-1+\alpha} \Psi(n_2,n_6,n_7)\Bigr)^2
\end{equs}

{\bf Case 7.} For the final term $B_1 = I_1 ( (h^{(1),t,x,N}\otimes_1 (\cJ(h^{(3),t,x,N}) \otimes_2 \cJ(h^{(3),t,x,N})))$, we have
\begin{equs}\label{B1}
{}&\| \|B_1\|_{H^{s-1}_x} \|^2_{L^p(\Omega)} 
\\&
\lesssim p^1T^5 \sum_{n_7 \in \Z^2_N} \jpb{n_6}^{2s-2}   \Big( \sum_{n_1,n_2,n_3\in \Z^2_N} \Psi(n_1,n_2,n_3) \jpb{n_1}^{-1+\alpha} \Psi(n_2,n_3,n_7)\Bigr)^2. 
\end{equs}

This finishes the preliminary estimates on all severn terms. 

To apply Lemma~\ref{lem:quintic-non-resonance}, we introduce the notion of pairing, with which all terms can be expressed in one unified formula. All cases (except the empty paring case) is depicted  in the
 figure below.
In a nutshell, by Lemma~\ref{lem:quintic-non-resonance},  there exists $\epsilon, C>0$ such that the dyadic scale version of the right hand sides of (\ref{A7}-\ref{B1})
is bounded by 
$$C\maxcurly{N_{1234567},N_1,\cdots,N_7}^{-\epsilon}.$$
 the preliminary bounds  obtained in each case is controlled by
$$\maxcurly{N_{1234567},N_1,\cdots,N_7}^{-\epsilon},$$
where $\epsilon$ is a suitable small positive constant.
Hence, for any $s<1-\alpha$ and $\f 14 \le \alpha < \f 3 8$, we have
\begin{equ}
\| \| \<7>^{\mathrm{dyadic}}_N\|_{H^{s-1}_x} \|^2_{L^p(\Omega)} \lesssim p^7T^5 \maxcurly{N_{1234567},N_{123},N_{4},N_{567}}^{-\epsilon}.
\end{equ}

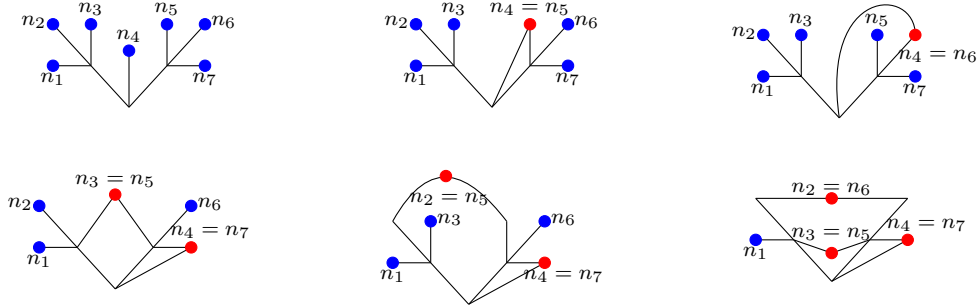
\begin{center}
    \begin{minipage}{0.3\textwidth}
        \centering
   \begin{tikzpicture}[scale=0.25]  
    \draw(6,0)--(6,3);         
    \draw[fill,blue](6,3) circle(.3);     
    \scriptsize\node at (6,3.8) {$n_4$};          
    \draw(6,0)--(10,4.4);       
    \draw[fill,blue](10,4.4) circle(.3);    
    \scriptsize\node at (11,4.5) {$n_6$};        
    \draw(8,2.2)--(8,4.4);    
    \draw[fill,blue](8,4.4) circle(.3);    
    \scriptsize\node at (8,5.2) {$n_5$};       
    \draw(8,2.2)--(10,2.2);    
    \draw[fill,blue](10,2.2) circle(.3);    
    \scriptsize\node at (10,1.5) {$n_7$}; 
    \draw(6,0)--(2,4.4);      
    \draw[fill,blue](2,4.4) circle(.3);   
    \scriptsize\node at (1,4.5) {$n_2$};   
    \draw(4,2.2)--(2,2.2);      
    \draw[fill,blue](2,2.2) circle(.3);   
    \scriptsize\node at (2,1.5) {$n_1$}; 
    \draw(4,2.2)--(4,4.4);      
    \draw[fill,blue](4,4.4) circle(.3);   
    \scriptsize\node at (4,5.2) {$n_3$}; 
\end{tikzpicture}
    \end{minipage}
    \begin{minipage}{0.3\textwidth}
        \centering
        \begin{tikzpicture}[scale=0.25] 
    \draw(6,0)--(8,4.4);                  
    \draw(6,0)--(10,4.4);       
    \draw[fill,blue](10,4.4) circle(.3);    
    \scriptsize\node at (11,4.5) {$n_6$};        
    \draw(8,2.2)--(8,4.4);    
    \draw[fill,red](8,4.4) circle(.3);    
    \scriptsize\node at (8,5.2) {$n_4=n_5$};       
    \draw(8,2.2)--(10,2.2);    
    \draw[fill,blue](10,2.2) circle(.3);    
    \scriptsize\node at (10,1.5) {$n_7$}; 
    \draw(6,0)--(2,4.4);      
    \draw[fill,blue](2,4.4) circle(.3);   
    \scriptsize\node at (1,4.5) {$n_2$};   
    \draw(4,2.2)--(2,2.2);      
    \draw[fill,blue](2,2.2) circle(.3);   
    \scriptsize\node at (2,1.5) {$n_1$}; 
    \draw(4,2.2)--(4,4.4);      
    \draw[fill,blue](4,4.4) circle(.3);   
    \scriptsize\node at (4,5.2) {$n_3$}; 
\end{tikzpicture}
    \end{minipage}
    \begin{minipage}{0.3\textwidth}
     \centering
    \begin{tikzpicture}[scale=0.25]  
\draw (6,0) .. controls (5,7) and (10,7) .. (10,4.4);                  
    \draw(6,0)--(10,4.4);       
    \draw[fill,red](10,4.4) circle(.3);    
    \scriptsize\node at (11.2,3.3) {$n_4=n_6$};        
    \draw(8,2.2)--(8,4.4);    
    \draw[fill,blue](8,4.4) circle(.3);    
    \scriptsize\node at (8,5.2) {$n_5$};       
    \draw(8,2.2)--(10,2.2);    
    \draw[fill,blue](10,2.2) circle(.3);    
    \scriptsize\node at (10,1.5) {$n_7$}; 
    \draw(6,0)--(2,4.4);      
    \draw[fill,blue](2,4.4) circle(.3);   
    \scriptsize\node at (1,4.5) {$n_2$};   
    \draw(4,2.2)--(2,2.2);      
    \draw[fill,blue](2,2.2) circle(.3);   
    \scriptsize\node at (2,1.5) {$n_1$}; 
    \draw(4,2.2)--(4,4.4);      
    \draw[fill,blue](4,4.4) circle(.3);   
    \scriptsize\node at (4,5.2) {$n_3$}; 
\end{tikzpicture}
\end{minipage}
\end{center}
\begin{center}
    \begin{minipage}{0.3\textwidth}
        \centering
 \begin{tikzpicture}[scale=0.25] 
   \draw(6,0)--(10,2.2);             
    \draw(6,0)--(10,4.4);       
    \draw[fill,blue](10,4.4) circle(.3);    
    \scriptsize\node at (11,4.5) {$n_6$};        
    \draw(8,2.2)--(6,5);    
    \scriptsize\node at (6,5.7) {$n_3=n_5$};       
    \draw(8,2.2)--(10,2.2);    
    \draw[fill,red](10,2.2) circle(.3);    
    \scriptsize\node at (11,3) {$n_4=n_7$}; 
    \draw(6,0)--(2,4.4);      
    \draw[fill,blue](2,4.4) circle(.3);   
    \scriptsize\node at (1,4.5) {$n_2$};   
    \draw(4,2.2)--(2,2.2);      
    \draw[fill,blue](2,2.2) circle(.3);   
    \scriptsize\node at (2,1.5) {$n_1$}; 
    \draw(4,2.2)--(6,5);      
    \draw[fill,red](6,5) circle(.3);   
\end{tikzpicture}
    \end{minipage}
    \begin{minipage}{0.3\textwidth}
        \centering
        \begin{tikzpicture}[scale=0.25] 
 \draw (2,4.4) .. controls (4,8) and (6,7) .. (8,4.4);
     \draw[fill,red](4.8,6.8) circle(.3); 
        \scriptsize\node at (5,5.6) {$n_2=n_5$};
    \draw(6,0)--(10,2.2);                   
    \draw(6,0)--(10,4.4);       
    \draw[fill,blue](10,4.4) circle(.3);    
    \scriptsize\node at (11,4.5) {$n_6$};        
    \draw(8,2.2)--(8,4.4);               
    \draw(8,2.2)--(10,2.2);    
    \draw[fill,red](10,2.2) circle(.3);    
    \scriptsize\node at (11,1.5) {$n_4=n_7$}; 
    \draw(6,0)--(2,4.4);          
    \draw(4,2.2)--(2,2.2);      
    \draw[fill,blue](2,2.2) circle(.3);   
    \scriptsize\node at (2,1.5) {$n_1$}; 
    \draw(4,2.2)--(4,4.4);      
    \draw[fill,blue](4,4.4) circle(.3);   
    \scriptsize\node at (5,4.5) {$n_3$}; 
\end{tikzpicture}
    \end{minipage}
    \begin{minipage}{0.3\textwidth}
     \centering
   \begin{tikzpicture}[scale=0.25]  
    \draw(2,4.4)--(10,4.4);
    \draw(6,0)--(10,2.2);                
    \draw(6,0)--(10,4.4);               
    \draw(8,2.2)--(6,1.5);    
\scriptsize\node at (6,2.5) {$n_3=n_5$};       
    \draw(8,2.2)--(10,2.2);    
    \draw[fill,red](10,2.2) circle(.3);    
    \scriptsize\node at (11,3) {$n_4=n_7$}; 
    \draw(6,0)--(2,4.4);        
    \draw(4,2.2)--(2,2.2);      
    \draw[fill,blue](2,2.2) circle(.3);   
    \scriptsize\node at (2,1.5) {$n_1$}; 
    \draw(4,2.2)--(6,1.5);      
    \draw[fill,red](6,1.5) circle(.3);   
    \draw[fill,red](6,4.4) circle(.3); 
\scriptsize\node at (6,4.9) {$n_2=n_6$}; 
\end{tikzpicture}
\end{minipage}
\end{center}
\begin{center}
\centering
Pairs in $\<30>_N\<1>_N\<30>_N$.
\end{center}

Finally, we explain how to related to Lemma~\ref{lem:quintic-non-resonance} and thus complete the proof of Proposition~\ref{prop:reg-septic}.  To thsi end, we first introduce the notion of \emph{pairing} \cite[pp74]{Bri20b}
and  convert the right hand side of  the bounds given in (\ref{A7}-\ref{B1}) in the language of Lemma~\ref{lem:quintic-non-resonance}, with which we may conclude the proof of the proposition.

 \begin{definition}\label{pairing}
Let $J\ge 1$ and consider $ M=\{1,2,\cdots, J \}$.  A relation $\cP$, i.e. a subset of $M\times M$, is called a pairing if the following conditions hold
\begin{itemize}
\item $\cP$ is anti-reflexive: $(j,j)\notin \cP$ for all $1 \le j \le J$;
\item $\cP$ is symmetric: $(i,j)\in \cP$ if and only if $(j,i)\in \cP$;
\item $\cP$ is univalent: for every $1\le i \le J$, there is at most one $1 \le j \le J$ such that $(i,j)\in \cP$.
\end{itemize}
If $(i,j)\in \cP$, the tuple $(i,j)$ is called a pair. If  $j\in M$  is contained in a pair, we say it is paired and otherwise unpaired.
Furthermore, let $Q=(Q_\ell)_{1\le \ell \le L}$ be a partition of the set $\{1,2,\cdots, J \}$. We say that the pairing $\cP$ is with respect to $Q$ if no pairing in $\cP$ belongs to the same block of the partition.
\end{definition}

 \begin{notation}
 If $\cP\subset M\times M$ is a pairing, we denote by $\pi (\cP)\subset \{1,\cdots,7\}$ the collection of unaired numbers (it is literally the projection of the pairing $\cP$, considered as a set, to one of its first coordinate). Furthermore let $\pi^c(\cP)=\{1,\cdots, 7\}\setminus \pi(\cP)$. In case $\cP$ is fixed we simply use the notation $\pi$ and $\pi^c$.
 
 The notation $\sum_{n_j :, j \in \pi^c}$ denotes summation over all $n_j \in \mathbb{Z}^2$ with indices $j$ ranging through the set $\pi^c$. On the other hand,  the notation $\sum_{j\in \pi}$ popularly used in literature represents summation runs through all pairs in $\cP$: If $i\in \pi$ and $j$ is the unique element such that $(i,j)\in \cP$, we set $n_j=-n_i$ in the summand and sum over $n_i$. So 
$$\sum_{n_j: j\in \pi}F(n_1, \dots, n_7):=\sum_{n_j: (i,j)\in \cP, n_i=-n_j}F_{i,j}(n_1, \dots, n_7)$$
 where $F_{i,j}$ is obtained with $n_j$ replaced by $-n_i$.
\end{notation}

We now translate the terms in the above proposition in the language of summing over paring and over unpaired indices. For this, we set $J=7$ and consider the partition $$Q=\{ 1, 2, 3\}, \; \{4\},\; \{5,6,7\}.$$
Let 
\begin{equation}\label{sum-over-paring}
n_{\pi^c(\cP)} = \sum_{n_j :  j \in \pi^c(\cP)} n_j.
\end{equation}
Then it is easy to check that the sum of the right hand side of the seven terms in question is bounded by the term below 
\begin{equ}\label{xx}
\begin{aligned}
& \sum_{\cP}  \sum_{n_j:  j \in\pi^c(\cP)} \jpb{n_{\pi^c(\cP)}}^{2s-2} \Bigl( \sum_{n_i : (i,j)\in \cP, n_j=-n_i} \1_{|n_{1234567}| \sim N_{1234567}}  \prod_{j=1}^7 \1_{|n_j| \sim N_j} \Psi(n_1,n_2,n_3) \jpb{n_4}^{-1+\alpha} \Psi(n_5,n_6,n_7) \Bigr)^2.
\end{aligned}
\end{equ}
We illustrate some cases here for clarification:
The term $A_7$ corresponding to $\cP=\emptyset$, i.e. there is no contraction of indices. The summation in this case is over all indices, the tsummand in question is
\begin{equs}
{}&A_7 \lesssim \sum_{n_1,\cdots, n_7 \in \Z^2_N} \jpb{n_{1234567}}^{2s-2} \mathbf{1}_{|n_{1234567}| \sim N_{1234567}}\prod_{j=1}^7 \1_{|n_j| \sim N_j}  \Psi^2(n_1,n_2,n_3) \jpb{n_4}^{-2+2\alpha}  \Psi^2(n_5,n_6,n_7). 
\end{equs}
An upper bound for \ref{xx} is given in Lemma~\ref{lem:septic-sum} of Section \ref{sec:counting}, which allows us to conclude the lemma. The proof for Lemma~\ref{lem:septic-sum} has been divided into 7 types of parings, which corresponds to the severn cases considered above.
With this, we have completed the entire proof of Proposition~\ref{prop:reg-septic}. 
\end{proof}

We  may conclude this section.

\subsection{Cubic-cubic-cubic symbol $\<9>_N$ and $\<90>_N$}
If we attempt to use trilinear estimates to control $\cJ(\<30>^3)$, it becomes difficult to achieve this directly using trilinear estimate when $\alpha \ge \f 3 {10}$. This difficulty arises because the regularity of $v$ in \eqref{eq:2nd-v} is $2\alpha-\f 14 +$, which is higher than the regularity, $\f 54 -3\alpha-$, of $\<30>$. Consequently, we must instead define a new cubic-cubic-cubic symbol. Similar to the septic symbol, we do not seek the optimal regularity for this symbol. To obtain a local solution to equation \eqref{eq:2nd-v},  it suffices to work with regularity $2\alpha-\f 14+$. 
Set
$$\<9>_N = (\<30>_N)^3.$$
A naive regularity count suggests that  $\<9>_N$ has regularity $\f 94-3\alpha-$, which is already significantly larger than the threshold $2\alpha-\f 14 $ whenever $\alpha < \f 12$.

\begin{proposition}\label{prop:reg-cubic-cubic-cubic}
Assume $s<\f 94-3\alpha$ and $\f 14 \le \alpha < \f 5{12}$. 
Then, for any $b > \f 12$, we have
\begin{equ}
\| \| \<9>_N\|_{\X^{s-1,b-1}_T} \|_{L^p(\Omega)} \lesssim p^{\f 9 2} T.
\end{equ}
\end{proposition}
\begin{proof}
By inhomogeneous Strichartz estimate \ref{lem:inhomo-Stri} and Gaussian hyper-contractivity, we have 
\begin{equ}
\| \| \<9>_N\|_{\X^{s-1,b-1}_T} \|_{L^p(\Omega)} \lesssim p^{\f92} \sup_{t \in [0,T]} \| \| \<9>_N\|_{H^{s-1}_x} \|_{L^2(\Omega)}.
\end{equ}
Denote by $\cP$  a pairing of $ \{ 1,2,3,4,5,6,7,8,9 \}$ with respect to the partition 
\begin{equation}
\{ 1, 2, 3\}, \quad \{4,5,6\}, \quad \{7,8,9\},
\end{equation}
which means that no pairing in $\cP$ belongs to the same block of the partition, see Definition~\ref{pairing}Define:
\begin{equation}
\Psi(n_i,n_j,n_k) :=  \sum_{m \in \Z} \jpb{m}^{-1} \jpb{n_{ijk}}^{-1} \Pi_{\ell =i,j,k} \jpb{n_\ell}^{\alpha-1} \mathbf{1}_{|\phi-m| \leq 1}.
\end{equation}

Repeat the argument in Proposition \ref{prop:reg-septic}, denoting $\pi(\cP)$ the projection of $\CP$ to one of its factors, we have 
\begin{equs}
\| \| \<9>_N\|_{H^{s-1}_x} \|^2_{L^2(\Omega)} \lesssim & T^\theta  \sum_{n_j, j\in\pi^c(\cP)} \jpb{(n_{\pi^c})^{2s-2}} \Bigl( \sum_{n_i : (i,j)\in \cP, n_j=-n_i} \1_{|n_j| \sim N_j} \Psi(n_1,n_2,n_3) \Psi(n_4,n_5,n_6) \Psi(n_7,n_8,n_9) \Bigr)^2,
\end{equs}
for some $\theta >0$. 

Consider the dyadic piece by inserting Littlewood-Paley decomposition with respect to frequencies $n_{1234567},n_1,\cdots,n_9$, we have 
\begin{equs}
\| \| \<9>^{\mathrm{dyadic}}_N\|_{H^{s-1}_x} \|^2_{L^2(\Omega)} \lesssim  
& T^\theta N^{2s-2}_{123456789}  \sum_{n_j, j \in\pi^c(\cP)}
\Bigl(\sum_{n_i : (i,j)\in \cP, n_j=-n_i} \1_{|n_{123456789}| \sim N_{123456789}}\\
& \quad \prod_{j=1}^9 \1_{|n_j| \sim N_j} \Psi(n_1,n_2,n_3) \Psi(n_4,n_5,n_6) \Psi(n_7,n_8,n_9) \Bigr)^2
\end{equs} 
where we used the fact $|n_{\pi^c(\cP)}| \sim N_{123456789}$. 

Note the $(1,2,3),(4,5,6)$ and $(7,8,9)$ are symmetric in the sense as the case where $4$ is not paired in the proof of Lemma \ref{lem:septic-sum}. Hence, by applying Cauchy-Schwarz inequality and relabelling the indices, we have 
\begin{equs}
\| \| \<9>^{\mathrm{dyadic}}_N\|_{H^{s-1}_x} \|^2_{L^2(\Omega)}
& \lesssim  T^\theta 
(N_{123456789})^{2s-2}   \Bigl( \sum_{n_1,n_2,n_3 \in \Z^2} \prod_{j=1}^3 \1_{|n_j| \sim N_j}  \Psi^2(n_1,n_2,n_3) \Big)\\
 & \times \Bigl( \sum_{n_4,n_5,n_6 \ in \Z^2} \prod_{j=4}^6 \1_{|n_j| \sim N_j}  \Psi^2(n_4,n_5,n_6) \Big)  \Bigl( \sum_{n_7,n_8,n_9 \in \Z^2} \prod_{j=7}^9 \1_{|n_j| \sim N_j}  \Psi^2(n_7,n_8,n_9) \Bigr).
\end{equs}

By cubic sum estimate \eqref{cubic-repeat}, we have 
\begin{equs}
{}&\| \| \<9>^{\mathrm{dyadic}}_N\|_{H^{s-1}_x} \|^2_{L^2(\Omega)} \\
\lesssim & T^\theta (N_{123456789})^{2s-2}  \maxcurly{N_1,N_2,N_3}^{-2s_\alpha+\gamma}\maxcurly{N_4,N_5,N_6}^{-2s_\alpha+\gamma} \maxcurly{N_7,N_8,N_9}^{-2s_\alpha+\gamma},
\end{equs}
where $s_\alpha$ defined as \eqref{eq:s-alpha}, and for some $0<\gamma  \ll 1$ which comes from absorbing factor $\log^2(2+\maxcurly{N_i,N_j,N_k})$.

Let $s<\f 94-3\alpha$, and choose $\epsilon = \mincurly{ \f 94-3\alpha - s, 2s_\alpha-\gamma}$, such that 
\begin{equ}
\| \| \<9>^{\mathrm{dyadic}}_N\|_{H^{s-1}_x} \|_{L^2(\Omega)} \lesssim  T^\theta \maxcurly{N_{123456789},N_1,\cdots,N_9}^{-2\epsilon},
\end{equ}
which is summable on dyadic scale $N_1,\cdots,N_9,N_{123456789}$.
\end{proof}

\subsection{Quartic symbol $\<31>$}
This section addresses the quartic term, for which we again employ the product formula for multiple Wiener integrals, as we did for quintic term. Recall the definitions of $\<1>_N$ and $\<30>_N$ from \eqref{eq:chicken-foot}. Define
$$\<31>_N := \<1>_N \<30>_N. $$
By applying the product formula \eqref{product}, we obtain
\begin{equs}
\<31>_N  &= I_1[h^{(1),t,x,N}] I_3[h^{(3),t,x,N,\cJ}]\\
&= I_4[h^{(1),t,x,N} \otimes_0 h^{(3),t,x,N,\cJ}] + 3I_2[h^{(1),t,x,N} \otimes_1 h^{(3),t,x,N,\cJ}],
\end{equs}
where the coefficients arise from the combinatorial structure of the product formula. We refer to the two terms above as the non-resonant term and the single-resonant term, respectively, since they correspond to distinct types of contractions. We denote these components by$\<31>^{(4)},\<31>^{(2)}$. 

More explicitly, we have
\begin{equs}
\<31>^{(4)}_N(t) =&\sum_{n_1,\cdots, n_4 \in \Z^2_N} e^{\i (n_1+\cdots n_4)\cdot x} \I_4[k^{(4),t}_{n_1,\cdots,n_4}] \\
\<31>^{(2)}_N(t) =&\sum_{n_3,n_4  \in \Z^2_N} e^{\i (n_3+ n_4)\cdot x} \I_2[k^{(2),t}_{n_3,n_4}] \\
\end{equs}
where
\begin{equs}
k^{(4),t}_{n_1,\cdots,n_4}(s_1, \cdots,s_4) = & \1_{[0,t]}(s_1) \f {\sin((t-s_1) \jpb{n_1})}{\jpb{n_1}^{1-\alpha}} \\
&\quad \times \int_0^t \f {\sin((t-s)\jpb{n_{234}})}{\jpb{n_{234}} }\prod_{j=2}^4 \1_{[0,s]}(s_j) \f {\sin((s-s_j) \jpb{n_j})}{\jpb{n_j}^{1-\alpha}} \, \d s, \label{k4} \\
k^{(2),t}_{n_3,n_4}(s_3, s_4) = &   \sum_{n_2 \in \Z^2_N} \int_0^t \f {\sin((t-s)\jpb{n_{234}})}{\jpb{n_{234}} }\\
&\quad \times \prod_{j=3}^4 \1_{[0,s]}(s_j) \Bigl( \int_0^s  \frac{\sin ((t-s_2)\jpb{n_2}) \sin ((s-s_2)\jpb{n_2})}{ \jpb{n_2}^{2-2\alpha}}  \, \d s_2 \Bigr) \, \d s. \label{k2}  \\
\end{equs}

\begin{proposition}\label{prop:reg-of-quartic}
Let $\f 14 \le \alpha <\f 38$. For any  $s <-\alpha$,
 we have
\begin{equ}
\| \| \<31>_N\|_{L^\infty_tW^{s,\infty}_x([0,T] \times \T^2)} \|_{L^p(\Omega)}  \lesssim p^2 T^2.
\end{equ}
Furthermore,  $\{ \<31>_N \}$ is a Cauchy sequence in  $L^p(\Omega; L^\infty_tW^{s,\infty}_x([0,T] \times \T^2))$.
\end{proposition}

\begin{proof}
We divide our proof in several steps as follow.
\begin{itemize}
\item[Step 1.] Reduce $L^\infty_tW^{s,\infty}_x$ to $L^\infty_tH^{s'}_x$-norm. To compute the stochastic symbols in $C_tW^{s,\infty}_x$-space, it is straightforward to involve the Kolmogorov's type of argument. However, we encounter 
some issue near some t$\alpha$. Let $f(t,x,\omega)$ be a random field from $[0,T] \times \T^2 \times \Omega$ to $\C$. Take $2 \le q \le \infty$ to be sufficient large, by Sobolev embedding both in $t$ and $x$ and Bernstein's inequality, we have
\begin{equ}
\| P_N f \|_{L^\infty_tW^{s,\infty}_x([0,T] \times \T^2)} \lesssim N^{\f 3 q} \| \jpb{\nabla}^s P_N f\|_{L^q_tL^q_x([0,T] \times \T^2))}.
\end{equ}
where $P_N$ is the Littlewood projector at frequency scale $N$. If $f$ is a Gaussian process which is spatial stationary, then $ \jpb{\nabla_x}^{s} f$ remains stationary which leads to 
\begin{equs}
\| \| P_N \jpb{\nabla_x}^s f \|_{L^\infty_tW^{s,\infty}_x([0,T] \times \T^2)} \|_{L^p(\Omega)} \le & C_p \|  \|P_N \jpb{\nabla_x}^sf\|_{L^2(\Omega)L^q_x(\T^2)} \|_{L^q_t([0,T] )} \\
= & C_p \|  \|P_N \jpb{\nabla_x}^s f\|_{L^2(\Omega)L^2_x(\T^2)} \|_{L^q_t([0,T] )} \\
\le & C_p T^{\f 1 q} \|  \|P_N \jpb{\nabla_x}^s f\|_{L^2_x(\T^2)L^2(\Omega)} \|_{L^\infty_t([0,T])} 
\end{equs}
where $p \ge q$ and we applied Minkowski's integral inequality, $\|P_N f\|_{L^2(\Omega)}$ is constant in $x$, as well as H\"{o}lder's inequality in $t$. As a conclusion, we obtain
\begin{equ}\label{eq:Wsinfty-to-Hs}
\| \| P_N f \|_{L^\infty_tW^{s,\infty}_x([0,T] \times \T^2)} \|_{L^p(\Omega)}\lesssim C_pN^{\f 3 q} T^{\f 1q} \|  \|P_N f\|_{H^s_x(\T^2)L^2(\Omega)} \|_{L^\infty_t([0,T])}.
\end{equ}

\item[Step 2.] Analysis of non-resonance term $\<31>^{(4)}$. By \eqref{eq:Wsinfty-to-Hs}, for fixed $t \in [0,T]$, we have 
\begin{align}\label{eq:k4-in-Sobolev}
\| \|\<31>^{(4)}_N \|_{H^{s}_x}\|^2_{L^2(\Omega)} =  &\Big\| \E\Big[ \sum_{\substack{n=n_{1234},n_1,\cdots,n_4\in\Z^2_N  }} \I_4[k^{(4),t}_{n_1,\cdots,n_4}] \overline{\sum_{\substack{n=n'_{1234},n'_1,\cdots,n'_4 \in\Z^2_N  }} \I_4[k^{(4),t}_{n'_1,\cdots,n'_4}] } \Big] \Big\|^2_{\ell^2_n} \\
=&\sum_{n_1,\cdots,n_4 \in\Z^2_N  } \jpb{n_{1234}}^{2s} \int_{[0,\infty)^4} |k^{(4),t}_{n_1,\cdots,n_4}|^2 \, \d s_1\cdots \d s_4 
\end{align}
where we used the Fubini's theorem to interchange the order of $L^2(\Omega)$ and $\ell^2_n$ norms and It\^{o} isometry. 

If we apply the Littlewood-Paley decomposition to the frequencies $n_1,n_2,n_3$ and $n_4$, and denote the dyadic piece by $\<31>^{(4),\mathrm{dyadic}}_N$, we have 
\begin{equ}
\| \|\<31>^{(4)}_N \|_{H^{s-1}_x}\|^2_{L^2(\Omega)} \le \sum_{N_1,\cdots,N_4 \ge 1} \| \|\<31>^{(4),\mathrm{dyadic}}_N \|_{H^{s-1}_x}\|^2_{L^2(\Omega)} 
\end{equ}
where 
\begin{equ}\label{eq:31-4-dyadic}
\| \|\<31>^{(4),\mathrm{dyadic}}_N \|_{H^{s}_x}\|^2_{L^2(\Omega)} = \sum_{n_1,\cdots,n_4 \in \Z^2_N} \jpb{n_{1234}}^{2s} \int_{[0,\infty)^4} |k^{(4),t,\mathrm{dyadic}}_{n_1,\cdots,n_4}|^2 \, \d s_1\cdots \d s_4 
\end{equ}
with 
\begin{equ}
k^{(4),t,\mathrm{dyadic}}_{n_1,\cdots,n_4} = \prod_{j=1}^4 \1_{|n_j| \sim N_j} k^{(4),t}_{n_1,\cdots,n_4}.
\end{equ}

Recall the definition of $k^{(4),t}_{n_1,\cdots,n_4}$, and use \eqref{eq:sin-to-e}, we obtain that
\begin{equs}
|k^{(4),t,\mathrm{dyadic}}_{n_1,\cdots,n_4}| \le \sum_{\sigma,\sigma_1,\cdots,\sigma_4 \in \{ \pm\} }\prod_{j=1}^4 \1_{|n_j| \sim N_j} \jpb{n_j}^{-1+\alpha} \jpb{n_{234}}^{-1} \f 1{\jpb{\phi(\bn)}}
\end{equs}
where $\phi(\bn) = \sigma \jpb{n_{234}} + \sigma_2\jpb{n_2} +\sigma_3\jpb{n_3} + \sigma_4 \jpb{n_4}$.

Put it back to \eqref{eq:31-4-dyadic}, we have
\begin{equs}
\| \|\<31>^{(4),\mathrm{dyadic}}_N \|_{H^{s}_x}\|^2_{L^2(\Omega)} \lesssim T^4 \sum_{n_1,\cdots,n_4 \in \Z^2_N} \jpb{n_{1234}}^{2s} \prod_{j=1}^4 \1_{|n_j| \sim N_j} \jpb{n_j}^{-2+2\alpha} \jpb{n_{234}}^{-2} \f 1{\jpb{\phi(\bn)}^2} \\
\lesssim T^4  \sup_{m \in \Z} \sum_{n_1,\cdots,n_4 \in \Z^2_N} \jpb{n_{1234}}^{2s} \prod_{j=1}^4 \1_{|n_j| \sim N_j} \jpb{n_j}^{-2+2\alpha} \jpb{n_{234}}^{-2} \1_{|\phi-m| \le 1}
\end{equs}

Set $s=-\alpha-\epsilon$, by quartic non-resonance sum estimate Lemma \ref{lem:quartic-non-resonance}, we have
\begin{equs}
\| \|\<31>^{(4),\mathrm{dyadic}}_N \|_{H^{s-1}_x}\|^2_{L^2(\Omega)} \lesssim T^4 N^{-2\epsilon}_4 \maxcurly{N_1,N_2,N_3,N_4}^{-2s_\alpha},
\end{equs}
for $0<\alpha<\f 5{12}$ and $s_\alpha$ defined as \eqref{eq:s-alpha}.

\item[Step 3.] Analysis of single-resonance term $\<31>^{(2)}_N$. By \eqref{eq:Wsinfty-to-Hs}, for fixed $t \in [0,T]$, using we have 
\begin{align}\label{eq:k2-in-Sobolev}
\| \|\<31>^{(2)}_N \|_{H^{s}_x}\|^2_{L^2(\Omega)} =&\sum_{n_3, n_4 \in\Z^2_N  } \jpb{n_{34}}^{2s} \int_{[0,\infty)^2} |k^{(2),t}_{n_3,n_4}|^2 \, \d s_3 \d s_4.
\end{align}

By introducing the Littlewood-Paley decomposition again and denote the dyadic piece by $\<31>^{(2),\mathrm{dyadic}}_N$, we have 
\begin{equ}
\| \|\<31>^{(2),\mathrm{dyadic}}_N \|_{H^{s}_x}\|^2_{L^2(\Omega)} =\sum_{n_3, n_4 \in\Z^2_N  } \jpb{n_{34}}^{2s} \int_{[0,\infty)^2} |k^{(2),t,\mathrm{dyadic}}_{n_3,n_4}|^2 \, \d s_3 \d s_4, \label{eq:k2-dyadic}
\end{equ}
where
\begin{equ}
k^{(2),t,\mathrm{dyadic} }=  \prod_{j=2}^4 \1_{|n_j| \sim N_j} k^{(2),t}_{n_3,n_4}.
\end{equ}

However, due the single resonance, it takes more efforts to bound $|k^{(2),t,\mathrm{dyadic} }|$.  Applying \eqref{eq:sin-to-e}, we have
\begin{equs}
|k^{(2),t,\mathrm{dyadic} }_{n_3,n_4}| \le \sum_{\sigma,\sigma_1,\cdots,\sigma_4 \in \{ \pm\} } \prod_{j=3}^4 \1_{|n_j| \sim N_j} \jpb{n_j}^{-1+\alpha} \sum_{|n_2| \sim N_2} \jpb{n_{234}}^{-1} \jpb{n_2}^{-2+2\alpha}  \f {1}{\jpb{\phi(\bn)}}
\end{equs}
where $\phi(\bn) = \sigma \jpb{n_{234}} + \sigma_2\jpb{n_2} +\sigma_3\jpb{n_3} + \sigma_4 \jpb{n_4}$.

By inserting the extra $m \in \Z$ and applying the basic resonance estimate, we have 
\begin{equs}
|k^{(2),t,\mathrm{dyadic} }_{n_3,n_4}| \le &\sum_{\sigma,\sigma_1,\cdots,\sigma_4 \in \{ \pm\} } \prod_{j=3}^4 \1_{|n_j| \sim N_j} \jpb{n_j}^{-1+\alpha} \sum_{m \in \Z} \f 1{\jpb{m}} \sum_{|n_2| \sim N_2} \jpb{n_{234}}^{-1} \jpb{n_2}^{-2+2\alpha}  \\
\lesssim& \log(2+\maxcurly{N_2,N_3,N_4}) \prod_{j=3}^4 \1_{|n_j| \sim N_j} \jpb{n_j}^{-1+\alpha} \jpb{n_{34}}^{-\f 32+2\alpha}.
\end{equs}
Hence, it yields
\begin{equs}
\| \|\<31>^{(2),\mathrm{dyadic}}_N \|_{H^{s}_x}\|^2_{L^2(\Omega)}  \lesssim & T^2  \log(2+\maxcurly{N_2,N_3,N_4})  \sum_{n_3, n_4 \in\Z^2_N  }  \jpb{n_{34}}^{2s-3 +4\alpha} \prod_{j=3}^4 \1_{|n_j| \sim N_j} \jpb{n_j}^{-2+2\alpha} 
\end{equs}

Using the discrete convolution inequality, we have 
\begin{equs}
\| \|\<31>^{(2),\mathrm{dyadic}}_N \|_{H^{s}_x}\|^2_{L^2(\Omega)}  \lesssim T^2 \log(2+\maxcurly{N_2,N_3,N_4})  \maxcurly{N_3,N_4}^{2s-3+8\alpha}. \label{eq:k2-dyadic-dispersion}
\end{equs}

To remove the logarithmic divergence in $N_2$, we need to utilize the sine-cancellation again. It suffices to consider the case $N_2 \gg \maxcurly{N_3,N_4}$. Let us go back to \eqref{eq:k2-dyadic}. By setting
\begin{equ}
f(t,s,n_2;s_3,s_4,n_3,n_4) := \jpb{n_{234}}^{-1} \jpb{n_3}^{-2} \sin((t-s_3)\jpb{n_3}) \sin((t-s_4)\jpb{n_4})  \prod_{j=3}^4 \1_{[0,s]}(s_j),
\end{equ}
we have
\begin{equs}
{} & \| \|\<31>^{(2),\mathrm{dyadic}}_N \|_{H^{s}_x}\|^2_{L^2(\Omega)} \\
\le & 3 N^{2\alpha}_2 \sum_{n_3, n_4 \in \Z^2_N} \prod_{j=3}^4 \1_{|n_j|\sim N_j} \jpb{n_j}^{-2+2\alpha} \cdot \jpb{n_{34}}^{2s}    \\
& \quad \times \int_{[0,T]^2}  \Big| \sum_{|n_2| \sim N_2} \sin((t-s)\jpb{n_2+n_{34}}) \cos((t-s)\jpb{n_2}) f(t,s,n_2;\cdots) \, \d s \Big|^2   \, \d s_3 \d s_4.
\end{equs}
It is easy to check $f(t,s,n_2;\cdots)$ satisfies the condition in Lemma \ref{lem:sine-cancellation}, hence, we have 
\begin{equs}
 \Big| \sum_{|n_2| \sim N_2} \sin((t-s)\jpb{n_2+n_{34}}) \cos((t-s)\jpb{n_2}) f(t,s,n_2;\cdots) \, \d s \Big| \le T^2 N^2_{34} \log(2+N_2) N^{-1}_2.
 \end{equs}
 Put all together and by discrete convolution sum inequality, we have 
 \begin{equs}
 \| \|\<31>^{(2),\mathrm{dyadic}}_N \|_{H^{s}_x}\|^2_{L^2(\Omega)} \lesssim  T^6 \maxcurly{N_3,N_4}^{2s-4\alpha+4} \log^2(2+N_2) N^{-2}_2. \label{eq:k2-dyadic-sine}
 \end{equs}
 when $s> -1$. 
 
 By interpolating \eqref{eq:k2-dyadic-dispersion} and \eqref{eq:k2-dyadic-sine} and set $s=-\alpha-\epsilon$ with $\epsilon< 1-\alpha$, there exists $\epsilon$ such that 
  \begin{equs}
 \| \|\<31>^{(2),\mathrm{dyadic}}_N \|_{H^{s}_x}\|^2_{L^2(\Omega)} \lesssim  T^\theta \maxcurly{N_2,N_3,N_4}^{-2\epsilon},
 \end{equs}
 which is summable at dyadic scale $N_2,N_3,N_4$. For $s \le -1$, by embedding, we can conclude the same result.
\end{itemize}
\end{proof}

\section{Deterministic and random tensor estimates}\label{sec:deterministic-tensor}
In this subsection, we study the regularity of the stochastic object $\cJ(\<2>_N v)$. We begin by introducing notation for random tensors and recalling some basic results. Note that a tensor defined on an index set $C$ equipped with a partition induces a linear operator, while testing such a tensor with iterated Wiener integrals of a symmetric function (indexed by a subset of $C$) produces a random tensor. These concepts allow us to exploit resonance before working with the $L^p(\Omega)$ norm of the object.

\subsection {Tensors and random tensors}
Let $\X$ be a discrete set,  which we later take to be $\Z^2$. Let  $B, C$ be a partition of a finite set $A$. Let us consider a tensor $h$ on $A$, which is an $L^2$ function from $\X^A$ to $\C$. In canonical basis coordinates, we write $h=(h_{n_A}, n_A \in \X^A)\in \ell^2(\X^A)$. Given a function $n_B$ in $\X^B$ and a function $n_C$ in $\X^C$, we denote the concatenation function by $n_{B\sqcup C}$:
$$n_B \sqcup n_C=\left\{ \begin{aligned}
n_B, \qquad \hbox{ on } B\\
n_C, \qquad \hbox{ on } C.
\end{aligned}\right.$$

It is convenient to view $h$  as a linear map from $\ell^2(\X^B)$ to $\ell^2(\X^C)$, denoted as $h_{B\to C}$. Using the identical identification,  $\ell^2(\Omega_1\times \Omega_2)=\ell^2(\Omega_1) \otimes \ell^2(\Omega_2)$, for any $u\in \ell^2(\X^B)$ and $v\in \ell^2(\X^C)$, 
$$\langle  h_{B\to C}(u), v\rangle= \langle h, u\otimes v\rangle.$$
Let $n_B\in \X^B$, $n_C\in \X^C$ and  $u\in \ell_2(\X^B)$, then the tensor products of the delta functions  sitting at the single elements at $n_B$ and $n_C$ is:
$\delta_{n_B}\otimes \delta_{n_C}(x,y)=\delta_{n_B}(x)\delta_{n_C}(y)=\delta_{n_B \sqcup n_C}(x,y)$ where $x\in \X^B$ and $y\in \X^C$.

The operator norm of the operator is given by
\begin{equ}\label{tensor-norm}
\|h_{B\to C}\|:=\|h_{B\to C}\|_{\mathcal{L} (\ell^2(\X^{B});\ell^2(\X^{C}))}
\end{equ}
Indeed,  in local coordinates we have:
$$ h_{B\to C}(u)(n_C) =   \sum_{n_B \in \X^B} h(n_B \sqcup n_C) u(n_B),$$
where $h$ is viewed as an element of $ \ell^2(\X^B)\otimes \ell^2(\X^C)$, and
$$\|h_{B\to C}\|^2=  \sup_{\|u\|_{\ell_2(\X^B)}=1}\sum_{n_C\in \X^C} | \sum_{n_B \in \X^B} h(n_B \sqcup n_C) u(n_B)|^2.$$

\begin{remark}
This operator norm defined above in \eqref{tensor-norm} agrees with the norm  $\|h\|_{n_B\to n_C}$ used in the stochastic wave equation literature, see e.g. \cite{DNY22, OWZ22, Bri20b}, 
where for $h$ a tensor on $(\Z^2)^A$ and  $(B,C)$ be a partition of $A$, they used the following formalism:
\begin{equ}
(\|h\|_{n_B\to n_C})^2 = \sup_{\|z_{n_B}\|_{\ell^2_{n_B}}=1} \{ \sum_{n_C \in (\Z^2)^C} \Big| \sum_{n_B \in (\Z^2)^B} h_{n_A}z_{n_B} \Big|^2 \}. 
\end{equ}
To emphasize:
\begin{equ}\label{coordinate-tensor-norm}
\|h_{B\to C}\|=\|h\|_{n_B\to n_C}.
\end{equ}
\end{remark}

We also need the concept of the contracted random tensor and its moment estimate.
\begin{definition}\label{c-tensor}
Let $S\subset A$ be finite sets and let $h$ be a tensor on $(\Z^2)^A$.  Assume that
$h_{n_A}$ vanishes whenever $|n_A|\ge M$ where $M$ is a given number. Let us consider a family of symmetric functions
$$f_{n_S}: L^2(\R^S, \C), \qquad n_S \in (\Z^2)^S$$
such that their $L^2$ norm is uniformly bounded: $\| \|f_{n_S}(t_S)\|_{L^2_{t_S}} \|_{\ell^\infty_{n_S}}<\infty$.
 The tensor $h$ tested on $f$ defines  a random tensor on $(\Z^2)^{A\setminus S}$ as follows: for any $n_{A\setminus S}\in (\Z^2)^{A\setminus S}$,
 $$h^{\mathrm {con}}_{n_{A\setminus S}}= \sum_{n_{S} \in (\Z^2)^S} \I_k [h_{n_A}f_{n_S}(t_S)],$$
where $n_A=n_A\sqcup n_{A\setminus S}$. We refer $h^{\mathrm {con}}$ as the contracted random tensor by $f$ with respect to the partition $\{S, A\setminus S\}$ of $A$ (in short w.r.t. $S\subset A$).
\end{definition}

We denote by $\part(A)$ the set of non-trivial partitions of $A$.

\begin{proposition}\label{contracted-tensor-est}
Let $A,S, f, h$, $M$  as in Definition \ref{c-tensor}. Denote by $h^{\mathrm{con}}$ the contraction of $h$ by $f$ with respect to $S\subset A$.
Let $B,C$ be a partition of $A\setminus S$, and consider the collection $\P$ of partitions 
of $A$ that extends the partition $(B,C)$ of $A\setminus S$ in the sense that 
$$\part_{B,C}=\{ (X,Y)\in \part(A): B \subset X, C\subset Y\}.$$
Then, for all $p \ge 2$ and any $\theta > 0$, we have
\begin{equation}\label{eq:c-tensor-est}
\| \| h^{\mathrm{con}} \|_{n_B\to n_C} \|_{L^p(\Omega)} 
\lesssim p^{k/2} N^{\theta}
 \max_{(X,Y)\in \part_{B,C}} \|h\|_{n_X \to n_Y} \| \|f_{n_S}(t_S)\|_{L^2_{t_S}} \|_{\ell^\infty_{n_S}}.
\end{equation}
\end{proposition}
The proof can be found in  \cite{DNY22}  or \cite{OWZ22}. 

In our application, we set $n_S = (n_1,n_2)$, $n_{A\backslash S} = (n,n_3)$, and $n_B = n_3, n_C= n$.

\subsection{Proof of boundedness of the quadratic random operator}
Given a function $z$, we define 
\begin{equ}\label{frequency-sign}
\F_t(z^+)(\lambda) = \1_{\lambda \ge 0} \F_t (z)(\lambda), \quad \F_t(z^-)(\lambda) = \1_{\lambda < 0} \F_t(z)(\lambda).
\end{equ}
Then $z=z^++z^-$. Note the convenient decomposition for the norm:
\begin{equ}\label{sum-xsbnorm}
 \| z\|^2_{\X^{s',b'}}=\| z^+\|^2_{\X^{s',b'}}+ \| z^-\|^2_{\X^{s',b'}}.
\end{equ}
In addition since the frequencies of  $z^{\pm}$ lie entirely on $\R+$ or $\R_-$, we have,  for any $s',b'  \in \R$, 
\begin{equ}
\| z^\pm\|^2_{\X^{s',b'}}
=\int_\R\|  \jpb{m}^{s'} \jpb{\pm \mu-\jpb{m}}^{b'} \F_{t,x}(z^{\pm})(\mu,m)\|^2_{\ell^2_{m}} \, \d \mu.
\end{equ}
Therefore, by the invariance of the Lebesgue measure, we have
\begin{equ}
\| z\|^2_{\X^{s',b'}}
=\sum_{\tau\in \{\pm 1\}} 
\Bigl \|\jpb{m}^{s'}  \int_\R \jpb{\lambda}^{b'}\F_{t,x}(z^{\tau})(\lambda+\tau  \jpb{m},m)\, \d \lambda\Bigr\|^2_{\ell^2_m}. \label{z-norm}
\end{equ}
 The following well known lemma holds:
\begin{lemma}\label{le:xsb-product}
For any $s,b\in \R$,
\begin{equs}
\|fg\|_{\X^{s,b}}
\lesssim \sum_{\sigma, \tau \in \{\pm 1\}}
\|\jpb{n}^s\jpb{\lambda}^b\sum_m \int_\R \F_{t,x} (f)(\lambda -\sigma(n)-\mu, n-m) \F_{t,x} (g^\tau) (\mu, m) d\mu \|_{\ell^2_n L^2_\lambda}
\end{equs}
\end{lemma}
{\begin{proof}
This follows from \eqref{eq:Xsb-upper},
\begin{equs}
\| fg \|_{\X^{s,b}} &\le \sum_{\sigma \in \{ \pm 1\}} \| \jpb{n}^{s} \jpb{\lambda}^b \cF_{t,x}(fg)(\lambda+\sigma \jpb{n},n)\|_{\ell^2_nL^2_\lambda}\\
& \le    \sum_{\sigma \in \{ \pm 1\}} \| \jpb{n}^{s} \jpb{\lambda}^b \sum_{n_3 \in \Z^2} \int_{\R} \cF_{t,x}(f)(\lambda+\sigma \jpb{n}-\lambda_3,n-n_3)\cF_{t,x}(g)(\lambda_3,n_3) \, \d \lambda_3 \|_{\ell^2_nL^2_\lambda}
\end{equs}
Introducing the decomposition $g = g^+ + g^-$ to conclude.\end{proof}

\begin{proposition}\label{tensors}
Let $T\le 1$, $p \ge 2$, and $b>\f 12$. 
\begin{itemize}
\item For $0<\alpha<\f 14$, we work with the space $\X^{\f 14+\epsilon,b}$  with the restriction~$0< \epsilon<\f 12$,
\item
For $\f 14 \le \alpha < \f 38$,  we work with $\X^{2\alpha-\f 14+\epsilon, b}$ where $0<\epsilon < \f 12(3-8\alpha)$.
\end{itemize} 
For such $s,b$, we have
\begin{equ}\label{eq:1st-tensor}
\sup_{N \ge 1}\| \sup_{I \subset [0,T]} \sup_{\|z\|_{\X^{s,b}_I} = 1}
 \| \<2>_N z \|_{\X^{s-1,b-1}_I} \|_{L^p(\Omega} \lesssim 1.
\end{equ}
Consequently,  for almost surely all $\omega$, the operator $\{ \cJ^{\<2>_N} \}$, where
$$\cJ^{\<2>^N}:  z \in \X^{s,b}_I \mapsto  \cJ(\<2>_N z) \in \X^{s,b}_I,$$
is a Cauchy sequence in  $\mathcal{L}(\X^{s,b}_I;\X^{s,b}_I)$, the space of bounded linear operator.
We denote its limit by $\cJ^{\<2>}$. Furthermore, its operaator norm belongs to $L^p$.
\end{proposition}

\begin{proof}
Let $s,b$ satisfy the conditions of the statement and let $I\subset [0,T]$. 
 We divide the proof  into two steps.

\textbf{Step 1}. 
We first reduce the estimate in $L^p(\Omega;X^{s-1,b-1}_I)$ to an operator-norm bound.
By the definition of the restricted norm, for any $s',b' \in \R$, 
\begin{align}
 \| \sup_{I \subset [0,T]} \sup_{\|z\|_{\cX^{s',b'}_I} = 1} \| \<2>_N z \|_{\cX^{s-1,b-1}_I} \|_{L^p(\Omega)}
\lesssim & \| \sup_{\|z\|_{\cX^{s',b'}_T} = 1} \| \<2>_N z \|_{X^{s-1,b-1}_T} \|_{L^p(\Omega)} \\
\lesssim & \| \sup_{\|z\|_{\cX^{s',b'}} =1} \| \1_{[0,T]} \<2>_N z \|_{\cX^{s-1,b-1}} \|_{L^p(\Omega)}.\label{eq:operator-norm}
\end{align}

To encode the interactions between $\<2>_N$ and $z$, we decompose $z$,  as in~\eqref{frequency-sign}):
\begin{equation}
    z= z^+ + z^-.
\end{equation}
To estimate $\| \1_{[0,T]} \<2>_N z \|_{\cX^{s-1,b-1}}$, we work with its Fourier transform. By discrete convolution, 
\begin{equs}
{}&\cF_{t,x}(\1_{[0,T]} \<2>_N z)(\lambda-\sigma \jpb{n}, n)\\
=& \sum_{\sigma_3 \in \{\pm\}} \sum_{n_3 \in \bZ^2} \int_{\bR} \cF_{t,x}(\1_{[0,T]}\<2>_N)(\Lambda,n-n_3) \cF_{t,x}(z^{\sigma_3})(\lambda_3+\sigma_3\jpb{n_3},n_3)\d \lambda_3.
\end{equs}
where we have introduced the shorthand:
$$\Lambda=\lambda-\lambda_3-\sigma\jpb{n}-\sigma_3\jpb{n_3}.$$ 

We next compute the Fourier transform of $\1_{[0,T]}\<2>_N$. Recall the tensor defined in~\eqref{define-h}:
$$h^{(2),N}_{n_1,n_2}(t)= \prod_{j=1}^2 \1_{[0,t]}(s_j) \1_{|n_j| \le N}\f {\sin((t-s_j) \jpb{n_j})}{\jpb{n_j}^{1-\alpha}}
$$
With this, we may write $\<2>_N(t) = \sum_{n_1,n_2 \in \Z^2} e^{\i (n_1+n_2)\cdot x}\cI_2[h^{(2),N}_{n_1,n_2}](t)$. 
We also introduce the phase function $\phi(\boldsymbol{n}) )$
  \begin{equ}
\phi(\boldsymbol{n}) =\sigma\jpb{n}+\sigma_1\jpb{n_1}+\sigma_2\jpb{n_2}+\sigma_3\jpb{n_3}.
\end{equ}

With these notations in place, setting $n_1+n_2=n-n_3$, we obtain the frequency component of $\cF_{t,x}(\1_{[0,T]} \<2>_N z)$ at $n-n_3$:
\begin{align}
 \cF_{t,x}(\1_{[0,T]}\<2>_N)(\Lambda,n-n_3) 
&=  \int_0^T e^{-\i t \Lambda} 
\cF_x(\<2>_N(t)) (n-n_3) \, \d t\\
&=  \int_0^T e^{-\i t \Lambda}  \sum_{n_1,n_2 \in \Z^2} \1_{n-n_3=n_1+n_2}\cI_2[h^{(2),N}_{n_1,n_2}](t)\, \d t\\
&= \cI_2[\hat h_{n,n_3}(\lambda,\lambda_3)],
\end{align}
where we have set 
\begin{align}
{}&\hat h_{n,n_3}(\lambda,\lambda_3)= \sum_{n_1,n_2\in\Z^2} \hat h_{n,n_1,n_2, n_3}(\lambda,\lambda_3) \\
 {}&  \hat h_{n,n_1,n_2n_3}(\lambda,\lambda_3),\\
   &=  \1_{n=n_{123}} \prod_{j=1}^2 \jpb{n_j}^{-1+\alpha} \1_{|n_j| \le N} \int_0^T e^{-\i t (\lambda-\lambda_3-\sigma\jpb{n}-\sigma_3\jpb{n_3})} \Bigl( \prod_{j=1}^2 \sin((t-s_j)\jpb{n_j})) \1_{[0,t]}(s_j) \Bigr) \, \d t\\
    &= -\frac{1}{4}  \1_{n=n_{123}} \sum_{\sigma_1,\sigma_2 \in \{ \pm 1\}} \prod_{j=1}^2 \jpb{n_j}^{-1+\alpha} \1_{|n_j| \le N}  \int_0^T e^{-\i t (\lambda-\lambda_3-\phi(\boldsymbol{n}) )} e^{-\i \sum_{j=1}^2 \sigma_j s_j \jpb{n_j}} \, \d t.
\end{align}
We used the notation $n_{123}=n_1+n_2+n_3$. Inserting this back, we obtain
\begin{equ}
\cF_{t,x}(\1_{[0,T]} \<2>_N z)(\lambda - \sigma\jpb{n},n)
= \sum_{\sigma_3 \in \{\pm \}} \int_{\bR} \sum_{n_3 \in \bZ^2}  \cI_2[\hat h_{n,n_3}(\lambda,\lambda_3)] \cF_{t,x}(z^{\sigma_3})(\lambda_3+\sigma_3\jpb{n_3},n_3)\d \lambda_3.
\end{equ}
We recall the inequality~\eqref{eq:Xsb-upper}:
$$
\|u\|_{\X^{s-1,b-1}} \lesssim \sum_{\sigma \in \{ \pm \}} \|\jpb{n}^{s-1} \jpb{\lambda}^{b-1} \cF_{t,x}(u)(\lambda - \sigma \jpb{n}, n)\|_{\ell^2_nL^2_\lambda}$$
and hence
\begin{align}
{}& \| \1_{[0,T]}\<2>_N z \|_{\cX^{s-1,b-1}} \\
&\lesssim \sum_{\sigma \in \{\pm\}} \Bigl \| \jpb{n}^{s-1}\jpb{\lambda}^{b-1}\!\!\!\sum_{\sigma_3 \in \{\pm\}} \sum_{n_3 \in \bZ^2} \int_{\bR}  \cI_2[\hat h_{n,n_3}(\lambda,\lambda_3)] \cF_{t,x}(z^{\sigma_3})(\lambda_3+\sigma_3\jpb{n_3},n_3)\d \lambda_3 \Bigr\|_{\ell^2_nL^2_\lambda} .
\end{align}

By interchanging the $\ell^2_n$-norm and  $L^1_{\lambda_3}$-norm,  and applying the Cauchy-Schwartz inequality, we further obtain the bound
\begin{align}
{}& \| \1_{[0,T]}\<2>_N z \|_{\cX^{s-1,b-1}} \\
&\lesssim\sum_{\sigma,\sigma_3 \in \{ \pm \}} \| \jpb{n}^{s-1}\jpb{\lambda}^{b-1}
 \int_{\bR}  \Big\| \sum_{n_3 \in \bZ^2} \cI_2[\hat h_{n,n_3}(\lambda,\lambda_3)] \cF_{t,x}(z^{\sigma_3})(\lambda_3+\sigma_3\jpb{n_3},n_3)  \Big\|_{\ell^2_n}  \, \d \lambda_3 \|_{L^2_\lambda}.
\end{align}

We are not able to conclude that the Sobolev norm of the above object admits a finite $L^p(\Omega)$-bound. Instead, fixing the $z$-variable, we interpret the expression as the contraction of two tensors.
The resulting object is a random tensor which acts as an operator on $\ell^2(\Z^2))$. 

More precisely, fixing $\lambda, \lambda_3, \sigma_3$, and contracting in the $n_3$ variable, we obtain the random operator:
$$z\mapsto \sum_{n_3 \in \bZ^2}  \cI_2[\hat h_{n,n_3}(\lambda,\lambda_3)] \cF_{t,x}(z^{\sigma_3})(\lambda_3+\sigma_3\jpb{n_3},n_3).$$
Fixing $z$, we interpret the tensor product on the right-hand side as a linear operator on $\ell^2(\Z^2)$. To this end, set 
$$v^{\sigma_3}(z)(\lambda_3+\sigma_3\jpb{n_3},n_3)=\jpb{\lambda_3}^{b'} \jpb{n_3}^{s'}\cF_{t,x}(z^{\sigma_3})(\lambda_3+\sigma_3\jpb{n_3},n_3).$$ 
We then consider the linear map
$$v_3^\sigma(z) \mapsto K_{n_3,n}(\lambda,\lambda_3)(v^\sigma_3(\lambda_3+\sigma_3\jpb{n_3}, n_3)),$$ which we denote by $K(\lambda,\lambda_3)$. Explicitly,
\begin{equ} \label{eq:kernel-K}
K(\lambda,\lambda_3)(v^\sigma_3(z))(n)
 =\sum_{n_3 \in \bZ^2}  \jpb{n}^{s-1}\jpb{n_3}^{-s'}\cI_2[ \hat h_{n,n_3}(\lambda,\lambda_3)]v^{\sigma_3}(\lambda_3+\sigma_3\jpb{n_3},n_3).
\end{equ}

Substituting the representation \eqref{eq:kernel-K} back into \eqref{eq:operator-norm} and applying Lemma~\ref{le:xsb-product}, we obtain
\begin{align}
{}&\| \sup_{\|z\|_{\cX^{s',b'}} = 1} \| \1_{[0,T]}\<2>_N z \|_{\cX^{s-1,b-1}} \|_{L^p(\Omega)}\\
&\lesssim \Bigl\|  \sup_{\|z\|_{\cX^{s',b'}} = 1}
\sum_{\sigma,\sigma_3 \in \{ \pm \}}    \| \jpb{\lambda}^{b-1}
 \int_{\bR}  \Big\| \jpb{n}^{s-1}\sum_{n_3 \in \bZ^2} \cI_2[\hat h_{n,n_3}(\lambda,\lambda_3)] \cF_{t,x}(z^{\sigma_3})(\lambda_3+\sigma_3\jpb{n_3},n_3)  \Big\|_{\ell^2_n}  \, \d \lambda_3 \|_{L^2_\lambda}\Bigr\|_{L^p(\Omega)}\\
& \lesssim \Bigl\|  \sup_{\|z\|_{\cX^{s',b'}} = 1}
\sum_{\sigma,\sigma_3 \in \{ \pm \}}    \| \jpb{\lambda}^{b-1}
 \int_{\bR}  \jpb{\lambda_3}^{-b'}\Big\| K(\lambda,\lambda_3)(v^\sigma_3(z))(n)  \Big\|_{\ell^2_n}  \, \d \lambda_3 \|_{L^2_\lambda}\Bigr\|_{L^p(\Omega)}\\
& \lesssim \sum_{\sigma \in \{ \pm \}} 
\Bigl \|  \jpb{\lambda}^{b-1}\jpb{\lambda_3}^{-b'} 
\| K(\lambda,\lambda_3)\| \; 
\Big\|_{L^2(\lambda_3), L^2_\lambda,L^p(\Omega)}
  \sup_{\|z\|_{\cX^{s',b'}} = 1}\sum_{\sigma_3 \in \{ \pm \}} 
  \|v^{\sigma_3}(z) (\lambda_3+\sigma_3\jpb{n_3},n_3)\big\|_{L_{\lambda_3}^2\ell^2_{n_3}}.
  \end{align}
 In the final step, we applied the Cauchy-Schwartz inequality. We now take advantage of the decomposition into $z^\pm$, for which
\begin{equ}
\| z^{\sigma_3} \|^2_{\X^{s',b'}}
=\Bigl \| \int_\R \jpb{n_3}^{2s'}\jpb{\lambda_3}^{2b'} |\F_{t,x}(z^{\sigma_3})(\lambda_3+\sigma_3 \jpb{n_3},n_3)|^2\, \d \lambda_3\Bigr\|^2_{\ell^1_{n_3}}\end{equ}
and we observe that
\begin{equation}\label{eq:norm-transformation}
\|v^{\sigma_3}(z)(\lambda_3+\sigma_3 \jpb{n_3},n_3)\|_{L^2_{\lambda_3}\ell^2_{n_3}} = \| z^{\sigma_3} \|_{\cX^{s',b'}}.\end{equation}
It remains to show that
$$\Bigl \|  \jpb{\lambda}^{b-1}\jpb{\lambda_3}^{-b'} 
\| K(\lambda,\lambda_3)\| \; 
\Big\|_{L^2(\lambda_3), L^2_\lambda,L^p(\Omega)}<\infty
$$
where $\| K(\lambda,\lambda_3)\|$ denotes the operator norm on $\ell^2(\Z^2)$ which coincides with the tensor norm introduced used in \cite{DNY22,OWZ22}. 

Note $b'>\frac{1}{2}$ and assume $b<\frac{1}{2}$. (This is consistent with our approach in most $\cX^{s,b}$ norm estimates where we first treat $b<\frac{1}{2}$ and then upgrade to $b>\frac{1}{2}$ by interpolation; 
see Remark \ref{rmk:b-vs-b_}.) Under these assumptions, both $\jpb{\lambda}^{b-1}$ and $\jpb{\lambda_3}^{-b'}$ are integrable in $L^2$. Therefore, it remains only to show that
$$\sup_{\lambda, \lambda_3} \|  \| K(\lambda,\lambda_3)\| \|_{L^p(\Omega)}<\infty.$$

\textbf{Step 2}. We estimate the operator norm of $K$ where
\begin{equ}
K(\lambda,\lambda_3)(v^\sigma_3(z))(n)
 =\sum_{n_3 \in \bZ^2}  \jpb{n}^{s-1}\jpb{n_3}^{-s'}\cI_2[ \hat h_{n,n_3}(\lambda,\lambda_3)]v^{\sigma_3}(\lambda_3+\sigma_3\jpb{n_3},n_3).
\end{equ}
For $m \in \Z$, set
\begin{align}
\tilde h^m_{n,n_1,n_2,n_3} = \1_{n=n_{123}} \1_{|\phi-m| \le 1} \jpb{n}^{s-1} \jpb{n_3}^{-s'} \prod_{j=1}^2 \jpb{n_j}^{-1+\alpha} \1_{|n_j| \le N} \\
\tilde k^m_{n,n_1,n_2,n_3}(\lambda,\lambda_3) =  -\f1 4\int_0^T \1_{|\phi-m| \le 1}e^{-\i t (\lambda-\lambda_3-\phi(\boldsymbol{n}) )} e^{-\i \sum_{j=1}^2 \sigma_j s_j \jpb{n_j}} \, \d t.
\end{align}

{\it In the remaining of the proof, we simplify notation by writing $h^m$ for $\tilde h^m$ and $k^m$ for $\tilde k^m$.}
Then the contracted operator takes the following form, 
\begin{equ}
K_{n,n_3}(\lambda,\lambda_3) = \sum_{n_1,n_2 \in \Z^2} \sum_{m \in \Z} \I_2[h^m_{n,n_1,n_2,n_3} k^m_{n,n_1,n_2,n_3}(\lambda,\lambda_3) ].
\end{equ}

It suffices to bound the operator on a dyadic scale, where we insert the dyadic decomposition with respect to frequencies $n_1,n_2$ and $n_3$ and denote the corresponding dyadic pieces by 
\begin{align}
\tilde h^{m,\mathrm{dyadic}}_{n,n_1,n_2,n_3} = \1_{n=n_{123}} \1_{|\phi-m| \le 1} \jpb{n}^{s-1} \prod_{j=1}^3 \1_{|n_j| \sim N_j} \jpb{n_3}^{-s'} \prod_{j=1}^2 \jpb{n_j}^{-1+\alpha} \1_{|n_j| \le N} \\
\tilde k^{m,\mathrm{dyadic}}_{n_1,n_2}(\lambda,\lambda_3) =  -\f1 4\prod_{j=1}^2 \1_{|n_j| \sim N_j} \int_0^T \1_{|\phi-m| \le 1}e^{-\i t (\lambda-\lambda_3-\phi(\boldsymbol{n}) )} e^{-\i \sum_{j=1}^2 \sigma_j s_j \jpb{n_j}} \, \d t.
\end{align}

Applying Proposition~\ref{contracted-tensor-est}, we obtain, for any $\theta>0$,
\begin{align}
\| \| K_{n,n_3}(\lambda,\lambda_3)\|_{} \|_{L^p(\Omega)} \lesssim & p \maxcurly{N_1,N_2,N_3}^\theta \\
&\quad \sum_{m \in \Z} \max_{(X,Y) \in \P} \Big(\|h^{m,\text{dyadic}}_{n,n_1,n_2,n_3} \|_{n_X \to n_Y}\Big) \|k^{m,\mathrm{dyadic}}_{n_1,n_2}(\lambda,\lambda_3) \|_{\ell^\infty_{n_1,n_2}L^2_{s_1,s_2}}
\end{align}
where $n_X$ contains $n_3$ and $n_Y$ contains $n$.
Note that 
\begin{equ}
\|k^{m,\mathrm{dyadic}}_{n_1,n_2}(\lambda,\lambda_3) \|_{\ell^\infty_{n_1,n_2}L^2_{s_1,s_2}} \lesssim  \frac{1}{\jpb{\lambda-\lambda_3-m}},
\end{equ}
and 
\begin{equs}
\max_{(X,Y) \in \P} \|h^{m,\mathrm{dyadic}}_{n,n_1,n_2,n_3} \|_{n_X \to n_Y} = &\max \Bigl( \|h^{m,\mathrm{dyadic}}\|_{n_1n_2n_3 \to n}, \|h^{m,\mathrm{dyadic}}\|_{n_1n_3\to nn_2},  \\
&\quad  \|h^{m,\mathrm{dyadic}}\|_{n_2n_3\to nn_1},  \|h^{m,\mathrm{dyadic}}\|_{n_3\to nn_1n_2}\Bigr)
\end{equs}

Thus, we obtain, for any $\theta,\kappa >0$
\begin{equs}
{} &\sup_{\lambda,\lambda_3 \in \R} \| \| K^{\mathrm{dyadic}}_{n,n_3}(\lambda,\lambda_3)\|_{n_3 \to n} \|_{L^p(\Omega)} \\
\lesssim & \maxcurly{N_1,N_2,N_3}^\theta \sum_{m \in \Z} \max \Bigl( \|h^{m,\mathrm{dyadic}}\|_{n_1n_2n_3 \to n}, \|h^{m,\mathrm{dyadic}}\|_{n_1n_3\to nn_2},  \\
&\quad  \|h^{m,\mathrm{dyadic}}\|_{n_2n_3\to nn_1},  \|h^{m,\mathrm{dyadic}}\|_{n_3\to nn_1n_2}\Bigr)  \frac{1}{\jpb{\lambda-\lambda_3-m}} \\
\le &\maxcurly{N_1,N_2,N_3}^\theta \sup_{m \in \Z} \Bigl( \|h^{m,\mathrm{dyadic}}\|_{n_1n_2n_3 \to n}, \|h^{m,\mathrm{dyadic}}\|_{n_1n_3\to nn_2}, \|h^{m,\mathrm{dyadic}}\|_{n_2n_3\to nn_1},  \\
&\quad   \|h^{m,\mathrm{dyadic}}\|_{n_3\to nn_1n_2}\Bigr) \Big(\sup_{\lambda,\lambda_3 \in \R} \sum_{m \in \Z}  \frac{1}{\jpb{\lambda-\lambda_3-m}} \Big). \\
\lesssim & \maxcurly{N_1,N_2,N_3}^{(\theta+\kappa)} \sup_{m \in \Z} \Bigl( \|h^{m,\mathrm{dyadic}}\|_{n_1n_2n_3 \to n}, \|h^{m,\mathrm{dyadic}}\|_{n_1n_3\to nn_2}, \\
&\quad \|h^{m,\mathrm{dyadic}}\|_{n_2n_3\to nn_1},  \|h^{m,\mathrm{dyadic}}\|_{n_3\to nn_1n_2}\Bigr) 
\end{equs}
 where in the last step we use the following fact 
\begin{equ}
\sup_{\lambda,\lambda_3 \in \R} \sum_{m \in \Z}  \frac{1}{\jpb{\lambda-\lambda_3-m}} \lesssim \log(2+\maxcurly{N_1,N_2,N_3}),
\end{equ}
due to $|m| \lesssim \maxcurly{N_1,N_2,N_3}$.

If $s=s' = \f 14 +\epsilon$ or $ s= s'=2\alpha -\f 14 +\epsilon$, by the tensor estimate from Lemma \ref{tensor}, there exists some appropriate $\epsilon=\epsilon(s,\alpha)>0$, we have
\begin{align}
\sup_{\lambda,\lambda_3 \in \R} \| \| K^{\mathrm{dyadic}}_{n,n_3}(\lambda,\lambda_3)\|_{n_3 \to n} \|_{L^p(\Omega)} \lesssim \maxcurly{N_1,N_2,N_3}^{-\epsilon+\theta+\kappa}.
\end{align}
Taking $\epsilon>\theta+\kappa$ and summing over all the dyadic scales, we get \eqref{eq:1st-tensor} and complete the proof.
\end{proof}

\begin{remark}
By Lemma \ref{lem:nlXsb} and \eqref{eq:1st-tensor}, it yields that 
\begin{equ}
\|  \| \cJ(\<2>_N z) \|_{\X^{s,b}_I} \|_{L^p(\Omega)} \lesssim \|z\|_{\X^{s',b'}_I}.
\end{equ}
Setting  $s'=s$, $b'=b$, we see that, for any $I \subset [0,T]$ almost surely, $\{\cJ^{\<2>^N}\}$ is Cauchy  in $\mathcal{L}(\X^{s,b}_I;\X^{s,b}_I)$.
Furthermore,  the operator norm of $\cJ^{\<2>_N}$ satisfies
\begin{equ}
\Bigl\| \|\cJ^{\<2>_N}\| \Bigr\|_{L^p(\Omega)} =\Bigl\|  \sup_{\|z\|_{\X^{s,b}_I} =1}  \| \cJ(\<2>_N z) \|_{\X^{s,b}_I}\Bigr\|_{L^p(\Omega)}
\end{equ}
is bounded as claimed.
\end{remark}

Now we show a tensor estimate, which is used in the proof of Proposition \ref{tensors}.

\begin{lemma}[Tensor estimate]\label{tensor}
Let $N_1,N_2,N_3,N_{123} \in 2^{\N}$, $m \in \Z$, and let 
\begin{align}
h^m_{nn_1n_2n_3} &= \1_{n=n_{123}}  \cdot \1_{|n_{123}| \sim N_{123}} \prod_{j=1}^3 \1_{|n_j| \sim N_j} \cdot  \1_{|\phi-m| \le 1} \f {\jpb{n}^{s-1}}{ \jpb{n_1}^{1-\alpha} \jpb{n_2}^{1-\alpha} \jpb{n_3}^{s'}},
\end{align}
where the phase is given by
\begin{equ}
\phi(n_1,n_2,n_3) = \sigma\jpb{n} + \sigma_1\jpb{n_1}+\sigma_2\jpb{n_2}+\sigma_3\jpb{n_3}
\end{equ}
with $\sigma, \sigma_1, \sigma_2, \sigma_3 \in \{ \pm 1\}$ given.

Then the following statements hold:
\begin{itemize}
\item For $0 < \alpha < \f 14$ and $\f 14 < s = s' < \f 34$, there exists $\epsilon = \epsilon(s,\alpha)>0$ such that 
\begin{align}
{}& \sup_{m \in \Z} \max \Bigl( \|h^{m,\mathrm{dyadic}}_{nn_1n_2n_3}\|_{n_1n_2n_3 \to n},\|h^{m,\mathrm{dyadic}}_{nn_1n_2n_3}\|_{n_1n_3\to nn_2},  \|h^{m,\mathrm{dyadic}}_{nn_1n_2n_3}\|_{n_2n_3\to nn_1},  \|h^{m,\mathrm{dyadic}}_{nn_1n_2n_3}\|_{n_3\to nn_1n_2}\Bigr)\\
\lesssim &\maxcurly{N_1,N_2,N_3}^{-\epsilon}.
\end{align}
\item For $\f 14 \le  \alpha < \f 3 8$ and $2\alpha - \f 14 < s = s' < \f 34$, there exists $\epsilon = \epsilon(s,\alpha)>0$ such that 
\begin{align}
{}&\sup_{m \in \Z} \max \Bigl( \|h^{m,\mathrm{dyadic}}_{nn_1n_2n_3}\|_{n_1n_2n_3 \to n},\|h^{m,\mathrm{dyadic}}_{nn_1n_2n_3}\|_{n_1n_3\to nn_2},  \|h^{m,\mathrm{dyadic}}_{nn_1n_2n_3}\|_{n_2n_3\to nn_1},  \|h^{m,\mathrm{dyadic}}_{nn_1n_2n_3}\|_{n_3\to nn_1n_2}\Bigr)\\
&\lesssim \maxcurly{N_1,N_2,N_3}^{-\epsilon}.
\end{align}
\end{itemize}
\end{lemma}

\begin{proof}
We will estimate term by term.
\begin{itemize}
\item $\|h^{m,\mathrm{dyadic}}_{nn_1n_2n_3}\|_{n_1n_2n_3 \to n}$. Noticing the symmetry in $n_1,n_2$, we can assume $N_1 \geq N_2$ without loss of generality. 
Recall the Schur test is given by 
\begin{equ}
\|T\|^2_{L^2 \to L^2} \le \Big(\sup_{x \in X} \int_Y |K(x,y)|\, \d y \Big) \cdot \Big( \sup_{y \in Y} \int_x |K(x,y)|\, \d x\Big)
\end{equ}
where $K$ is the non-negative kernel of operator $T$.

Adapting this test to $\ell^2 \to \ell^2$ setting, we have
\begin{align}
\|h^{m,\mathrm{dyadic}}_{nn_1n_2n_3}\|^2_{n_1n_2n_3 \to n}  \lesssim & N^{2s-2}_{123}N^{-2+2\alpha}_1N^{-2+2\alpha}_2N^{-2s'}_3 \\
& \quad \times \sup_{n \in \Z^2} \sum_{n_1,n_2,n_3 \in \Z^2}  \prod_{j=1}^3 \mathbf{1}_{|n_j| \sim N_j} \mathbf{1}_{|n| \sim N_{123}} \mathbf{1}_{n=n_{123}}\mathbf{1}_{|\phi -m| \leq 1}\\
& \quad \times \sup_{n_1,n_2,n_3 \in \Z^2} \sum_{n \in \Z^2}  \prod_{j=1}^3 \mathbf{1}_{|n_j| \sim N_j} \mathbf{1}_{|n| \sim N_{123}} \mathbf{1}_{n=n_{123}}\mathbf{1}_{|\phi -m| \leq 1}.
\end{align}
Once $n_1,n_2,n_3$ is fixed, $n$ is also fixed. Hence the last term can be simply bounded by $1$. For the term in second line,  we apply the sup-counting estimate \eqref{eq:cubic-sup-count1} (i) to get
\begin{align}
\|h^{m,\mathrm{dyadic}}_{nn_1n_2n_3}\|^2_{n_1n_2n_3 \to n}   \lesssim & N^{2s-2}_{123}N^{-2+2\alpha}_1N^{-2+2\alpha}_2N^{-2s'}_3 \\
\times &\medcurly{N_1,N_2,N_3}^2 \mincurly{N_1,N_2,N_3}^{\f 32}.
\end{align}

We discuss case by case.
\begin{itemize}
\item[(i)] $N_1 \ge N_2 \ge N_3$. We have 
\begin{align}
\|h^{m,\mathrm{dyadic}}_{nn_1n_2n_3}\|^2_{n_1n_2n_3 \to n}   \lesssim & N^{2s-2}_{123} N^{-2+2\alpha}_1 N^{2\alpha}_2 N^{-2s'+\f 32}_3,
\end{align}
If $s'<\f 34$, then one has 
\begin{equ}
\|h^{m,\mathrm{dyadic}}_{nn_1n_2n_3}\|^2_{n_1n_3\to nn_2} \lesssim N^{2s-2}_{123}N^{-2+4\alpha}_1 \le N^{-2+4\alpha}_1.
\end{equ}
Hence, for $0<\alpha\le \f 14$, by setting both $s,s'$ to be $\f 14 +\epsilon$  where $0<\epsilon < \f 12$, we have the desired estimate since $-2+2\alpha+4\alpha$ is negative. If $\f 14 \le \alpha < \f 12$, by setting both $s,s'$ to be $2\alpha - \f 14 +\epsilon$ where $\epsilon < 1-2\alpha$, we also have the estimate.

\item[(ii)] $N_1 \ge N_3 \ge N_2$. Then, we get
\begin{align}
\|h^{m,\mathrm{dyadic}}_{nn_1n_2n_3}\|^2_{n_1n_2n_3 \to n}  \lesssim &N^{2s-2}_{123} N^{-2+2\alpha}_1N^{-2+2\alpha}_2N^{-2s'}_3 \cdot N^2_3N^{\f 32}_2,
\end{align}
As case (i), if $0<\alpha\le \f 14$ and $s=s'=\f 14 +\epsilon$ with $0<\epsilon<\f 34$, we have 
\begin{equ}
\|h^{m,\mathrm{dyadic}}_{nn_1n_2n_3}\|^2_{n_1n_3\to nn_2} \lesssim N^{-\f 32+2\epsilon}_{123} N^{-2+2\alpha}_1 N^{-\f 1 2+2\alpha}_2 \le N^{-2+2\alpha}_1;
\end{equ}
If $\f 14 \le \alpha < \f 12$ and $s=s'=2\alpha-\f 14 +\epsilon$ with $0<\epsilon<\f 54-2\alpha$, we have
\begin{equ}
\|h^{m,\mathrm{dyadic}}_{nn_1n_2n_3}\|^2_{n_1n_3\to nn_2} \lesssim N^{-\f 52+4\alpha+2\epsilon}_{123} N^{-\f 52+4\alpha}_1 \le N^{-\f 52+4\alpha}_1.
\end{equ}

\item[(iii)] $N_3 \ge N_1 \ge N_2$. We have
\begin{align}
\|h^{m,\mathrm{dyadic}}_{nn_1n_2n_3}\|^2_{n_1n_2n_3 \to n} \lesssim & N^{2s-2}_{123} N^{-2+2\alpha}_1 N^{-2+2\alpha}_2 N^{-2s'}_3  \cdot N^2_1N^{\f 32}_2
\end{align}
if $0<\alpha\le \f 14$ and $s=s'=\f 14 +\epsilon$ with $0<\epsilon<\f 34$, we have 
\begin{equ}
\|h^{m,\mathrm{dyadic}}_{nn_1n_2n_3}\|^2_{n_1n_3\to nn_2} \lesssim N^{-\f 32+2\epsilon}_{123}  N^{-\f 12+2\alpha}_2 N^{-\f 1 2+2\alpha-2\epsilon}_3 ;
\end{equ}
If $\f 14 \le \alpha < \f 12$ and $s=s'=2\alpha-\f 14 +\epsilon$ with $0<\epsilon<\f 54-2\alpha$, we have
\begin{equ}
\|h^{m,\mathrm{dyadic}}_{nn_1n_2n_3}\|^2_{n_1n_3\to nn_2} \lesssim N^{-\f 52+4\alpha+2\epsilon}_{123} N^{-2\epsilon}_3 \le N^{-2\epsilon}_3.
\end{equ}

\end{itemize}

\item $\|h^{m,\mathrm{dyadic}}_{nn_1n_2n_3}\|_{n_3\to nn_1n_2}$. In fact, this is similar to the case above with applying the sup-counting estimate (ii) instead of (i) in Lemma \ref{lem:cubic-sup-count}. We point out this case will give a more restrictive upper bound for $\alpha$. Assume $N_1 \geq N_2$.

By Schur test, we deduce $\|h^{m,\mathrm{dyadic}}_{nn_1n_2n_3}\|_{n_1n_3\to nn_2}$ to 
\begin{align}
\|h^{m,\mathrm{dyadic}}_{nn_1n_2n_3}\|^2_{n_1n_3\to nn_2} \lesssim & N^{2s-2}_{123} N^{-2+2\alpha}_1 N^{-2+2\alpha}_2 N^{-2s'}_3 \\
\quad \times &\medcurly{N_1,N_2,N_{123}}^2 \mincurly{N_1,N_2,N_{123}}^{\frac{3}{2}}.
\end{align}
Again, we do it case by case.
\begin{itemize}
\item[(i)] $N_1 \ge N_2 \ge N_{123}$. It yields that
\begin{align}
\|h^{m,\mathrm{dyadic}}_{nn_1n_2n_3}\|^2_{n_1n_3\to nn_2} \lesssim & N^{2s-2}_{123} N^{-2+2\alpha}_1N^{-2+2\alpha}_2 N^{-2s'}_3 \cdot N^2_2N^{\frac{3}{2}}_{123}\\
= & N^{2s- \f 12}_{123} N^{-2+2\alpha}_1N^{2\alpha}_2 N^{-2s'}_3.
\end{align}

If $N_3 = \maxcurly{N_1,N_3}$, then $N_{123} \lesssim N_3$ and we have 
\begin{align}
\|h^{m,\mathrm{dyadic}}_{nn_1n_2n_3}\|^2_{n_1n_3\to nn_2} \lesssim &  N^{-2+4\alpha}_1 N^{-2s'+2s- \f 12}_3 \le N^{-2s'+2s- \f 12}_3,
\end{align}
when $0<\alpha<\f 12$. Hence, for any $s=s'$, we have the desired estimate.

If $N_1 = \maxcurly{N_1,N_3}$, then $N_{123} \lesssim N_1$ and we have 
\begin{align}
\|h^{m,\mathrm{dyadic}}_{nn_1n_2n_3}\|^2_{n_1n_3\to nn_2} \lesssim   N^{-\f 52+2s+4\alpha}_1,
\end{align}
due to $s'>0$.

Hence, for $0< \alpha <\f 14$, and $s=s' =\f 14+\epsilon$, we have $\|h^{m,\mathrm{dyadic}}_{nn_1n_2n_3}\|^2_{n_1n_3\to nn_2} \lesssim   N^{-2+2\epsilon+2\alpha}_1$ which is accpetable when $\epsilon < 1-\alpha$. For $\alpha \ge \f 14$, by setting $s=s'=2\alpha-\f 14 +\epsilon$, we can see 
\begin{align}
\|h^{m,\mathrm{dyadic}}_{nn_1n_2n_3}\|^2_{n_1n_3\to nn_2} \lesssim   N^{-3+8\alpha+2\epsilon}_1.
\end{align}
Therefore, we need $\alpha <\f 38$ and $\epsilon < \f12(3-8\alpha)$.

\item[(ii)] $N_1 \ge N_{123} \ge N_2$. We have 
\begin{align}
\|h^{m,\mathrm{dyadic}}_{nn_1n_2n_3}\|^2_{n_1n_3\to nn_2} \lesssim & N^{2s-2}_{123} N^{-2+2\alpha}_1N^{-2+2\alpha}_2 N^{-2s'}_3 \cdot N^2_{123}N^{\frac{3}{2}}_2 \\
= & N^{2s}_{123} N^{-2+2\alpha}_1N^{-\f 12 +2\alpha}_2 N^{-2s'}_3.
\end{align}
Hence, by the fact $N_{123} \le N_1$, we have 
\begin{equ}
\|h^{m,\mathrm{dyadic}}_{nn_1n_2n_3}\|^2_{n_1n_3\to nn_2} \lesssim N^{-2+2\alpha+2s}_1 N^{-\f 12 +2\alpha}_2 N^{-2s'}_3.
\end{equ}
What follows is the same as case $N_1 \ge N_2 \ge N_{123}$, we omit the details.

\item[(iii)] $N_{123} \ge N_1 \ge N_2$. Note that $N_{123} \lesssim \maxcurly{N_1,N_3}$, we have the similar analysis as case $N_1 \ge N_2 \ge N_{123}$, we omit the details as well.
\end{itemize}

\item $\|h^m_{nn_1n_2n_3}\|_{n_1n_3 \to n_2n}$. In this case, we can directly abandon the effect of the dispersion by using the crude bound $\1_{|\phi -m | \leq 1} \leq 1$.
Again, by Schur's test,
\begin{align}
\|h^m_{nn_1n_2n_3}\|^2_{n_1n_3 \to n_2n} \lesssim & N^{2s-2}_{123} N^{-2+2\alpha}_1N^{-2+2\alpha}_2N^{-2s'}_3\\
\times & \sup_{n,n_2 \in \Z^2} \sum_{n_1,n_3 \in \Z^2}  \prod_{j=1}^3 \mathbf{1}_{|n_j| \sim N_j} \1_{|n| \sim N_{123}} \mathbf{1}_{n=n_{123}}\\
\times &  \sup_{n_1,n_3 \in \Z^2} \sum_{n_2, n \in \Z^2} \prod_{j=1}^3 \mathbf{1}_{|n_j| \sim N_j} \1_{|n| \sim N_{123}} \mathbf{1}_{n=n_{123}}\\
\lesssim & N^{2s-2}_{123} N^{-2+2\alpha}_1N^{-2+2\alpha}_2N^{-2s'}_3 \mincurly{N_1,N_3}^2 \mincurly{N_2,N_{123}}^2.
\end{align}
This leads to the following four subcases.

\begin{itemize}
\item[(i)] $N_1 \le N_3$ and $N_2 \le N_{123}$. In this case, we have $\maxcurly{N_1,N_2,N_3} = N_3$. Then
\begin{equs}
\|h^m_{nn_1n_2n_3}\|^2_{n_1n_3 \to n_2n} \lesssim N^{2s-2}_{123} N^{2\alpha}_1 N^{2\alpha}_2 N^{-2s'}_3 \le N^{2s-2+2\alpha}_{123}N^{-2s'+2\alpha}_3
\end{equs}

If $0 < \alpha < \f 14$ and $s=s'=\f 14+\epsilon$, we need $\f 14 +\epsilon < 1-\alpha$ so that 
\begin{equs}
\|h^m_{nn_1n_2n_3}\|^2_{n_1n_3 \to n_2n} \lesssim N^{2s-2+2\alpha}_{123}N^{-2s'+2\alpha}_3 \le N^{-\f 12 +2\alpha -2\epsilon}_3
\end{equs}
which implies $\epsilon < \f 34 -\alpha$.

If $\f 14 \le \alpha < \f 5{12}$ and $s=s'=2\alpha - \f 14+\epsilon$, we need $2\alpha -\f 14 +\epsilon < 1-\alpha$ so that 
\begin{equs}
\|h^m_{nn_1n_2n_3}\|^2_{n_1n_3 \to n_2n} \lesssim N^{-\f 52+6\alpha+2\epsilon}_{123}N^{\f 12-2\alpha-2\epsilon}_3 \le N^{-2\epsilon}_3
\end{equs}
which requires $\epsilon < \f 54 -3\alpha$.

\item[(ii)] $N_1 \le N_3$ and $N_{123} \le N_2$. In this case, we don't know whether $N_3$ is larger than $N_2$ or not. However, we have 
\begin{equs}
\|h^m_{nn_1n_2n_3}\|^2_{n_1n_3 \to n_2n} \lesssim N^{2s-2+2\alpha}_2  N^{-2s'+2\alpha}_3.
\end{equs}
If $0 < \alpha < \f 14$ and $s=s'=\f 14+\epsilon$ with $0<\epsilon < \f 34 -\alpha$, we obtain 
\begin{equs}
\|h^m_{nn_1n_2n_3}\|^2_{n_1n_3 \to n_2n} \lesssim N^{-\f 32+2\alpha+2\epsilon}_2  N^{-\f 12+2\alpha-2\epsilon}_3.
\end{equs}
Note by taking a smaller $\epsilon< \f 14 -\alpha$ and we have $\|h^m_{nn_1n_2n_3}\|^2_{n_1n_3 \to n_2n} \lesssim (N_2N_3)^{-2\epsilon}$. 

If $\f 14 \le \alpha < \f 5{12}$ and $s=s'=2\alpha - \f 14+\epsilon$ with $0<\epsilon < \f 54 -3\alpha$ , we obtain 
\begin{equs}
\|h^m_{nn_1n_2n_3}\|^2_{n_1n_3 \to n_2n} \lesssim N^{-\f 52+6\alpha+2\epsilon}_2 N^{\f 12-2\alpha-2\epsilon}_3 \le (N_2N_3)^{-2\epsilon}
\end{equs}
where we relabel $\epsilon = \mincurly{\f 54-3\alpha-\epsilon, -\f 14+\alpha+\epsilon}$.

\item[(iii)] $N_3 \le N_1$ and $N_{123} \le N_2$. Then we have 
\begin{equs}
\|h^m_{nn_1n_2n_3}\|^2_{n_1n_3 \to n_2n} \lesssim N^{-2s'+2\alpha}_1 N^{2s-2+2\alpha}_2
\end{equs}
which is similar to case (i).
\item[(iv)] $N_3 \le N_1$ and $N_2 \le N_{123}$. In this case, we have $\maxcurly{N_1,N_2,N_3} = N_1$. 
\begin{equs}
\|h^m_{nn_1n_2n_3}\|^2_{n_1n_3 \to n_2n} \lesssim N^{-2s'+2\alpha}_1 N^{2s-2+2\alpha}_2
\end{equs}
which is similar to case (i).
\end{itemize}

\item $\|h^m_{nn_1n_2n_3}\|_{n_1n_3 \to n_2n}$. This can be done similarly as $\|h^m_{nn_1n_2n_3}\|_{n_1n_3 \to n_2n}$.
\end{itemize}
\end{proof}

\subsection{Sine-cancellation}
\begin{lemma}\label{lem:sine-cancellation}
Let  $T \geq 1$,  $ J \subset [0,T]$ an interval, and $A,N \in \N$. Furthermore, assume $|a| \lesssim A \ll N$. Let $f \colon \R \times \R \times \Z^2 \to \C$ satisfies
\begin{align}\label{eqn:sine-cancellation}
\ab{f(t',t,n)} \leq A \jpb{n}^{-3}, \quad \ab{f(t',t,n)-f(t',t,-n)}\leq A\jpb{n}^{-4}, \quad \ab{\partial_{t'}f(t',t,n)} \leq A\jpb{n}^{-3}.
\end{align}
Then, we have
\begin{align}
\sup_{|t| \leq T} \ab{\sum_{|n| \sim N} \int_0^t \1_J \sin\pa{(t-t')\jpb{a+n}} \cos\pa{(t-t')\jpb{n}}f(t',t,n)\, \d t'} \lesssim T^2 A^2 \log(2+N) N^{-1}.
\end{align}
\end{lemma}

\begin{proof}
We first split the integral:	
\begin{align}
		&2\sum_{|n| \sim N} \int_{0}^{t} \1_J\sin((t-t')\jpb{a+n})\cos((t-t')\jpb{n})f(t',t,n)\,\d t'\\
		= & \sum_{|n| \sim N} \int_{0}^{t} \1_J \Bigl( \sin((t-t')(\jpb{a+n}-\jpb{n})) f(t',t,n)\,\d t'
		+\sin((t-t')(\jpb{a+n}+\jpb{n})) f(t',t,n)\Bigr)\,\d t'.
	\end{align}
We estimate the two summands separately.

A. We first estimate
		$$I=\sum_{|n| \sim N} \int_{0}^{t} \1_J\sin((t-t')(\jpb{a+n}-\jpb{n})) f(t',t,n)\,\d t'.$$
To obtain a decay in $N$, we observe that sum in invariant under the map $n\mapsto -n$. Matching this symmetry with the anti-symmetry in  the $\sin$ function, we see that
		\begin{align}
			& \ab{\sum_{|n| \sim N} \sin((t-t')(\jpb{a+n}-\jpb{n})) f(t',t,n)}\\
			=& \frac{1}{2} \ab{\sum_{|n| \sim N} \brt{\sin((t-t')(\jpb{a+n}-\jpb{n}))  f(t',t,n) +\sin((t-t')(\jpb{a-n}-\jpb{n}))  f(t',t,-n)  } }\\
			\leq & \sum_{|n| \sim N}\ab{\sin\pa{(t-t')(\jpb{a+n}-\jpb{n})} +\sin\pa{(t-t')(\jpb{a-n}-\jpb{n})} } \cdot |f(t',t,n)|\\
			&+  \sum_{|n| \sim N} \ab{f(t',t,n) -f(t',t,-n) }
		\end{align}
By the assumption on $f$, the last term is bounded above by $AN^{-2}$. On the other hand,		\begin{align}
			& \ab{\sin\pa{(t-t')(\jpb{a+n}-\jpb{n})} +\sin\pa{(t-t')(\jpb{a-n}-\jpb{n})}}\\
			= & \ab{\sin\pa{(t-t')(\jpb{a+n}-\jpb{n})} -\sin\pa{(t-t')(-\jpb{a-n}+\jpb{n})}}\\
			\le & T\ab{\jpb{a+n}-\jpb{n}+\jpb{a-n}-\jpb{n}}\\
			\lesssim &TA^2N^{-1}, 
		\end{align}
		where we used the fact  that by Taylor expansion
		\begin{equs}
		 \jpb{a\pm n}-\jpb{n} = \f {a \cdot n}{\jpb{n}} + O(\f {A^2}{N})
		\end{equs}	
		Now, integrating with $t'$, one can see $I$ can be controlled by $T^2A^2N^{-1}$.
		
B.  Different from the firs term, we can expect to gain a decay, $N^{-1}$, in 
				$$II:= \sum_{|n| \sim N} \int_{0}^{t} \1_J\sin((t-t')(\jpb{a+n}+\jpb{n})) f(t',t,n)\,\d t'$$
		from integrating in time.  since $N \lesssim \jpb{a+n}+\jpb{n}$. That is 
		\begin{align}
			& \ab{\sum_{|n| \sim N} \int_{0}^{t} \1_J\sin((t-t')(\jpb{a+n}+\jpb{n})) f(t',t,n)\,\d t'}\\
			\lesssim & \max_{\pm} \ab{\sum_{|n| \sim N} \int_{0}^{t} \1_J \exp{\i(\jpb{a+n} \pm \jpb{n})t'} f(t',t,n)\,\d t'} \\
			\lesssim & \max_{\pm} \sum_{|n| \sim N} \frac{1}{\jpb{ \jpb{a+n}+\jpb{n}}}\pa{\sup_{0 \leq t' \leq T} |f(t',t,n)| + T\sup_{0 \leq t' \leq T}|\partial_{t'}f(t',t,n)| }\\
			\lesssim & \sum_{|n| \sim N} \frac{1}{\jpb{ \jpb{a+n}+\jpb{n}} }\cdot TAN^{-3}.
		\end{align}
			
	Since $x \mapsto \jpb{x}$ is globally $1-$Lipschitz, we may transfer the discrete sum to volume computation:
		\begin{align}
			\sum_{|n| \sim N} \frac{1}{1+\ab{\jpb{a+n}\pm \jpb{n}}} \lesssim \int_{\R^2} \mathbf{1}_{|\xi| \sim N} \frac{\d \xi}{1+\ab{\jpb{a + \xi}\pm\jpb{\xi}}}.
		\end{align}
		
		By rotation, without loss of generality, we can assume $a=(0,|a|)$. Consider $(\xi=\xi^{(1)},\xi^{(2)})$ in the polar coordinates $(r,\theta)$. For fixed $\theta$, $r \mapsto \jpb{a + \xi}+\jpb{\xi}$ is bi-Lipschitz on $r \sim N$. By a further change of variable, one has
		\begin{equation}
			\int_{\R^2} \mathbf{1}_{|\xi| \sim N} \frac{\d \xi}{1+\ab{\jpb{a + \xi}\pm \jpb{\xi}}} \lesssim N \int_0^\infty \mathbf{1}_{r \sim N} \frac{1}{1+\ab{r \pm \lambda}}\, \d r \lesssim N \log(2+N),
		\end{equation}
		which seems a little better than we expected.
		
This can also be obtained with the basic counting estimate. 
Define $\phi(a,n)=\jpb{a+n}+\jpb{n}$. If we insert the indicator function $\1_{|\phi-m| \le 1}$, then we have
		\begin{align}
			\sum_{|n| \sim N} \frac{1}{\jpb{\jpb{a+n}+\jpb{n}}} & \lesssim \sum_{|n| \sim N} \sum_{m \in \Z} \frac{1}{\jpb{m}} \1_{|\phi-m| \leq 1} = \sum_{m \in \Z} \frac{1}{\jpb{m}}  \sum_{|n| \sim N} \1_{|\phi-m| \leq 1} \\
			&=\sum_{m \in \Z^2} \frac{1}{\jpb{m}} \#\{n \in \Z^2 \colon |n| \sim N, \, |\jpb{a+n}+\jpb{n}-m| \leq 1\}\\
			&\lesssim \log(2+N) \sup_{m \in \Z}\#\{n \in \Z^2 \colon |n| \sim N, \, |\jpb{a+n}+\jpb{n}-m| \leq 1\} \\
			&\lesssim \log(2+N)N,
		\end{align}
	where in the last line we have applied the improved bound in Remark \ref{elliptic-vs-hyperbolic}.
\end{proof}

\section{Counting estimates}\label{sec:counting}
We collect counting estimates essential in describing the wave interactions and make a convention on the notations. 

Let the capital letters $N_1,N_2,\cdots$ in this section denote dyadic numbers $2^k, 2^{k_1}, 2^{k_2}, \cdots$ for some $\ell, k, k_1, k_2 \in \N$. The symbol $A \sim B$ refers to that quantity $A$ and $B$ are comparable, specifically $\f12 B \le A \le 2B$. This is a symmetric relation. The symbols $A \ll B$ and $A \gg B$ refer to $A \le \f 12 B$ and $A \ge 2B$ respectively. The letters $n,n_1,n_2, \cdots \in \Z^2$  denote frequencies. We will frequently use the notation $|n_j| \sim N_j$ where $j\in\N$ instead of $|n_j| \sim N_{k_j} =2^{k_j}$ for some $k_j \in \N$.
In addition, we denote $\medcurly{N_1, N_2, N_3}$ the middle number, so if $N_1<N_2<N_3$, then 
$$\medcurly{N_1, N_2, N_3}=N_2.$$

\subsection{Discrete convolution sum} 
We recall two basic counting estimates that have been instrumental in the study of quadratic nonlinearities.
\begin{lemma}[Discrete convolution sum]\label{lem:convolution-sum}~\\
\begin{itemize}
\item Let $d \geq 1$ and $\alpha,\beta \in \R$ satisfy 
\begin{equation}
\alpha + \beta > d, \quad \mbox{and}\quad \alpha,\beta < d.
\end{equation}
Then, for any $n \in \Z^{d}$, we have
\begin{equation}
\sum_{n=n_{1}+n_{2}} \frac{1}{\jpb{n_{1}}^{\alpha}\jpb{n_{2}}^{\beta}} \lesssim \jpb{n}^{d-\alpha-\beta}.
\end{equation}
\item Let $d \geq 1$ and $\alpha,\beta \in \R$ satisfy $\alpha + \beta > d$. Then, we have
\begin{equation}
\sum_{\substack{n=n_{1}+n_{2} \\ |n_{1}| \sim |n_{2}|}} \frac{1}{\jpb{n_{1}}^{\alpha}\jpb{n_{2}}^{\beta}} \lesssim \jpb{n}^{d-\alpha-\beta}.
\end{equation}
\end{itemize}
\end{lemma}
The  Lemma is derived from basic calculations. a proof
can be found in \cite[Lemmas 4.1 \& 4.2]{MWX16}.

\subsection{Basic counting estimate}
We now derive the basic counting estimates for cubic nonlinearities that are used in this article. We first quote the following well known frequency counting estimates.
\begin{lemma}[Basic counting estimate.]\label{lem:basic-count}
Let $A,N \in 2^{\N}$, $a \in \Z^{2}$ satisfying $|a| \sim A$. Then,
\begin{equs}
\sup_{m \in \Z} \#\{ n \in \Z^{2} \colon |n| \sim N, |a| \sim A, |\jpb{a+n}-\jpb{n}-m| \le 1 \}
\lesssim \f{N^{2}}{\min(A,N)^{\f 12}},\label{eq:basic-count-hyperbolic}\\
\sup_{m \in \Z} \#\{ n \in \Z^{2} \colon |n| \sim N, |a| \sim A, |\jpb{a+n}+\jpb{n}-m| \le 1 \}
\lesssim \f{N^{2}}{N^{1/2}} .\label{eq:basic-count-elliptic}
\end{equs}
\end{lemma}
Such lattice point estimates are fundamental in the study of wave equations. The above lemma is directly quoted from \cite[Lemma 3.12]{Bri23}, see also the references therein.

\begin{remark}\label{elliptic-vs-hyperbolic}
The two estimates in Lemma \ref{lem:basic-count} are respectively referred to as elliptic and hyperbolic.  we refer to the positive sign case $|\jpb{a+n}+ \jpb{n}-m|$,   as elliptic,  and to the negative sign case $|\jpb{a+n}- \jpb{n}-m|$ as hyperbolic. 
The elliptic estimates is better than the hyperbolic estimates. Furthermore, by a geometric argument, the elliptic estimate can be even improved to $N$.
\end{remark}

We adapt these basic estimates to include one more constraint $B$.
\begin{lemma}[Two ball counting estimate.]\label{lem:two-ball-count}
Let $N,A,B\in 2^{\N_0}$ and $a \in \Z^{2}$ satisfying $|a| \sim A$. Then
\begin{equ}\label{eq:two-ball-count}
\sup_{m \in \Z} \# \{ n \in \Z^{2} \colon |n| \sim N, |n+a| \sim B, |\jpb{a+n}\pm \jpb{n}-m| \le 1 \} \lesssim\f{ \mincurly{B,N}^{2}}{ \mincurly{A,B,N}^{1/2}}.
\end{equ}
\end{lemma}

\begin{proof} 
By Lemma \ref{lem:basic-count},  if $|a|\sim A$, denoting $$\D_{N,m,a, A}= \{ n \in \Z^{2} \colon |n| \sim N, |a| \sim A, |\jpb{a+n} \pm \jpb{n}-m| \le 1 \},$$
 \begin{equ}
\sup_{m \in \Z} \#\D_{N,m,a, A}
 \lesssim \f{N^{2}}{\min(A,N)^{\f 12}}.\label{eq:basic-count}
\end{equ}
Now the set on the left hand side of (\ref{eq:two-ball-count}) is included in $\D_{N,m,a, A}$, and
\begin{align}
\mathrm{LHS}  \le  \sup_{m \in \Z} \# \D_{N,m,a, A} \lesssim \f{N^{2}}{\mincurly{A,N}^{1/2}}.
\end{align}

On the other hand, by a change of variable $n'=n+a$, we have another bound on  the left hand side of (\ref{eq:two-ball-count}):
\begin{align}
\mathrm {LHS}
\le  \sup_{m \in \Z} \#\{ n' \in \Z^{2} \colon |n'| \sim B, |\jpb{n'} \pm \jpb{n'-a}-m| \le 1 \} \lesssim  \f{B^2} {\mincurly{A,B}^{1/2}}.
\end{align}
Combining the two, the desired estimate follows.
\end{proof}

\subsection{Cubic counting estimate}
Given $\sigma, \sigma_1, \sigma_2,\sigma_3 \in \{\pm 1\}$, we define the phase function
\begin{equ}\label{eq:cubic-phase}
\phi(n_1,n_2,n_3) = \sigma\jpb{n_{123}}+\sigma_1\jpb{n_1} +\sigma_2\jpb{n_2} +\sigma_3\jpb{n_3}.
\end{equ}
For brevity, we will frequently omit the argument $(n_1,n_2,n_3)$.
\begin{lemma}[Cubic counting estimate.]\label{lem:cubic-count}
Let the signs $\sigma, \sigma_1, \sigma_2,\sigma_3 \in \{\pm 1\}$ be given. The following  estimates hold.
\begin{itemize}
\item[(i)] For any fixed dyadic numbers $N_1,N_2,N_3$, \begin{equs}
& \sup_{m \in \Z} \#\{ (n_1,n_2,n_3) \colon |\phi-m|\le 1,  |n_i| \sim N_i, i=1,2,3\}
 \lesssim \f{(N_1 N_2 N_3)^{2}}{\medcurly{N_1,N_2,N_3}^{\f 12}}.
\label{eq:cubic-count1}
\end{equs}
\item[(ii)] For any dyadic numbers $N_{123},N_1,N_2$, 
\begin{align}
& \sup_{m \in \Z} \#\{ (n_1,n_2,n_{123}) \colon  |\phi-m|\le 1, |n_1| \sim N_1, |n_2| \sim N_2, |n_{123}| \sim N_{123}\}\\
& \lesssim \f{(N_{123} N_1 N_2)^{2}}{\medcurly{N_{123},N_1,N_2}^{1/2}}.
\label{eq:cubic-count2}
\end{align}
\item[(iii)] For any dyadic numbers $N_2,N_{3},N_{34}$ and $\psi(n_2,n_3,n_4) = \sigma\jpb{n_{234}}+\sigma_2\jpb{n_2} +\sigma_3\jpb{n_3} +\sigma_4\jpb{n_4}$,
\begin{align}
& \sup_{m \in \Z}\#\{ (n_2,n_3,n_{34}) \colon  |\psi-m|\le 1, |n_2| \sim N_2, |n_3| \sim N_3, |n_{34}| \sim N_{34}\}\\
& \lesssim \f{(N_2N_3N_{34})^{2}} { \mincurly{N_{34},\maxcurly{N_2,N_3}} ^{1/2}}.
\label{eq:cubic-count3}
\end{align}
\item[(iv)] For any $N_{3},N_{34},N_{234} \ge 1$ and $\psi$ as in (iii),
\begin{align}
& \sup_{m \in \Z} \#\{ (n_3,n_{34},n_{234}) \colon  |\psi-m|\le 1, |n_3| \sim N_3, |n_{34}| \sim N_{34}, |n_{234}| \sim N_{234}\}\\
& \lesssim \f{(N_3N_{34}N_{234})^{2}} { \mincurly{N_{34},\maxcurly{N_3,N_{234}}} ^{1/2}}.
\label{eq:cubic-count4}
\end{align}
\end{itemize}
\end{lemma}

\begin{proof}

(i) By symmetry, without loss of generality, we may assume $N_1 \ge N_2 \ge N_3$. We rewrite
$$ |\phi-m|=|\jpb{n_2+n_{13}} \pm \jpb{n_2}- (m \mp \jpb{n_1} \mp \jpb{n_3})|.$$
Applying the basic counting estimate \eqref{eq:basic-count-hyperbolic}, with$N=|n_2|$ and $a=n_{13}$, we obtain:
\begin{equs}
{} & \#\{ (n_1,n_2,n_3) \colon |\phi- m|\le 1,  |n_1| \sim N_1, |n_2| \sim N_2 , |n_3| \sim N_3 \}\\ 
\lesssim  & \sum_{n_1,n_3 \in \Z^2} \prod_{j=1,3} \1_{|n_j| \sim N_j} \Big( \sum_{n_2 \in \Z^2} \1_{|\phi- m|\le 1,  |n_1| \sim N_1, |n_2| \sim N_2 , |n_3| \sim N_3 } \Big)\\ 
\lesssim  & \sum_{n_1,n_3 \in \Z^2} \prod_{j=1,3} \1_{|n_j| \sim N_j}  \f{(N_2)^2}  {\min(\jpb{n_{13}},N_2)^{1/2}} \\
\le  & \sum_{n_1,n_3 \in \Z^2} \prod_{j=1,3} \1_{|n_j| \sim N_j}  \f{(N_2)^2}  {\jpb{n_{13}}^{1/2}} \1_{\jpb{n_{13}\le N_2}} + \sum_{n_1,n_3 \in \Z^2} \prod_{j=1,3} \1_{|n_j| \sim N_j} (N_2)^{\f 32}\\
\lesssim &(N_1)^{\f 32} (N_2)^2 (N_3)^2 + (N_1)^2  (N_2)^{\f 32} (N_3)^2
\lesssim \f{ (N_1 N_2 N_3)^2}{(N_2)^{\f 12}},
\end{equs}
proving (i).

(ii) By viewing $n_{123}$ as a free variable, the phase function will be replaced by
\begin{equ}
\phi(n_1,n_2,n_{123}) = \sigma \jpb{n_{123}}  + \sigma_1 \jpb{n_1} + \sigma_2 \jpb{n_2}  + \sigma_3 \jpb{n_{123}-n_1-n_2}.
\end{equ}
Taking the transformation $(n_1,n_2) \mapsto (-n_1,-n_2)$, we obtain
\begin{equ}
\phi(n_1,n_2,n_{123}) = \sigma \jpb{n_{123}}  + \sigma_1 \jpb{n_1} + \sigma_2 \jpb{n_2}  + \sigma_3 \jpb{n_{123}+n_1+n_2}.
\end{equ}
If we relabel $n_{123}$ as $n_3$, the counting problem reduces to the same one as in case (i), except that we now impose $|n_3| \sim N_{123}$.

(iii) Note that the phase function can be rewritten as
\begin{equ}
\psi(n_2,n_3,n_{34}) = \sigma \jpb{n_2+n_{34}}  + \sigma_2 \jpb{n_2} + \sigma_3 \jpb{n_3}  + \sigma_4 \jpb{n_{34}-n_3}.
\end{equ}
where we view $n_3$, $n_{34}$ and $n_{234}$ as independent indices.
By first summing over $n_3 \in \Z^2$ with basic counting estimate \eqref{eq:basic-count-hyperbolic}, and then $n_{34},n_2 \in \Z^2$, we obtain
\begin{align}
\#\{ (n_2,n_3,n_{34}) \colon  |\psi-m|\le 1, |n_j| \sim N_j, j=2,3, |n_{34}| \sim N_{34}\} \lesssim  \mincurly{N_{34},N_3} ^{-\f 12} N^2_2N^2_3N^2_{34}.
\end{align}
However, if we first sum in $n_2 \in \Z^2$, we have
\begin{align}
\#\{ (n_2,n_3,n_{34}) \colon  |\psi-m|\le 1, |n_j| \sim N_j, j=2,3, |n_{34}| \sim N_{34}\} \lesssim  \mincurly{N_{34},N_2} ^{-\f 12} N^2_2N^2_3N^2_{34}.
\end{align}
Since $\maxcurly{\mincurly{N_{34},N_3},\mincurly{N_{34},N_2} } = \mincurly{N_{34},\maxcurly{N_2,N_3}}$, combining the above two inequalities give the desired estimate.

(iv)  Compared with case (iii),  there is the additional restriction $|n_{234}|\sim N_{234}$.
The proof is almost the same as (iii) by viewing the phase function as 
\begin{equ}
\psi(n_3,n_{34},n_{234}) = \sigma \jpb{n_{234}}  + \sigma_2 \jpb{n_{234}-n_{34}} + \sigma_3 \jpb{n_3}  + \sigma_4 \jpb{n_{34}-n_3}.
\end{equ}
where we again view $n_3$, $n_{34}$ and $n_{234}$ as independent indices.
We first sum over $n_{234} \in \Z^2$, viewing $ \sigma_3 \jpb{n_3}  + \sigma_4 \jpb{n_{34}-n_3}$ as the auxiliary $m$, applying the basic counting estimate \eqref{eq:basic-count-hyperbolic}, and then bound the summation over $n_{34},n_{234} \in \Z^2$ by $N^2$, we obtain
\begin{align}
&\#\{ (n_3,n_{34},n_{234}) \colon  |\psi-m|\le 1, |n_j| \sim N_j, |n_3|\sim N_3, |n_{34}| \sim N_{34}, |n_{234}| \sim N_{234}\}\\
& \lesssim  \mincurly{N_{34},N_3} ^{-\f 12} N^2_3N^2_{34}N^2_{234}.
\end{align} 
We then sum over $n_3 \in \Z^2$, applying again the basic counting estimate \eqref{eq:basic-count-hyperbolic}  and setting $m=\sigma \jpb{n_{234}}  + \sigma_2 \jpb{n_{234}-n_{34}} $, the summand does not involve $n_3$, 
and finally bound the number of  $n_{34}$ and $n_{234}$ by $N^2$ to obtain
\begin{align}
&\#\{ (n_3,n_{34},n_{234}) \colon  |\psi-m|\le 1, |n_j| \sim N_j, |n_3|\sim N_3, |n_{34}| \sim N_{34}, |n_{234}| \sim N_{234}\}\\
& \lesssim  \mincurly{N_{34},N_3} ^{-\f 12} N^2_3N^2_{34}N^2_{234}.
\end{align} 
Finally, by interpolation, we have the desired result, c.f. the proof of (iii).
\end{proof}

Below $\phi(n_1,n_2,n_3)$ is as in \eqref{eq:cubic-phase} and recall $s_\alpha$ from
 \eqref{eq:s-alpha}:
\begin{equ}
s_\alpha = \begin{cases}
1-2\alpha,\quad &\mbox{if} \quad 0<\alpha<\frac{1}{4},\\
\frac{5}{4}-3\alpha,\quad &\mbox{if} \quad \frac{1}{4}\leq \alpha<\frac{5}{12}.
\end{cases}
\end{equ}

\begin{convention}\label{convention-dya-sum}
When introducing dyadic summations,we often use abbreviated notation. For example,
$$  \sum_{N_{123} =1 }^{\max{(N_1, N_2, N_3)}} N_{123}^d N_1^aN_2^bN_3^c
= \sum_{k_{123} =0 }^{\log (\max{(N_1, N_2, N_3)})} 2^{d\; k_{123}} \;N_1^a N_2^bN_3^c.  $$
\end{convention}

\begin{proposition}[Cubic sum estimate.]\label{prop:cubic-sum}
  Let   $\sigma, \sigma_1, \sigma_2,\sigma_3 \in \{\pm 1\}$. Then for any $s > 0$, and dyadic numbers $N_1,N_2,N_3$,
 the following holds:
\begin{equ}\label{eq:cubic-sum}
\sup_{m \in \Z} \sum_{n_1,n_2,n_3 \in \Z^2} \jpb{n_{123}}^{2(s-1)}  \Bigl[ \prod_{j=1}^3 \1_{|n_j|\sim N_j}  \jpb{n_{j}}^{-2+2\alpha} \1_{|\phi -m|\le 1} \Bigr] \lesssim  \maxcurly{N_1,N_2,N_3}^{2(s-s_{\alpha})}.
\end{equ}
\end{proposition}

\begin{proof}
By symmetry in $n_1,n_2,n_3$, without loss of generality that  $N_1 \ge N_2 \ge N_3$. 
Define
\begin{equs}
\mathcal{C}(m) 
&= \{(n_1,n_2,n_3) \in (\Z^{2})^{3} \colon |n_j| \sim N_j, j=1,2,3, |n_{123}| \sim N_{123}, |\phi -m| \le 1 \}.
\end{equs}
Then from (\ref{eq:cubic-count2}) in Lemma \ref{lem:cubic-count}, 
$$ \#\mathcal{C}(m)\lesssim  \medcurly{N_{123},N_2,N_3}^{-1/2} (N_{123} N_2 N_3)^2.$$

By inserting a dyadic decomposition with respect to frequency $n_{123}$, and noting that $|n_{123}|\lesssim N_1+N_2+N_3$, we obtain:
\begin{equs}
{} & \sum_{n_1,n_2,n_3 \in \Z^2} \jpb{n_{123}}^{2(s-1)}  \Bigl[ \prod_{j=1}^3 \1_{|n_j|\sim N_j}  \jpb{n_{j}}^{-2+2\alpha} \1_{|\phi -m|\le 1} \Bigr] \\
\lesssim  & \sum_{k_{123} =0 }^{\log (N_1)}\#\mathcal{C}(m) \;2^{(2s-2)k_{123}} (N_1N_2N_3)^{-2+2\alpha}\\
\lesssim & \sum_{N_{123} =1 }^{N_{1} } (N_{123})^{2s} (N_1)^{-2+2\alpha}(N_2)^{2\alpha}(N_3)^{2\alpha} \medcurly{N_{123},N_2,N_3}^{-\f 12}.
\label{eq:cubic-sum-est}
\end{equs}
Here $N_{123}$ is of the form $2^{k_{123}}$ where $k_{123}$ takes values in $\N$. In the last step we applied (\ref{eq:cubic-count2}) and adopted Convention \ref{convention-dya-sum} on dyadic summations introduced earlier.

It remains to show that $$ \sum_{N_{123} =1 }^{N_{1} } (N_{123})^{2s} (N_1)^{-2+2\alpha}(N_2)^{2\alpha}(N_3)^{2\alpha}<\infty.$$
Note that under the ordering $N_1 \ge N_2 \ge N_3$, we must have $\medcurly{N_{123},N_2,N_3} \le N_2$ which leads to two further subcases.

\textbf{Case 1}.  If $\medcurly{N_{123},N_{2},N_{3}} \sim N_{2}$,  we have
\begin{align}
\eqref{eq:cubic-sum-est} \lesssim \sum_{N_{123} = 1}^{N_1} (N_{123})^{2s} (N_1)^{-2+2\alpha}(N_2)^{-\f 12 + 2\alpha}(N_3)^{2\alpha}.
\end{align}
Although we can control $N_3$ with $N_2$, it is sufficient to control it with $N_1$, then
the exponent $-\f 12 + 2\alpha$ on $N_2$ distinguishes the regularities by $0<\alpha<\f 14$ and $\alpha \ge \f 14$. Hence, if  $2\alpha -\f 12<0$, 
we bounded $(N_2)^{-\f 12 + 2\alpha}$ with $1$ 
 and bound the resulting geometric series 
 $$ \sum_{N_{123} =1 }^{N_{1} } (N_{123})^{2s} = \sum_{k_{123} = 0}^{\log(N_1)} 2^{2sk_{123}} \lesssim (N_1)^{2s}.$$
Consequently,
\begin{equ}
\eqref{eq:cubic-sum-est} \lesssim   (N_1)^{2s-2+4\alpha}(N_2)^{-\f 12 + 2\alpha} \lesssim (N_1)^{2s-2+6\alpha}
\end{equ}
This is summable over the dyadic scale $N_1=2^{k_{123}}$ (i.e. over $k_{123}$), precisely when $s<1-2\alpha$. 

If $\alpha \ge \f 14$, we should bound all frequencies by $N_1$, which leads to 
\begin{equ}
\eqref{eq:cubic-sum-est} \lesssim (N_1)^{2(s-\f 54+3\alpha)}.\end{equ}
which is summable when $s<\f 54 -3\alpha$. 

\textbf{Case 2}. If $\medcurly{N_{123},N_{2},N_{3}} \ll N_{2}$, then the decay will only occur on $N_3$ or $N_{123}$. This gives rise to two further subcases.

\textbf{Case 2.1} If $\medcurly{N_{123},N_{2},N_{3}} = N_3$, then
\begin{equ}
\eqref{eq:cubic-sum-est}
 \lesssim \sum_{N_{123} \ge 1} (N_{123})^{2s} (N_1)^{-2+2\alpha}(N_2)^{2\alpha}(N_3)^{2\alpha-\f 12}.
\end{equ}
As in case~1, we distinguish between the two regimes $\alpha<\f 14$, or $\alpha \ge \f 14$. Accordingly \begin{align}
\sum_{N_{123} \ge 1} (N_{123})^{2s} (N_1)^{-2+2\alpha}(N_2)^{2\alpha} (N_3)^{2\alpha-\f 12 }
 \lesssim \begin{cases}
N^{2(s-1+2\alpha)}_{1},\quad &\mbox{if}\, 0<\alpha <\frac{1}{4}\\
N^{2(s-\frac{5}{4}+3\alpha)}_{1},\quad &\mbox{if}\, \alpha \ge \frac{1}{4}
\end{cases}.
\end{align}
In both regimes the exponents are negative for $s<s_\alpha$, and the sums over the dyadic scale $k_{123}$ are therefore convergent for such $s$.

\textbf{Case 2.2} If $\medcurly{N_{123},N_{2},N_{3}} = N_{123}$, we have 
\begin{equ}
\eqref{eq:cubic-sum-est}
 \lesssim \sum_{N_{123} \ge 1} (N_{123})^{2s} (N_1)^{-2+2\alpha}(N_2)^{2\alpha}(N_3)^{2\alpha-\f 12},
\end{equ}
where we bounded $(N_{123})^{-\f 12}$ with $(N_3)^{-\f 12}$. If $s$ is negative, we can bound $\sum_{N_{123}\ge 1}N_{123}^{2s}$ with $1$ and bound $(N_2)^{2\alpha}$ with $(N_1)^{2\alpha}$, and obtain the following estimates :\begin{equ}
\eqref{eq:cubic-sum-est}
 \lesssim  (N_1)^{-2+4\alpha}(N_3)^{2\alpha-\f 1 2}.
\end{equ}
I f $\alpha< \f 14$, all exponents are negative, and the right hand side is summable proving the assertion. 
If on the other hand $\alpha \ge\f 14$,  $ (N_1)^{-\f 5 2+6\alpha}$ which is negative when $\alpha<\f 5 {12}$.

When the exponent $s$ is positive we have: $\sum_{N_{123} \ge 1} (N_{123})^{2s} \lesssim (N_1)^{2s}$. We used the fact that $N_{123}\le N_1$ and the summation of the geometric series yield the bound.
As before splitting into two cases. If $0<\alpha<\f 14$, we bound $N_2$ with $N_1$ and throw away $N_3$ to obtain: \begin{equ}
\eqref{eq:cubic-sum-est} \lesssim (N_1)^{2(s-1+2\alpha)},
\end{equ}
which has negative exponent when $s<s_\alpha$, proving the claim in this case.
For the other case $\alpha\ge \f 14$, the exponents of $N_2$ and $N_3$ are positive and so we have
\begin{equ}
\eqref{eq:cubic-sum-est} \lesssim (N_1)^{2s-\f 52+6\alpha},
\end{equ}
which is summable on dyadic scale for $s<s_\alpha$. This completes the proof of the proposition.
\end{proof}

\begin{lemma}[Cubic sup-counting estimate.]\label{lem:cubic-sup-count}
Let $N_{123},N_{1},N_{2},N_{3} \in 2^{\N}$, $m \in \Z$, and $\sigma, \sigma_1, \sigma_2,\sigma_3 \in \{\pm 1\}$. Then, we have the following counting estimates:
\begin{equs}
{}&\sup_{m \in \Z} \sup_{n \in \Z^{2}} \#\{ (n_1,n_2,n_3) \colon |n_j| \sim N_j, j=1,2,3, n=n_{123}, |\phi-m|\le 1\} \\
& \quad \lesssim  \medcurly{N_1,N_2,N_3}^2 \mincurly{N_1,N_2,N_3}^{\f 32};
\label{eq:cubic-sup-count1}\\
{}& \sup_{m \in \Z} \sup_{n_3 \in \Z^2 } \#\{ (n,n_1,n_2) \colon |n| \sim N_{123}, |n_j| \sim N_j, j=1,2, n=n_{123}, |\phi-m|\le 1\} \\
&\quad \lesssim   \medcurly{N_{123},N_1,N_2}^2 \mincurly{N_{123},N_1,N_2}^{\f 32}.
\label{eq:cubic-sup-count2}
\end{equs}
\end{lemma}

\begin{proof}
 By symmetry, we may assume without loss of generality that $N_{1} \ge N_{2} \ge N_{3}$. 
 For any $n \in \Z^2$, since $n=n_1+n_2+n_3$, the value of $n_1$ is fixed once $n_2,n_3$ are chosen. Thus, 
\begin{align}
{}&\#\{ (n_{1},n_{2},n_{3}) \colon |n_j| \sim N_j, j=1,2,3,\, n=n_{123}, |\varphi-m|\leq 1\}\\
\le  &  \#\{ (n_{2},n_{3}) \colon |n_j| \sim N_j, j=2,3,\, n=n_{123}, \\
& \qquad |\sigma\jpb{n}+ \sigma_1\jpb{n-n_2-n_3} + \sigma_2\jpb{n_2}+\sigma_3\jpb{n_3}-m|\le1 \}\\
\lesssim & \sum_{n_3 \in \Z^{2}} \1_{|n_3 |\sim N_3} \mincurly{\jpb{n-n_3}^{-1},N_{2}}^{-\f 12}N^2_2\\
\lesssim  &\sum_{n_3 \in \Z^{2}} \1_{|n_3 |\sim N_3}\jpb{n-n_3}^{-\f 12}N^2_2 + N^{\f 32}_2 N^2_3 \lesssim  (N_2)^2(N_3)^{\f 32},
\end{align}
where in the second line we applied the basic counting estimate Lemma \ref{lem:basic-count} in $n_2$-variable. This proves the first estimate.
By viewing $n_{123}$ as a free variable ((in the same role as $n_3$ in the first inequality), we obtain the desired result.
\end{proof}

\begin{lemma}
Let $\f 1 4 \le \alpha < \f 5 {12}$. Let $N_1,N_2,N_3 \in 2^{\N}$ and the phase is given by
\begin{equ}
\phi(n_2,n_3,n_4) = \sigma \jpb{n_{234}} + \sigma_2 \jpb{n_2} + \sigma_3 \jpb{n_3} + \sigma_4 \jpb{n_4}.
\end{equ}
 Then, there exists $\epsilon =\epsilon(\alpha) >0$ such that 
\begin{equ}\label{eq:special-cubic-sum}
\sup_{m \in \Z} \sum_{n_2,n_3,n_4 \in \Z^2} \jpb{n_{234}}^{-1} \jpb{n_{34}}^{-1}  \Bigl[\prod_{j=2}^4 \1_{|n_j|\sim N_j}  \jpb{n_{j}}^{-2+2\alpha} \1_{|\phi(\bn) -m|\le 1} \Bigr] \lesssim  \maxcurly{N_2,N_3,N_4}^{-2\epsilon}.
\end{equ}
\end{lemma}
\begin{proof}
Define the set
\begin{equ}
\mathcal{C}(m) = \{(n_2,n_3,n_4) \in (\Z^{2})^{3} \colon |n_j| \sim N_j, j=2,3,4, |n_{234}| \sim N_{234}, |n_{34}| \sim |N_{34}|, |\phi -m| \le 1 \}.
\end{equ}
Then we have 
\begin{equs}
&\sum_{n_2,n_3,n_4 \in \Z^2} \Bigl[ \1_{|n_{234}| \sim N_{234}} \jpb{n_{234}}^{-1} \1_{|n_{34}| \sim N_{34}} \jpb{n_{34}}^{-1} \prod_{j=2}^4 \1_{|n_j|\sim N_j}\jpb{n_{j}}^{-2+2\alpha}   \1_{|\phi -m|\le 1} \Bigr]  \\
\lesssim & \sum_{N_{234},N_{34}} (N_{234})^{-1} (N_{34})^{-1} \prod_{j=2}^4 (N_j)^{-2+2\alpha}  \#\mathcal{C}(m). \label{eq:aux-cubic-sum-for-quntic}
\end{equs}

It is sufficient to prove that the right hand side is bounded by $(N_j)^{-c})$ for some negative number $c$ and some $j=1,2,3,4$. We split the study in two cases: wither $N_3 \sim N_4$ or $N_3 \not \sim N_4$.

Suppose that $N_3 \not \sim N_4$, we assume $N_3 \ge N_4$ by symmetry, which leads to two further subcases: either $N_3 \sim N_4$ or $N_3 \gg N_4$. 

If $N_3 \gg N_4$. We have the following subcases:
$N_2 \sim N_3$ or $N_2 \ll N_3, \, N_2 \gg N_3$.
\begin{item}
\item[Case 1.1] If $N_3 \gg N_4$ and $N_2 \sim N_3$, then $N_{34} \sim N_3$. 
In this case, we must have either $N_3 \gg N_{234} \sim N_4$ or $N_{234} \sim N_3 \gg N_4$. Hence, by applying the cubic counting estimate \eqref{eq:cubic-count2}, we obtain
\begin{equs}
\eqref{eq:aux-cubic-sum-for-quntic} \lesssim & \sum_{N_{234},N_{34}} N^{-1}_{234} N^{-1}_{34} \prod_{j=2}^4 N^{-2+2\alpha}_j  \medcurly{N_{234},N_3,N_4}^{-\f 12} N^2_3N^2_4N^2_{234} \\
\lesssim & \sum_{N_{234},N_{34}}  N_{234} N^{-3+4\alpha}_3 N^{2\alpha}_4 \medcurly{N_{234},N_3,N_4} ^{-\f 12}
 \lesssim  \log(2+N_3) N^{-\f 5 2 +6\alpha}_3.
\end{equs}
If $0<\alpha <\f 5{12}$, for any $0<\epsilon< \f 5 2 - 6\alpha$,  $(N_3)^{-\f 5 2 +6\alpha+\epsilon} $ is summable on dyadic scale, proving the claim.

\item[Case 1.2]  $N_3 \gg N_4$ and $N_2 \ll N_3$. In this case, we have $N_3 \sim N_{234}$. By the cubic counting estimate \eqref{eq:cubic-count1}, we obtain
\begin{equs}
\eqref{eq:aux-cubic-sum-for-quntic} \lesssim & \sum_{N_{234},N_{34}} N^{-1}_3 N^{-1}_{34} \prod_{j=2}^4 N^{-2+2\alpha}_j  \medcurly{N_2,N_3,N_4}^{-\f 12} N^2_2 N^2_3 N^2_4\\
=& \sum_{N_{234},N_{34}}  N^{-2}_3  N^{2\alpha}_2 N^{2\alpha}_3 N^{2\alpha}_4  \medcurly{N_2,N_3,N_4}^{-\f 12}  \lesssim  N^{-\f 5 2 +6\alpha}_3,
\end{equs}
where in the first line we used the fact that $N_3 \sim N_{234}$ and $N_3 \sim N_{34}$. If $0<\alpha <\f 5{12}$, we have the desired estimate.

\item[Case 1.3]  $N_3 \gg N_4$ and $N_2 \gg N_3$. This leads to $N_{234} \sim \maxcurly{N_2,N_3}$. By the cubic counting estimate \eqref{eq:cubic-count1}, we obtain
\begin{equs}
\eqref{eq:aux-cubic-sum-for-quntic} \lesssim &\sum_{N_{234},N_{34}}  N^{-1}_{234} N^{-1}_3 \prod_{j=2}^4 N^{-2+2\alpha}_j   \medcurly{N_{234},N_3,N_4}^{-\f 12} N^2_{234} N^2_3 N^2_4\\
\lesssim &\sum_{N_{234},N_{34}}  N_{234} N^{-2+2\alpha}_2 N^{-\f 32 +2\alpha}_3 N^{2\alpha}_4 \lesssim N^{-1+2\alpha}_2 N^{-\f 32 +4\alpha}_3.
\end{equs}
Hence, if $\alpha < \f 38$, then $\eqref{eq:aux-cubic-sum-for-quntic} \lesssim (N_2)^{-1+2\alpha}$.
If $\alpha \ge \f 38$, then $\eqref{eq:aux-cubic-sum-for-quntic} \lesssim N^{-\f 52+6\alpha}_2$, which in turn requires $\alpha < \f 5{12}$.
\end{item}

Now we consider the case where $N_3 \sim N_4$, it yields two subcases.
\begin{item}
\item[Case 2.1] If $N_3 \sim N_4$ and $N_3 \ll N_2$, then $N_{234} \sim N_2$. By the cubic counting estimate \eqref{eq:cubic-count3}, we obtain
\begin{equs}
\eqref{eq:aux-cubic-sum-for-quntic} \lesssim &\sum_{N_{234},N_{34}}  N^{-1}_2 N^{-1}_{34} \prod_{j=2}^4 N^{-2+2\alpha}_j \mincurly{N_{34},\maxcurly{N_2,N_3}}^{-\f 12} N^2_2 N^2_3 N^2_{34} \\
\lesssim & \sum_{N_{234},N_{34}}  N^{-1+2\alpha}_2 N^{\f 12}_{34} N^{-2+4\alpha}_3 \lesssim \sum_{N_{234},N_{34}} N^{-1+2\alpha}_2 N^{-\f 32 +4\alpha}_3.
\end{equs}
where in the last step we use $N_{34} \lesssim N_3$. What follows is identical to the second subcase of
Case 1.2.
\item[Case 2.2] Suppose that $N_3 \sim N_4$ and $N_3 \gtrsim N_2$. 

If $N_3 \sim N_4, N_3 \gtrsim N_2$ and $N_2 \ll N_{34}$, we must have $N_{234} \sim N_{34}$. 
Assume first $N_{34} \ll N_3$. By  \eqref{eq:cubic-count3}, we have
\begin{equs}
\eqref{eq:aux-cubic-sum-for-quntic} \lesssim &\sum_{N_{234},N_{34}}  N^{-1}_{234} N^{-1}_{34} \prod_{j=2}^4 N^{-2+2\alpha}_j  \mincurly{N_{34},\maxcurly{N_2,N_3}}^{-\f 12} N^2_2 N^2_3 N^2_{34}\\
\lesssim & \sum_{N_{234},N_{34}}  N^{2\alpha}_2  N^{-2+4\alpha}_3 \mincurly{N_{34},N_3}^{-\f 12} 
\lesssim N^{2\alpha-\f 12}_2  N^{-2+4\alpha}_3 \lesssim N^{-\f 5 2+6\alpha}_3.
\end{equs}
We have used $N_2 \ll N_{34}$ in the last step. The remaining case $N_{34} \sim N_3$ is straightforward and yields the same bound.

If $N_3 \sim N_4, N_3 \gtrsim N_2$ and $N_2 \gtrsim N_{34}$, then we require a stronger positive power on
 $N_{234}$, since $$N_{234} \lesssim \maxcurly{N_2,N_{34}} \lesssim N_2.$$
Applying the cubic counting estimate \eqref{eq:cubic-count4}, we obtain
\begin{equs}
\eqref{eq:aux-cubic-sum-for-quntic} \lesssim &\sum_{N_{234},N_{34}}  N^{-1}_{234} N^{-1}_{34} \prod_{j=2}^4 N^{-2+2\alpha}_j  \mincurly{N_{34},\maxcurly{N_3,N_{234}}}^{-\f 12} N^2_3 N^2_{34} N^2_{234}\\
\lesssim & \sum_{N_{234},N_{34}}  N_{234} N_{34} N^{-2+2\alpha}_2  N^{-2+4\alpha}_3 \mincurly{N_{34},N_3}^{-\f 12}  
\end{equs}

If $N_3 \le N_{34}$,  then automatically $N_{34} \sim N_{3}$, and hence
\begin{equ}
\eqref{eq:aux-cubic-sum-for-quntic} \lesssim N^{-\f 12+2\alpha}_2  N^{-2+4\alpha}_3 \lesssim N^{-\f 52 +6\alpha}_3.
\end{equ}
If $N_{34} \le N_3$, it suffices to consider $N_{34} \ll N_3$. Since $N_2 \gtrsim N_{34}$, the same bound follows, which completes the proof.
\end{item}
\end{proof}

In the lemma below we follow the notation that $n_{12}=n_1+n_2$ and $n_{123}=n_{12}+n_3$.

\begin{lemma}[Basic resonance estimate.]\label{lem:basic-resonant} 
Let $n_{1},n_{2} \in \Z^{2}$ be fixed and $N_{3}  \in 2^{\N}$. Then, for all $0< \alpha <\frac{1}{2}$, it holds that
\begin{equ}\label{eq:basic-resonant}
\sum_{m \in \Z} \f 1 {\jpb{m}}  \sum_{n_3 \in \Z^2} \1_{|n_3| \sim N_3} \jpb{n_{123}}^{-1} \jpb{n_3}^{-2+2\alpha} \1_{ |\phi-m|\leq 1}  \lesssim  \log(2+N_3)\jpb{n_{12}}^{-\tilde{s}_{\alpha}}
\end{equ}
where
\begin{equation}\label{eq:tilde-s-alpha}
\tilde{s}_{\alpha} = 
\begin{cases}
1, & 0 < \alpha <\frac{1}{4}\\
\frac{3}{2}-2\alpha, & \frac{1}{4} \le \alpha < \frac{1}{2} 
\end{cases}.
\end{equation}
\end{lemma}

\begin{proof}
Since $n_{1},n_{2}\in \Z^{2}$ are fixed and the phase function $\phi$ is globally Lipschitz, there is at most $\sim N_3$ number of $m$ which gives $\f 1 {\jpb{m}} \lesssim \log(2+N_3)$,
 it suffices to show that 
\begin{align}
\sup_{m \in \Z} \sum_{n_{3} \in \Z^{2}} \mathbf{1}_{|n_3|\sim N_3} \jpb{n_{123}}^{-1} \jpb{n_3}^{-2+2\alpha} \1_{|\phi-m|\leq 1} \lesssim  \jpb{n_{12}}^{-\tilde{s}_\alpha}.
\end{align}
Decomposing the  frequency $n_{123}$ dyadically, we have
\begin{align}
& \sum_{n_{3} \in \Z^{2}} \1_{|n_3|\sim N_3} \jpb{n_{123}}^{-1} \jpb{n_3}^{-2+2\alpha} \1_{|\phi-m|\leq 1}\\
\lesssim & N^{-2+2\alpha}_{3} \sum_{N_{123} \ge 1} N^{-1}_{123} \sum_{n_3 \in \Z^2}\1_{|n_3| \sim N_3} \1_{|\phi-m|\leq 1}
\end{align}
	
	 Assuming that  $|n_{12}|\sim N_{12}$, and counting the number of $n_3$, the right-hand side of the above inequality can be controlled by the counting, 
\begin{equation}
\lesssim \sum_{N_{123} \geq 1} N^{-2+2\alpha}_{3}  N^{-1}_{123} \mincurly{N_{123},N_{12},N_{3}}^{-\frac{1}{2}}\mincurly{N_{123},N_{3}}^{2},
\end{equation}
We have applied Lemma \ref{lem:two-ball-count}, by setting $a=n_{12}$ and $n=n_3$.

	We analyze case by case.
	\begin{itemize}
		\item If $N_{123} \ll N_{3}$, then $N_{123} \sim N_{12}$. Consequently,		\begin{align}
			& \sum_{N_{123} \ll N_{3}}N^{-2+2\alpha}_{3}  N^{-1}_{123} \mincurly{N_{123},N_{12},N_{3}}^{-\frac{1}{2}}\mincurly{N_{123},N_{3}}^{2}\\
			\lesssim  &\sum_{N_{123} \ll N_{3}} N^{\frac{1}{2}}_{123}  N^{-2+2\alpha}_{3} \lesssim N^{-\frac{3}{2}+2\alpha}_{3} \sim N^{-\frac{3}{2}+2\alpha}_{12}.
		\end{align}
		\item If $N_{123} \sim N_{3}$, then $N_{12} \lesssim N_{123} \sim N_{3}$ and we have
		\begin{align}
			& \sum_{N_{123} \sim N_{3}} N^{-2+2\alpha}_{3}  N^{-1}_{123} \mincurly{N_{123},N_{12},N_{3}}^{-\frac{1}{2}}\mincurly{N_{123},N_{3}}^{2}\\
			\lesssim  &\sum_{N_{123} \sim N_{3}} N^{-1+2\alpha}_{3}N^{-\frac{1}{2}}_{12} \lesssim N^{-\frac{3}{2}+2\alpha}_{12},
		\end{align}
	when $\alpha < \f 12$.	
	
		\item If $N_{123} \gg N_{3}$, we have $N_{123} \sim N_{12}$ and it follows that
		\begin{align}
			& \sum_{N_{123} \gg N_{3}}N^{-2+2\alpha}_{3}  N^{-1}_{123} \mincurly{N_{123},N_{12},N_{3}}^{-\frac{1}{2}}\mincurly{N_{123},N_{3}}^{2}\\
			\lesssim  &\sum_{N_{123} \gg N_{3}} N^{-1}_{123}  N^{-\frac{1}{2}+2\alpha}_{3} \lesssim 
			\begin{cases}
				N^{-1}_{12}, & 0 < \alpha <\frac{1}{4}\\
				N^{-\frac{3}{2}+2\alpha}_{12}, &  \frac{1}{4} \leq \alpha < \frac{1}{2} 
			\end{cases}.
		\end{align}
	\end{itemize}
	This concludes the proof.
	\end{proof}

\begin{lemma}[Quartic non-resonance  sum estimate.]\label{lem:quartic-non-resonance}
Let $N_1,N_2,N_3$ be dyadic numbers and $\sigma, \sigma_1, \sigma_2,\sigma_3 \in \{\pm 1\}$ be given.
Let $\phi(n_1,n_2,n_3)$ be as defined by \eqref{eq:cubic-phase}. Then, for all $-1< s<-\alpha$, there exists $\epsilon>0$ which may depend on $\alpha,s$, such that that 
\begin{align}
{} & \sup_{m \in \Z}\sum_{n_1,n_2,n_3,n_4  \in \Z^2}  \Bigl( \prod_{j=1}^4 \1_{|n_j| \sim N_j} \cdot \jpb{n_j}^{-2+2\alpha}\Bigr) \jpb{n_{1234}}^{2s} \jpb{n_{123}}^{-2} \1_{|\phi-m| \le 1}\\
\lesssim & N^{-2\epsilon}_4 \maxcurly{N_1,N_2,N_3}^{-2s_\alpha}.
\label{eq:quartic-non-resonant}
\end{align}
\end{lemma}
\begin{proof}
 Since $ \jpb{n_4}\sim N_4$ we may extract the factor $(N_4)^{-2\epsilon}$ before summing over $n_4$. More precisely,
\begin{align}
{} & \sum_{n_1,n_2,n_3,n_4  \in \Z^2}  \Bigl( \prod_{j=1}^4 \mathbf{1}_{|n_j| \sim N_j} \cdot \jpb{n_j}^{-2+2\alpha} \Bigr) \jpb{n_{1234}}^{2s} \jpb{n_{123}}^{-2}  \1_{|\phi-m| \le 1}\\
\lesssim &  \sum_{n_1,n_2,n_3  \in \Z^2} \Bigl(\prod_{j=1}^3 \mathbf{1}_{|n_j| \sim N_j}  \jpb{n_j}^{-2+2\alpha}\Bigr)  \jpb{n_{123}}^{-2}  \cdot  \1_{|\phi-m| \le 1} N^{-2\epsilon}_4 \Bigl( \sum_{n_4 \in \Z^2} \jpb{n_{1234}}^{2s} \jpb{n_4}^{-2+2\alpha+2\epsilon} \Bigr)\\
\lesssim & N^{-2\epsilon}_4 \sum_{n_1,n_2,n_3  \in \Z^2} \Bigl( \prod_{j=1}^3 \mathbf{1}_{|n_j| \sim N_j}  \jpb{n_j}^{-2+2\alpha}\Bigr)  \jpb{n_{123}}^{-2} \cdot  \jpb{n_{123}}^{2s+2\alpha+2\epsilon}  \1_{|\phi-m| \le 1}\\
\le & (N_4)^{-2\epsilon} \sum_{n_1,n_2,n_3  \in \Z^2} \prod_{j=1}^3 \mathbf{1}_{|n_j| \sim N_j}  \jpb{n_j}^{-2+2\alpha} \jpb{n_{123}}^{-2}   \1_{|\phi-m| \le 1}.
\end{align}
The penultimate line follows from discrete convolution inequality from Lemma \ref{lem:convolution-sum}.
The last line follows by choosing $\epsilon>0$ sufficiently small so that $2s+2\alpha+2\epsilon<0$.
The cubic sum estimate from Proposition \ref{prop:cubic-sum}, allow us to complete the proof.
\end{proof}

For the Quartic non-resonance sum estimate,  we need the following lemma, \cite[Lemma A.1.] {Bri20b}:
\begin{lemma}\label{dyadic-reduced}
Let $N_1,N_5,N_{1234},N_{12345}$ be dyadic frequency scales such that there are frequencies $n_1, n_5, m_{1234}$, and $n_{12345}$ which realizes the respective magnitude. We denote this by the following expression:
\begin{equ}
\1_{|n_1| \sim N_1}\cdot  \1_{|n_5| \sim N_5}\cdot \1_{|n_{12345}| \sim N_{12345}} \cdot \1_{|n_{1234}| \sim N_{1234}}  \ne 0. 
\end{equ}
Then, it holds that  
\begin{equ}
\frac{\mincurly{N_5,N_{12345}}  \mincurly{N_{1234},N_1}^{\f 12}}{\mincurly{N_{12345},N_{1234},N_5}^{\frac{1}{2}}}  \le N^{\f 12}_5 N^{\f 1 2}_{12345}.
\end{equ}
\end{lemma}

In the lemma, $n_{12345}$ denotes $n_1+n_2+n_3+n_4+n_5$. Similarly, $n_{234}$, $n_{1234}$, $n_{1345}$ are defined.

\begin{lemma}[Quintic non-resonance sum estimate] \label{lem:quintic-non-resonance}
Let $N_1,N_2,N_3,N_4,N_5\in 2^{\N}$, $\sigma_{12345},\sigma_{234}$ and $\sigma_1,\sigma_2,\sigma_3, \sigma_4,\sigma_5 \in \{\pm\}$ be given. Define:\begin{align}
& \phi(\bn)= \sigma_{234} \jpb{n_{234}} + \sum_{j=2}^4 \jpb{n_j},\\
& \widetilde{\phi}(\bn) = \sigma_{12345}\jpb{n_{12345}} + \sigma_{234} \jpb{n_{234}} + \sigma_1 \jpb{n_1} + \sigma_5 \jpb{n_5} ,\\
& \widetilde{\widetilde{\phi}}(\bn) = \sigma_{12345}\jpb{n_{12345}} + \sigma_{234} \jpb{n_{234}} + \sum_{j=1}^5 \jpb{n_j}.
\end{align}
Suppose that $\alpha , s$ is a pairs of numbers satisfying one of the following constraints:
\begin{equs}
&\alpha \in (0,\frac{1}{4}),  \quad s<\frac{3}{4}\\
&\alpha \in (\frac{1}{4} ,\frac{5}{12}), \quad  s<\frac{5}{6}-\alpha.
\end{equs}
 Then, there exists $\epsilon =\epsilon(\alpha,s)> 0$ such that
\begin{align}\label{eq:quintic-non-resonance}
&\sup_{m, m' \in \Z} \sum_{n_1,\cdots, n_5 \in \Z^{2}} \jpb{n_{12345}}^{2(s-1)} \jpb{n_{234}}^{-2}  \prod_{j=1}^{5} \1_{|n_j| \sim N_j} \jpb{n_{j}}^{-2+2\alpha} \1_{|\phi-m|\leq 1} ( \1_{|\tilde{\varphi}-m'|\leq 1} + \1_{|\widetilde{\widetilde{\phi}}-m'|\leq 1}) \\
& \lesssim \maxcurly{N_1, N_2, N_3, N_4, N_5}^{-\epsilon}.
\end{align}
\end{lemma}
\begin{proof}
We fix $m,m' \in \Z$. Note that
\begin{align*} \widetilde{\phi}(\bn)-m'&= \sigma_{12345}\jpb{n_{12345}} + \sigma_5 \jpb{n_5}
 -(m'-\sigma_{234}\jpb{n_{234}} -\sigma_1 \jpb{n_1})).
\end{align*}
We introduce $\tilde{m}=-\sigma_{234} \jpb{n_{234}} - \sigma_1 \jpb{n_1} +m'$. 
 We also introduce the dyadic decomposition for frequencies $n_{12345}$ and $n_{1234}$ and work with the following sum over $n_5$:
\begin{align*}
\sum_{N_{12345}, N_{1234}}
\sum_{n_5 \in \Z^{2}} \1_{|n_{12345}| \sim N_{12345}} \1_{|n_{1234}| \sim N_{1234}} \jpb{n_{12345}}^{2(s-1)} \1_{|n_5| \sim N_5} \jpb{n_{5}}^{-2+2\alpha}   \1_{|\tilde{\varphi}-m'|\leq 1} 
\end{align*}
Note that on each dyadic piece,  the size of $n_5$, $n_{12345}$, and $n_{1234}$ are determined.
\begin{align}
&\sum_{N_{12345}, N_{1234}}\sum_{n_5 \in \Z^{2}} \1_{|n_{12345}| \sim N_{12345}} \1_{|n_{1234}| \sim N_{1234}} \jpb{n_{12345}}^{2(s-1)} \1_{|n_5| \sim N_5} \jpb{n_{5}}^{-2+2\alpha} 
  \1_{|\tilde{\varphi}-m'|\leq 1} \\
&=\sum_{N_{12345}, N_{1234}} N_{12345}^{2(s-1)}  N_5^{-2+2\alpha} 
\sum_{n_5 \in \Z^{2}} \1_{|n_{12345}|\sim N_{12345}} \1_{|n_{1234}| \sim N_{1234}} \1_{|n_5| \sim N_5}   \1_{|\tilde{\varphi}-m'|\leq 1}.
\end{align}
These turns the estimate into counting. We set $a=n_{1234}$, $A=N_{1234}$, $n=n_5$, $N=N_5$, $B=N_{12345}$, 
and then sum over $n_{5} \in \Z^2$. With these notations, we apply the two-ball counting (\ref{eq:two-ball-count}) to the above expression to obtain:
\begin{align}
&\sum_{n_5 \in \Z^{2}} \1_{|n_{12345}|\sim N_{12345}} \1_{|n_{1234}| \sim N_{1234}} \1_{|n_5| \sim N_5}   \1_{|\tilde{\varphi}-m'|\leq 1} \\
&\lesssim
  \mincurly{N_{12345},N_{1234},N_5}^{-1/2} \mincurly{N_5,N_{12345}}^{2}. 
\end{align}
We have used the fact that $|\tilde{\varphi}-m'|=|\jpb{a+n}\pm \jpb{n}-m'|$. To summarize:
\begin{align}
&\sum_{N_{12345}, N_{1234}}\sum_{n_5 \in \Z^{2}} \1_{|n_{12345}| \sim N_{12345}} \1_{|n_{1234}| \sim N_{1234}} \jpb{n_{12345}}^{2(s-1)} \1_{|n_5| \sim N_5} \jpb{n_{5}}^{-2+2\alpha} 
  \1_{|\tilde{\varphi}-m'|\leq 1} \\
&=\sum_{N_{12345}, N_{1234}} N_{12345}^{2(s-1)}  N_5^{-2+2\alpha} 
 \mincurly{N_{12345},N_{1234},N_5}^{-1/2} \mincurly{N_5,N_{12345}}^{2}.
\end{align}

The term with $ \1_{|\tilde{\tilde{\varphi}}-m'|\leq 1}$ can be proved similarly be rewriting:
\begin{align*} 
 \widetilde{\widetilde{\phi}}(\bn) -m'= \sigma_{12345}\jpb{n_{12345}} +\sigma_5\jpb{n_5}
 -(m'- \sum_{j=1}^4 \jpb{n_j}- \sigma_{234} \jpb{n_{234}}).
\end{align*}
For the remaining factor we use a crude bound on each dyadic piece:
\begin{align}
&\sum_{n_1, \cdots, n_4 \in \Z^2} \1_{|n_{1234}| \sim N_{1234}} \jpb{n_{234}}^{-2}
 \prod_{j=1}^{4} \1_{|n_j| \sim N_j} \jpb{n_{j}}^{-2+2\alpha} 
  \1_{|\phi-m|\leq 1} \\
 &\lesssim
 \sum_{n_2, n_3, n_4 \in \Z^2}\jpb{n_{234}}^{-2}
 \prod_{j=2}^{4} \1_{|n_j| \sim N_j} \jpb{n_{j}}^{-2+2\alpha}  \1_{|\phi-m|\leq 1} \sum_{n_1 \in \Z^2} (N_1)^{-2+2\alpha} \1_{|n_1| \sim N_1} 
  \1_{|n_{1234}| \sim N_{1234}} \\
 &\lesssim (N_1)^{-2+2\alpha} 
 \sum_{n_2, n_3, n_4 \in \Z^2}\jpb{n_{234}}^{-2}
 \prod_{j=2}^{4} \1_{|n_j| \sim N_j} \jpb{n_{j}}^{-2+2\alpha}  \1_{|\phi-m|\leq 1}   \min\{N_1, N_{1234}\}^2.
  \end{align}
  Puttign all estimates together we have:
  \begin{align}
 & \sum_{n_1,\cdots, n_5 \in \Z^{2}} \jpb{n_{12345}}^{2(s-1)} \jpb{n_{234}}^{-2}  \prod_{j=1}^{5} \1_{|n_j| \sim N_j} \jpb{n_{j}}^{-2+2\alpha} \1_{|\phi-m|\leq 1} ( \1_{|\tilde{\varphi}-m'|\leq 1} + \1_{|\widetilde{\widetilde{\phi}}-m'|\leq 1})\\
\lesssim & \sum_{N_{12345},N_{1234}} N^{2s-2}_{12345}  \mincurly{N_{12345},N_{1234},N_5}^{-1/2} \mincurly{N_5,N_{12345}}^{2}  N^{-2+2\alpha}_5 \\
&\quad  \times N^{-2+2\alpha}_1 N^{-2+2\alpha}_5 \sum_{n_2,n_3,n_4 \in \Z^2} \prod_{j=2}^4 \1_{|n_j| \sim N_j} \jpb{n_j}^{-2+2\alpha} \cdot \jpb{n_{234}}^{-2} \1_{|\phi-m|\leq 1}. \label{eq:rhs-quintic-non-resonance}
\end{align}
Finally, applying the cubic sum estimate \eqref{prop:cubic-sum} with $s=0$ to the expression begins with $\prod_{j=2}^4 \1_{|n_j| \sim N_j} \jpb{n_j}^{-2+2\alpha} $, we have 
\begin{align}
\eqref{eq:rhs-quintic-non-resonance}
&\lesssim \sum_{N_{12345},N_{1234}} N^{2s-2}_{12345} (N_1N_5)^{-2+2\alpha} \frac{\mincurly{N_5,N_{12345}}^{2} \mincurly{N_{1234},N_1}^2}{\mincurly{N_{12345},N_{1234},N_5}^{\frac{1}{2}}} \maxcurly{N_2,N_3,N_4}^{-2s_\alpha}\\
&\lesssim \sum_{N_{12345},N_{1234}} N^{2s-\frac{3}{2}}_{12345}N^{-2+2\alpha}_1 N^{-\frac{3}{2}+2\alpha}_5\mincurly{N_5,N_{12345}} \mincurly{N_{1234},N_1}^{\frac{3}{2}}\maxcurly{N_2,N_3,N_4}^{-2s_\alpha},
\label{eq:rhs-quintic-non-resonance-2}
\end{align}
where in the last step we used Lemma \ref{dyadic-reduced}. 

Next, we prove that the right hand side is bounded with the required decay rate 
$$\maxcurly{N_1, N_2, N_3, N_4, N_5}^{-\epsilon}.$$
 Note that the summation is over $N_{12345}$ and $N_{1234}$, and we observe that $N_{1234} \lesssim \maxcurly{N_5,N_{12345}}$. We
 consider the following cases (1) $N_{12345} \le N_5$ and (2) $N_{12345} > N_5$.

(1) First we suppose that $N_{12345}\leq N_5$, then $N_{1234}\le N_5$ and for any $\theta\in [0,1]$, we have
\begin{align}
\hbox{ RHS of  } \eqref{eq:rhs-quintic-non-resonance-2}
\le  & \sum_{  N_{12345}, N_{1234} } N^{2s-\frac{3}{2}}_{12345} N^{-2+2\alpha}_1 N^{-\frac{3}{2}+2\alpha}_5 N_{12345} N^{\frac{3}{2}\theta}_1 N^{\frac{3}{2}(1-\theta)}_{1234}\maxcurly{N_2,N_3,N_4}^{-2s_\alpha}\\
&\lesssim \sum_{N_{12345}} (N_{12345})^{2s-\frac{1}{2}} (N_1)^{-2+2\alpha+\frac{3}{2}\theta}( N_5)^{-\frac{3}{2}+2\alpha+\frac{3}{2}(1-\theta)} \maxcurly{N_2,N_3,N_4}^{-2s_\alpha}.
\end{align} 
This follows from the fact that  $N_{1234} \lesssim N_5$ and so the geometric sum below has the indicated bound:
$$\sum_{N_{1234}} (N_{1234})^{\frac{3}{2}(1-\theta)}\lesssim (N_5)^{\frac{3}{2}(1-\theta)}.$$
Below we choose $\theta$ to obtain the desired control.  The aim is to produce negative exponents of $N_i$. 
For this, we bound the unwanted terms by $1$ if it has negative exponent, otherwise by a larger frequencies. 
In more detail, we consider  two further subcases. 

(1a)  Suppose that $s < \f 1 4$, we can throw away the term involving $N_{12345}$. We further require the exponents of $N_1$ and $N_5$ to be negative, i.e.
\begin{equ}
-2+2\alpha +\frac{3}{2}\theta <0, \quad -\frac{3}{2}+2\alpha+\frac{3}{2}(1-\theta) =2\alpha -\f 32 \theta <0.
\end{equ}
Equivalently,
\begin{equ}\label{choose-theta-1}
\frac{4}{3}\alpha <\theta <\frac{2}{3}(2-2\alpha)
\end{equ}
Let us take $\theta =\frac{2}{3}$. Since $0<\alpha<\f 5 {12}$, the two exponents listed above are indeed. \begin{equs}
\hbox{ RHS of  } \eqref{eq:rhs-quintic-non-resonance-2} 
&\lesssim N^{2\alpha-1}_1 N^{2\alpha-1}_5 \maxcurly{N_2,N_3,N_4}^{-2s_\alpha}  \\
&\lesssim \maxcurly{N_1,N_2,N_3,N_4,N_5}^{-\epsilon},
\end{equs}
where $\epsilon = \mincurly{1-2\alpha, 2s_\alpha}>0$. In this case, the expression is summable on the dyadic scales.

(1b)  Suppose that $s\ge \f 1 4$. In addition to (\ref{choose-theta-1}), since $N_{12345}<N_5$, we request a third constraint, the sum of their exponent is negative:
\begin{equ}
2s -\frac{1}{2}  -\frac{3}{2}+2\alpha+\frac{3}{2}(1-\theta) < 0.
\end{equ}
The overall constraint is:
\begin{equ}
\frac{4}{3}\alpha <\theta <\frac{2}{3}(2-2\alpha) , \quad \mbox{and} \quad \frac{1}{4}  \le s <\frac{1}{4}-\alpha+\frac{3}{4}\theta.
\end{equ}
If $0 < \alpha < \frac{1}{4}$,  we set $\theta = 1$, and assume that $s\in [\f 14, 1-\alpha)$. Observe that  $[\f 14, 1-\alpha)\subset [\f 14, \f 34]$. Then, 
\begin{equ}
\hbox{ RHS of  } \eqref{eq:rhs-quintic-non-resonance-2} \lesssim N^{-\frac{1}{2}+2\alpha}_1N^{2s-2+2\alpha}_5\maxcurly{N_2,N_3,N_4}^{-2s_\alpha}.
\end{equ}
Recall that $s_\alpha =1-2\alpha>\f 12-2\alpha$, we set $\epsilon = \mincurly{\f 12 - 2\alpha, 2-2s-2\alpha}$ to conclude.

If $\frac{1}{4} \le \alpha <\frac{5}{12}$,  we let $\theta = \frac{7}{9}$. Then,  we require that $s<\frac{1}{4}-\alpha+\frac{3}{4}\theta=\f 56 -\alpha$. WE also have,
\begin{equ}
\hbox{ RHS of  } \eqref{eq:rhs-quintic-non-resonance-2} \lesssim N^{-\f 5 6 +2\alpha}_1N^{2s-\f 5 6 +2\alpha}_5\maxcurly{N_2,N_3,N_4}^{-2s_\alpha}.
\end{equ}
Recall that $s_\alpha=\f 54-3\alpha>\f 56-2s-2\alpha$ when $ \f14 \le \alpha<\f 5{12}$. Set $\epsilon = \mincurly{\f 5 6 - 2\alpha, \f 56-2s-2\alpha}$ to conclude. 

(ii) Next we consider the case  $ N_5 \le N_{12345}$. In \eqref{eq:rhs-quintic-non-resonance-2} we replace $\mincurly{N_5,N_{12345}} $ with $N_5$ and obtain:
\begin{align}
\hbox{ RHS of  } \eqref{eq:rhs-quintic-non-resonance-2}
\le & \sum_{N_{12345}, N_{1234}} N^{2s-\frac{3}{2}}_{12345} N^{-2+2\alpha}_1 N^{-\frac{3}{2}+2\alpha}_5 N_5 N^{\frac{3}{2}\theta}_1 N^{\frac{3}{2}(1-\theta)}_{1234} \maxcurly{N_2,N_3,N_4}^{-2s_\alpha}\\
\lesssim &\sum_{N_{12345}  } (N_{12345})^{2s-\frac{3}{2}+\frac{3}{2}(1-\theta)} (N_1)^{-2+2\alpha+\frac{3}{2}\theta}
 (N_5)^{-\frac{1}{2}+2\alpha} \maxcurly{N_2,N_3,N_4}^{-2s_\alpha}
\end{align} 
where $\theta\in [0,1]$.

 (iia) Let $0<\alpha < \f 1 4$. We require the exponents of $N_{12345}$ and $N_1$ to be positive:
 \begin{equ}
-2+2\alpha +\frac{3}{2}\theta <0, \quad -\frac{3}{2}+2s+\frac{3}{2}(1-\theta) < 0.
\end{equ}
Equivalently,
\begin{equ}\label{choose-theta-2}
\frac{4}{3}s <\theta <\frac{4}{3}(1-\alpha).
\end{equ}
Take $\theta =1$, then $s<\f 34$. We have
\begin{equ}
\hbox{ RHS of  } \eqref{eq:rhs-quintic-non-resonance-2} \lesssim N^{-\frac{1}{2}+2\alpha}_1N^{-\f 12+2\alpha}_5 \maxcurly{N_2,N_3,N_4}^{-2s_\alpha}.
\end{equ} 
Let
 $\epsilon=\min{\{ \f 12-2\alpha, 2s_\alpha \}}$, we conclude the proof.

(iib) If  $\alpha \ge  \f 1 4$, we require not only \eqref{choose-theta-2}, but also the sum of the exponents of $N_5$ and $N_{12345}$:
\begin{equ}
2s-\frac{3}{2}+\frac{3}{2}(1-\theta) -\frac{1}{2}+2\alpha< 0.
\end{equ}
Since $s<\f 56 -\alpha$, the inequalities are satisfied by  $\theta=\f 7 9$. Letting $\theta=\f 79$, we have
\begin{equ}
\hbox{ RHS of  } \eqref{eq:rhs-quintic-non-resonance-2} \lesssim N^{-\f 5 6+2\alpha}_1 \maxcurly{N_2,N_3,N_4}^{-2s_\alpha}.
\end{equ} 
Set $\epsilon = \mincurly{\f 5 6-2\alpha,2s_\alpha}$, also using $N_5 \le N_{12345} \lesssim \maxcurly{N_1,N_2,N_3,N_4}$, we complete the proof.
\end{proof}

In the lemma below,  the phase function is
\begin{equation}
\phi\equiv \phi(n_i,n_j,n_k) := \sigma_{ijk} \jpb{n_{ijk}}+\sigma_i \jpb{n_i}+\sigma_j\jpb{n_j}+\sigma_k \jpb{n_k},
\end{equation}
where $\sigma_{ijk}, \sigma_i, \sigma_j, \sigma_k \in \{ \pm 1\}$  are given.

\medskip

 Before stating the next estimate, we recall the notation concerning pairings. Let $\cP\subset \{1, \cdots, 7\}^2$ be a pairing ( c.f. Definition \ref{pairing}) we denote by $\pi\subset \{1,2,3,4,5, 6,7\}$ the collection of paired numbers. Furthermore let $\pi^c=\{1,\cdots, 7\}\setminus \pi$ denotes  its complement. 

\begin{lemma}[The Septic sum estimate]
\label{lem:septic-sum}
Let $N_{1234567}, N_1,N_2,\cdots,N_7 \in 2^{\N}$. Define
\begin{align}
\Psi(n_i,n_j,n_k) = \sum_{\substack{\sigma_{ijk},\sigma_i,\sigma_j \\ \sigma_k \in \{\pm 1\} }} \sum_{m \in \Z} \jpb{m}^{-1} \jpb{n_{ijk}}^{-1}  \Bigl(\jpb{n_i}\jpb{n_j}\jpb{n_k}\Bigr)^{-1+\alpha} \1_{|\phi-m| \leq 1}.
\end{align}
Let $\cP$ be a partition of $\{1,\cdots ,7\}$ which respects the partition $Q:=\{\{1,2,3\}, \{4\}, \{5,6,7\}\}$ (this includes the empty partition).Then, for each $ \f 1 4 \le \alpha < \f 3 8$ and given $0 < \kappa \ll \f1{12}$, for any $s < 1-\alpha-\kappa$ and , there exist $\epsilon=\epsilon(s,\alpha,\kappa)>0$ such that
\begin{equs}
&\sum_{n_j \in \Z^2: j \in \pi^c}
 \jpb{  {\textstyle \sum}_{j\in \pi^c }n_j }^{2(s-1)} \cdot \\
&\Bigl(\sum_{n_i : (i,j)\in \cP, n_j=-n_i}
\mathbf{1}_{|n_{1234567}| \sim N_{1234567}} \prod_{j=1}^7 \1_{|n_j| \sim N_j}  
\Psi(n_1,n_2,n_3) \times \frac{1}{\jpb{n_4}^{1-\alpha}}\Psi(n_5,n_6,n_7) \Bigr)^2\\
 &\lesssim \maxcurly{N_{1234567},N_1,N_2,N_3,N_4,N_5,N_6,N_7}^{-\epsilon}.\label{eq:septic-sum}
\end{equs}
The proportionality constant is independent of the partition.
\end{lemma}
\begin{proof}
We note that the function $$\tilde \Psi=\Psi(n_1,n_2,n_3) \times \frac{1}{\jpb{n_4}^{1-\alpha}}\Psi(n_5,n_6,n_7).$$
is invariant under the permutations of $\{n_1,n_2, n_3\}$ and invariant under permutations of $\{5,6,7\}$. 
We also note that 
$$I=\mathbf{1}_{|n_{1234567}| \sim N_{1234567}}\prod_{j=1}^7 \1_{|n_j| \sim N_j} 
$$
is also invariant under the permutations within the partition. 

Due to the invariance of the permutation between the two groups $\{1,2,3\}$ and $\{5,6,7\}$, and also the invariance within the partition in the permutation of the index, the problem reduces to the following cases:

Case 1. The index $4$ is not paired, then we may have the following $4$ distinct cases:\\
(1a)  $\cP=\emptyset$;\\ 
(1b) $\cP=\{(1,7)\}$;\\
(1c) $\cP=\{(1,7), (2,6)\}$;\\
(1d) $\cP=\{(1,7), (2,6), (3,5)\}$.

In all four cases, the main idea is first sum over $n_4$ and then sum over the remaining indices.
We also use the notation  $n_{\pi^c}:=\sum_{j\in \pi^c}n_j$.

\begin{equs}
{}&\sum_{n_j, j \in\pi^c} \jpb{n_{\pi^c}}^{2(s-1)} \Bigl( \sum_{n_i: (i,j)\in \cP, n_j=-n_i} \mathbf{1}_{|n_{1234567}| \sim N_{1234567}}\prod_{j=1}^7 \1_{|n_j| \sim N_j} \Psi(n_1,n_2,n_3) \\
& \quad \times \frac{1}{\jpb{n_4}^{1-\alpha}}\Psi(n_5,n_6,n_7) \Bigr)^2 \\
= &\sum_{n_j, j \in\pi^c \backslash\{4\}}  \sum_{n_4} \mathbf{1}_{|n_{\pi^c}| \sim N_{1234567}} \1_{|n_4|\sim N_4}\jpb{n_{\pi^c}}^{2(s-1)}  \jpb{n_4}^{-2+2\alpha} \\
& \quad  \times \Bigl( \sum_{n_j: (i,j)\in \cP, n_j=-n_i} \prod_{1,2,3,5,6,7} \1_{|n_j| \sim N_j} \Psi(n_1,n_2,n_3) \Psi(n_5,n_6,n_7) \Bigr)^2
\end{equs}
Here $N_{1234567}$ denotes the corresponding reduced sum, after pairing have been taken into account.  For example if $(1,7)$ are paired, $N_{1234567}=N_{23456}$.
Note that by the discrete convolution inequality, we obtain:
\begin{equs}
{}&\sum_{n_4} \mathbf{1}_{|n_{\pi^c}| \sim N_{1234567}} \1_{|n_4|\sim N_4}\jpb{n_{\pi^c}}^{2(s-1)}  \jpb{n_4}^{-2+2\alpha} \\
\lesssim & N^{-2\kappa}_4 \sum_{n_4} \jpb{n_{\pi^c}}^{2(s-1)} \jpb{n_4}^{-2+2\alpha+2\kappa}  \1_{|n_{\pi^c \backslash \{4\}}| \sim N_{1234567}} \lesssim N^{-2\kappa}_4 |n_{\pi^c \backslash \{4\}}|^{2s-2+2\alpha+2\kappa}
\end{equs}
where  $0<s<1-\alpha-\kappa$. 

Observing that $N_{1234567} \sim |n_{\pi^c \backslash \{4\}}| \le |n_{\pi^c}| + |n_4| \lesssim \max(|n_{\pi^c \backslash \{4\}}|,N_4)$, we have the following two cases: 

\begin{itemize}
\item \textbf{Case I}. \, Suppose that $|n_{\pi^c \backslash \{4\}}| \le N_4$, we have $N_{1234567} \lesssim N_4$. Since $0<s<1-\alpha-\kappa$, 
then $|n_{2356}|^{2s-2+2\alpha+2\kappa} $, with negative exponent,  is bounded by $1$.
Similarly, $$|n_{\pi^c \backslash \{4\}}|^{2s-2+2\alpha+2\kappa} \lesssim 1 
\lesssim N^{-\kappa}_{1234567}N^\kappa_4.$$
\item 
\textbf{Case II}. \,  If on the other hand $|n_{\pi^c \backslash \{4\}}| \ge N_4$, we have $N_{1234567} \lesssim |n_{\pi^c \backslash \{4\}}|$.  Since the exponent $2s-2+2\alpha+2\kappa$ is negative,   $|n_{\pi^c \backslash \{4\}}|^{2s-2+2\alpha+2\kappa}\lesssim N_{1234567}^{2s-2+2\alpha+2\kappa}$.

\end{itemize}
\medskip

Put the two cases  together, we have the bound
\begin{equs}
{}&\sum_{n_4} \mathbf{1}_{|n_{\pi^c}| \sim N_{1234567}} \1_{|n_4|\sim N_4}\jpb{n_{\pi^c}}^{2(s-1)}  \jpb{n_4}^{-2+2\alpha} \lesssim (N_{1234567}^{-\kappa} \vee N_{1234567}^{2s-2+2\alpha+2\kappa}) N^{-\kappa}_4.
\end{equs}

To obtain a suitable bound to the remaining  sum,   we apply the Cauchy-Schwarz and relabel the indices and study case by case. We show shortly that,  in each case,  we have:
 \begin{equs}
 {}&\sum_{n_j, j \in\pi^c \backslash\{4\}}  \Bigl( \sum_{n_i: (i,j)\in \cP, n_j=-n_i} \prod_{1,2,3,5,6,7} \1_{|n_j| \sim N_j} \Psi(n_1,n_2,n_3) \Psi(n_5,n_6,n_7) \Bigr)^2\\
&\lesssim  \Bigl( \sum_{n_1,n_2,n_3}\prod_{j=1}^3  \1_{|n_j| \sim N_j} \Psi^2(n_1,n_2,n_3) \Bigr) \Bigl( \sum_{n_5,n_6,n_7}\prod_{j=5}^7 \1_{|n_j| \sim N_j} \Psi^2(n_5,n_6,n_7) \Bigr)
\\& \lesssim \max(N_1,N_2,N_3)^{-2s_\alpha+\kappa}\max(N_5,N_6,N_7)^{-2s_\alpha+\kappa},
 \end{equs}
 which together with the bound for the other factor allows us to conclude.  The estimates for the two factors are completely identical, below we bound the term $\Bigl( \sum_{n_1,n_2,n_3}\prod_{j=1}^3 \Psi^2(n_1,n_2,n_3) \Bigr)$ justifying the last step.
 
  We observe that, by definition
\begin{equs}
\sum_{n_1,n_2,n_3} \prod_{j=1}^3 \1_{|n_j| \sim N_j}\Psi^2(n_1,n_2,n_3) 
&= \sum_{\substack{\sigma,\sigma_1,\sigma_2,\sigma_3,\\ \sigma',\sigma'_1,\sigma'_2,\sigma'_3}}  \sum_{n_1,n_2,n_3} \jpb{n_{123}}^{-2} \prod_{j=1}^3 \1_{|n_j| \sim N_j} \jpb{n_j}^{-2+2\alpha} \\
& \quad \times   \sum_{m \in \Z} \sum_{m' \in \Z} \frac{1}{\jpb{m}} \frac{1}{\jpb{m'}}  \1_{|\phi-m| \le 1} \1_{|\phi-m'| \le 1} 
\end{equs} 
where for a set of signs fixed, $\phi \equiv \phi(n_1,n_2,n_3) = \sigma\jpb{n_{123}} +\sigma_1 \jpb{n_1} + \sigma_2 \jpb{n_2} + \sigma_3 \jpb{n_3}$.

Note for fixed $N_1,N_2,N_3$, we have at most $\lesssim \max(N_1,N_2,N_3)$ choices of $m$. Hence by taking $\ell^\infty$ norm on factor $\1_{|\phi-m| \le 1}$ and $\ell^1$ norm on $\frac{1}{\jpb{m}}$, and note there are only a finite number of choices of signs,
 we have 
\begin{equs}
{}&\sum_{n_1,n_2,n_3} \prod_{j=1}^3 \1_{|n_j| \sim N_j}\Psi^2(n_1,n_2,n_3) \\
\lesssim & \log^2(\max(N_1,N_2,N_3)) \sup_{m,m' \in \Z} \sum_{n_1,n_2,n_3} \jpb{n_{123}}^{-2} \prod_{j=1}^3 \1_{|n_j| \sim N_j} \jpb{n_j}^{-2+2\alpha}  \1_{|\phi-m| \le 1} \1_{|\phi-m'| \le 1}\\
\le & \log^2(\max(N_1,N_2,N_3)) \sup_{m \in \Z} \sum_{n_1,n_2,n_3} \jpb{n_{123}}^{-2} \prod_{j=1}^3 \1_{|n_j| \sim N_j} \jpb{n_j}^{-2+2\alpha}  \1_{|\phi-m| \le 1} 
\end{equs} 
Now we apply the following cubic sum estimate  from  \eqref{eq:cubic-sum}:
\begin{equ}
\sup_{m \in \Z} \sum_{n_1,n_2,n_3 \in \Z^2} \jpb{n_{123}}^{2(s-1)}  \Bigl[ \prod_{j=1}^3 \1_{|n_j|\sim N_j}  \jpb{n_{j}}^{-2+2\alpha} \1_{|\phi -m|\le 1} \Bigr] 
\lesssim  \maxcurly{N_1,N_2,N_3}^{2(s-s_{\alpha})}
\end{equ}
to obtain:
\begin{equs}
\sum_{n_1,n_2,n_3} \prod_{j=1}^3 \1_{|n_j| \sim N_j}\Psi^2(n_1,n_2,n_3)
& \lesssim  \log^2(\max(N_1,N_2,N_3))  \max(N_1,N_2,N_3)^{-2s_\alpha}\\
&\lesssim \max(N_1,N_2,N_3)^{-2s_\alpha+},
\label{cubic-repeat}
\end{equs}
 where $\kappa$ is any small positive number. 
 
 Put everything together, we have 
\begin{equs}
&\sum_{n_j, j \in \pi^c } \jpb{n_{\pi^c}}^{2(s-1)} \Bigl( \sum_{n_i: (i,j)\in \cP, n_j=-n_i} \mathbf{1}_{|n_{1234567}| \sim N_{1234567}}\prod_{j=1}^7 \1_{|n_j| \sim N_j} \Psi(n_1,n_2,n_3) \\
& \qquad \times \frac{1}{\jpb{n_4}^{1-\alpha}}\Psi(n_5,n_6,n_7) \Bigr)^2 \\
&\lesssim (N_{1234567}^{-\kappa} \vee N_{1234567}^{2s-2+2\alpha+2\kappa}) N^{-\kappa}_4 \max(N_1,N_2,N_3)^{-2s_\alpha+\gamma}\max(N_5,N_6,N_7)^{-2s_\alpha+\gamma} 
\end{equs}
where $\gamma$ comes from the $\log$ factor. By setting $s<1-\alpha-\kappa$, and $\epsilon = \min(1-\alpha+\kappa-s,\kappa,2s_\alpha-\gamma)>0$ such that
\begin{equs}
&\sum_{n_j, j \in \pi^c } \jpb{n_{\pi^c}}^{2(s-1)} \Bigl( \sum_{n_i: (i,j)\in \cP, n_j=-n_i} \mathbf{1}_{|n_{1234567}| \sim N_{1234567}}\prod_{j=1}^7 \1_{|n_j| \sim N_j} \Psi(n_1,n_2,n_3) \\
& \qquad \times \frac{1}{\jpb{n_4}^{1-\alpha}}\Psi(n_5,n_6,n_7) \Bigr)^2 \\
&\lesssim \max(N_{1234567},N_1,N_2,N_3,N_4,N_5,N_6,N_7)^{-\epsilon}
\end{equs}
completing the proof for the case where $4$ is not paired.

We now conduct the case study obtaining the bound  $\Bigl( \sum_{n_1,n_2,n_3}\prod_{j=1}^3  \1_{|n_j| \sim N_j} \Psi^2(n_1,n_2,n_3) \Bigr) ^2$, to which we have apply cubic estimates in the paragraph above.

\textbf{(1a)} Suppose that $\cP=\emptyset$. We have
\begin{equs}
{}&\sum_{n_1,n_2,n_3,n_5,n_6,n_7}  \prod_{1,2,3,5,6,7} \1_{|n_j| \sim N_j} \Psi^2(n_1,n_2,n_3) \Psi^2(n_5,n_6,n_7) \\
& = \Big(\sum_{n_1,n_2,n_3} \prod_{j=1}^3 \1_{|n_j| \sim N_j} \Psi^2(n_1,n_2,n_3) \Big)\Big(\sum_{n_4,n_5,n_6} \prod_{j=4}^6 \1_{|n_j| \sim N_j}\Psi^2(n_4,n_5,n_6) \Big).\end{equs}

\textbf{(1b)} Assume that there is one pair, set $\cP=\{(1,7)\}$. Since $n_7=-n_1$, applying Cauchy-Schwartz to $n_1$, we can reduce the sum to a product of two sums:
\begin{equs}
{}& \sum_{n_j, j \in\pi^c \backslash\{4\}}  \Bigl( \sum_{n_i: (i,j)\in \cP, n_j=-n_i} \prod_{1,2,3,5,6,7} \1_{|n_j| \sim N_j} \Psi(n_1,n_2,n_3) \Psi(n_5,n_6,n_7) \Bigr)^2\\
= &\sum_{n_2,n_3,n_5,n_6} \prod_{2,3,5,6} \1_{|n_j| \sim N_j}  \Bigl( \sum_{n_1}\1_{|n_1| \sim N_1} \1_{|-n_1| \sim N_7} \Psi(n_1,n_2,n_3) \Psi(n_5,n_6,-n_1) \Bigr)^2 \\
\le & \sum_{n_2,n_3,n_5,n_6} \prod_{2,3,5,6} \1_{|n_j| \sim N_j} \Bigl( \sum_{n_1} \1_{|n_1| \sim N_1} \Psi^2(n_1,n_2,n_3) \Bigr) \Bigl( \sum_{n_1} \1_{|-n_1| \sim N_7}  \Psi^2(n_5,n_6,-n_1) \Bigr) \\
= &\Bigl( \sum_{n_1,n_2,n_3}\prod_{j=1}^3  \1_{|n_j| \sim N_j} \Psi^2(n_1,n_2,n_3) \Bigr) \Bigl( \sum_{n_5,n_6,n_7}\prod_{j=5}^7 \1_{|n_j| \sim N_j} \Psi^2(n_5,n_6,n_7) \Bigr).
\end{equs}

\textbf{(1c)}Assuming that there are two pairs, set $\cP=\{(1,7), (2,6)\}$. We  apply Cauchy-Schwarz with respect to $n_1,n_2$,
\begin{equs}
{}& \sum_{n_j, j \in\pi^c \backslash\{4\}}  \Bigl( \sum_{n_i: (i,j)\in \cP, n_j=-n_i} \prod_{1,2,3,5,6,7} \1_{|n_j| \sim N_j} \Psi(n_1,n_2,n_3) \Psi(n_5,n_6,n_7) \Bigr)^2\\
= &\sum_{n_3,n_5} \prod_{3,5} \1_{|n_j| \sim N_j}  \Bigl( \sum_{n_1,n_2} \1_{|n_1| \sim N_1}  \1_{|n_2| \sim N_2}  \Psi(n_1,n_2,n_3) \1_{|-n_2| \sim N_6}\1_{|-n_1| \sim N_7}\Psi(n_5,-n_2,-n_1) \Bigr)^2 \\
\le & \sum_{n_3,n_5} \prod_{3,5} \1_{|n_j| \sim N_j} \Bigl( \sum_{n_1,n_2} \1_{|n_1| \sim N_1} \1_{|n_2| \sim N_2} \Psi^2(n_1,n_2,n_3) \Bigr) \Bigl( \sum_{n_1,n_2} \1_{|-n_2| \sim N_5} \1_{|-n_1| \sim N_7}  \Psi^2(n_5,-n_2,-n_1) \Bigr) \\
= &\Bigl( \sum_{n_1,n_2,n_3}\prod_{j=1}^3 \Psi^2(n_1,n_2,n_3) \Bigr) \Bigl( \sum_{n_5,n_6,n_7}\prod_{j=5}^7 \Psi^2(n_5,n_6,n_7) \Bigr).
\end{equs}

\textbf{(1d)}\, The three pair case  is analogous to case (1c).

Case 2.  The index $4$ is paired with one index from $\{1,2,3\}\cup \{5,6,7\}$. By the symmetry, we have the following cases:\\
(2a) $\cP=\{(1,4)\}$;\\
(2b) $\cP=\{(1,4), (2,5)\}$;\\
(2c) $\cP=\{(1,4), (2,5), (3,6)\}$.

 Denote the pairing after removing $(3,4)$ by $\mathcal{P}\backslash\{3,4\}$. Inserting dyadic decomposition with respect to frequencies $n_1,n_2$ and $n_3$. Then by applying the basic resonance estimate Lemma \ref{lem:basic-resonant}, with summing over $n_4$, we have

The main idea now is first sum over $n_4$ using basic resonant estimate instead of discrete convolution inequality. Indeed, the basic resonance estimate Lemma \ref{lem:basic-resonant}. For $\f 14 \le \alpha < \f 12$ 
the exponent collected in the estimate is $-\tilde{s}_{\alpha} =-( \frac{3}{2}-2\alpha)$. This leads to
\begin{equs}
{}&\sum_{n_4} \1_{|n_4| \sim N_4} \Psi(-n_4,n_2,n_3) \jpb{n_4}^{-1+\alpha}\\
= & \sum_{n_4} \1_{|n_4| \sim N_4} \sum_{m \in \Z} \frac{1}{\jpb{m}} \jpb{-n_4+n_2+n_3}^{-1} \jpb{n_4}^{-2+2\alpha} \jpb{n_2}^{-1+\alpha} \jpb{n_3}^{-1+\alpha} \\
\lesssim &  N^{-\kappa}_4 \mathrm{log}(N_4) \jpb{n_{23}}^{-\f 32+2\alpha+\kappa} \jpb{n_2}^{-1+\alpha} \jpb{n_3}^{-1+\alpha}.
\end{equs}
Hence, we have 
\begin{equs}
{} &\sum_{n_j, j \in \pi^c } \jpb{n_{\pi^c}}^{2(s-1)} \Bigl( \sum_{n_i: (i,j)\in \cP, n_j=-n_i} \mathbf{1}_{|n_{1234567}| \sim N_{1234567}} \prod_{j=1}^7 \1_{|n_j| \sim N_j} \Psi(n_1,n_2,n_3) \\
& \qquad \times \frac{1}{\jpb{n_4}^{1-\alpha}}\Psi(n_5,n_6,n_7) \Bigr)^2  \\
\lesssim & N^{2s-2}_{1234567} N^{-\kappa}_4 \mathrm{log}(N_4)  \sum_{n_j, j \in \pi^c }  \Bigl( \sum_{n_i: (i,j)\in \cP\backslash \{(1,4)\}, n_j=-n_i} \prod_{j=2,3,5,6,7} \1_{|n_j| \sim N_j} \\
& \qquad \times \jpb{n_{23}}^{-\f 32+2\alpha+\kappa} \jpb{n_2}^{-1+\alpha} \jpb{n_3}^{-1+\alpha} \Psi(n_5,n_6,n_7) \Bigr)^2 ,
\label{4-paired-eq-1}
\end{equs}
where we used the factor after paired, $n_1=-n_4$, $$ \jpb{n_{\pi^c}}^{2(s-1)} \mathbf{1}_{|n_{1234567}| \sim N_{1234567}}  \sim N^{2s-2}_{1234567}.$$
Now we consider the sub-cases.

\textbf{(2a) } Suppose that $\cP=\{(1,4)\}$, so here is no other pair. The summation  in \eqref{4-paired-eq-1}reduces to:
\begin{equs}
{}& \sum_{n_j, j \in \pi^c }  \Bigl( \sum_{n_i: (i,j)\in \cP\backslash \{(1,4)\}, n_j=-n_i} \prod_{j=2,3,5,6,7} \1_{|n_j| \sim N_j} \\
& \qquad \times \jpb{n_{23}}^{-\f 32+2\alpha+\kappa} \jpb{n_2}^{-1+\alpha} \jpb{n_3}^{-1+\alpha} \Psi(n_5,n_6,n_7) \Bigr)^2 \\
=&\Big( \sum_{n_2,n_3} \prod_{j=2}^3 \1_{|n_j| \sim N_j} \jpb{n_{23}}^{-3+4\alpha+2\kappa} \jpb{n_2}^{-2+\alpha} \jpb{n_3}^{-2+\alpha} \Big) \Big(\sum_{n_5,n_6,n_7} \prod_{j=5}^7 \1_{|n_j| \sim N_7} \Psi^2(n_5,n_6,n_7) \Big).
\end{equs}

\textbf{(2b) Let $\cP=\{(1,4), (2,5)\}$}. We have 
\begin{equs}
{}& \sum_{n_j, j \in \pi^c }  \Bigl( \sum_{n_i: (i,j)\in \cP\backslash \{(1,4)\}, n_j=-n_i} \prod_{j=2,3,5,6,7} \1_{|n_j| \sim N_j} \\
& \qquad \times \jpb{n_{23}}^{-\f 32+2\alpha+\kappa} \jpb{n_2}^{-1+\alpha} \jpb{n_3}^{-1+\alpha} \Psi(n_5,n_6,n_7) \Bigr)^2 \\
= &\sum_{n_3,n_6,n_7} \prod_{j=3,6,7} \1_{|n_j| \sim N_j} \jpb{n_3}^{-2+\alpha} \Big(\sum_{n_2} \jpb{n_{23}}^{-\f 32+2\alpha+\kappa} \1_{|n_2| \sim N_2}\jpb{n_2}^{-1+\alpha} \1_{|-n_2| \sim N_5} \Psi(-n_2,n_6,n_7) \Big)^2\\
\le & \sum_{n_3,n_6,n_7} \prod_{j=3,6,7} \1_{|n_j| \sim N_j} \jpb{n_3}^{-2+\alpha}   \Big(\sum_{n_2} \jpb{n_{23}}^{-3+4\alpha+2\kappa} \jpb{n_2}^{-2+\alpha} \Big) \Big(\sum_{n_2} \1_{|-n_2| \sim N_5} \Psi^2(-n_2,n_6,n_7)  \Big)
\end{equs}
where in the last step we applied the Cauchy Schwarz inequality.

By relabelling $-n_2$ to be $n_5$, we arrive at
\begin{equs}
\Big( \sum_{n_2,n_3} \prod_{j=2}^3 \1_{|n_j| \sim N_j} \jpb{n_{23}}^{-3+4\alpha+2\kappa} \jpb{n_2}^{-2+\alpha} \jpb{n_3}^{-2+\alpha} \Big) \Big(\sum_{n_5,n_6,n_7} \prod_{j=5}^7 \1_{|n_j| \sim N_7} \Psi^2(n_5,n_6,n_7) \Big).
\end{equs}

\textbf{(2c) The case $\cP=\{(1,4), (2,5), (3,6)\}$} can be treated similar to (2b) and we omit the repetition.

It reduces to control the bound of $\Big( \sum_{n_2,n_3} \prod_{j=2}^3 \1_{|n_j| \sim N_j} \jpb{n_{23}}^{-3+4\alpha+2\kappa} \jpb{n_2}^{-2+\alpha} \jpb{n_3}^{-2+\alpha} \Big)$. Thanks to the discrete convolution inequality, we have 
\begin{equs}
{} & \sum_{n_2,n_3} \prod_{j=2}^3 \1_{|n_j| \sim N_j} \jpb{n_{23}}^{-3+4\alpha+2\kappa} \jpb{n_2}^{-2+2\alpha} \jpb{n_3}^{-2+2\alpha} \lesssim \max(N_2,N_3)^{-3+8\alpha+2\kappa}
\end{equs}
where we require $3-4\alpha-2\kappa<2$ which yields $\alpha > \f 1 4 -\f \kappa 2$. 

By \eqref{cubic-repeat} and put all together, we obtain
\begin{equs}
{} &\sum_{n_j, j \in \pi^c } \jpb{n_{\pi^c}}^{2(s-1)} \Bigl( \sum_{n_i: (i,j)\in \cP, n_j=-n_i} \mathbf{1}_{|n_{1234567}| \sim N_{1234567}} \prod_{j=1}^7 \1_{|n_j| \sim N_j} \Psi(n_1,n_2,n_3) \\
& \qquad \times \frac{1}{\jpb{n_4}^{1-\alpha}}\Psi(n_5,n_6,n_7) \Bigr)^2  \\
\lesssim &N^{2s-2}_{1234567} N^{-2\kappa}_4 \mathrm{log}^2(N_4) \mathrm{log}^2(\max(N_5,N_6,N_7))  \max(N_5,N_6,N_7)^{-2s_\alpha} \max(N_2,N_3)^{-3+8\alpha+2\kappa} \\
\lesssim &N^{2s-2}_{1234567} N^{-2\kappa+2\gamma}_4  \max(N_5,N_6,N_7)^{-2s_\alpha+2\gamma} \max(N_2,N_3)^{-3+8\alpha+2\kappa},
\end{equs}
where we absorb the log factor by assuming $\gamma \ll \kappa$.

Note if $s<1-\alpha-\kappa$ and $\alpha < \f 3 8$ and also use the factor in this case $N_1 \sim N_4$, we have 
\begin{equs}
{} &\sum_{n_j, j \in \pi^c } \jpb{n_{\pi^c}}^{2(s-1)} \Bigl( \sum_{n_i: (i,j)\in \cP, n_j=-n_i} \mathbf{1}_{|n_{1234567}| \sim N_{1234567}} \prod_{j=1}^7 \1_{|n_j| \sim N_j} \Psi(n_1,n_2,n_3) \\
& \qquad \times \frac{1}{\jpb{n_4}^{1-\alpha}}\Psi(n_5,n_6,n_7) \Bigr)^2  \\
\lesssim &\max(N_{1234567},N_1,N_2,N_3,N_4,N_5,N_6,N_7)^{-\epsilon}
\end{equs}
where we choose $\epsilon = \min(2\alpha+2\kappa,\kappa - \epsilon,3-8\alpha-2\kappa,s_\alpha-\gamma)>0$ which completes the proof.
\end{proof}

\end{document}